\documentclass{amsart}
\usepackage{graphicx}
\usepackage[all]{xy}
\usepackage{mathtools}%

\xyoption{arc}

\usepackage{txfonts, amsmath,amstext,amsthm,amscd,amsopn,verbatim,amssymb, amsfonts}
\usepackage{fullpage}
\usepackage[makeroom]{cancel}

\usepackage{float}
\restylefloat{table}

\usepackage[bbgreekl]{mathbbol}

\usepackage{tikz}
\usetikzlibrary{matrix}
\usetikzlibrary{shapes}
\usetikzlibrary{arrows}
\usetikzlibrary{calc,3d}
\usetikzlibrary{decorations,decorations.pathmorphing,decorations.pathreplacing,decorations.markings}
\usetikzlibrary{through}

\tikzset{ext/.style={circle, draw,inner sep=1pt},int/.style={circle,draw,fill,inner sep=1.4pt},nil/.style={inner sep=1pt}}
\tikzset{cy/.style={circle,draw,fill,inner sep=2pt},scy/.style={circle,draw,inner sep=2pt},scyx/.style={draw,cross out,inner sep=2pt},scyt/.style={draw,regular polygon,regular polygon sides=3,inner sep=0.95pt}}
\tikzset{exte/.style={circle, draw,inner sep=3pt},inte/.style={circle,draw,fill,inner sep=3pt}}
\tikzset{diagram/.style={matrix of math nodes, row sep=3em, column sep=2.5em, text height=1.5ex, text depth=0.25ex}}
\tikzset{diagram2/.style={matrix of math nodes, row sep=0.5em, column sep=0.5em, text height=1.5ex, text depth=0.25ex}}

\tikzset{
  rightorange/.style={
    decoration={markings,mark=at position .8 with {\arrow[scale=1.5,orange]{latex}}},
    postaction={decorate},
    shorten >=0.4pt}}
\tikzset{
  leftorange/.style={
    decoration={markings,mark=at position .55 with {\arrowreversed[scale=1.5,orange]{latex}}},
    postaction={decorate},
    shorten >=0.4pt}}

\tikzset{
  rightyellow/.style={
    decoration={markings,mark=at position .8 with {\arrow[scale=1.6,yellow]{latex}}},
    postaction={decorate},
    shorten >=0.4pt}}
\tikzset{
  leftyellow/.style={
    decoration={markings,mark=at position .55 with {\arrowreversed[scale=1.6,yellow]{latex}}},
    postaction={decorate},
    shorten >=0.4pt}}

 \tikzset{
  rightpurple/.style={
    decoration={markings,mark=at position .8 with {\arrow[scale=1.6,purple]{latex}}},
    postaction={decorate},
    shorten >=0.4pt}}
\tikzset{
  leftpurple/.style={
    decoration={markings,mark=at position .55 with {\arrowreversed[scale=1.6,purple]{latex}}},
    postaction={decorate},
    shorten >=0.4pt}}

\tikzset{
  rightgreen/.style={
    decoration={markings,mark=at position .8 with {\arrow[scale=1.6,green]{latex}}},
    postaction={decorate},
    shorten >=0.4pt}}
\tikzset{
  leftgreen/.style={
    decoration={markings,mark=at position .55 with {\arrowreversed[scale=1.6,green]{latex}}},
    postaction={decorate},
    shorten >=0.4pt}}

\tikzset{
  rightblue/.style={
    decoration={markings,mark=at position .8 with {\arrow[scale=1.7,blue]{latex}}},
    postaction={decorate},
    shorten >=0.4pt}}
\tikzset{
  leftblue/.style={
    decoration={markings,mark=at position .55 with {\arrowreversed[scale=1.7,blue]{latex}}},
    postaction={decorate},
    shorten >=0.4pt}}
\tikzset{
  rightred/.style={
    decoration={markings,mark=at position .45 with {\arrow[scale=1.3,red]{latex}}},
    postaction={decorate},
    shorten >=0.4pt}}
\tikzset{
  leftred/.style={
    decoration={markings,mark=at position .2 with {\arrowreversed[scale=1.3,red]{latex}}},
    postaction={decorate},
    shorten >=0.4pt}}

\def\id{{\mbox{1 \hskip -7pt 1}}}
\newcommand{\sgn}{{sgn}}
 \newcommand{\lon}{\longrightarrow}
 \newcommand{\bu}{\bullet}
 
 \newcommand{\rar}{\rightarrow}

\newcommand{\Der}{\mathrm{Der}}

 \newcommand{\Z}{{\mathbb Z}}
 \newcommand{\bS}{{\mathbb S}}

 \newcommand{\R}{{\mathbb{R}}}
 \newcommand{\N}{{\mathbb N}}
 \newcommand{\K}{{\mathbb K}}

 \newcommand{\ot}{\otimes}

\newcommand{\Def}{\mbox{\sf Def}}
\newcommand{\fGC}{\mbox{\sf fGC}}
\newcommand{\fcGC}{\mbox{\sf fcGC}}
\newcommand{\GC}{\mbox{\sf GC}}

 \newcommand{\LB}{\mathcal{L}{ieb}}

\newcommand{\LBcd}{\mathcal{L}{ieb}_{c,d}}

\newcommand{\HoLBcd}{\mathcal{H}{olieb }_{c,d}}

\newcommand{\wHoLBcd}{\widehat{\mathcal{H}{olieb}}_{c,d}}

\newcommand{\HoLB}{\mathcal{H}{olieb}}

 \newcommand{\Beq}{\begin{equation}}
 \newcommand{\Eeq}{\end{equation}}
 \newcommand{\Beqr}{\begin{eqnarray}}
 \newcommand{\Eeqr}{\end{eqnarray}}
 \newcommand{\Beqrn}{\begin{eqnarray*}}
 \newcommand{\Eeqrn}{\end{eqnarray*}}
 \newcommand{\Ba}{\begin{array}}
 \newcommand{\Ea}{\end{array}}
 \newcommand{\Bi}{\begin{itemize}}
 \newcommand{\Ei}{\end{itemize}}
 \newcommand{\Bc}{\begin{center}}
 \newcommand{\Ec}{\end{center}}
\newcommand{\fa}{{\mathfrak a}}

 \newcommand{\fg}{{\mathfrak g}}
 
 \newcommand{\fm}{{\mathfrak m}}

\newcommand{\fr}{{\mathfrak r}}
\newcommand{\fs}{{\mathfrak s}}

\newcommand{\ft}{{\mathfrak t}}

\newcommand{\ii}{{\mathfrak i}}
  \newcommand{\f}{{\mathcal O}r}
 
 \newcommand{\cA}{{\mathcal A}}
 \newcommand{\cB}{{\mathcal B}}
 \newcommand{\cC}{{\mathcal C}}
 
 \newcommand{\cE}{{\mathcal E}}
 \newcommand{\cF}{{\mathcal F}}
 \newcommand{\cG}{{\mathcal G}}
 
 \newcommand{\cI}{{\mathcal I}}

 \newcommand{\caL}{{\mathcal L}}

 \newcommand{\cP}{{\mathcal P}}
 
 \newcommand{\cR}{{\mathcal R}}
 \newcommand{\cS}{{\mathcal S}}

 \newcommand{\al}{\alpha}
 \newcommand{\be}{\beta}
 \newcommand{\ga}{\gamma}
 
 \newcommand{\Ga}{\Gamma}

 \newcommand{\om}{\omega}

 \newcommand{\Hom}{\text{Hom}}

 \newcommand{\sip}{\smallskip}
 \newcommand{\bip}{\bigskip}
 \newcommand{\mip}{\vspace{2.5mm}}

\theoremstyle{plain}
\swapnumbers

\newtheorem{prop-def}[theorem]{Proposition-definition}

\newtheorem{f-theorem}{Formality Theorem}[section]
\newtheorem{main-theorem}{Main~Theorem}[section]
\newtheorem{section-theorem}{Theorem}[section]

\theoremstyle{definition}

\begin{document}

 \sloppy

 \newenvironment{proo}{\begin{trivlist} \item{\sc {Proof.}}}
  {\hfill $\square$ \end{trivlist}}

\long\def\symbolfootnote[#1]#2{\begingroup%
\def\thefootnote{\fnsymbol{footnote}}\footnote[#1]{#2}\endgroup}

\title{ Multi-oriented props and  homotopy algebras with branes}

\author{Sergei~Merkulov}
\address{Sergei~Merkulov:  Mathematics Research Unit, Luxembourg University, 
Grand Duchy of Luxembourg }
\email{sergei.merkulov@uni.lu}

\begin{abstract} We introduce a new category of differential graded {\em multi-oriented}\, props whose representations (called homotopy algebras with branes) in a graded vector space require a choice of a collection of $k$ linear subspaces
in that space,  $k$ being the number of extra orientations (if $k=0$ this structure recovers an ordinary prop); symplectic vector spaces equipped  with $k$ Lagrangian subspaces play a distinguished role in this theory.
 Manin triples is a classical example of an algebraic structure (concretely, a Lie bialgebra structure) given in terms of a  vector space and its subspace; in the context of this paper Manin triples are precisely symplectic Lagrangian representations of the {\em 2-oriented} generalization of the classical operad  of Lie algebras. In a sense, the theory of multi-oriented props provides us with a far reaching strong homotopy generalization
of Manin triples type constructions.

\sip

The homotopy theory of multi-oriented props can be quite non-trivial (and different from that of ordinary props).
The famous Grothendieck-Teichm\"uller group acts faithfully as homotopy
non-trivial automorphisms on infinitely many multi-oriented props, a fact which motivated much the present work as it gives us a hint
to a non-trivial deformation quantization theory in every geometric dimension $d\geq 4$ generalizing to higher dimensions Drinfeld-Etingof-Kazhdan's quantizations of Lie bialgebras (the case $d=3$)
and  Kontsevich's quantizations of Poisson structures (the case $d=2$).


\end{abstract}
\maketitle
\markboth{}{}


{\Large
\section{\bf Introduction}
}

\mip

\subsection{Why bother with multi-oriented props?} A short answer to this question: ``Because of the Grothendieck-Teichm\"uller group $GRT_1$". It is the latter beautiful and mysterious structure which is the {main} motivation for introducing and studying  a new category of multi-oriented props as well as their representations (``homotopy algebras with branes"). In geometric dimensions $2$ and $3$ the group $GRT_1$ acts   on some ordinary props of odd/even  strong homotopy Lie bialgebras \cite{MW1} and plays thereby the classifying role in the associated transcendental deformation quantizations of Poisson and, respectively, ordinary Lie bialgebra structures. In geometric dimension $d\geq 4$ the group $GRT_1$ survives in the form of symmetries of some {\em multi-oriented}\, props of even/odd homotopy Lie bialgebras so that  deformation quantizations in higher dimensions (in which $GRT_1$ retains its  fundamental classifying role) should involve
a really new class of algebro-geometric structures --- the  homotopy algebras {\em with branes}. It is an attempt to understand what could be a higher ($d\geq 4$) analogue of two famous formality theorems, one for Poisson structures \cite{Ko2} (the case $d=2$), and another for ordinary Lie bialgebras \cite{EK,Me4} (the case $d=3$),   that lead the author to the  category of multi-oriented props after reading the paper \cite{Z} by Marko \v Zivkovi\' c and its predecessor  \cite{Wi2} by Thomas Willwacher (see \S 5 for a brief but self-contained description of their remarkable results).

\sip

 It is not that hard to define multi-oriented props in general, and multi-oriented generalizations of some concrete classical operads and props in particular, at the purely combinatorial level: the rules of the game with multi-oriented decorated graphs are  more or less standard (see \S 2) --- for any given $k\geq 1$ one just adds $k$ extra directions (=orientations) to each edge/leg of a 1-oriented graph,
 $$
\begin{tikzpicture}[baseline=-1ex]
\node[] (a) at (0,-0.05) {};
\node[] (b) at (1.5,-0.05) {};
\draw (a) edge[->] (b);
\end{tikzpicture}
\rightsquigarrow
\Ba{c}
\begin{tikzpicture}[baseline=-2ex]
\draw (0,0) edge[rightblue] node[above] {\ \ \ $\scriptstyle 1$}(0.9,0);
\draw (0.9,0) edge[leftred] node[above] {$\scriptstyle 2$\ \ }(1.5,0);
\draw (1.5,0) edge[] node[above] {$...$\ \ \ }(1.9,0);
\draw (1.9,0) edge[leftpurple] node[above] {\ \ \ $\scriptstyle k$ }(2.3,0);
\draw (2.2,0) edge[-latex](2.7,0);
\end{tikzpicture}
\Ea
$$
and defines rules for multi-oriented prop compositions via contractions along admissible multi-oriented subgraphs. However, it is much less evident how to transform
 that more or less standard rules into non-trivial and interesting representations (i.e.\ examples) --- the
 intuition from the theory of ordinary (wheeled) props does not help much. Adding new $k$ directions to each edge of a decorated graph of an ordinary prop can be naively understood
 as extending that ordinary prop into a $2^{k}$-coloured one, but then the requirement that the new directions on graph edges  do not create ``wheels"
 (that is, closed directed paths of edges with respect to {\em each}\, of the new orientations) kills that
 naive picture immediately ---  the elements of the set of $2^{k}$ new colours start interacting with each other in a non-trivial way. We know which structure distinguishes ordinary props (the ones with no wheels
 in the given single orientation, i.e.\ the ones which are {\em 1-oriented}\, in the terminology of this paper) from the ordinary wheeled props (that is, {\em 0-oriented 1-directed}\, props in the terminology of this paper) in terms of representations in, say, a graded vector space $V$ --- it is the {\em dimension of $V$}. In general, a wheeled prop can have well-defined representations in $V$ only in the case  $\dim V<\infty$ as graphs with wheels generate the trace operation $V\times V^*\rar \K$ which explodes in the case
 $\dim V=\infty$; this phenomenon explains the need for {\em 1-oriented}\, props.
 \sip

 How to explain the need for all (or some) of the {\em extra}\, $k$ directions to be oriented?
 Which structure on  a graded vector space $V$ can be used to separate (in the sense of representations) $0$-oriented $(k+1)$-directed  props  from $(k+1)$-oriented ones? Perhaps the main result of this paper is a rather
 surprising answer to that question: one has to work again
with a certain class of infinite-dimensional vector spaces $V$, but now equipped with $k$ linear subspaces $W_1,\ldots, W_k \subset V$ together with their complements, and interpret a {\em single}\,  element of a $(k+1)$-oriented prop $\cP$
as a {\em collection of $k$  linear maps}\,  from various intersections of subspaces $W_\bu$ and their complements  and their duals to themselves; then indeed graphs with wheels in one or another extra
``coloured direction" get exploded under generic representations and hence must be prohibited. Representations when $V$ has a symplectic structure and the
subspaces $W_1,\ldots, W_k$ are Lagrangian play a special role in this story,
a fact which becomes obvious once all the general definitions are given (see \S 4).

 \sip

 A well-known example of such a ``brane" algebraic structure in finite dimensions is provided by Manin triples \cite{D}. In the context of this paper Manin triples construction  emerges as a (reduced symplectic Lagrangian) representation of the 2-oriented {\em operad}\, of Lie algebras (note that one can not describe
 Manin triples using  {\em ordinary operads} --- one needs an   {\em ordinary prop}\, of Lie bialgebras for that purpose).
 In a sense, multi-oriented  props provide us with a far reaching strong homotopy generalization of Manin triples  type constructions; they are really a new kind of ``Cheshire cat smiles"  controlling (via representations) homotopy algebras {\em with branes} and admitting in some interesting cases a highly non-trivial action of the Grothendieck-Teichm\"uller group \cite{A,MW1}.

\sip

This is the first of a sequence of papers on multi-oriented props. In the following paper \cite{Me5} we study several transcendental constructions with multi-oriented props (elucidating their role
 as the construction material for building new highly non-trivial representations of {\em ordinary props}) and use them to prove several concrete deformation quantization theorems.
 This  paper attempts to be as simple as possible and aims for more  general audience: we explain here the main notion, illustrate it with  examples, prove some theorems on multi-oriented resolutions, and, most importantly, discuss in full details the most surprising part of the story ---
 the representation theory of multi-oriented props in the category of dg  vector spaces with branes.

\subsection{Finite dimensionality versus infinite one in the context of ordinary props}
The theory of (wheeled) operads and props originated in 60s and 70s in algebraic topology, and has seen since an explosive development
(see, e.g., the books \cite{LV,MSS}  or the articles \cite{Ma, MMS,Va}  for details and references). Operads and props provide us with effective tools to discover surprisingly deep
and unexpected links between  different theories and even branches of mathematics.

\sip

A building block of a prop(erad) $\cP=\{\cP(m,n)\}_{m,n\geq 0}$ is a graph (often called {\em corolla})
$$
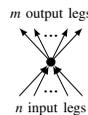
\Ba{c}\resizebox{14mm}{!}  {
\begin{tikzpicture}[baseline=-1ex]
\node[] at (0,0.8) {$\scriptstyle m\ \mathrm{output\ legs}$};
\node[] at (0,-0.8) {$\scriptstyle n\ \mathrm{input\ legs}$};
\node[int] (a) at (0,0) {};
\node[] (u1) at (-0.65,0.7) {};
\node[] (u2) at (-0.35,0.7) {};
\node[] (u3) at (0.35,0.7) {};
\node[] (u4) at (0.65,0.7) {};
\node[] at (0.0,0.45) {...};
\node[] (d1) at (-0.65,-0.7) {};
\node[] (d2) at (-0.35,-0.7) {};
\node[] (d3) at (0.35,-0.7) {};
\node[] (d4) at (0.65,-0.7) {};
\node[] at (0.0,-0.5) {...};
\draw (a) edge[-latex] (u1);
\draw (a) edge[-latex] (u2);
\draw (a) edge[-latex] (u3);
\draw (a) edge[-latex] (u4);
\draw (a) edge[latex-] (d1);
\draw (a) edge[latex-] (d2);
\draw (a) edge[latex-] (d3);
\draw (a) edge[latex-] (d4);
\end{tikzpicture}
}
\Ea
$$

consisting of one vertex  which has $n$ incoming legs and
$m$ outgoing legs and which is decorated with an element of some $\bS_m^{op}\times \bS_n$ module $\cP(m,n)$. Upon a representation of $\cP$ in a graded vector space $V$ this $(m,n)$-corolla gets transformed into a
linear map from $V^{\ot n}\rar V^{\ot m}$, i.e.\ every leg
corresponds, roughly speaking, to a copy of $V$,
$$
\Ba{c}
\begin{tikzpicture}[baseline=-1ex]
\node[] (a) at (0,-0.08) {};
\node[] (u1) at (1,-0.08) {};
\draw (a) edge[-latex] (u1);
\end{tikzpicture}\Ea \ \ \Leftrightarrow \ \ V\ .
$$
Such linear maps can be composed which leads us to the idea of considering all possible graphs,
for example these ones
$$
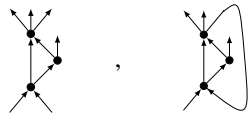
\Ba{c}\resizebox{10mm}{!}  {
\begin{tikzpicture}[baseline=-1ex]
\node[int] (a) at (0,0) {};
\node[int] (u1) at (0.5,0.5) {};
\node[int] (u2) at (0.0,1.0) {};
\node[] (uu3) at (0.0,1.6) {};
\node[] (uu4) at (-0.5,1.6) {};
\node[] (uu5) at (0.5,1.6) {};
\node[] (u11) at (0.5,1.1) {};
\node[] (d1) at (-0.5,-0.6) {};
\node[] (d2) at (-0.0,-0.6) {};
\node[] (d3) at (0.5,-0.6) {};
\draw (a) edge[-latex] (u1);
\draw (a) edge[-latex] (u2);
\draw (a) edge[latex-] (d1);
\draw (a) edge[latex-] (d3);
\draw (u2) edge[-latex] (uu3);
\draw (u2) edge[-latex] (uu4);
\draw (u2) edge[-latex] (uu5);
\draw (u1) edge[-latex] (u2);
\draw (u1) edge[-latex] (u11);
\end{tikzpicture}
}\Ea
\ \ \ \ \ , \ \ \ \ \
\Ba{c}\resizebox{11mm}{!}  {\begin{tikzpicture}[baseline=-1ex]
\node[int] (a) at (0,0) {};
\node[int] (u1) at (0.5,0.5) {};
\node[int] (u2) at (0.0,1.0) {};
\node[] (uu3) at (0.0,1.6) {};
\node[] (uu4) at (-0.5,1.6) {};
\node[] (uu5) at (0.5,1.6) {};
\node[] (u11) at (0.5,1.1) {};
\node[] (d1) at (-0.5,-0.6) {};
\node[] (d2) at (-0.0,-0.6) {};
\node[] (d3) at (0.26,-0.2) {};
\draw (a) edge[-latex] (u1);
\draw (a) edge[-latex] (u2);
\draw (a) edge[latex-] (d1);
\draw (a) edge[latex-] (d3);
\draw (u2) edge[-latex] (uu3);
\draw (u2) edge[-latex] (uu4);
\draw (u2) edge[] (uu5);
\draw (u1) edge[-latex] (u2);
\draw (u1) edge[-latex] (u11);
\draw[-] (0.38,1.45) .. controls (0.7, 1.7) .. (0.8, 0.5);
\draw[-] (0.0,0.0) .. controls (0.9, -0.7) .. (0.8, 0.5);
\end{tikzpicture}
}\Ea
$$
 composed
from corollas by connecting some output legs of one corolla with  input legs of another corolla and so on. The graphs shown above --- when translated into linear maps upon some representation of $\cP$ in $V$ --- give us two very different situations: if the left graph makes sense for representations in both finite- and infinite-dimensional vector spaces $V$, the right graph gives us a well-defined linear map only for {\em finite-dimensional}\, vector spaces $V$ as it contains a closed path of directed edges (``wheel") and hence involves a trace map $V\ot V^*\rar \K$ which is not a well-defined operation in infinite dimensions in general. Hence in order  to be able to work in infinite dimensions\footnote{As the symmetric monoidal category of {\em infinite-dimensional}\, vector spaces is not closed, one must be careful about the definition of the {\em endomorphism prop}\, $\cE nd_V$ in this category, see {\ref{App: subsec on inf dim vector spaces}} for details.}  one has to prohibit certain graphs --- the graphs with wheels --- and work solely with {\em oriented}\, (from bottom to the top) graphs. Similarly, in order to be able to work with certain  completions (defined in {\ref{App: subsubsec on tensor algebra of Ws}}) of various intersections
of linear subspaces $W_1\subset V, \ldots, W_k\subset V$ (``branes"), one has to prohibit certain divergent  multi-directed graphs (which already have no wheels with respect to the basic direction!) --- and this leads us to the new notion of $(k+1)$-{\em oriented}\, prop which takes care about more sophisticated divergences associated with branes (the case $k=0$ recovers the ordinary props); this important ``divergency handling"  part of our story is discussed in detail in \S 4.

\sip

There is a nice generalization of the notion of prop which takes care about collections of vector spaces $W_1, \ldots, W_N$. The corresponding  props are called {\em coloured props}\ and, say, $N$-coloured (wheeled) prop $\cP$ is generated by corollas
$$
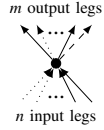
\Ba{c}\resizebox{18mm}{!}  {
\begin{tikzpicture}[baseline=-1ex]
\node[] at (0,0.8) {$\scriptstyle m\ \mathrm{output\ legs}$};
\node[] at (0,-0.8) {$\scriptstyle n\ \mathrm{input\ legs}$};
\node[int] (a) at (0,0) {};
\node[] (u1) at (-0.65,0.7) {};
\node[] (u2) at (-0.35,0.7) {};
\node[] (u3) at (0.35,0.7) {};
\node[] (u4) at (0.65,0.7) {};
\node[] at (0.0,0.45) {...};
\node[] (d1) at (-0.65,-0.7) {};
\node[] (d2) at (-0.35,-0.7) {};
\node[] (d3) at (0.35,-0.7) {};
\node[] (d4) at (0.65,-0.7) {};
\node[] at (0.0,-0.5) {...};
\draw (a) edge[-latex] (u1);
\draw (a) edge[-latex, dotted] (u2);
\draw (a) edge[-latex] (u3);
\draw (a) edge[-latex, dashed] (u4);
\draw (a) edge[latex-,dotted] (d1);
\draw (a) edge[latex-,dotted] (d2);
\draw (a) edge[latex-] (d3);
\draw (a) edge[latex-] (d4);
\end{tikzpicture}
}\Ea
$$
whose input and output legs are ``colored" (say, the unique vertex has $a_1$ input legs in ``straight colour" ,
$a_2$ input legs in ``dotted colour", etc) and correspond to linear maps of the form
$$
W_1^{\ot a_1}\ot W_2^{\ot a_2} \ot \ldots \ot W_N^{\ot a_N} \lon
W_1^{\ot b_1}\ot W_2^{\ot b_2} \ot \ldots \ot W_N^{\ot b_N}, \ \ \ a_1+\ldots+ a_N=n,
 \ \ b_1+\ldots+ b_N=m.
$$
Again it makes sense to talk about wheeled (i.e.\ $0$-{\em oriented $1$-directed}) and ordinary (i.e.\ 1-{\em oriented})
coloured props. In this theory an oriented leg in ``colour" $i\in \{1,\ldots,N\}$ corresponds to the  $i$-th vector space
$$
\Ba{c}
\begin{tikzpicture}[baseline=-1ex]
\node[] (a) at (0,-0.08) {};
\node[] (u1) at (1,-0.08) {};
\draw (a) edge[-latex,dotted] (u1);
\end{tikzpicture}\Ea \ \ \Leftrightarrow \ \ W_i
$$
in the collection $\{W_1, \ldots, W_i, \ldots, W_N\}$.

\subsection{From branes to multidirected props} A $(k+1)$-directed prop $\cP^{k+1}=\{\cP^{k+1}(m,n)\}$ is generated (modulo, in general, some relations) by corollas
$
\Ba{c}\resizebox{21mm}{!}  {
\begin{tikzpicture}[baseline=-1ex]
\node[] at (0,1.1) {$\scriptstyle m\ \mathrm{output\ legs}$};
\node[] at (0,-1.1) {$\scriptstyle n\ \mathrm{input\ legs}$};
\node[int] (a) at (0,0) {};
\node[] (u1) at (-1.2,1) {};
\node[] (u2) at (-0.45,1) {};
\node[] (u3) at (0.45,1) {};
\node[] (u4) at (1.2,1) {};
\node[] at (0.0,0.7) {...};
\node[int] (a) at (0,0) {};
\node[] (d1) at (-1.2,-1) {};
\node[] (d2) at (-0.45,-1) {};
\node[] (d3) at (0.45,-1) {};
\node[] (d4) at (1.2,-1) {};
\node[] at (0.0,-0.7) {...};
\draw (a) edge[-latex, rightblue, leftred] (u1);
\draw (a) edge[-latex, rightblue, rightred] (u2);
\draw (a) edge[-latex, leftblue, leftred] (u3);
\draw (a) edge[-latex,rightblue, rightred] (u4);
\draw (a) edge[latex-, rightblue, leftred] (d1);
\draw (a) edge[latex-, rightblue, rightred] (d2);
\draw (a) edge[latex-, leftblue, leftred] (d3);
\draw (a) edge[latex-,rightblue, rightred] (d4);
\end{tikzpicture}
}\Ea
$
whose vertex is decorated with an element of some module
$\cP^{k+1}(m,n)$ (see \S 2 for details) and
whose input and output legs are decorated with extra (labelled by integers from $1$ to $k$ or by some colours -- blue, red, etc --- as in the picture above) directions. The ``original" (or basic) direction is always shown in pictures in black colour as in the case of ordinary props; it is this basic direction which permits us to call this creature an $(m,n)$-corolla (it can have different numbers of input and output legs with respect to directions in other colours). Comparing this picture to the definition
of an $N$-coloured prop, one can immediately see  that a $(k+1)$-directed prop
is just a special case of a coloured prop when the number of colours is a power of $2$,
\Beq\label{1: formula for number of colours}
N=2^k.
\Eeq
If we allow all possible compositions of such multi-oriented corollas along legs
with identical extra directions, then we get nothing but a $2^k$-coloured prop indeed (called $0$-{\em oriented $k+1$-directed prop}). Keeping in mind the key distinction between ordinary and wheeled props, one might contemplate the possibility of  prohibiting compositions along graphs
which have closed wheels along any of the extra orientations (and no wheels along the basic one), i.e.\ prohibiting compositions of the form
$$
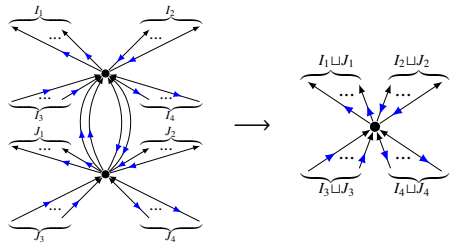
\Ba{c}\resizebox{35mm}{!}  {
\begin{tikzpicture}[baseline=-.65ex]
\node[] at (-1.3,2.2) {$^{\scriptstyle I_1}$};
\node[] at (-1.3,2.05) {$\overbrace{\ \ \ \ \ \ \ \ \ \ \ \ \ \ }$};
\node[] at (-1,1.7) {$...$};
\node[] at (1.3,2.2) {$^{\scriptstyle I_2}$};
\node[] at (1.3,2.05) {$\overbrace{\ \ \ \ \ \ \ \ \ \ \ \ \ \ }$};
\node[] at (1,1.7) {$...$};
\node[] at (-1.3,0.15) {$_{\scriptstyle I_3}$};
\node[] at (-1.3,0.35) {$\underbrace{\ \ \ \ \ \ \ \ \ \ \ \ \ \ }$};
\node[] at (-1.2,0.5) {$...$};
\node[] at (1.3,0.15) {$_{\scriptstyle I_4}$};
\node[] at (1.3,0.35) {$\underbrace{\ \ \ \ \ \ \ \ \ \ \ \ \ \ }$};
\node[] at (1.2,0.5) {$...$};
%
\node[] at (-1.3,-0.20) {$^{\scriptstyle J_1}$};
\node[] at (-1.3,-0.35) {$\overbrace{\ \ \ \ \ \ \ \ \ \ \ \ \ \ }$};
\node[] at (-1.2,-0.5) {$...$};
\node[] at (1.3,-0.20) {$^{\scriptstyle J_2}$};
\node[] at (1.3,-0.35) {$\overbrace{\ \ \ \ \ \ \ \ \ \ \ \ \ \ }$};
\node[] at (1.2,-0.5) {$...$};
\node[] at (-1.3,-2.25) {$_{\scriptstyle J_3}$};
\node[] at (-1.3,-2.05) {$\underbrace{\ \ \ \ \ \ \ \ \ \ \ \ \ \ }$};
\node[] at (-1,-1.7) {$...$};
\node[] at (1.3,-2.25) {$_{\scriptstyle J_4}$};
\node[] at (1.3,-2.05) {$\underbrace{\ \ \ \ \ \ \ \ \ \ \ \ \ \ }$};
\node[] at (1,-1.7) {$...$};
\node[] (uu_1) at (-2,2.0) {};
\node[] (uu_2) at (-1,2.0) {};
\node[] (uu_3) at (2,2.0) {};
\node[] (uu_4) at (1,2.0) {};
\node[] (ud_1) at (-2,0.4) {};
\node[] (ud_2) at (-1,0.4) {};
\node[] (ud_3) at (2,0.4) {};
\node[] (ud_4) at (1,0.4) {};
\node[] (dd_1) at (-2,-2.0) {};
\node[] (dd_2) at (-1,-2.0) {};
\node[] (dd_3) at (2,-2.0) {};
\node[] (dd_4) at (1,-2.0) {};
\node[] (du_1) at (-2,-0.4) {};
\node[] (du_2) at (-1,-0.4) {};
\node[] (du_3) at (2,-0.4) {};
\node[] (du_4) at (1,-0.4) {};
 \node[int] (a) at (0,-1.0) {};
 \node[int] (b) at (0,1.0) {};
\node[] (d) at (0,-2.4) {};
\draw (uu_1) edge[latex-, leftblue] (b);
\draw (uu_2) edge[latex-, leftblue] (b);
\draw (uu_3) edge[latex-, rightblue] (b);
\draw (uu_4) edge[latex-, rightblue] (b);
\draw (b) edge[latex-, leftblue] (ud_1);
\draw (b) edge[latex-, leftblue] (ud_2);
\draw (b) edge[latex-, rightblue] (ud_3);
\draw (b) edge[latex-, rightblue] (ud_4);
\draw (du_1) edge[latex-, leftblue] (a);
\draw (du_2) edge[latex-, leftblue] (a);
\draw (du_3) edge[latex-, rightblue] (a);
\draw (du_4) edge[latex-, rightblue] (a);
\draw (a) edge[latex-, leftblue] (dd_1);
\draw (a) edge[latex-, leftblue] (dd_2);
\draw (a) edge[latex-, rightblue] (dd_3);
\draw (a) edge[latex-, rightblue] (dd_4);
 \draw (b) edge[latex-, rightblue, bend left=50] (a);
 \draw (b) edge[latex-, leftblue, bend right=50] (a);
  \draw (b) edge[latex-, rightblue, bend left=30] (a);
 \draw (b) edge[latex-, leftblue, bend right=30] (a);
\end{tikzpicture}}\Ea
\lon
\Ba{c}\resizebox{24mm}{!}
{\begin{tikzpicture}[baseline=-1ex]
\node[] at (-0.6,1.0) {$^{\scriptstyle I_1\sqcup J_1}$};
\node[] at (-0.7,0.8) {$\overbrace{\  \ \ \ \ \ \ \ \ \ \ }$};
\node[] at (-0.45,0.5) {$...$};
\node[] at (0.6,1.0) {$^{\scriptstyle I_2\sqcup J_2}$};
\node[] at (0.7,0.8) {$\overbrace{\  \ \ \ \ \ \ \ \ \ \ }$};
\node[] at (0.45,0.5) {$...$};
\node[] at (-0.6,-1.0) {$_{\scriptstyle I_3\sqcup J_3}$};
\node[] at (-0.7,-0.8) {$\underbrace{\  \ \ \ \ \ \ \ \ \ \ }$};
\node[] at (-0.45,-0.5) {$...$};
\node[] at (0.6,-1.0) {$_{\scriptstyle I_4\sqcup J_4}$};
\node[] at (0.7,-0.8) {$\underbrace{\  \ \ \ \ \ \ \ \ \ \ }$};
\node[] at (0.45,-0.5) {$...$};
\node[int] (0) at (0,0) {};
\node[] (1) at (-1.2,0.8) {};
\node[] (2) at (-0.3,0.8) {};
\node[] (3) at (0.3,0.8) {};
\node[] (4) at (1.2,0.8) {};
\node[int] (0) at (0,0) {};
\node[] (5) at (-1.2,-0.8) {};
\node[] (6) at (-0.3,-0.8) {};
\node[] (7) at (0.3,-0.8) {};
\node[] (8) at (1.2,-0.8) {};
\draw (1) edge[latex-,leftblue] (0);
\draw (2) edge[latex-,leftblue] (0);
\draw (3) edge[latex-,rightblue] (0);
\draw (4) edge[latex-,rightblue] (0);
\draw (0) edge[latex-,leftblue] (5);
\draw (0) edge[latex-,leftblue] (6);
\draw (0) edge[latex-,rightblue] (7);
\draw (0) edge[latex-,rightblue] (8);
\end{tikzpicture}}\Ea
$$
and allowing only compositions
$$
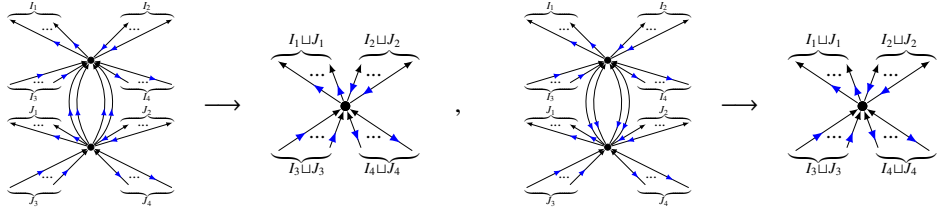
\Ba{c}\resizebox{34mm}{!}  {
\begin{tikzpicture}[baseline=-.65ex]
\node[] at (-1.3,2.2) {$^{\scriptstyle I_1}$};
\node[] at (-1.3,2.05) {$\overbrace{\ \ \ \ \ \ \ \ \ \ \ \ \ \ }$};
\node[] at (-1,1.7) {$...$};
\node[] at (1.3,2.2) {$^{\scriptstyle I_2}$};
\node[] at (1.3,2.05) {$\overbrace{\ \ \ \ \ \ \ \ \ \ \ \ \ \ }$};
\node[] at (1,1.7) {$...$};
\node[] at (-1.3,0.15) {$_{\scriptstyle I_3}$};
\node[] at (-1.3,0.35) {$\underbrace{\ \ \ \ \ \ \ \ \ \ \ \ \ \ }$};
\node[] at (-1.2,0.5) {$...$};
\node[] at (1.3,0.15) {$_{\scriptstyle I_4}$};
\node[] at (1.3,0.35) {$\underbrace{\ \ \ \ \ \ \ \ \ \ \ \ \ \ }$};
\node[] at (1.2,0.5) {$...$};
%
\node[] at (-1.3,-0.20) {$^{\scriptstyle J_1}$};
\node[] at (-1.3,-0.35) {$\overbrace{\ \ \ \ \ \ \ \ \ \ \ \ \ \ }$};
\node[] at (-1.2,-0.5) {$...$};
\node[] at (1.3,-0.20) {$^{\scriptstyle J_2}$};
\node[] at (1.3,-0.35) {$\overbrace{\ \ \ \ \ \ \ \ \ \ \ \ \ \ }$};
\node[] at (1.2,-0.5) {$...$};
\node[] at (-1.3,-2.25) {$_{\scriptstyle J_3}$};
\node[] at (-1.3,-2.05) {$\underbrace{\ \ \ \ \ \ \ \ \ \ \ \ \ \ }$};
\node[] at (-1,-1.7) {$...$};
\node[] at (1.3,-2.25) {$_{\scriptstyle J_4}$};
\node[] at (1.3,-2.05) {$\underbrace{\ \ \ \ \ \ \ \ \ \ \ \ \ \ }$};
\node[] at (1,-1.7) {$...$};
\node[] (uu_1) at (-2,2.0) {};
\node[] (uu_2) at (-1,2.0) {};
\node[] (uu_3) at (2,2.0) {};
\node[] (uu_4) at (1,2.0) {};
\node[] (ud_1) at (-2,0.4) {};
\node[] (ud_2) at (-1,0.4) {};
\node[] (ud_3) at (2,0.4) {};
\node[] (ud_4) at (1,0.4) {};
\node[] (dd_1) at (-2,-2.0) {};
\node[] (dd_2) at (-1,-2.0) {};
\node[] (dd_3) at (2,-2.0) {};
\node[] (dd_4) at (1,-2.0) {};
\node[] (du_1) at (-2,-0.4) {};
\node[] (du_2) at (-1,-0.4) {};
\node[] (du_3) at (2,-0.4) {};
\node[] (du_4) at (1,-0.4) {};
 \node[int] (a) at (0,-1.0) {};
 \node[int] (b) at (0,1.0) {};
\node[] (d) at (0,-2.4) {};
\draw (uu_1) edge[latex-, leftblue] (b);
\draw (uu_2) edge[latex-, leftblue] (b);
\draw (uu_3) edge[latex-, rightblue] (b);
\draw (uu_4) edge[latex-, rightblue] (b);
\draw (b) edge[latex-, leftblue] (ud_1);
\draw (b) edge[latex-, leftblue] (ud_2);
\draw (b) edge[latex-, rightblue] (ud_3);
\draw (b) edge[latex-, rightblue] (ud_4);
\draw (du_1) edge[latex-, leftblue] (a);
\draw (du_2) edge[latex-, leftblue] (a);
\draw (du_3) edge[latex-, rightblue] (a);
\draw (du_4) edge[latex-, rightblue] (a);
\draw (a) edge[latex-, leftblue] (dd_1);
\draw (a) edge[latex-, leftblue] (dd_2);
\draw (a) edge[latex-, rightblue] (dd_3);
\draw (a) edge[latex-, rightblue] (dd_4);
 \draw (b) edge[latex-, leftblue, bend left=50] (a);
 \draw (b) edge[latex-, leftblue, bend right=50] (a);
  \draw (b) edge[latex-, leftblue, bend left=30] (a);
 \draw (b) edge[latex-, leftblue, bend right=30] (a);
\end{tikzpicture}}\Ea
\lon
\Ba{c}\resizebox{24mm}{!}
{\begin{tikzpicture}[baseline=-1ex]
\node[] at (-0.6,1.0) {$^{\scriptstyle I_1\sqcup J_1}$};
\node[] at (-0.7,0.8) {$\overbrace{\  \ \ \ \ \ \ \ \ \ \ }$};
\node[] at (-0.45,0.5) {$...$};
\node[] at (0.6,1.0) {$^{\scriptstyle I_2\sqcup J_2}$};
\node[] at (0.7,0.8) {$\overbrace{\  \ \ \ \ \ \ \ \ \ \ }$};
\node[] at (0.45,0.5) {$...$};
\node[] at (-0.6,-1.0) {$_{\scriptstyle I_3\sqcup J_3}$};
\node[] at (-0.7,-0.8) {$\underbrace{\  \ \ \ \ \ \ \ \ \ \ }$};
\node[] at (-0.45,-0.5) {$...$};
\node[] at (0.6,-1.0) {$_{\scriptstyle I_4\sqcup J_4}$};
\node[] at (0.7,-0.8) {$\underbrace{\  \ \ \ \ \ \ \ \ \ \ }$};
\node[] at (0.45,-0.5) {$...$};
\node[int] (0) at (0,0) {};
\node[] (1) at (-1.2,0.8) {};
\node[] (2) at (-0.3,0.8) {};
\node[] (3) at (0.3,0.8) {};
\node[] (4) at (1.2,0.8) {};
\node[int] (0) at (0,0) {};
\node[] (5) at (-1.2,-0.8) {};
\node[] (6) at (-0.3,-0.8) {};
\node[] (7) at (0.3,-0.8) {};
\node[] (8) at (1.2,-0.8) {};
\draw (1) edge[latex-,leftblue] (0);
\draw (2) edge[latex-,leftblue] (0);
\draw (3) edge[latex-,rightblue] (0);
\draw (4) edge[latex-,rightblue] (0);
\draw (0) edge[latex-,leftblue] (5);
\draw (0) edge[latex-,leftblue] (6);
\draw (0) edge[latex-,rightblue] (7);
\draw (0) edge[latex-,rightblue] (8);
\end{tikzpicture}}\Ea\ \ , \ \ \
\Ba{c}\resizebox{34mm}{!}  {
\begin{tikzpicture}[baseline=-.65ex]
\node[] at (-1.3,2.2) {$^{\scriptstyle I_1}$};
\node[] at (-1.3,2.05) {$\overbrace{\ \ \ \ \ \ \ \ \ \ \ \ \ \ }$};
\node[] at (-1,1.7) {$...$};
\node[] at (1.3,2.2) {$^{\scriptstyle I_2}$};
\node[] at (1.3,2.05) {$\overbrace{\ \ \ \ \ \ \ \ \ \ \ \ \ \ }$};
\node[] at (1,1.7) {$...$};
\node[] at (-1.3,0.15) {$_{\scriptstyle I_3}$};
\node[] at (-1.3,0.35) {$\underbrace{\ \ \ \ \ \ \ \ \ \ \ \ \ \ }$};
\node[] at (-1.2,0.5) {$...$};
\node[] at (1.3,0.15) {$_{\scriptstyle I_4}$};
\node[] at (1.3,0.35) {$\underbrace{\ \ \ \ \ \ \ \ \ \ \ \ \ \ }$};
\node[] at (1.2,0.5) {$...$};
%
\node[] at (-1.3,-0.20) {$^{\scriptstyle J_1}$};
\node[] at (-1.3,-0.35) {$\overbrace{\ \ \ \ \ \ \ \ \ \ \ \ \ \ }$};
\node[] at (-1.2,-0.5) {$...$};
\node[] at (1.3,-0.20) {$^{\scriptstyle J_2}$};
\node[] at (1.3,-0.35) {$\overbrace{\ \ \ \ \ \ \ \ \ \ \ \ \ \ }$};
\node[] at (1.2,-0.5) {$...$};
\node[] at (-1.3,-2.25) {$_{\scriptstyle J_3}$};
\node[] at (-1.3,-2.05) {$\underbrace{\ \ \ \ \ \ \ \ \ \ \ \ \ \ }$};
\node[] at (-1,-1.7) {$...$};
\node[] at (1.3,-2.25) {$_{\scriptstyle J_4}$};
\node[] at (1.3,-2.05) {$\underbrace{\ \ \ \ \ \ \ \ \ \ \ \ \ \ }$};
\node[] at (1,-1.7) {$...$};
\node[] (uu_1) at (-2,2.0) {};
\node[] (uu_2) at (-1,2.0) {};
\node[] (uu_3) at (2,2.0) {};
\node[] (uu_4) at (1,2.0) {};
\node[] (ud_1) at (-2,0.4) {};
\node[] (ud_2) at (-1,0.4) {};
\node[] (ud_3) at (2,0.4) {};
\node[] (ud_4) at (1,0.4) {};
\node[] (dd_1) at (-2,-2.0) {};
\node[] (dd_2) at (-1,-2.0) {};
\node[] (dd_3) at (2,-2.0) {};
\node[] (dd_4) at (1,-2.0) {};
\node[] (du_1) at (-2,-0.4) {};
\node[] (du_2) at (-1,-0.4) {};
\node[] (du_3) at (2,-0.4) {};
\node[] (du_4) at (1,-0.4) {};
 \node[int] (a) at (0,-1.0) {};
 \node[int] (b) at (0,1.0) {};
\node[] (d) at (0,-2.4) {};
\draw (uu_1) edge[latex-, leftblue] (b);
\draw (uu_2) edge[latex-, leftblue] (b);
\draw (uu_3) edge[latex-, rightblue] (b);
\draw (uu_4) edge[latex-, rightblue] (b);
\draw (b) edge[latex-, leftblue] (ud_1);
\draw (b) edge[latex-, leftblue] (ud_2);
\draw (b) edge[latex-, rightblue] (ud_3);
\draw (b) edge[latex-, rightblue] (ud_4);
\draw (du_1) edge[latex-, leftblue] (a);
\draw (du_2) edge[latex-, leftblue] (a);
\draw (du_3) edge[latex-, rightblue] (a);
\draw (du_4) edge[latex-, rightblue] (a);
\draw (a) edge[latex-, leftblue] (dd_1);
\draw (a) edge[latex-, leftblue] (dd_2);
\draw (a) edge[latex-, rightblue] (dd_3);
\draw (a) edge[latex-, rightblue] (dd_4);
 \draw (b) edge[latex-, rightblue, bend left=50] (a);
 \draw (b) edge[latex-, rightblue, bend right=50] (a);
  \draw (b) edge[latex-, rightblue, bend left=30] (a);
 \draw (b) edge[latex-, rightblue, bend right=30] (a);
\end{tikzpicture}}\Ea
\lon
\Ba{c}\resizebox{24mm}{!}
{\begin{tikzpicture}[baseline=-1ex]
\node[] at (-0.6,1.0) {$^{\scriptstyle I_1\sqcup J_1}$};
\node[] at (-0.7,0.8) {$\overbrace{\  \ \ \ \ \ \ \ \ \ \ }$};
\node[] at (-0.45,0.5) {$...$};
\node[] at (0.6,1.0) {$^{\scriptstyle I_2\sqcup J_2}$};
\node[] at (0.7,0.8) {$\overbrace{\  \ \ \ \ \ \ \ \ \ \ }$};
\node[] at (0.45,0.5) {$...$};
\node[] at (-0.6,-1.0) {$_{\scriptstyle I_3\sqcup J_3}$};
\node[] at (-0.7,-0.8) {$\underbrace{\  \ \ \ \ \ \ \ \ \ \ }$};
\node[] at (-0.45,-0.5) {$...$};
\node[] at (0.6,-1.0) {$_{\scriptstyle I_4\sqcup J_4}$};
\node[] at (0.7,-0.8) {$\underbrace{\  \ \ \ \ \ \ \ \ \ \ }$};
\node[] at (0.45,-0.5) {$...$};
\node[int] (0) at (0,0) {};
\node[] (1) at (-1.2,0.8) {};
\node[] (2) at (-0.3,0.8) {};
\node[] (3) at (0.3,0.8) {};
\node[] (4) at (1.2,0.8) {};
\node[int] (0) at (0,0) {};
\node[] (5) at (-1.2,-0.8) {};
\node[] (6) at (-0.3,-0.8) {};
\node[] (7) at (0.3,-0.8) {};
\node[] (8) at (1.2,-0.8) {};
\draw (1) edge[latex-,leftblue] (0);
\draw (2) edge[latex-,leftblue] (0);
\draw (3) edge[latex-,rightblue] (0);
\draw (4) edge[latex-,rightblue] (0);
\draw (0) edge[latex-,leftblue] (5);
\draw (0) edge[latex-,leftblue] (6);
\draw (0) edge[latex-,rightblue] (7);
\draw (0) edge[latex-,rightblue] (8);
\end{tikzpicture}}\Ea
$$
along the subgraphs with {\em no wheels with respect to any of the directions}; let us call a prop generated by such multi-oriented corollas and equipped with such compositions laws (see \S 2 for the full list of axioms) a $(k+1)$-{\em oriented}\, one.

 \sip

 Which structure on   graded vector spaces $V$ can be used to separate (in the sense of representations) $0$-oriented $(k+1)$-directed props from $(k+1)$-oriented ones (or, more generally, $(l+1)$ oriented $(k+1)$-directed with $l\geq 0$)? Note that the compositions prohibited in the $(k+1)$-oriented prop are still nicely oriented with respect to the basic direction, so  the answer can not be {\em the dimension of $V$}\, only.

\sip

To make sense of these new restrictions (which have no analogue in the theory of coloured props) we suggest to define a {\em representation}\, of a $(k+1)$-oriented prop
in a graded vector space $V$ as follows. Assume $V$ contains a collection of linear subspaces (satisfying certain restrictions in the infinite-dimensional case, see \S 4)
\Beq\label{1: k subspaces in V}
W_1^+\subset V,\  W_2^+\subset V, \ \ldots \ , \  W_k^+\subset V
\Eeq
with chosen complements
$$
V/W_1^+\simeq W_1^-\subset V,\ \ \ V/W_2^+\simeq W_2^-\subset V, \ \ \ldots \ \ ,
V/W_k^+ \simeq W_k^-\subset V.
$$
Then to a $(k+1)$-directed outgoing leg we associate (roughly) an
intersection\footnote{Strictly speaking, this is true only in finite dimensions. In infinite dimensions the subspaces $W_i^+$ are defined as {\em direct limits}\, of systems of finite-dimensional spaces  while their complements $W^-_i$ always come as {\em projective limits}, so their intersection makes sense only
at the level of finite-dimensional systems  first (it is here where the interpretation of $W^+$ and $W^-$ as subspaces of {\em one and the same}\, vector space plays its role), and then taking either  the direct or projective limit in accordance with the rule explained in \S 4.}

$$
\Ba{c}
\begin{tikzpicture}[baseline=-2ex]
\draw (0,0) edge[rightblue] node[above] {\ \ \ $\scriptstyle 1$}(0.9,0);
\draw (0.9,0) edge[leftred] node[above] {$\scriptstyle 2$\ }(1.5,0);
\draw (1.5,0) edge[] node[above] {$...$\ \ \ }(1.9,0);
\draw (1.9,0) edge[rightgreen] node[above] {$\scriptstyle k$\ }(2.3,0);
\draw (2.2,0) edge[-latex](2.7,0);
\end{tikzpicture}
\Ea\ \ \ \ \ \Leftrightarrow  \ \ \  W^+_1\cap W_2^- \cap \ldots \cap W^+_k
$$
obtained by the intersection of ``branes" according to the rule:
\Bi
\item to the basic direction we always associate the ``full" vector space $V$;
\item to the $i$-th direction we associate the vector subspace $W^+_i$ if that direction is in agreement with the basic one, or its complement $W_i^-$ if it is not.
\Ei
Then any multi-oriented corolla gets interpreted as a collection of $k$ linear maps, one map for each coloured orientation. For example, a 3-directed corolla
$
\Ba{c}\resizebox{17mm}{!}  {
\begin{tikzpicture}[baseline=-1ex]
\node[int] (a) at (0,0) {};
\node[] (u1) at (-1,1) {};
\node[] (u2) at (-0.65,0.9) {};
\node[] (u3) at (0.65,0.9) {};
\node[] (u4) at (1,1) {};
\node[int] (a) at (0,0) {};
\node[] (d2) at (-0.65,-0.9) {};
\node[] (d3) at (0.65,-0.9) {};
\draw (a) edge[-latex, rightblue, rightred] (u2);
\draw (a) edge[-latex, rightblue, leftred] (u3);
\draw (a) edge[latex-, leftblue, leftred] (d2);
\draw (a) edge[latex-, rightblue, leftred] (d3);
\end{tikzpicture}
}\Ea
$
 gets represented in a graded vector space $V$ equipped with two branes $ W_{\color{blue} 1}^\pm,  W_{\color{red} 2}^\pm \subset V$ as two linear maps, one corresponding to one blue input and three blue outputs of the corolla,
 $$
 W_{\color{blue} 1}^+\cap  W_{\color{red} 2}^+ \ \lon \  ( W_{\color{blue} 1}^+\cap  W_{\color{red} 2}^+)\ot ( W_{\color{blue} 1}^+\cap  W_{\color{red} 2}^-)\ot ( W_{\color{blue} 1}^-\cap  W_{\color{red} 2}^+)^*,
 $$
and another to three red inputs and one red output of the same corolla,
 $$
 (W_{\color{blue} 1}^+\cap  W_{\color{red} 2}^+)\ot ( W_{\color{blue} 1}^-\cap  W_{\color{red} 2}^+)\ot ( W_{\color{blue} 1}^+\cap  W_{\color{red} 2}^-)^* \lon     W_{\color{blue} 1}^+\cap  W_{\color{red} 2}^+.
 $$
 In finite dimensions both maps are just re-incarnations of one and the same linear map
 $$
 (W_{\color{blue} 1}^+\cap  W_{\color{red} 2}^+)\ot ( W_{\color{blue} 1}^-\cap  W_{\color{red} 2}^+)  \lon    ( W_{\color{blue} 1}^+\cap  W_{\color{red} 2}^+) \ot ( W_{\color{blue} 1}^+\cap  W_{\color{red} 2}^-),
 $$
 which is far from being the case in {\em infinite}\, dimensions.
 Most importantly,  this approach to the representation theory of
 multi-directed props explains nicely why compositions along graphs
 with wheels in at least one extra orientation must be prohibited
 (we show explicit examples of the associated divergences in \S 4). This approach also explains the formula (\ref{1: formula for number of colours}) for the associated number of ``colours" on legs.

 \subsection{Structure of the paper} In \S 2 we give a detailed (combinatorial type)  definition of multi-oriented props. In \S 3 we
 consider concrete examples. In particular, we introduce and study multi-oriented analogues,
 $\cA ss^{(k+1)}$ and $\caL ie^{(k+1)}$,  of the classical operads of associative algebras and Lie algebras, and explicitly describe their minimal resolutions $\cA ss^{(k+1)}_\infty$ and $\caL ie^{(k+1)}_\infty$; we also construct surprising ``forgetting the basic direction" maps from $\cA ss^{(2)}$ to the dioperad of infinitesimal
 bialgebras, and from $\caL ie^{(2)}$ to the dioperad of Lie bialgebras
 (proving that among representations of multi-oriented props we can recover sometimes classical structures); we also introduce a family of $(k+1)$-oriented props
 of homotopy Lie bialgebras $\HoLBcd^{(c+d-1)\text{-or}}$ on which the Grothendieck-Teichm\"uller group acts faithfully (see \S 5). In \S 4, the main section of this paper, we define the notion of a representation of
 a  multi-oriented prop in the category of graded vector spaces {\em with branes}, and show, as an illustration,  that Manin triples give us a class of symplectic Lagrangian representations of $\caL ie^{(2)}$.

\mip

{\bf Acknowledgement}. It is a pleasure to thank Assar Andersson, Anton Khoroshkin, Thomas Willwacher  and  Marko \v
Zivkovi\' c  for valuable discussions. I am also grateful to the referee for several very useful comments and suggestions.

\bip

{\Large
\section{\bf Multi-oriented props}
}

\mip

\subsection{$\bS$-bimodules reinterpreted} For a finite set $I$ let $S_I^{(1)}$ be
 the set of  maps
$$
\fs: I \rar \{out, in\}
$$
from $I$ to the set consisting of two elements
called $out$ and $in$.
A finite set $I$ together with a fixed function $\fs\in S_I^{(1)}$ is called {\em 1-oriented}. The collection of 1-oriented sets forms a groupoid $\cS^{(1)}$ with isomorphisms
$$
(I,\fs) \lon (I', \fs')
$$
being bijections $\sigma: I\rar I'$  of finite sets
such that $\fs'= \fs\circ \sigma^{-1}$. The latter condition says that the groupoid $\cS^{(1)}$ can be identified with the groupoid of cartesian products,
$\{I_{in}:=\fs^{-1}(in)\ \times \ I_{out}:=\fs^{-1}(out)\}$, of finite sets.

\sip

Let $\cC$ be a symmetric monoidal category.
 A functor
$$
\Ba{rccc}
\cP: & \cS^{(1)} & \lon & \cC\\
     & (I,\fs) & \lon &  \cP(I,\fs)
\Ea
$$
is called an $\cS^{(1)}$-{\em module}. An element  $a\in \cP(I,\fs)$
can be represented pictorially as a corolla with $\# I$ legs labelled by elements of
$I$ and oriented via the rule: if $\fs(i)=out$ (resp.,  $\fs(i)=in$) we orient the $i$-labelled leg by putting the direction $``>"$ {\em from}\, (resp., {\em towards}) the vertex; the vertex itself is decorated with $a$. For example, an element $a\in \cP([6],\fs)$ can have a pictorial representation of the form
$$
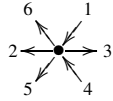
\xy
(6.6,0)*{_3},
(4,5)*{^1},
(-4,5)*{^6}="3",
(-6,0)*{_2}="4",
(4,-5)*{_4}="5",
(-4,-5)*{_5}="6",
 (0,0)*{\bullet}="a",
(5,0)*{}="1",
(3,4)*{}="2",
(-3,4)*{}="3",
(-5,0)*{}="4",
(3,-4)*{}="5",
(-3,-4)*{}="6",
\ar @{->} "a";"1" <0pt>
\ar @{<-} "a";"2" <0pt>
\ar @{->} "a";"3" <0pt>
\ar @{->} "a";"4" <0pt>
\ar @{<-} "a";"5" <0pt>
\ar @{->} "a";"6" <0pt>
\endxy
$$
The category of finite sets has a skeleton whose objects are sets $[N]=\{1,2,\ldots,N\}$ for some $N\geq 0$ (with $[0]=\emptyset$). For $I=[N]$, we often abbreviate $\cP_\fs(N,\fs):=\cP_\fs([N],\fs)$. Note that the above corolla is not assumed to be planar  so that it can be equivalently represented in a more standard way,
$$
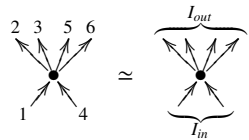
\xy
(5,6)*{^6},
(-5,6)*{^2},
(-2,6)*{^3},
(2,6)*{^5},
(4,-5)*{_4},
(-4,-5)*{_1},
 (0,0)*{\bullet}="a",
(5,5)*{}="1",
(2,5)*{}="2",
(-2,5)*{}="3",
(-5,5)*{}="4",
(3,-4)*{}="5",
(-3,-4)*{}="6",
\ar @{->} "a";"1" <0pt>
\ar @{->} "a";"2" <0pt>
\ar @{->} "a";"3" <0pt>
\ar @{->} "a";"4" <0pt>
\ar @{<-} "a";"5" <0pt>
\ar @{<-} "a";"6" <0pt>
\endxy\ \  \simeq \ \
\xy
(0,6.2)*{\overbrace{\ \ \ \ \ \ \ \ \ \ \ \ \ \ }},
(0,8.0)*{^{I_{out}}},
(0,-5.2)*{\underbrace{\ \ \ \ \ \ \ \ \ \ }},
(0,-7.2)*{_{I_{in}}},
 (0,0)*{\bullet}="a",
(5,5)*{}="1",
(2,5)*{}="2",
(-2,5)*{}="3",
(-5,5)*{}="4",
(3,-4)*{}="5",
(-3,-4)*{}="6",
\ar @{->} "a";"1" <0pt>
\ar @{->} "a";"2" <0pt>
\ar @{->} "a";"3" <0pt>
\ar @{->} "a";"4" <0pt>
\ar @{<-} "a";"5" <0pt>
\ar @{<-} "a";"6" <0pt>
\endxy
$$
which respects the flow of orientations going from the bottom to the top.

\sip

Any $\bS$-bimodule $E=\{E(m,n)\}_{m,n\geq 0}$ (with each $E(m,n)$ being an $\bS_m^{op}\times \bS_n$-module), gives rise to a $\cS^{(1)}$-module in the obvious way (and vice versa).

\subsection{Multi-oriented modules}\label{2new: subsec on multi-oriented modules}  For a natural number $k\geq 0$ let
$
\f_{\scriptstyle k^+}$ be the set of all maps  $\fm: [ k^+]:=\{{{0},{{1}},\ldots,{k}}\} \rar \{out,in\}$; the value $\fm({\tau})\subset \{out,in\}$ on ${\tau}\in [{k^+}]$ is called ${\tau}$-th {\em orientation}; the zero-th orientation $\fm({ 0})$ is called the {\em basic}\, one; the map $\fm$ is called itself a {\em multi-direction}.
The elements of $[k^+]$ are often called {\em colours}.
One can represent
pictorially a multi-direction $\fm\in \f_{k^+}$  as an ``outgoing leg" if $\fm({ 0})=out$,
$$
\fm \leftrightarrow \Ba{c}
\begin{tikzpicture}[baseline=-2ex]
\node[int] at (0,0) {};
\node[] at (2.75,0) {};
\draw (0,0) edge[leftgreen] node[above] { \ \ \ $\scriptstyle \fm( 1)$}(0.9,0);
\draw (0.9,0) edge[rightblue] node[above] {\ $\scriptstyle \fm( 2)$ }(1.5,0);
\draw (1.5,0) edge[] node[above] {$...$\ \ \ }(1.9,0);
\draw (1.9,0) edge[leftred] node[above] { $\scriptstyle \fm( k)$ }(2.3,0);
\draw (2.2,0) edge[-latex] node[above] {\ \ \ \ \ $\scriptstyle \fm(0)$ }(2.8,0);
\end{tikzpicture}
\Ea
$$
or ``ingoing leg" if $\fm({ 0})=in$
$$
\fm \leftrightarrow \Ba{c}
\begin{tikzpicture}[baseline=-2ex]
\node[] at (0,0) {};
\node[int] at (2.85,0) {};
\draw (0,0) edge[leftgreen] node[above] { \ \ \ $\scriptstyle \fm(1)$}(0.9,0);
\draw (0.9,0) edge[rightblue] node[above] {\ $\scriptstyle \fm(2)$ }(1.5,0);
\draw (1.5,0) edge[] node[above] {$...$\ \ \ }(1.9,0);
\draw (1.9,0) edge[leftred] node[above] { $\scriptstyle \fm(k)$ }(2.3,0);
\draw (2.2,0) edge[-latex] node[above] {\ \ \ \ \ $\scriptstyle \fm(0)$ }(2.8,0);
\end{tikzpicture}
\Ea
$$
using the obvious rule: for any $\tau\in [k]$ the value  $\fm(\tau)$ is represented by  the $\tau$-coloured symbol $``>_{\tau}"$ oriented
in the same direction as $\fm({ 0})=``>"$ if  $\fm(\tau)=\fm({  0})$ , or in the opposite direction, $``<_\tau"$,
if  $\fm(\tau)\neq \fm( { 0})$.

\sip

For a finite set $I$ consider the associated set $S_I^{(k+1)}$
 of  maps
$$
\Ba{rccc}
\fs: & I & \lon & \f_{k^+} \\
     & i & \lon & \fs_i:=\fs(i): [k^+]\rar \{out,in\}.
     \Ea
$$
 For $i\in I$ the value $\fs_i$
on $\tau\in [{ k^+}]$ is called $\tau$-{\em th orientation
 (or $\tau$-th direction) of the element $i$}. For any such a function $\fs$ there is associated the {\em opposite}\, function
$\fs^{opp}:  I \rar \f_{k^+}$ which is uniquely defined by the following condition: for each $i\in I$ and  each $\tau\in [{ k^+}]$ the value of $\fs^{opp}_i$ on $\tau$ is different from that of $\fs_i$ on $\tau$. Thus the set $S_I^{(k+1)}$ comes equipped with an involution. The restriction of the function $\fs_i$ to the subset
$[k]\subset [k^+]$ is denoted by $\bar{\fs}_i$; hence we can write
$$
\bar{\fs}_i \in \f_k:=\left\{[k]\rar \{out,in\}\right\}, \ \ \ \forall i\in I.
$$
This function takes care about {\em extra}\, (i.e.\ non-basic) orientation assigned to an element $i\in I$.
In some pictorial representations of multi-oriented sets $(I,\fs)$   we show explicitly only the basic orientation while compressing all the extra ones into this symbol
$\bar{\fs}_i$ (see below).
\sip

Note that for any given  multi-oriented set $(I,\fs)$ and any fixed colour
$\tau\in [k^+]$ there is an associated map
\Beq\label{2new: fs_tau map}
\Ba{rccc}
\check{\fs}_\tau: & I & \lon & \{out,in\}\\
            & i & \lon & \check{\fs}_\tau(i):= \fs_i(\tau)
\Ea
\Eeq
which we use in several constructions below.

\sip
Using the above pictorial  interpretation of elements of $\f_{k^+}$ as multi-oriented legs,
one can uniquely represent any element $\fs\in S_I^{(k+1)}$ as a {\em multi-directed (or multi-oriented) corolla}, that is, as a (non-planar) graph with one vertex $\bu$ and $\# I$ legs such that  each leg  is (i) distinguished by an element $i$ from $I$ and
(ii) decorated with the multi-direction $\fs_i\in \f_{k^+}$ as explained just above. For example, a corolla
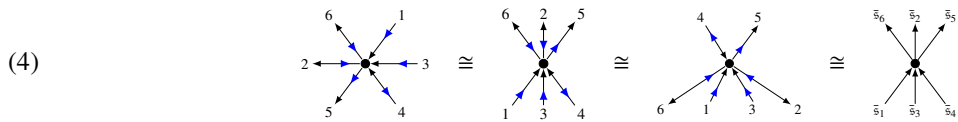
\Beq\label{4: maultidir corolla}
\Ba{c}\resizebox{21mm}{!}
{\begin{tikzpicture}[baseline=-1ex]
\node[] at (-0.6,0.8) {$\scriptstyle 6$};
\node[] (2) at (0.6,0.8) {$\scriptstyle 1$};
\node[] (3) at (-1,0) {$\scriptstyle 2$};
\node[] (4) at (1,0) {$\scriptstyle 3$};
\node[] (5) at (-0.6,-0.8) {$\scriptstyle 5$};
\node[] (6) at (0.6,-0.8) {$\scriptstyle 4$};
\node[int] (0) at (0,0) {};
\node[] (1) at (-0.6,0.8) {};
\node[] (2) at (0.6,0.8) {};
\node[] (3) at (-1,0) {};
\node[] (4) at (1,0) {};
\node[] (5) at (-0.6,-0.8) {};
\node[] (6) at (0.6,-0.8) {};
\draw (1) edge[latex-,rightblue] (0);
\draw (0) edge[latex-,leftblue] (2);
\draw (3) edge[latex-,rightblue] (0);
\draw (5) edge[latex-,leftblue] (0);
\draw (0) edge[latex-,leftblue] (4);
\draw (0) edge[latex-,rightblue] (6);
\end{tikzpicture}}\Ea
\cong
\Ba{c}\resizebox{15mm}{!}
{\begin{tikzpicture}[baseline=-1ex]
\node[] at (-0.6,0.8) {$\scriptstyle 6$};
\node[] (2) at (0.6,0.8) {$\scriptstyle 5$};
\node[] (3) at (0,0.8) {$\scriptstyle 2$};
\node[] (4) at (0,-0.8) {$\scriptstyle 3$};
\node[] (5) at (-0.6,-0.8) {$\scriptstyle 1$};
\node[] (6) at (0.6,-0.8) {$\scriptstyle 4$};
\node[int] (0) at (0,0) {};
\node[] (1) at (-0.6,0.8) {};
\node[] (2) at (0.6,0.8) {};
\node[] (3) at (0,0.8) {};
\node[] (4) at (0,-0.8) {};
\node[] (5) at (-0.6,-0.8) {};
\node[] (6) at (0.6,-0.8) {};
\draw (1) edge[latex-,rightblue] (0);
\draw (2) edge[latex-,leftblue] (0);
\draw (3) edge[latex-,rightblue] (0);
\draw (0) edge[latex-,leftblue] (5);
\draw (0) edge[latex-,leftblue] (4);
\draw (0) edge[latex-,rightblue] (6);
\end{tikzpicture}}\Ea
\cong
\Ba{c}\resizebox{25mm}{!}
{\begin{tikzpicture}[baseline=-1ex]
\node[] at (-1.2,-0.8) {$\scriptstyle 6$};
\node[] (2) at (0.5,0.8) {$\scriptstyle 5$};
\node[] (3) at (-0.5,0.8) {$\scriptstyle 4$};
\node[] (4) at (0.4,-0.8) {$\scriptstyle 3$};
\node[] (5) at (-0.4,-0.8) {$\scriptstyle 1$};
\node[] (6) at (1.2,-0.8) {$\scriptstyle 2$};
\node[int] (0) at (0,0) {};
\node[] (1) at (-1.2,-0.8) {};
\node[] (2) at (0.6,0.8) {};
\node[] (3) at (-0.5,0.8) {};
\node[] (4) at (0.5,-0.8) {};
\node[] (5) at (-0.4,-0.8) {};
\node[] (6) at (1.2,-0.8) {};
\draw (1) edge[latex-,rightblue] (0);
\draw (2) edge[latex-,leftblue] (0);
\draw (0) edge[latex-,rightblue] (3);
\draw (0) edge[latex-,leftblue] (5);
\draw (0) edge[latex-,leftblue] (4);
\draw (6) edge[latex-,rightblue] (0);
\end{tikzpicture}}\Ea
\cong
\Ba{c}\resizebox{15mm}{!}
{\begin{tikzpicture}[baseline=-1ex]
\node[] at (-0.6,0.8) {$\scriptstyle \bar{\fs}_6$};
\node[] (2) at (0.6,0.8) {$\scriptstyle \bar{\fs}_5$};
\node[] (3) at (0,0.8) {$\scriptstyle \bar{\fs}_2$};
\node[] (4) at (0,-0.8) {$\scriptstyle \bar{\fs}_3$};
\node[] (5) at (-0.6,-0.8) {$\scriptstyle \bar{\fs}_1$};
\node[] (6) at (0.6,-0.8) {$\scriptstyle \bar{\fs}_4$};
\node[int] (0) at (0,0) {};
\node[] (1) at (-0.6,0.8) {};
\node[] (2) at (0.6,0.8) {};
\node[] (3) at (0,0.8) {};
\node[] (4) at (0,-0.8) {};
\node[] (5) at (-0.6,-0.8) {};
\node[] (6) at (0.6,-0.8) {};
\draw (1) edge[latex-] (0);
\draw (2) edge[latex-] (0);
\draw (3) edge[latex-] (0);
\draw (0) edge[latex-] (5);
\draw (0) edge[latex-] (4);
\draw (0) edge[latex-] (6);
\end{tikzpicture}}\Ea
\Eeq
represents non-ambiguously some element
$\fs\in \cS^{(2)}_{[6]}$. In the theory of ordinary props corollas are often depicted in such a way that the orientation flow runs from the bottom to the top. In the multi-directed case such a respecting flow representation (now non-unique --- one for each coloured direction) also makes sense in applications.

\sip

A finite set $I$ together with a fixed function $\fs\in S_I^{(k+1)}$ is called {\em $(k+1)$-oriented}. The collection of $(k+1)$-oriented sets forms a groupoid $\cS^{(k+1)}$ with isomorphisms
$$
(I,\fs) \lon (I', \fs')
$$
being isomorphisms $\sigma: I\rar I'$  of finite sets
such that $\fs'= \fs\circ \sigma^{-1}$. For example, the automorphism group of the object $([6],\fs)$ given by corolla (\ref{4: maultidir corolla}) is $\bS_2\times \bS_2$ as
we can permute only labels $(1,3)$ and independently $(5,6)$ using morphisms in the category $\cS^{(2)}$.

\sip

Let $\cC$ be a symmetric monoidal category.
 A functor
$$
\Ba{rccc}
\cP^{(k+1)}: & \cS^{(k+1)} & \lon & \cC\\
     & (I,\fs) & \lon &  \cP^{(k+1)}(I,\fs)
\Ea
$$
is called an $\cS^{(k+1)}$-{\em module}. Thus an element of  $\cP^{(k+1)}(I,\fs)$ is a of pair of the form
$$
\left(c=
\Ba{c}\resizebox{21mm}{!}
{\begin{tikzpicture}[baseline=-1ex]
\node[] at (-0.6,0.8) {$\scriptstyle 6$};
\node[] (2) at (0.6,0.8) {$\scriptstyle 1$};
\node[] (3) at (-1,0) {$\scriptstyle 2$};
\node[] (4) at (1,0) {$\scriptstyle 3$};
\node[] (5) at (-0.6,-0.8) {$\scriptstyle 5$};
\node[] (6) at (0.6,-0.8) {$\scriptstyle 4$};
\node[int] (0) at (0,0) {};
\node[] (1) at (-0.6,0.8) {};
\node[] (2) at (0.6,0.8) {};
\node[] (3) at (-1,0) {};
\node[] (4) at (1,0) {};
\node[] (5) at (-0.6,-0.8) {};
\node[] (6) at (0.6,-0.8) {};
\draw (1) edge[latex-,rightblue] (0);
\draw (0) edge[latex-,leftblue] (2);
\draw (3) edge[latex-,rightblue] (0);
\draw (5) edge[latex-,leftblue] (0);
\draw (0) edge[latex-,leftblue] (4);
\draw (0) edge[latex-,rightblue] (6);
\end{tikzpicture}}\Ea, \ \ \cP^{(k+1)}(c)\in Objects(\cC)\right)
$$
Note that $\cP^{(k+1)}(c)$ carries a representation of the group $Aut(c)$ (in this particular case, of $\bS_2\times \bS_2$). We shall work in this paper in category of topological vector spaces so that it make sense to talk about elements $v$ of $\cP^{(k+1)}(c)$. The pairs $(c, v)$ are called {\em decorated corollas}\, and are often represented pictorially by the corolla $c$ with its vertex decorated (often tacitly) by the vector $v$. Such decorated corollas span  $\cP^{(k+1)}(I,\fs)$.

\sip

When $k$ is clear, we often abbreviate $\cP=\cP^{(k+1)}$. The case $k=0$ corresponds to the ordinary $\bS$-bimodule.

\sip

Given $\cS^{(k+1)}$-modules $\cP$ and $\cP'$.
A natural transformation of functors
$$
f: \cP \lon \cP' \ \ \ \Leftrightarrow\ \ \ \ \left\{f_{I,\fs}: \cP(I,\fs) \rar \cP'(I,\fs)\mid \ f_{I',\fs'}\circ \cP(\al)=\cP'(\al)\circ f_{I,\fs}\ \ \forall \al: (I,\fs)\rar (I,\fs')\right\}
$$
is called a {\em morphism of $\cS^{(k+1)}$-modules}.

\subsection{Directed and multidirected graphs}\label{2new: subsec on multid graphs}
A {\em graph with legs}\, $\Ga$ is, by definition, a finite set $H(\Ga)$ (whose elements
are called {\em half-edges}) equipped with
\Bi
\item[(a)] a partition  into the disjoint union of subsets,
$
H(\Ga)=\coprod_{v\in V(\Ga)} H(v),
$
parameterized by a set $V(\Ga)$  whose elements are called {\em vertices of $\Ga$}; the subset $H(v)\subset H(\Ga)$ is called the {\em set of half-edges attached to the vertex $v$} and
its cardinality $\# H(v)$ is called the {\em valency}\, of $v$,
\item[(b)] an involution $\ii : H(\Ga)\rar H(\Ga)$ whose fixed points are called {\em legs}\, of $\Ga$ and whose orbits $(h,\ii(h))\subset H(\Ga)$ of cardinality two are called  {\em internal edges}\, or simply {\em edges}; the set of edges is denoted by  $E(\Ga)$
    and the set of legs is denoted by $L(\Ga)$.
\Ei
 Any graph with legs $\Ga$  can be identified with its {\em geometric realization}\, which is a topological space constructed as follows: (i) for each vertex $v$ consider the disjoint union  of ${\# H(v)}$ copies of the unit interval $[0,1]\subset \R$
and identify all the copies of the end point $0$, the result is a topological space (equipped with the quotient topology) which is
called the {\em corolla of $v$}; (ii) consider the union of all stars and, for each internal edge $(h,\ii(h))$, identify the end-points $1$ of the intervals $[0,1]$ labelled by $h$ and $\ii(h)$.

\sip

A graph with legs is called {\em directed}\, if every leg or internal edge  in its geometric realization  comes equipped with a fixed (one of the two possible) orientations. It is convenient
to identify this orientation with a flow on a geometric edge making the latter into an arrow. Here is an example of a directed graph
$$
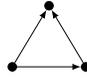
\Ba{c}\resizebox{11mm}{!}{\begin{tikzpicture}[baseline=-1ex]
\node[int] (b) at (0,1) {};
\node[int] (c) at (-0.6,0) {};
\node[int] (d) at (0.6,0) {};
\draw (b) edge[latex-] (c);
\draw (b) edge[latex-] (d);
\draw (d) edge[latex-] (c);
\end{tikzpicture}}\Ea
$$
with three vertices and three edges.

\sip

 If the involution $\ii$ has no fixed points, i.e.\ $L(\Ga)=\emptyset$, the graph with legs is called simply a {\em graph}. We shall use directed graphs with legs below when defining multi-oriented props while in this subsection we continue on only with graphs without legs
 (i.e.\ simply with graphs).

\sip

By a {\em multi-directed}, more precisely, $(k+1)$-{\em directed}\, graph
we understand a pair $\displaystyle \left(\Ga, s\in \cS^{(k+1)}_{E(\Ga)}\right)$ consisting of a directed graph and a function $\fs: E(\Ga)\rar \f_{k^+}$ such that
for each $e\in E(\Ga)$ the value of the associated function $\fs_e:[k^+]\rar \{out,in\}$ takes value ``out" (or, equivalently, ``in")  at the ``zero-th colour" $0$   and is identified pictorially with the original (basic) direction of $e$,
$$
e=\Ba{c}\resizebox{11mm}{!}{\begin{tikzpicture}[baseline=-1ex]
\node[int] (c) at (-0.6,-0.05) {};
\node[int] (d) at (0.6,-0.05) {};
\
\draw (d) edge[latex-] (c);
\end{tikzpicture}}\Ea
\rightsquigarrow
\Ba{c}
\begin{tikzpicture}[baseline=-2ex]
\node[int] at (0,0) {};
\node[int] at (3.15,0) {};
\draw (0,0) edge[leftgreen] node[above] { \ \ \ $\scriptstyle \fs_e(1)$}(0.9,0);
\draw (0.9,0) edge[rightblue] node[above] {\ $\scriptstyle \fs_e(2)$ }(1.5,0);
\draw (1.5,0) edge[] node[above] {$...$\ \ \ }(1.9,0);
\draw (1.9,0) edge[leftred] node[above] { $\scriptstyle \fs_e(k)$ }(2.5,0);
\draw (2.5,0) edge[-latex] node[above] {\ \ \ \ \ \ \ \  $\scriptstyle \fs_e(0)=out$ }(3.10,0);
\end{tikzpicture}
\Ea
$$
 Thus the data $\displaystyle \left(\Ga, s\in \cS^{(k+1)}_{E(\Ga)}\right)$  can be represented pictorially as a graph whose edges are equipped with
$(k+1)$-directions.
 Here is an example of a $3$-directed graph.
$$
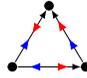
\Ba{c}\resizebox{11mm}{!}{\begin{tikzpicture}[baseline=-1ex]
\node[int] (b) at (0,1) {};
\node[int] (c) at (-0.6,0) {};
\node[int] (d) at (0.6,0) {};
\draw (b) edge[latex-,leftblue, rightred] (c);
\draw (b) edge[latex-,leftblue, leftred] (d);
\draw (d) edge[latex-,rightblue, leftred] (c);
\end{tikzpicture}}\Ea
$$
Let $G^{0\uparrow k+1}$ denote set of all such $(k+1)$-directed graphs.
The permutation group $\bS_{k+1}$ acts on this set via its canonical action on the set of colours $[k^+]$.

\sip

Let $A$ be a subset of $[k^+]$. A $(k+1)$-directed graph $\Ga \in G^{0\uparrow k+1}$ is called $A$-{\em oriented}\, if $\Ga$ contains no closed directed paths of edges (``wheels" or ``loops") in every colour
$c\in A$. The subset of $A$-oriented graphs is denoted by $G^{A\uparrow k+1}$. If $A$ is non-empty, then applying a suitable element of the automorphism group $\bS_{k+1}$ we can (and will) assume without loss of generality that $A=\{0,1,2, \ldots, l\}$ for some $l\geq 0$ and re-denote
 $G^{l+1\uparrow k+1}:= G^{A\uparrow k+1}$. If $l=k$, we further abbreviate $G^{(k+1)\text{-or}}:=G^{k+1\uparrow k+1}$ and call its elements {\em multi-oriented}\, graphs.

\subsection{From multi-directed graphs to endofunctors on $\cS^{(k+1)}$-modules} Fix an integer $k\geq 0$ and an integer $l$ in the range $-1\leq l\leq k$.
For a finite set $I$ define $G^{l+1\uparrow k+1}(I)$ to be the set of $(k+1)$-directed $(l+1)$-oriented graphs $\Ga$ equipped with an injection  $\ii: I\rar V_1(\Ga)$, where $V_1(\Ga)\subset V(\Ga)$ is the subset of univalent vertices. The univalent vertices lying in the image $L(\Ga):=\ii(I)$ of this map are called $(k+1)$-directed {\em legs}\, of $\Ga$; each such leg is labelled therefore by an element $i$ of $I$ and is called
    an $i$-leg (in pictures we show it as a leg indeed with the 1-valent vertex {\em erased}\, and the index $i$ put on its place); vertices in $V_{int}(\Ga):=V(\Ga)\setminus L(\Ga)$ are called
    {\em internal}. Edges connecting internal vertices are called {\em internal}; there is a decomposition $E(\Ga)=L(\Ga)\sqcup E_{int}(\Ga)$. Here are examples of 2-directed graphs, one with 4 internal vertices and  4 legs, the other with two internal vertices and 3 legs,
$$
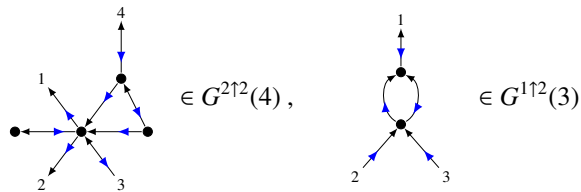
\Ba{c}\resizebox{21mm}{!}
{\begin{tikzpicture}[baseline=-1ex]
\node[] at (-0.6,0.8) {$\scriptstyle 1$};
\node[] (4) at (0.6,1.8) {$\scriptstyle 4$};
\node[] (5) at (-0.6,-0.8) {$\scriptstyle 2$};
\node[] (6) at (0.6,-0.8) {$\scriptstyle 3$};
\node[int] (0) at (0,0) {};
\node[] (1) at (-0.6,0.8) {};
\node[int] (2) at (0.6,0.8) {};
\node[int] (3) at (-1,0) {};
\node[int] (4) at (1,0) {};
\node[] (5) at (-0.6,-0.8) {};
\node[] (6) at (0.6,-0.8) {};
\node[] (7) at (0.6,1.8) {};
\draw (1) edge[latex-,leftblue] (0);
\draw (0) edge[latex-,leftblue] (2);
\draw (3) edge[latex-,rightblue] (0);
\draw (5) edge[latex-,leftblue] (0);
\draw (0) edge[latex-,leftblue] (4);
\draw (0) edge[latex-,rightblue] (6);
\draw (2) edge[latex-,rightblue] (4);
\draw (7) edge[latex-,rightblue] (2);
\end{tikzpicture}}\Ea\ \in G^{2\uparrow 2}(4)\ ,
 \ \ \ \
 \Ba{c}\resizebox{18mm}{!}  {
\begin{tikzpicture}[baseline=-.65ex]
\node[] (4) at (0,1.2) {$\scriptstyle 1$};
\node[] (5) at (-0.7,-1.2) {$\scriptstyle 2$};
\node[] (6) at (0.7,-1.2) {$\scriptstyle 3$};
\node[] (u) at (0,1.2) {};
 \node[int] (a) at (0,-0.4) {};
 \node[int] (b) at (0,0.4) {};
\node[] (d') at (-0.7,-1.2) {};
\node[] (d'') at (0.7,-1.2) {};
\draw (u) edge[latex-, rightblue] (b);
\draw (a) edge[latex-, leftblue] (d');
\draw (a) edge[latex-, leftblue] (d'');
 \draw (b) edge[latex-, rightblue, bend left=60] (a);
 \draw (b) edge[latex-, leftblue, bend right=60] (a);
\end{tikzpicture}}\Ea \in G^{1\uparrow 2}(3)
$$
Given an internal vertex $v\in V_{int}(\Ga)$, there is an associated set $H_v$
of edges attached to $v$ and an obvious function (``the multi-oriented corolla at $v$")
$$
\fs_v: H_v \lon \f_{k^+}
$$
There is also an induced function
$$
\fs: I=L(\Ga) \lon \f_{k^+}
$$
on the set of legs defined uniquely by the pictorial rule explained in
the first paragraph of \S {\ref{2new: subsec on multi-oriented modules}}.
 Let $G^{l+1\uparrow k+1}(I,\fs) \subset G^{l+1\uparrow k+1}(I)$
be the subset of multi-directed (partially oriented, in general) graphs which have one and the same orientation function
$\fs$ on the set of legs $I$.

\sip

For an $\cS^{(k+1)}$-module $\cE=\{\cE(I,\fs)\}$ in a symmetric monoidal category $\cC$ with countable coproducts and a graph $\Ga\in G^{l+1\uparrow k+1}(I,\fs)$
consider the unordered tensor product\footnote{The (unordered) tensor product
 $\bigotimes_{i\in I} X_i$ of vector spaces $X_i$ labelled by elements $i$ of a finite set $I$ of cardinality, say, $n$ is defined as the space of $\bS_n$-coinvariants
$
\left(\bigoplus_{\sigma: [n]\stackrel{\simeq}{\lon} I}  X_{\sigma(1)}\ot X_{\sigma(2)} \ot \ldots \ot X_{\sigma(n)}\right)_{\bS_n}.$} (cf.\ \cite{Ma,MSS})

$$
\Ga\langle\cE\rangle(I,\fs):= \left(
\bigotimes_{v\in V_{int}(G)} \cE(H_v, \fs_v)\right)_{Aut (\Ga)}
$$
where $Aut(\Ga)$ stands for
the automorphism group of the graph $\Ga$, and define an $\cS^{(k+1)}$-module in $\cC$
$$
\Ba{rccc}
\cF\! ree^{l+1\uparrow k+1}\langle \cE \rangle: & \cS^{(k+1)} & \lon & \cC\\
     & (I,\fs) & \lon &  \cF\! ree^{l+1\uparrow k+1}\langle \cE \rangle(\fs, I):= \bigoplus_{\Ga \in
     G^{l+1\uparrow k+1}(I,\fs)}  \Ga\langle\cE\rangle(I,\fs)
     \Ea
$$
As we shall see below, the $\cS^{(k+1)}$-module  $\cF\! ree^{l+1\uparrow k+1}\langle \cE \rangle$ gives us an example of a $(l+1)$-oriented $(k+1)$-directed   prop (called the {\em free}\, prop generated by the $\cS^{(k+1)}$-module $\cE$. For $l=k=0$ this is precisely the ordinary free prop generated by the $\cS^{(1)}$-module $\cE$.
For $l=-1,k=0$ this is the free {\em wheeled prop}\, generated by $\cE$ \cite{Me1,MMS}.  If $l=k$, i.e.\ if all directions are oriented, we abbreviate
$\cF\! ree^{k+1\uparrow k+1}\langle \cE \rangle =:\cF\! ree^{(k+1)\text{-or}}\langle \cE \rangle$.

\subsection{Multi-oriented prop(erad)s}
A (possibly disconnected) subgraph $\ga$ of a (connected or disconnected) graph $\Ga\in G^{l+1\uparrow l+1}(I,\fs)$ is called {\em complete}\, if every edge of $\Ga$ connecting a pair of (not necessarily distinct) vertices of $\ga$ belongs to $\ga$.
  Let $\Ga/\ga$ be the graph obtained from $\Ga$ by contracting all internal vertices and all internal edges of $\ga$ to a single new vertex; note that the legs of $\Ga/\ga$
 are the same as in $\Ga$ so that $\Ga/\ga$ comes equipped with the same orientation function $\fs: L(\Ga/\ga)\rar \f_{k^+}$.   A complete subgraph $\ga\subset \Ga$ is called admissible if $\Ga/\ga$ belongs to  $G^{l+1\uparrow k+1}(I, \fs)$, i.e.\ the contraction procedure does not create
wheels in the first $l+1$ coloured directions. Note that by its very definition $\ga$ belongs to $\in G^{l+1\uparrow k+1}(I',\fs')$, where  $I'$ is the subset of $E(\Ga)$ consisting of (non-loop) those edges which are attached to precisely   {\em one}\, vertex $v$
of $\ga$, and  the function $\fs': I'\rar \f_{k^+}$ is determined by the corresponding functions $\fs_v$ in the obvious way.

\sip

An $(l+1)$-{\em oriented $(k+1)$-directed  prop in a symmetric monoidal category (with countable colimits)}\, $\cC$ is, by definition, an $\cS^{(k+1)}$-module $\cP=\{\cP(I,\fs)\}$ in $\cC$ together with a natural transformation of functors
$$
\Ba{rccc}
\mu: & \cF\! ree^{l+1\uparrow k+1}\langle \cP \rangle & \lon & \cP\\
\mu_{\Ga}:  &   \Ga\langle \cP \rangle (I,\fs) & \lon \cP(I,\fs)
\Ea
$$
such that for any graph $\Ga\in G^{l+1\uparrow k+1}(I,\fs)$ and any admissible subgraph
  $\ga\subset \Ga$ one has
\Beq\label{5: graph-associativity}
\mu_\Ga=\mu_{\Ga/\ga}\circ \mu_\ga',
\Eeq
where
  $\mu_\ga': \Ga\langle \cP \rangle(I,\fs) \rar (\Ga/\ga)\langle \cP\rangle(I,\fs)$ stands for the obvious map
which equals $\mu_\ga$ on the (decorated) subgraph $\ga$ and which is identity on all other vertices of $\Ga$.
\sip

The most interesting case for us is $k=l$. The associated props are called
{\em multi-oriented} (more precisely, $(k+1)$-{\em oriented}).
Thus a multi-oriented prop is an $\cS^{(k+1)}$-{module} $\cP$
equipped, in particular, with
\Bi
\item[(i)]
a {\em horizontal}\, composition
$$
\boxtimes: \cP(I_1,\fs_1) \ot \cP(I_2,\fs_2) \lon \cP(I_1\sqcup I_2, \fs_{1}\sqcup \fs_2),
$$
\item[(ii)] a  {\em reduced  vertical}\, composition:
 for any two injections of the same finite set $f_1:K\rar I_1$ and $f_2:K\rar I_2$ such that the compositions
$$
 K\stackrel{\fs_1\circ f_1}\lon \f_{k^+}\ \ \ , \ \ \
 K\stackrel{\fs_2\circ f_2}\lon \f_{k^+}
$$
 satisfy the condition $\fs_1\circ f_1= (\fs_2\circ f_2)^{opp}$ there is a {\em vertical}\, composition
$$
\cP(I_1,\fs_1)\circ_K \cP(I_2,\fs_2) \lon \cP_{\fs_{12}}\left((I_1\setminus f_1(K))\sqcup (I_2\setminus f_2(K)), \fs_{12}\right)
$$
where
$$
\fs_{12}: (I_1\setminus f_1(K)) \sqcup (I_2\setminus f_2(K))\lon  \f_{k^+}
$$
 is defined by $\fs_{12}(i)=\fs_1(i)$ for $i\in I_1\setminus f_1(K)$ and $\fs_{12}(i)=\fs_2(i)$ for $i\in I_2\setminus f_2(K)$.
\Ei
These compositions are required to satisfy the ``associativity" axioms which essentially say
that upon iterating  such compositions the order in which we do it does not matter. In terms of decorated corollas these compositions correspond to
contraction maps (for $k=1$)
$$
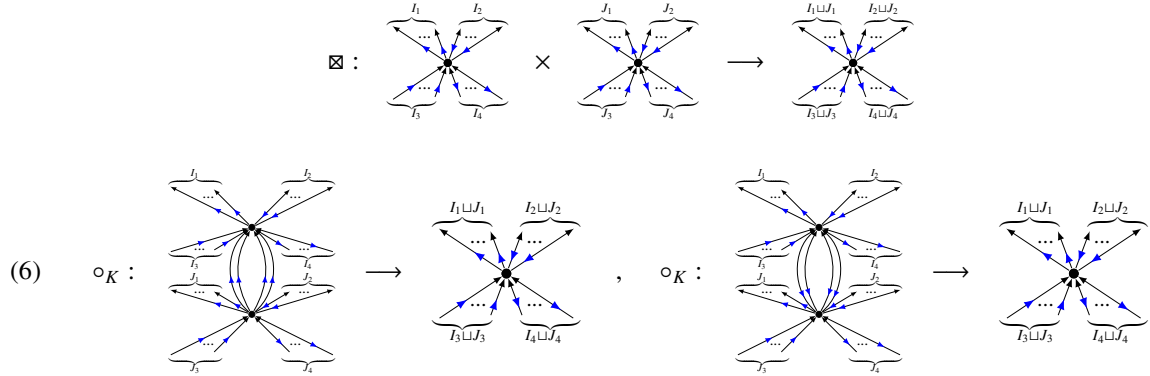
\boxtimes:
\Ba{c}\resizebox{24mm}{!}
{\begin{tikzpicture}[baseline=-1ex]
\node[] at (-0.6,1.0) {$^{\scriptstyle I_1}$};
\node[] at (-0.7,0.8) {$\overbrace{\  \ \ \ \ \ \ \ \ \ \ }$};
\node[] at (-0.45,0.5) {$...$};
\node[] at (0.6,1.0) {$^{\scriptstyle I_2}$};
\node[] at (0.7,0.8) {$\overbrace{\  \ \ \ \ \ \ \ \ \ \ }$};
\node[] at (0.45,0.5) {$...$};
\node[] at (-0.6,-1.0) {$_{\scriptstyle I_3}$};
\node[] at (-0.7,-0.8) {$\underbrace{\  \ \ \ \ \ \ \ \ \ \ }$};
\node[] at (-0.45,-0.5) {$...$};
\node[] at (0.6,-1.0) {$_{\scriptstyle I_4}$};
\node[] at (0.7,-0.8) {$\underbrace{\  \ \ \ \ \ \ \ \ \ \ }$};
\node[] at (0.45,-0.5) {$...$};
\node[int] (0) at (0,0) {};
\node[] (1) at (-1.2,0.8) {};
\node[] (2) at (-0.3,0.8) {};
\node[] (3) at (0.3,0.8) {};
\node[] (4) at (1.2,0.8) {};
\node[int] (0) at (0,0) {};
\node[] (5) at (-1.2,-0.8) {};
\node[] (6) at (-0.3,-0.8) {};
\node[] (7) at (0.3,-0.8) {};
\node[] (8) at (1.2,-0.8) {};
\draw (1) edge[latex-,leftblue] (0);
\draw (2) edge[latex-,leftblue] (0);
\draw (3) edge[latex-,rightblue] (0);
\draw (4) edge[latex-,rightblue] (0);
\draw (0) edge[latex-,leftblue] (5);
\draw (0) edge[latex-,leftblue] (6);
\draw (0) edge[latex-,rightblue] (7);
\draw (0) edge[latex-,rightblue] (8);
\end{tikzpicture}}\Ea
\times
\Ba{c}\resizebox{24mm}{!}
{\begin{tikzpicture}[baseline=-1ex]
\node[] at (-0.6,1.0) {$^{\scriptstyle J_1}$};
\node[] at (-0.7,0.8) {$\overbrace{\  \ \ \ \ \ \ \ \ \ \ }$};
\node[] at (-0.45,0.5) {$...$};
\node[] at (0.6,1.0) {$^{\scriptstyle J_2}$};
\node[] at (0.7,0.8) {$\overbrace{\  \ \ \ \ \ \ \ \ \ \ }$};
\node[] at (0.45,0.5) {$...$};
\node[] at (-0.6,-1.0) {$_{\scriptstyle J_3}$};
\node[] at (-0.7,-0.8) {$\underbrace{\  \ \ \ \ \ \ \ \ \ \ }$};
\node[] at (-0.45,-0.5) {$...$};
\node[] at (0.6,-1.0) {$_{\scriptstyle J_4}$};
\node[] at (0.7,-0.8) {$\underbrace{\  \ \ \ \ \ \ \ \ \ \ }$};
\node[] at (0.45,-0.5) {$...$};
\node[int] (0) at (0,0) {};
\node[] (1) at (-1.2,0.8) {};
\node[] (2) at (-0.3,0.8) {};
\node[] (3) at (0.3,0.8) {};
\node[] (4) at (1.2,0.8) {};
\node[int] (0) at (0,0) {};
\node[] (5) at (-1.2,-0.8) {};
\node[] (6) at (-0.3,-0.8) {};
\node[] (7) at (0.3,-0.8) {};
\node[] (8) at (1.2,-0.8) {};
\draw (1) edge[latex-,leftblue] (0);
\draw (2) edge[latex-,leftblue] (0);
\draw (3) edge[latex-,rightblue] (0);
\draw (4) edge[latex-,rightblue] (0);
\draw (0) edge[latex-,leftblue] (5);
\draw (0) edge[latex-,leftblue] (6);
\draw (0) edge[latex-,rightblue] (7);
\draw (0) edge[latex-,rightblue] (8);
\end{tikzpicture}}\Ea
\lon
\Ba{c}\resizebox{24mm}{!}
{\begin{tikzpicture}[baseline=-1ex]
\node[] at (-0.6,1.0) {$^{\scriptstyle I_1\sqcup J_1}$};
\node[] at (-0.7,0.8) {$\overbrace{\  \ \ \ \ \ \ \ \ \ \ }$};
\node[] at (-0.45,0.5) {$...$};
\node[] at (0.6,1.0) {$^{\scriptstyle I_2\sqcup J_2}$};
\node[] at (0.7,0.8) {$\overbrace{\  \ \ \ \ \ \ \ \ \ \ }$};
\node[] at (0.45,0.5) {$...$};
\node[] at (-0.6,-1.0) {$_{\scriptstyle I_3\sqcup J_3}$};
\node[] at (-0.7,-0.8) {$\underbrace{\  \ \ \ \ \ \ \ \ \ \ }$};
\node[] at (-0.45,-0.5) {$...$};
\node[] at (0.6,-1.0) {$_{\scriptstyle I_4\sqcup J_4}$};
\node[] at (0.7,-0.8) {$\underbrace{\  \ \ \ \ \ \ \ \ \ \ }$};
\node[] at (0.45,-0.5) {$...$};
\node[int] (0) at (0,0) {};
\node[] (1) at (-1.2,0.8) {};
\node[] (2) at (-0.3,0.8) {};
\node[] (3) at (0.3,0.8) {};
\node[] (4) at (1.2,0.8) {};
\node[int] (0) at (0,0) {};
\node[] (5) at (-1.2,-0.8) {};
\node[] (6) at (-0.3,-0.8) {};
\node[] (7) at (0.3,-0.8) {};
\node[] (8) at (1.2,-0.8) {};
\draw (1) edge[latex-,leftblue] (0);
\draw (2) edge[latex-,leftblue] (0);
\draw (3) edge[latex-,rightblue] (0);
\draw (4) edge[latex-,rightblue] (0);
\draw (0) edge[latex-,leftblue] (5);
\draw (0) edge[latex-,leftblue] (6);
\draw (0) edge[latex-,rightblue] (7);
\draw (0) edge[latex-,rightblue] (8);
\end{tikzpicture}}\Ea
$$

\Beq\label{2new: good compos in 2-or prop}
\circ_K:
\Ba{c}\resizebox{34mm}{!}  {
\begin{tikzpicture}[baseline=-.65ex]
\node[] at (-1.3,2.2) {$^{\scriptstyle I_1}$};
\node[] at (-1.3,2.05) {$\overbrace{\ \ \ \ \ \ \ \ \ \ \ \ \ \ }$};
\node[] at (-1,1.7) {$...$};
\node[] at (1.3,2.2) {$^{\scriptstyle I_2}$};
\node[] at (1.3,2.05) {$\overbrace{\ \ \ \ \ \ \ \ \ \ \ \ \ \ }$};
\node[] at (1,1.7) {$...$};
\node[] at (-1.3,0.15) {$_{\scriptstyle I_3}$};
\node[] at (-1.3,0.35) {$\underbrace{\ \ \ \ \ \ \ \ \ \ \ \ \ \ }$};
\node[] at (-1.2,0.5) {$...$};
\node[] at (1.3,0.15) {$_{\scriptstyle I_4}$};
\node[] at (1.3,0.35) {$\underbrace{\ \ \ \ \ \ \ \ \ \ \ \ \ \ }$};
\node[] at (1.2,0.5) {$...$};
%
\node[] at (-1.3,-0.20) {$^{\scriptstyle J_1}$};
\node[] at (-1.3,-0.35) {$\overbrace{\ \ \ \ \ \ \ \ \ \ \ \ \ \ }$};
\node[] at (-1.2,-0.5) {$...$};
\node[] at (1.3,-0.20) {$^{\scriptstyle J_2}$};
\node[] at (1.3,-0.35) {$\overbrace{\ \ \ \ \ \ \ \ \ \ \ \ \ \ }$};
\node[] at (1.2,-0.5) {$...$};
\node[] at (-1.3,-2.25) {$_{\scriptstyle J_3}$};
\node[] at (-1.3,-2.05) {$\underbrace{\ \ \ \ \ \ \ \ \ \ \ \ \ \ }$};
\node[] at (-1,-1.7) {$...$};
\node[] at (1.3,-2.25) {$_{\scriptstyle J_4}$};
\node[] at (1.3,-2.05) {$\underbrace{\ \ \ \ \ \ \ \ \ \ \ \ \ \ }$};
\node[] at (1,-1.7) {$...$};
\node[] (uu_1) at (-2,2.0) {};
\node[] (uu_2) at (-1,2.0) {};
\node[] (uu_3) at (2,2.0) {};
\node[] (uu_4) at (1,2.0) {};
\node[] (ud_1) at (-2,0.4) {};
\node[] (ud_2) at (-1,0.4) {};
\node[] (ud_3) at (2,0.4) {};
\node[] (ud_4) at (1,0.4) {};
\node[] (dd_1) at (-2,-2.0) {};
\node[] (dd_2) at (-1,-2.0) {};
\node[] (dd_3) at (2,-2.0) {};
\node[] (dd_4) at (1,-2.0) {};
\node[] (du_1) at (-2,-0.4) {};
\node[] (du_2) at (-1,-0.4) {};
\node[] (du_3) at (2,-0.4) {};
\node[] (du_4) at (1,-0.4) {};
 \node[int] (a) at (0,-1.0) {};
 \node[int] (b) at (0,1.0) {};
\node[] (d) at (0,-2.4) {};
\draw (uu_1) edge[latex-, leftblue] (b);
\draw (uu_2) edge[latex-, leftblue] (b);
\draw (uu_3) edge[latex-, rightblue] (b);
\draw (uu_4) edge[latex-, rightblue] (b);
\draw (b) edge[latex-, leftblue] (ud_1);
\draw (b) edge[latex-, leftblue] (ud_2);
\draw (b) edge[latex-, rightblue] (ud_3);
\draw (b) edge[latex-, rightblue] (ud_4);
\draw (du_1) edge[latex-, leftblue] (a);
\draw (du_2) edge[latex-, leftblue] (a);
\draw (du_3) edge[latex-, rightblue] (a);
\draw (du_4) edge[latex-, rightblue] (a);
\draw (a) edge[latex-, leftblue] (dd_1);
\draw (a) edge[latex-, leftblue] (dd_2);
\draw (a) edge[latex-, rightblue] (dd_3);
\draw (a) edge[latex-, rightblue] (dd_4);
 \draw (b) edge[latex-, leftblue, bend left=50] (a);
 \draw (b) edge[latex-, leftblue, bend right=50] (a);
  \draw (b) edge[latex-, leftblue, bend left=30] (a);
 \draw (b) edge[latex-, leftblue, bend right=30] (a);
\end{tikzpicture}}\Ea
\lon
\Ba{c}\resizebox{24mm}{!}
{\begin{tikzpicture}[baseline=-1ex]
\node[] at (-0.6,1.0) {$^{\scriptstyle I_1\sqcup J_1}$};
\node[] at (-0.7,0.8) {$\overbrace{\  \ \ \ \ \ \ \ \ \ \ }$};
\node[] at (-0.45,0.5) {$...$};
\node[] at (0.6,1.0) {$^{\scriptstyle I_2\sqcup J_2}$};
\node[] at (0.7,0.8) {$\overbrace{\  \ \ \ \ \ \ \ \ \ \ }$};
\node[] at (0.45,0.5) {$...$};
\node[] at (-0.6,-1.0) {$_{\scriptstyle I_3\sqcup J_3}$};
\node[] at (-0.7,-0.8) {$\underbrace{\  \ \ \ \ \ \ \ \ \ \ }$};
\node[] at (-0.45,-0.5) {$...$};
\node[] at (0.6,-1.0) {$_{\scriptstyle I_4\sqcup J_4}$};
\node[] at (0.7,-0.8) {$\underbrace{\  \ \ \ \ \ \ \ \ \ \ }$};
\node[] at (0.45,-0.5) {$...$};
\node[int] (0) at (0,0) {};
\node[] (1) at (-1.2,0.8) {};
\node[] (2) at (-0.3,0.8) {};
\node[] (3) at (0.3,0.8) {};
\node[] (4) at (1.2,0.8) {};
\node[int] (0) at (0,0) {};
\node[] (5) at (-1.2,-0.8) {};
\node[] (6) at (-0.3,-0.8) {};
\node[] (7) at (0.3,-0.8) {};
\node[] (8) at (1.2,-0.8) {};
\draw (1) edge[latex-,leftblue] (0);
\draw (2) edge[latex-,leftblue] (0);
\draw (3) edge[latex-,rightblue] (0);
\draw (4) edge[latex-,rightblue] (0);
\draw (0) edge[latex-,leftblue] (5);
\draw (0) edge[latex-,leftblue] (6);
\draw (0) edge[latex-,rightblue] (7);
\draw (0) edge[latex-,rightblue] (8);
\end{tikzpicture}}\Ea\ \ , \ \ \ \ \
\circ_K:
\Ba{c}\resizebox{34mm}{!}  {
\begin{tikzpicture}[baseline=-.65ex]
\node[] at (-1.3,2.2) {$^{\scriptstyle I_1}$};
\node[] at (-1.3,2.05) {$\overbrace{\ \ \ \ \ \ \ \ \ \ \ \ \ \ }$};
\node[] at (-1,1.7) {$...$};
\node[] at (1.3,2.2) {$^{\scriptstyle I_2}$};
\node[] at (1.3,2.05) {$\overbrace{\ \ \ \ \ \ \ \ \ \ \ \ \ \ }$};
\node[] at (1,1.7) {$...$};
\node[] at (-1.3,0.15) {$_{\scriptstyle I_3}$};
\node[] at (-1.3,0.35) {$\underbrace{\ \ \ \ \ \ \ \ \ \ \ \ \ \ }$};
\node[] at (-1.2,0.5) {$...$};
\node[] at (1.3,0.15) {$_{\scriptstyle I_4}$};
\node[] at (1.3,0.35) {$\underbrace{\ \ \ \ \ \ \ \ \ \ \ \ \ \ }$};
\node[] at (1.2,0.5) {$...$};
%
\node[] at (-1.3,-0.20) {$^{\scriptstyle J_1}$};
\node[] at (-1.3,-0.35) {$\overbrace{\ \ \ \ \ \ \ \ \ \ \ \ \ \ }$};
\node[] at (-1.2,-0.5) {$...$};
\node[] at (1.3,-0.20) {$^{\scriptstyle J_2}$};
\node[] at (1.3,-0.35) {$\overbrace{\ \ \ \ \ \ \ \ \ \ \ \ \ \ }$};
\node[] at (1.2,-0.5) {$...$};
\node[] at (-1.3,-2.25) {$_{\scriptstyle J_3}$};
\node[] at (-1.3,-2.05) {$\underbrace{\ \ \ \ \ \ \ \ \ \ \ \ \ \ }$};
\node[] at (-1,-1.7) {$...$};
\node[] at (1.3,-2.25) {$_{\scriptstyle J_4}$};
\node[] at (1.3,-2.05) {$\underbrace{\ \ \ \ \ \ \ \ \ \ \ \ \ \ }$};
\node[] at (1,-1.7) {$...$};
\node[] (uu_1) at (-2,2.0) {};
\node[] (uu_2) at (-1,2.0) {};
\node[] (uu_3) at (2,2.0) {};
\node[] (uu_4) at (1,2.0) {};
\node[] (ud_1) at (-2,0.4) {};
\node[] (ud_2) at (-1,0.4) {};
\node[] (ud_3) at (2,0.4) {};
\node[] (ud_4) at (1,0.4) {};
\node[] (dd_1) at (-2,-2.0) {};
\node[] (dd_2) at (-1,-2.0) {};
\node[] (dd_3) at (2,-2.0) {};
\node[] (dd_4) at (1,-2.0) {};
\node[] (du_1) at (-2,-0.4) {};
\node[] (du_2) at (-1,-0.4) {};
\node[] (du_3) at (2,-0.4) {};
\node[] (du_4) at (1,-0.4) {};
 \node[int] (a) at (0,-1.0) {};
 \node[int] (b) at (0,1.0) {};
\node[] (d) at (0,-2.4) {};
\draw (uu_1) edge[latex-, leftblue] (b);
\draw (uu_2) edge[latex-, leftblue] (b);
\draw (uu_3) edge[latex-, rightblue] (b);
\draw (uu_4) edge[latex-, rightblue] (b);
\draw (b) edge[latex-, leftblue] (ud_1);
\draw (b) edge[latex-, leftblue] (ud_2);
\draw (b) edge[latex-, rightblue] (ud_3);
\draw (b) edge[latex-, rightblue] (ud_4);
\draw (du_1) edge[latex-, leftblue] (a);
\draw (du_2) edge[latex-, leftblue] (a);
\draw (du_3) edge[latex-, rightblue] (a);
\draw (du_4) edge[latex-, rightblue] (a);
\draw (a) edge[latex-, leftblue] (dd_1);
\draw (a) edge[latex-, leftblue] (dd_2);
\draw (a) edge[latex-, rightblue] (dd_3);
\draw (a) edge[latex-, rightblue] (dd_4);
 \draw (b) edge[latex-, rightblue, bend left=50] (a);
 \draw (b) edge[latex-, rightblue, bend right=50] (a);
  \draw (b) edge[latex-, rightblue, bend left=30] (a);
 \draw (b) edge[latex-, rightblue, bend right=30] (a);
\end{tikzpicture}}\Ea
\lon
\Ba{c}\resizebox{24mm}{!}
{\begin{tikzpicture}[baseline=-1ex]
\node[] at (-0.6,1.0) {$^{\scriptstyle I_1\sqcup J_1}$};
\node[] at (-0.7,0.8) {$\overbrace{\  \ \ \ \ \ \ \ \ \ \ }$};
\node[] at (-0.45,0.5) {$...$};
\node[] at (0.6,1.0) {$^{\scriptstyle I_2\sqcup J_2}$};
\node[] at (0.7,0.8) {$\overbrace{\  \ \ \ \ \ \ \ \ \ \ }$};
\node[] at (0.45,0.5) {$...$};
\node[] at (-0.6,-1.0) {$_{\scriptstyle I_3\sqcup J_3}$};
\node[] at (-0.7,-0.8) {$\underbrace{\  \ \ \ \ \ \ \ \ \ \ }$};
\node[] at (-0.45,-0.5) {$...$};
\node[] at (0.6,-1.0) {$_{\scriptstyle I_4\sqcup J_4}$};
\node[] at (0.7,-0.8) {$\underbrace{\  \ \ \ \ \ \ \ \ \ \ }$};
\node[] at (0.45,-0.5) {$...$};
\node[int] (0) at (0,0) {};
\node[] (1) at (-1.2,0.8) {};
\node[] (2) at (-0.3,0.8) {};
\node[] (3) at (0.3,0.8) {};
\node[] (4) at (1.2,0.8) {};
\node[int] (0) at (0,0) {};
\node[] (5) at (-1.2,-0.8) {};
\node[] (6) at (-0.3,-0.8) {};
\node[] (7) at (0.3,-0.8) {};
\node[] (8) at (1.2,-0.8) {};
\draw (1) edge[latex-,leftblue] (0);
\draw (2) edge[latex-,leftblue] (0);
\draw (3) edge[latex-,rightblue] (0);
\draw (4) edge[latex-,rightblue] (0);
\draw (0) edge[latex-,leftblue] (5);
\draw (0) edge[latex-,leftblue] (6);
\draw (0) edge[latex-,rightblue] (7);
\draw (0) edge[latex-,rightblue] (8);
\end{tikzpicture}}\Ea
\Eeq

Note that compositions of the form
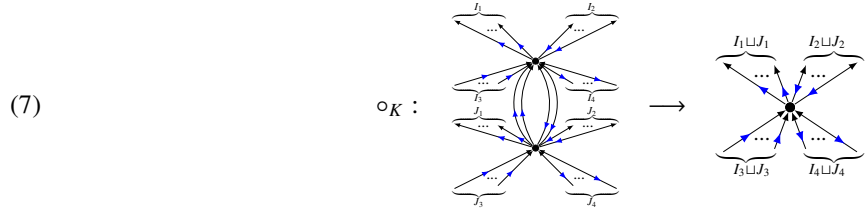
\Beq\label{2new: wrong compos in 2-or prop}
\circ_K:
\Ba{c}\resizebox{34mm}{!}  {
\begin{tikzpicture}[baseline=-.65ex]
\node[] at (-1.3,2.2) {$^{\scriptstyle I_1}$};
\node[] at (-1.3,2.05) {$\overbrace{\ \ \ \ \ \ \ \ \ \ \ \ \ \ }$};
\node[] at (-1,1.7) {$...$};
\node[] at (1.3,2.2) {$^{\scriptstyle I_2}$};
\node[] at (1.3,2.05) {$\overbrace{\ \ \ \ \ \ \ \ \ \ \ \ \ \ }$};
\node[] at (1,1.7) {$...$};
\node[] at (-1.3,0.15) {$_{\scriptstyle I_3}$};
\node[] at (-1.3,0.35) {$\underbrace{\ \ \ \ \ \ \ \ \ \ \ \ \ \ }$};
\node[] at (-1.2,0.5) {$...$};
\node[] at (1.3,0.15) {$_{\scriptstyle I_4}$};
\node[] at (1.3,0.35) {$\underbrace{\ \ \ \ \ \ \ \ \ \ \ \ \ \ }$};
\node[] at (1.2,0.5) {$...$};
%
\node[] at (-1.3,-0.20) {$^{\scriptstyle J_1}$};
\node[] at (-1.3,-0.35) {$\overbrace{\ \ \ \ \ \ \ \ \ \ \ \ \ \ }$};
\node[] at (-1.2,-0.5) {$...$};
\node[] at (1.3,-0.20) {$^{\scriptstyle J_2}$};
\node[] at (1.3,-0.35) {$\overbrace{\ \ \ \ \ \ \ \ \ \ \ \ \ \ }$};
\node[] at (1.2,-0.5) {$...$};
\node[] at (-1.3,-2.25) {$_{\scriptstyle J_3}$};
\node[] at (-1.3,-2.05) {$\underbrace{\ \ \ \ \ \ \ \ \ \ \ \ \ \ }$};
\node[] at (-1,-1.7) {$...$};
\node[] at (1.3,-2.25) {$_{\scriptstyle J_4}$};
\node[] at (1.3,-2.05) {$\underbrace{\ \ \ \ \ \ \ \ \ \ \ \ \ \ }$};
\node[] at (1,-1.7) {$...$};
\node[] (uu_1) at (-2,2.0) {};
\node[] (uu_2) at (-1,2.0) {};
\node[] (uu_3) at (2,2.0) {};
\node[] (uu_4) at (1,2.0) {};
\node[] (ud_1) at (-2,0.4) {};
\node[] (ud_2) at (-1,0.4) {};
\node[] (ud_3) at (2,0.4) {};
\node[] (ud_4) at (1,0.4) {};
\node[] (dd_1) at (-2,-2.0) {};
\node[] (dd_2) at (-1,-2.0) {};
\node[] (dd_3) at (2,-2.0) {};
\node[] (dd_4) at (1,-2.0) {};
\node[] (du_1) at (-2,-0.4) {};
\node[] (du_2) at (-1,-0.4) {};
\node[] (du_3) at (2,-0.4) {};
\node[] (du_4) at (1,-0.4) {};
 \node[int] (a) at (0,-1.0) {};
 \node[int] (b) at (0,1.0) {};
\node[] (d) at (0,-2.4) {};
\draw (uu_1) edge[latex-, leftblue] (b);
\draw (uu_2) edge[latex-, leftblue] (b);
\draw (uu_3) edge[latex-, rightblue] (b);
\draw (uu_4) edge[latex-, rightblue] (b);
\draw (b) edge[latex-, leftblue] (ud_1);
\draw (b) edge[latex-, leftblue] (ud_2);
\draw (b) edge[latex-, rightblue] (ud_3);
\draw (b) edge[latex-, rightblue] (ud_4);
\draw (du_1) edge[latex-, leftblue] (a);
\draw (du_2) edge[latex-, leftblue] (a);
\draw (du_3) edge[latex-, rightblue] (a);
\draw (du_4) edge[latex-, rightblue] (a);
\draw (a) edge[latex-, leftblue] (dd_1);
\draw (a) edge[latex-, leftblue] (dd_2);
\draw (a) edge[latex-, rightblue] (dd_3);
\draw (a) edge[latex-, rightblue] (dd_4);
 \draw (b) edge[latex-, rightblue, bend left=50] (a);
 \draw (b) edge[latex-, leftblue, bend right=50] (a);
  \draw (b) edge[latex-, rightblue, bend left=30] (a);
 \draw (b) edge[latex-, leftblue, bend right=30] (a);
\end{tikzpicture}}\Ea
\lon
\Ba{c}\resizebox{24mm}{!}
{\begin{tikzpicture}[baseline=-1ex]
\node[] at (-0.6,1.0) {$^{\scriptstyle I_1\sqcup J_1}$};
\node[] at (-0.7,0.8) {$\overbrace{\  \ \ \ \ \ \ \ \ \ \ }$};
\node[] at (-0.45,0.5) {$...$};
\node[] at (0.6,1.0) {$^{\scriptstyle I_2\sqcup J_2}$};
\node[] at (0.7,0.8) {$\overbrace{\  \ \ \ \ \ \ \ \ \ \ }$};
\node[] at (0.45,0.5) {$...$};
\node[] at (-0.6,-1.0) {$_{\scriptstyle I_3\sqcup J_3}$};
\node[] at (-0.7,-0.8) {$\underbrace{\  \ \ \ \ \ \ \ \ \ \ }$};
\node[] at (-0.45,-0.5) {$...$};
\node[] at (0.6,-1.0) {$_{\scriptstyle I_4\sqcup J_4}$};
\node[] at (0.7,-0.8) {$\underbrace{\  \ \ \ \ \ \ \ \ \ \ }$};
\node[] at (0.45,-0.5) {$...$};
\node[int] (0) at (0,0) {};
\node[] (1) at (-1.2,0.8) {};
\node[] (2) at (-0.3,0.8) {};
\node[] (3) at (0.3,0.8) {};
\node[] (4) at (1.2,0.8) {};
\node[int] (0) at (0,0) {};
\node[] (5) at (-1.2,-0.8) {};
\node[] (6) at (-0.3,-0.8) {};
\node[] (7) at (0.3,-0.8) {};
\node[] (8) at (1.2,-0.8) {};
\draw (1) edge[latex-,leftblue] (0);
\draw (2) edge[latex-,leftblue] (0);
\draw (3) edge[latex-,rightblue] (0);
\draw (4) edge[latex-,rightblue] (0);
\draw (0) edge[latex-,leftblue] (5);
\draw (0) edge[latex-,leftblue] (6);
\draw (0) edge[latex-,rightblue] (7);
\draw (0) edge[latex-,rightblue] (8);
\end{tikzpicture}}\Ea
\Eeq
are prohibited in 2-oriented props (as they contain at least one  wheel in blue
colour), but are allowed in 1- or 0-oriented 2-directed props.
In the case of 1-oriented (i.e.\ ordinary) prop the horizontal and reduced vertical compositions generate upon iteration {\em any}\, prop composition $\mu_\Ga:  \Ga\langle \cP \rangle (I,\fs) \rar \cP(I,\fs)$. However this is {\em not}\, true in the case of generic $(k+1)$-oriented props with $k\geq 2$. For example, a prop composition
$$
\mu_\Ga: \Ba{c}\resizebox{38mm}{!}{\begin{tikzpicture}[baseline=-1ex]
\node[] (u) at (-0.8,3.5) {};
\node[] (dl) at (-2.4,-2) {};
\node[] (dr) at (2.9,-2) {};
\node[] (b') at (1.1,1) {};
\node[] (b'') at (1.1,0.4) {};
\node[int] (b) at (1.1,1.8) {};
\node[int] (c) at (-0.6,0) {};
\node[] (d) at (0.6,0) {};
\node[] (d') at (0.3,0) {};
\node[int] (e) at (1.1,0) {};
\node[] (e') at (1.1,1.2) {};
\draw (b) edge[latex-,leftblue,rightred, rightgreen] (c);
\draw (b) edge[latex-,rightblue, rightred] (b'');
\draw (b') edge[rightgreen] (e);
\draw (d) edge[leftblue, leftred](c);
\draw (e) edge[latex-,rightgreen](d');
\draw (u) edge[latex-,leftblue,rightred, rightgreen] (b);
\draw (c) edge[latex-,leftblue,rightred, rightgreen] (dl);
\draw (e) edge[latex-,rightblue,leftred, leftgreen] (dr);
\end{tikzpicture}}\Ea
\lon
\Ba{c}\resizebox{31mm}{!}  {
\begin{tikzpicture}[baseline=-1ex]
\node[int] (a) at (0,0) {};
\node[] (u1) at (-1.8,1.8) {};
\node[int] (a) at (0,0) {};
\node[] (d1) at (-1.8,-1.8) {};
\node[] (d2) at (1.8,-1.8) {};
\draw (u1) edge[latex-,leftblue,rightred, rightgreen] (a);
\draw (a) edge[latex-,leftblue,rightred, rightgreen] (d1);
\draw (a) edge[latex-,rightblue,leftred, leftgreen] (d2);
\end{tikzpicture}
}\Ea
$$
can not  be represented as an iteration of more elementary (i.e.\ two vertex) compositions. It is important to notice, however, that such a representation is
{\em always}\, possible with respect to each coloured orientation, and different orientations
produce different (sometimes incompatible as in this particular case) decompositions of $\mu_\Ga$ into a sequence of horizontal and reduced vertical operations; moreover, it is this latter property which is really important when we consider representations of multi-oriented props (see \S 4 below).

\sip

If in the above definition of the natural transformation $\mu$ we restrict only to the subset $G_c^{l+1\uparrow k+1}(I,\fs) \subset G^{l+1\uparrow k+1}(I,\fs)$ of {\em connected}\, graphs, we get the notion of an {\em $(l+1)$-oriented $(k+1)$-directed  properad}\, $\cP$ (cf.\ \cite{Va}). In this case we do not have horizontal compositions in $\cP$, only vertical ones. There is an obvious exact functor from
$(k+1)$-directed properads to $(k+1)$-directed props.

\sip

For any $\cS^{(k+1)}$-module $\cE$ the associated $\cS^{(k+1)}$-module $\cF ree^{l+1\uparrow k+1}\langle \cE \rangle$ is a $(l+1)$-oriented $(k+1)$-directed  prop with contraction maps $\mu_\Ga$ being tautological. It is called the {\em free}\, multi-directed prop
generated by $\cE$. If $l=k$, it is called the free  {\em multi-oriented}\, (more precisely, $(k+1)$-{\em oriented}) prop
generated by $\cE$.

\subsection{Morphisms of multi-oriented props} Let $\cP$ and $\cP'$ be $(l+1)$-oriented $(k+1)$-directed props. A morphism
$$
f: \cP\lon \cP'
$$
of $\cS^{(k+1)}$-modules is called a {\em morphism of  $(l+1)$ oriented $(k+1)$-directed props}\, if for any $\Ga\in G^{l+1\uparrow k+1}(I,\fs)$ the associated diagram
$$
\xymatrix{
  \Ga\langle \cP \rangle (I,\fs) \  \ar[r]^-{\mu_\Ga}\ar[d]_{f^{\ot \# V(\Ga)}} &\  \cP(I,\fs) \ar[d]^-{f_{I,\fs}} \\
 \Ga\langle \cP' \rangle (I,\fs)\ \ar[r]^-{\mu_\Ga} &  \cP(I,\fs)' \\
}
$$
commutes.

\sip

As usual, any morphism $f:\cF ree^{l+1\uparrow k+1}\langle \cE \rangle \rar \cP'$
from a free prop is uniquely determined by its values on the generators $\cE$.

\subsection{Multi-oriented operads}
If a multi-oriented  properad  $\cP=\{ \cP(I,\fs)\}$  is such that $\cP(I,\fs)$ vanishes unless $\check{\fs}_{\scriptstyle 0}^{-1}(out)= 1$ (i.e.\ the functor $\cP$ is non-trivial only on multi-oriented corollas with precisely {\em one}\, outgoing leg with respect to the basic direction) it is called a {\em multi-oriented operad}. Note that there is no such a restriction on non-basic directions.

\subsection{Multi-oriented dioperads} Restricting the functor $\cF ree^{l+1\uparrow k+1}$ (and denoting it by
$\cF ree^{l+1\uparrow k+1}_0$) and the compositions $\mu_\Ga$  to connected graphs of genus zero only one gets the notion of a $(k+1)$-directed
$(l+1)$-oriented {\em dioperad}. Note that a multi-oriented operad is a special case of a multi-oriented dioperad, and the difference between the notions is rather small.

\subsection{Multidirected wheeled prop(erad)s} The above definitions make sense in the case $l=-1$ as well; the associated
$0$-oriented $(k+1)$-directed prop(erad)s $\cP$ are called  {\em $(k+1)$-directed wheeled prop(erad)s} (cf.\ \cite{Me1,MMS}). In this case  graphs $\Ga$ can have an internal edge connecting one and the same vertex, and any prop composition $\mu_\Ga$ can decomposed into an iteration of horizontal compositions and reduced vertical composition defined above, and also of a  trace map
$$
Tr_K: \cP(I,\fs) \lon \cP\left(I\setminus [f_1(K)\sqcup f_2(K), \fs']\right)
$$
which is defined for any two injections $f_{1}: K\rar I$, $f_{2}: K\rar I$ of a finite set $K$ such that
$f_1(K)\cap f_2(K)=\emptyset$ and $\fs\circ f_1= (\fs\circ f_2)^{opp}$; the orientation function $\fs'$ is obtained from $\fs$ by its restriction to the subset $I\setminus [f_1(K)\sqcup f_2(K)]$.
\mip

Thus in the the family of $(k+1)$-directed $(l+1)$ oriented props
the special case $(l=-1,k\geq 0)$ corresponds to the ordinary $2^k$-coloured wheeled prop
while the case $(l=0, k\geq 0)$ to the ordinary $2^k$ coloured  prop. Thus the  really new cases must have $k\geq l\geq 1$.

 \sip

 In the next section \S 3 we introduce multi-oriented versions of some classical operads and props
at the combinatorial level (which is straightforward). In \S 4 we discuss representations of multi-oriented props, i.e.\ explain what these multi-oriented graphical combinatorics give us
in practice (and this step is, perhaps, not that straightforward).

\subsection{Ordinary props as multi-oriented ones}  By the very definition, the category of ordinary props is precisely the the category of 1-oriented props. It is worth mentioning a canonical but very naive functor for any $k\geq 1$
$$
\Ba{rccc}
O^{(k+1)}: & \text{Category\ of\ ordinary\ props} & \lon &  \text{Category\ of}\ (k+1)\text{oriented props}\\
& P=\{P(I,\fs)\} & \lon & O^{(k+1)}(P)=\{P(I,\fs^{(k+1)})\}
\Ea
$$
which simply associates to an $\cS^{(1)}$-module $\{P(I,\fs)\}$ an $\cS^{(k+1)}$-module
 $\{P(I,\fs^{(k+1)})=P(I,\fs)\}$   which is non-trivial (and coincides with $\{P(I,\fs)\}$)
 only for those multi-orientations $\fs^{(k+1)}$ in which
 {\em all}\, extra directions are aligned coherently with the basic one.
More precisely, given an ordinary (1-oriented) prop
$$
P\left(I,\fs: I\rar \f_{0^+}\equiv\{out,in\}\right)
$$
we define
$$
P(I,\fs^{(k+1)}):=\left\{\Ba{rl}
P(I,\fs) & \text{if}\ \fs^{(k+1)}\ \text{satisfies}\ \fs^{(k+1)}_i(\tau):=\fs(i)\ \forall\ i\in I, \ \ \forall \tau\in [k^+],\\
0\ \text{or}\ \emptyset & \text{otherwise}.
\Ea
\right.
 $$
 We do not use this naive functor
in this paper (as it gives nothing new), but it is worth keeping in mind that all classical props can be ``embedded"
into the category of $(k+1)$-oriented props; at least nothing is lost.

\sip

Similarly one can interpret a $(k+1)$-oriented prop as a $(k+l+1)$-oriented prop for any $l\geq 1$. In the next section we consider much less naive extensions of classical operads and props to the multi-oriented setting.

\bip

\bip

{\Large
\section{\bf Multi-oriented versions of some classical operads and props}
}

\mip

\subsection{Multi-oriented operad of (strongly homotopy) associative algebras} Let us recall an explicit combinatorial description of the operad $\cA ss$ of associative algebras in terms of planar 1-oriented corollas. By definition  $\cA ss$ is the quotient,
$$
\cA ss:= {\cF\! ree^{1\text{-or}}\langle A\rangle}/(R)
$$
of the free operad $\cF\! ree^{1\text{-or}}\langle A \rangle$ generated
by  an  $\bS$-module $A=\{A(n)\}_{n\geq 0}$ with
\[
A(n):=\left\{
\Ba{rl}
\K[\bS_2]\equiv
\mbox{span}\left\langle
\Ba{c}\resizebox{8mm}{!}  {\begin{tikzpicture}[baseline=-1ex]
\node[] at (0,0.8) {$^{0}$};
\node[] at (-0.5,-0.71) {$_{1}$};
\node[] at (0.5,-0.71) {$_{2}$};
\node[int] (a) at (0,0) {};
\node[] (u1) at (0,0.8) {};
\node[int] (a) at (0,0) {};
\node[] (d1) at (-0.5,-0.7) {};
\node[] (d2) at (0.5,-0.7) {};
\draw (u1) edge[latex-] (a);
\draw (a) edge[latex-] (d1);
\draw (a) edge[latex-] (d2);
\end{tikzpicture}}\Ea
,
\Ba{c}\resizebox{8mm}{!}  {\begin{tikzpicture}[baseline=-1ex]
\node[] at (0,0.8) {$^{0}$};
\node[] at (-0.5,-0.71) {$_{2}$};
\node[] at (0.5,-0.71) {$_{1}$};
\node[int] (a) at (0,0) {};
\node[] (u1) at (0,0.8) {};
\node[int] (a) at (0,0) {};
\node[] (d1) at (-0.5,-0.7) {};
\node[] (d2) at (0.5,-0.7) {};
\draw (u1) edge[latex-] (a);
\draw (a) edge[latex-] (d1);
\draw (a) edge[latex-] (d2);
\end{tikzpicture}}\Ea
\right
\rangle
\ & \mbox{if}\ n=2, \vspace{3mm}\\
0 & \mbox{otherwise}
\Ea
\right.
\]
by the ideal generated by the relation (together with its $\bS_3$ permutations),
$$
\Ba{c}\resizebox{16mm}{!}  {
\begin{tikzpicture}[baseline=-1ex]
\node[] at (0.5,1.6) {$^{{0}}$};
\node[] at (-0.5,-0.72) {$_{{1}}$};
\node[] at (0.57,-0.72) {$_{{2}}$};
\node[] at (1.12,-0.02) {$_{{3}}$};
%
\node[int] (a) at (0,0) {};
\node[] (u2) at (0.5,1.6) {};
\node[] (ud2) at (1,0.0) {};
\node[int] (u1) at (0.5,0.8) {};
\node[int] (a) at (0,0) {};
\node[] (d1) at (-0.5,-0.7) {};
\node[] (d2) at (0.5,-0.7) {};
\draw (u1) edge[latex-] (a);
\draw (u2) edge[latex-] (u1);
\draw (u1) edge[latex-] (ud2);
\draw (a) edge[latex-] (d1);
\draw (a) edge[latex-] (d2);
\end{tikzpicture}
}\Ea
-
\Ba{c}\resizebox{12mm}{!}  {
\begin{tikzpicture}[baseline=-1ex]
\node[] at (0.5,1.6) {$^{0}$};
\node[] at (1.55,-0.72) {$_{3}$};
\node[] at (0.57,-0.72) {$_{2}$};
\node[] at (0.0,-0.02) {$_{1}$};
%
\node[] (u2) at (0.5,1.6) {};
\node[] (ud2) at (0,0.0) {};
\node[int] (u1) at (0.5,0.8) {};
\node[int] (a) at (1.0,0) {};
\node[] (d1) at (0.5,-0.7) {};
\node[] (d2) at (1.5,-0.7) {};
\draw (u1) edge[latex-] (a);
\draw (u2) edge[latex-] (u1);
\draw (u1) edge[latex-] (ud2);
\draw (a) edge[latex-] (d1);
\draw (a) edge[latex-] (d2);
\end{tikzpicture}
}\Ea =0
$$
Its minimal resolution  is a dg free operad,
$\cA ss_\infty:=(\cF\! ree\langle E\rangle, \delta)$ generated by the $\bS$-module
$E=\{E(n)\}_{\geq 2}$ (whose generators we represent pictorially as {\em planar}\, corollas of homological degree $2-n$)
$$
E(n):=\K[\bS_n][n-2]= \mbox{span}
\left(
\Ba{c}\resizebox{21mm}{!}
{\begin{tikzpicture}[baseline=-1ex]
\node[] at (-1.2,-0.9) {$_{\scriptstyle \sigma(1)}$};
\node[] at (-0.6,-0.9) {$_{\scriptstyle \sigma(2)}$};
\node[] at (-0.0,-0.5) {$...$};
\node[] at (1.3,-0.9) {$_{\scriptstyle \sigma(n)}$};
\node[] at (0.6,-0.9) {$_{\scriptstyle \sigma(n-1)}$};
\node[int] (0) at (0,0) {};
\node[] (1) at (0,0.8) {};
\node[int] (0) at (0,0) {};
\node[] (5) at (-1.2,-0.8) {};
\node[] (6) at (-0.6,-0.8) {};
\node[] (7) at (0.6,-0.8) {};
\node[] (8) at (1.2,-0.8) {};
\draw (1) edge[latex-] (0);
\draw (0) edge[latex-] (5);
\draw (0) edge[latex-] (6);
\draw (0) edge[latex-] (7);
\draw (0) edge[latex-] (8);
\end{tikzpicture}}\Ea
\right)_{\sigma\in \bS_n},
$$
and equipped with the differential given on the generators by
$$
\delta
\Ba{c}\resizebox{20mm}{!}
{\begin{tikzpicture}[baseline=-1ex]
\node[] at (-1.2,-0.9) {$_{\scriptstyle 1}$};
\node[] at (-0.6,-0.9) {$_{\scriptstyle 2}$};
\node[] at (-0.0,-0.5) {$...$};
\node[] at (1.3,-0.9) {$_{\scriptstyle n}$};
\node[] at (0.6,-0.9) {$_{\scriptstyle n-1}$};
\node[int] (0) at (0,0) {};
\node[] (1) at (0,0.8) {};
\node[int] (0) at (0,0) {};
\node[] (5) at (-1.2,-0.8) {};
\node[] (6) at (-0.6,-0.8) {};
\node[] (7) at (0.6,-0.8) {};
\node[] (8) at (1.2,-0.8) {};
\draw (1) edge[latex-] (0);
\draw (0) edge[latex-] (5);
\draw (0) edge[latex-] (6);
\draw (0) edge[latex-] (7);
\draw (0) edge[latex-] (8);
\end{tikzpicture}}\Ea
=\sum_{r=0}^{n-2}\sum_{l=2}^{n-r}
(-1)^{rl+n-r-l+1}
\Ba{c}\resizebox{24mm}{!}
{\begin{tikzpicture}[baseline=-1ex]
\node[] at (-1.4,0.0) {$_{\scriptstyle 1}$};
\node[] at (-0.6,0) {$_{\scriptstyle r}$};
\node[] at (1.4,0.0) {$_{\scriptstyle n}$};
\node[] at (0.6,0) {$_{\scriptstyle {r+l+1}}$};
\node[] at (0.2,-0.6) {$...$};
\node[] at (-0.75,0.2) {$...$};
\node[] at (0.75,0.2) {$...$};
\node[] at (-0.6,-0.85) {$_{\scriptstyle {r+1}}$};
\node[] at (-0.1,-0.85) {$_{\scriptstyle {r+2}}$};
\node[] at (0.7,-0.85) {$_{\scriptstyle {r+l}}$};
\node[int] (0) at (0,0) {};
\node[int] (1) at (0,0.8) {};
\node[] (u) at (0,1.6) {};
\node[] (l1) at (-1.4,0) {};
\node[] (l2) at (-0.7,0) {};
\node[] (r2) at (1.4,0) {};
\node[] (r1) at (0.7,0) {};
\node[int] (0) at (0,0) {};
\node[] (5) at (-0.6,-0.8) {};
\node[] (6) at (-0.1,-0.8) {};
\node[] (8) at (0.6,-0.8) {};
\draw (u) edge[latex-] (1);
\draw (1) edge[latex-] (l1);
\draw (1) edge[latex-] (l2);
\draw (1) edge[latex-] (r1);
\draw (1) edge[latex-] (r2);
\draw (1) edge[latex-] (0);
\draw (0) edge[latex-] (5);
\draw (0) edge[latex-] (6);
\draw (0) edge[latex-] (8);
\end{tikzpicture}}\Ea
$$

Let us first consider the most naive multi-oriented generalization of $\cA ss_\infty$ in which we enlarge
the set of generators by decorating each leg of each planar corolla with $k$ extra orientations in all possible ways while preserving its homological degree,
$$
\Ba{c}\resizebox{23mm}{!}
{\begin{tikzpicture}[baseline=-1ex]
\node[] at (0,0.8) {$^{\scriptstyle 0}$};
\node[] at (-1.2,-0.9) {$_{\scriptstyle 1}$};
\node[] at (-0.6,-0.9) {$_{\scriptstyle 2}$};
\node[] at (-0.0,-0.5) {$...$};
\node[] at (1.3,-0.9) {$_{\scriptstyle n}$};
\node[] at (0.6,-0.9) {$_{\scriptstyle {n-1}}$};
\node[int] (0) at (0,0) {};
\node[] (1) at (0,0.8) {};
\node[int] (0) at (0,0) {};
\node[] (5) at (-1.2,-0.8) {};
\node[] (6) at (-0.6,-0.8) {};
\node[] (7) at (0.6,-0.8) {};
\node[] (8) at (1.2,-0.8) {};
\draw (1) edge[latex-] (0);
\draw (0) edge[latex-] (5);
\draw (0) edge[latex-] (6);
\draw (0) edge[latex-] (7);
\draw (0) edge[latex-] (8);
\end{tikzpicture}}\Ea
\lon
\left\{
\Ba{c}\resizebox{23mm}{!}
{\begin{tikzpicture}[baseline=-1ex]
\node[] at (0,1.2) {$^{\scriptstyle 0}$};
\node[] at (-1.2,-1.25) {$_{\scriptstyle 1}$};
\node[] at (-0.6,-1.25) {$_{\scriptstyle 2}$};
\node[] at (-0.0,-0.6) {$...$};
\node[] at (1.2,-1.25) {$_{\scriptstyle {n}}$};
\node[] at (0.6,-1.25) {$_{\scriptstyle {n-1}}$};
\node[int] (0) at (0,0) {};
\node[] (1) at (0,1.2) {};
\node[int] (0) at (0,0) {};
\node[] (5) at (-1.2,-1.2) {};
\node[] (6) at (-0.6,-1.2) {};
\node[] (7) at (0.6,-1.2) {};
\node[] (8) at (1.2,-1.2) {};
\draw (1) edge[latex-,leftblue, leftred] (0);
\draw (0) edge[latex-,leftblue, leftred] (5);
\draw (0) edge[latex-,leftblue, rightred] (6);
\draw (0) edge[latex-,rightblue, rightred] (7);
\draw (0) edge[latex-,rightblue, leftred] (8);
\end{tikzpicture}}\Ea
,
\Ba{c}\resizebox{23mm}{!}
{\begin{tikzpicture}[baseline=-1ex]
\node[] at (0,1.2) {$^{\scriptstyle 0}$};
\node[] at (-1.2,-1.25) {$_{\scriptstyle 1}$};
\node[] at (-0.6,-1.25) {$_{\scriptstyle 2}$};
\node[] at (-0.0,-0.6) {$...$};
\node[] at (1.2,-1.25) {$_{\scriptstyle {n}}$};
\node[] at (0.6,-1.25) {$_{\scriptstyle {n-1}}$};
\node[int] (0) at (0,0) {};
\node[] (1) at (0,1.2) {};
\node[int] (0) at (0,0) {};
\node[] (5) at (-1.2,-1.2) {};
\node[] (6) at (-0.6,-1.2) {};
\node[] (7) at (0.6,-1.2) {};
\node[] (8) at (1.2,-1.2) {};
\draw (1) edge[latex-,rightblue, leftred] (0);
\draw (0) edge[latex-,leftblue, rightred] (5);
\draw (0) edge[latex-,rightblue, leftred] (6);
\draw (0) edge[latex-,rightblue, rightred] (7);
\draw (0) edge[latex-,rightblue, leftred] (8);
\end{tikzpicture}}\Ea,
\ \ \
\ldots\ \ \ \  \right\}
=
\left\{
\Ba{c}\resizebox{29mm}{!}
{\begin{tikzpicture}[baseline=-1ex]
\node[] at (0,0.8) {$^{\scriptstyle \bar{\fs}_0}$};
\node[] at (-1.2,-0.9) {$_{\scriptstyle \bar{\fs}_1}$};
\node[] at (-0.6,-0.9) {$_{\scriptstyle \bar{\fs}_2}$};
\node[] at (-0.0,-0.5) {$...$};
\node[] at (1.2,-0.9) {$_{\scriptstyle \bar{\fs}_{n}}$};
\node[] at (0.6,-0.9) {$_{\scriptstyle \bar{\fs}_{n-1}}$};
\node[int] (0) at (0,0) {};
\node[] (1) at (0,0.8) {};
\node[int] (0) at (0,0) {};
\node[] (5) at (-1.2,-0.8) {};
\node[] (6) at (-0.6,-0.8) {};
\node[] (7) at (0.6,-0.8) {};
\node[] (8) at (1.2,-0.8) {};
\draw (1) edge[latex-] (0);
\draw (0) edge[latex-] (5);
\draw (0) edge[latex-] (6);
\draw (0) edge[latex-] (7);
\draw (0) edge[latex-] (8);
\end{tikzpicture}}\Ea\right\}_{\forall \bar{\fs}_0, \bar{\fs}_1,\ldots \bar{\fs}_n\in \f_k}
$$
In some pictures we show explicitly only the basic direction while extra directions are indicated only by extra-orientation functions $\bar{\fs}_i$. Denote $\cA svbig_\infty^{(k+1)}:=\cF\! ree^{(k+1)\text{-or}}\left\langle A^{(k+1)}\right\rangle $ be the free (``very big")
  operad generated by these corollas, more precisely, by the associated $\cS^{(k+1)}$-module
   $A^{(k+1)}=\{A^{(k+1)}(I,\fs)\}$   defined formally as follows
$$
A^{(k+1)}(I,\fs)=\left\{\Ba{rl}
0 & \text{if}\ \# I\leq  2\\
0 & \text{if}\ \check{\fs}_{\scriptstyle 0}^{-1}(out)\neq 1\\
\text{span}\langle ord (I')\rangle[\# I-3] &   \text{otherwise}.
\Ea
\right.
$$
 Here $ord(I')$ is the set of total orderings on the finite set $I'=I\setminus \check{\fs}_{\scriptstyle 0}^{-1}(out)$.
The differential in $\cA ss_\infty$ can be extended to $\cA svbig_\infty^{(k+1)}$ by summing over all possible ways
of attaching extra directions $\bar{\fs}\in \f_k$ to the internal edge,
\Beq\label{3new: delta ass k+1 infty}
\delta \Ba{c}\resizebox{20mm}{!}
{\begin{tikzpicture}[baseline=-1ex]
\node[] at (0,0.8) {$^{\scriptstyle \bar{\fs}_0}$};
\node[] at (-1.2,-0.9) {$_{\scriptstyle \bar{\fs}_1}$};
\node[] at (-0.6,-0.9) {$_{\scriptstyle \bar{\fs}_2}$};
\node[] at (-0.0,-0.5) {$...$};
\node[] at (1.2,-0.9) {$_{\scriptstyle \bar{\fs}_{n}}$};
\node[] at (0.6,-0.9) {$_{\scriptstyle \bar{\fs}_{n-1}}$};
\node[int] (0) at (0,0) {};
\node[] (1) at (0,0.8) {};
\node[int] (0) at (0,0) {};
\node[] (5) at (-1.2,-0.8) {};
\node[] (6) at (-0.6,-0.8) {};
\node[] (7) at (0.6,-0.8) {};
\node[] (8) at (1.2,-0.8) {};
\draw (1) edge[latex-] (0);
\draw (0) edge[latex-] (5);
\draw (0) edge[latex-] (6);
\draw (0) edge[latex-] (7);
\draw (0) edge[latex-] (8);
\end{tikzpicture}}\Ea
=\sum_{r=0}^{n-2}\sum_{l=2}^{n-r}
\sum_{\bar{\fs}\in \f_k}
(-1)^{rl+n-r-l+1}
\Ba{c}\resizebox{23mm}{!}
{\begin{tikzpicture}[baseline=-1ex]
\node[] at (-0.1,0.2) {$^{\scriptstyle \bar{\fs}}$};
\node[] at (0,1.6) {$^{\scriptstyle \bar{\fs}_0}$};
\node[] at (-1.4,0.0) {$_{\scriptstyle \bar{\fs}_1}$};
\node[] at (-0.6,0) {$_{\scriptstyle \bar{\fs}_r}$};
\node[] at (1.4,0.0) {$_{\scriptstyle \bar{\fs}_n}$};
\node[] at (0.6,0) {$_{\scriptstyle \bar{\fs}_{r+l+1}}$};
\node[] at (0.2,-0.6) {$...$};
\node[] at (-0.75,0.2) {$...$};
\node[] at (0.75,0.2) {$...$};
\node[] at (-0.6,-0.85) {$_{\scriptstyle \bar{\fs}_{r+1}}$};
\node[] at (-0.1,-0.85) {$_{\scriptstyle \bar{\fs}_{r+2}}$};
\node[] at (0.7,-0.85) {$_{\scriptstyle \bar{\fs}_{r+l}}$};
\node[int] (0) at (0,0) {};
\node[int] (1) at (0,0.8) {};
\node[] (u) at (0,1.6) {};
\node[] (l1) at (-1.4,0) {};
\node[] (l2) at (-0.7,0) {};
\node[] (r2) at (1.4,0) {};
\node[] (r1) at (0.7,0) {};
\node[int] (0) at (0,0) {};
\node[] (5) at (-0.6,-0.8) {};
\node[] (6) at (-0.1,-0.8) {};
\node[] (8) at (0.6,-0.8) {};
\draw (u) edge[latex-] (1);
\draw (1) edge[latex-] (l1);
\draw (1) edge[latex-] (l2);
\draw (1) edge[latex-] (r1);
\draw (1) edge[latex-] (r2);
\draw (1) edge[latex-] (0);
\draw (0) edge[latex-] (5);
\draw (0) edge[latex-] (6);
\draw (0) edge[latex-] (8);
\end{tikzpicture}}\Ea
\Eeq
Hence $\cA svbig_\infty^{(k+1)}$ is just the $2^k$-coloured extension of $\cA ss_\infty$.
Note that the generating corollas in  $\cA svbig_\infty^{(k+1)}$ have at least one ingoing leg and at least one outgoing leg with respect to the basic direction (this condition kills ``curvature terms" in that direction). As we shall see in the next chapter (where we introduce representations of multi-oriented props), it is actually extra directions (if present) which play the genuine role of inputs and outputs. Hence to avoid ``curvature terms" with respect to {\em any}\, direction, we have to consider an ideal $I_1$ in the free operad  $\cA svbig_\infty^{(k+1)}$ generated by those corollas which have no output or no input leg(s) with respect to at least one extra orientation. It is easy to see that the above differential $\delta$ respects this ideal so that the quotient
$$
 \cA sbig_\infty^{(k+1)}:=  \cA svbig_\infty^{(k+1)}/I_1
$$
is a dg free operad again. It is generated by a ``smaller" set of generators,
but still that set can be further reduced. Note that once the basic direction is fixed, the set of extra orientations $\f_k$  can be identified with the set of words of length
$k$ in two letters, $>$ and $<$, and hence can be equipped with a lexicographic order $\leq$. Let us call a generating corolla
$$
\Ba{c}\resizebox{20mm}{!}
{\begin{tikzpicture}[baseline=-1ex]
\node[] at (0,0.8) {$^{\scriptstyle \bar{\fs}_0}$};
\node[] at (-1.2,-0.9) {$_{\scriptstyle \bar{\fs}_1}$};
\node[] at (-0.6,-0.9) {$_{\scriptstyle \bar{\fs}_2}$};
\node[] at (-0.0,-0.5) {$...$};
\node[] at (1.2,-0.9) {$_{\scriptstyle \bar{\fs}_{n}}$};
\node[] at (0.6,-0.9) {$_{\scriptstyle \bar{\fs}_{n-1}}$};
\node[int] (0) at (0,0) {};
\node[] (1) at (0,0.8) {};
\node[int] (0) at (0,0) {};
\node[] (5) at (-1.2,-0.8) {};
\node[] (6) at (-0.6,-0.8) {};
\node[] (7) at (0.6,-0.8) {};
\node[] (8) at (1.2,-0.8) {};
\draw (1) edge[latex-] (0);
\draw (0) edge[latex-] (5);
\draw (0) edge[latex-] (6);
\draw (0) edge[latex-] (7);
\draw (0) edge[latex-] (8);
\end{tikzpicture}}\Ea
$$
{\em special}\, if $\bar{\fs}_1\leq \bar{\fs}_2\leq \ldots \leq \bar{\fs}_n$
(i.e.\ if the planar order agrees with the lexicographic one), and let $I_2$
be the ideal in the free operad $\cA sbig_\infty^{(k+1)}$
generated by non-special corollas. It is again easy to see that the differential $\delta$ respects that second ideal so that the quotient
$$
\cA ss_\infty^{(k+1)}:=  \cA sbig_\infty^{(k+1)}/I_2
$$
is a dg free operad generated by the special corollas (essentially, the main point of this discussion is to motivate the claim that the derivation of  $\cA ss_\infty^{(k+1)}$ given on the generating special corollas by formula
(\ref{3new: delta ass k+1 infty}) is a {\em differential}). It is called the {\em multi-oriented operad
of strongly homotopy associative algebras}. Let $J$ be the differential closure
of the ideal in the free (viewed as a non-differential) operad $\cA ss_\infty^{(k+1)}$ generated by the above corollas with $n\geq 3$. The quotient
$$
\cA ss^{(k+1)}:= \cA ss_\infty^{(k+1)}/J
$$
is called a {\em multi-oriented operad of associative algebras}. We shall see below that this  multi-oriented {\em operad}\, controls  structures which
are governed, in some special case,  by ordinary {\em dioperads}. For example, a representation of $\cA ss^{(2)}$
in a symplectic vector space with one Lagrangian brane  can be identified with
an infinitesimal {\em bi}algebra structure on that brane.

\subsubsection{\bf The simplest non-trivial case $k=1$ in more detail}
The dg operad $\cA ss^{(2)}_\infty$ is generated by planar corollas of homological degree $2-\# I-\# J$
$$
\Ba{c}\resizebox{21mm}{!}
{\begin{tikzpicture}[baseline=-1ex]
\node[] at (-0.6,-1.0) {$_{\scriptstyle I}$};
\node[] at (-0.7,-0.8) {$\underbrace{\  \ \ \ \ \ \ \ \ \ \ }$};
\node[] at (-0.45,-0.5) {$...$};
\node[] at (0.6,-1.0) {$_{\scriptstyle J}$};
\node[] at (0.7,-0.8) {$\underbrace{\  \ \ \ \ \ \ \ \ \ \ }$};
\node[] at (0.45,-0.5) {$...$};
\node[int] (0) at (0,0) {};
\node[] (1) at (0,0.8) {};
\node[int] (0) at (0,0) {};
\node[] (5) at (-1.2,-0.8) {};
\node[] (6) at (-0.3,-0.8) {};
\node[] (7) at (0.3,-0.8) {};
\node[] (8) at (1.2,-0.8) {};
\draw (1) edge[latex-,leftblue] (0);
\draw (0) edge[latex-,leftblue] (5);
\draw (0) edge[latex-,leftblue] (6);
\draw (0) edge[latex-,rightblue] (7);
\draw (0) edge[latex-,rightblue] (8);
\end{tikzpicture}}\Ea
\ \
\text{and}\ \ \
\Ba{c}\resizebox{21mm}{!}
{\begin{tikzpicture}[baseline=-1ex]
\node[] at (-0.6,-1.0) {$_{\scriptstyle I}$};
\node[] at (-0.7,-0.8) {$\underbrace{\  \ \ \ \ \ \ \ \ \ \ }$};
\node[] at (-0.45,-0.5) {$...$};
\node[] at (0.6,-1.0) {$_{\scriptstyle J}$};
\node[] at (0.7,-0.8) {$\underbrace{\  \ \ \ \ \ \ \ \ \ \ }$};
\node[] at (0.45,-0.5) {$...$};
\node[int] (0) at (0,0) {};
\node[] (1) at (0,0.8) {};
\node[int] (0) at (0,0) {};
\node[] (5) at (-1.2,-0.8) {};
\node[] (6) at (-0.3,-0.8) {};
\node[] (7) at (0.3,-0.8) {};
\node[] (8) at (1.2,-0.8) {};
\draw (1) edge[latex-,rightblue] (0);
\draw (0) edge[latex-,leftblue] (5);
\draw (0) edge[latex-,leftblue] (6);
\draw (0) edge[latex-,rightblue] (7);
\draw (0) edge[latex-,rightblue] (8);
\end{tikzpicture}}\Ea, \ \ \ \ \ \ \# I + \#J \geq 2.
$$
where the finite sets of labels $I$ and $J$ are totally ordered (in agreement with the given planar structure of the corollas). The differential is given explicitly by
\Beqrn
\delta \hspace{-3mm}
\Ba{c}\resizebox{21mm}{!}
{\begin{tikzpicture}[baseline=-1ex]
\node[] at (-0.6,-1.0) {$_{\scriptstyle I}$};
\node[] at (-0.7,-0.8) {$\underbrace{\  \ \ \ \ \ \ \ \ \ \ }$};
\node[] at (-0.45,-0.5) {$...$};
\node[] at (0.6,-1.0) {$_{\scriptstyle J}$};
\node[] at (0.7,-0.8) {$\underbrace{\  \ \ \ \ \ \ \ \ \ \ }$};
\node[] at (0.45,-0.5) {$...$};
\node[int] (0) at (0,0) {};
\node[] (1) at (0,0.8) {};
\node[int] (0) at (0,0) {};
\node[] (5) at (-1.2,-0.8) {};
\node[] (6) at (-0.3,-0.8) {};
\node[] (7) at (0.3,-0.8) {};
\node[] (8) at (1.2,-0.8) {};
\draw (1) edge[latex-,leftblue] (0);
\draw (0) edge[latex-,leftblue] (5);
\draw (0) edge[latex-,leftblue] (6);
\draw (0) edge[latex-,rightblue] (7);
\draw (0) edge[latex-,rightblue] (8);
\end{tikzpicture}}\Ea
\hspace{-6mm}
&=&\hspace{-4mm}
\sum_{\ \ I=I_1\sqcup I_2 \sqcup I_3}\hspace{-1mm}
(-1)^{\# I_1\# I_2 +\#I_3 + \#J+1}
\hspace{-3mm}
\Ba{c}\resizebox{32mm}{!}
{\begin{tikzpicture}[baseline=-1ex]
\node[] at (-2.7,-0.15) {$_{\underbrace{\ \ \ }_{I_1}}$};
\node[] at (-0.7,-0.15) {$_{\underbrace{\ \ \ }_{I_3}}$};
\node[] at (0.7,-0.15) {$_{\underbrace{\ \ \ }_{J}}$};
\node[] at (-1.6,-0.6) {$...$};
\node[] at (0.65,0.2) {$...$};
\node[] at (-1.8,-0.9) {$_{\underbrace{\ \ \ \ \ \ \ \ }_{I_2}}$};
\node[int] (0) at (-1.7,0) {};
\node[int] (1) at (0,0.8) {};
\node[] (u) at (0,1.6) {};
\node[] (l1) at (-3.2,0) {};
\node[] (l2) at (-2.5,0) {};
\node[] (r2) at (-0.4,0) {};
\node[] (r1) at (-1.1,0) {};
\node[] (R2) at (1.2,0) {};
\node[] (R1) at (0.4,0) {};
%
\node[] (5) at (-2.4,-0.8) {};
\node[] (6) at (-1.9,-0.8) {};
\node[] (8) at (-1.2,-0.8) {};
\draw (u) edge[latex-,leftblue] (1);
\draw (1) edge[latex-,leftblue] (l1);
\draw (1) edge[latex-,leftblue] (l2);
\draw (1) edge[latex-,leftblue] (r1);
\draw (1) edge[latex-,leftblue] (r2);
\draw (1) edge[latex-,rightblue] (R1);
\draw (1) edge[latex-,rightblue] (R2);
\draw (1) edge[latex-,leftblue] (0);
\draw (0) edge[latex-,leftblue] (5);
\draw (0) edge[latex-,leftblue] (6);
\draw (0) edge[latex-,leftblue] (8);
\end{tikzpicture}}\Ea
+
\sum_{J=J_1\sqcup J_2 \sqcup J_3}\hspace{-1mm}
(-1)^{(\# I+ \#J_1)\# J_2 + \#J_3 +1}
\hspace{-3mm}
\Ba{c}\resizebox{32mm}{!}
{\begin{tikzpicture}[baseline=-1ex]
\node[] at (2.7,-0.15) {$_{\underbrace{\ \ \ }_{J_3}}$};
\node[] at (0.7,-0.15) {$_{\underbrace{\ \ \ }_{J_1}}$};
\node[] at (-0.7,-0.15) {$_{\underbrace{\ \ \ }_{I}}$};
\node[] at (1.6,-0.6) {$...$};
\node[] at (-0.65,0.2) {$...$};
\node[] at (1.8,-0.9) {$_{\underbrace{\ \ \ \ \ \ \ \ }_{J_2}}$};
\node[int] (0) at (1.7,0) {};
\node[int] (1) at (0,0.8) {};
\node[] (u) at (0,1.6) {};
\node[] (l1) at (3.2,0) {};
\node[] (l2) at (2.5,0) {};
\node[] (r2) at (0.4,0) {};
\node[] (r1) at (1.1,0) {};
\node[] (R2) at (-1.2,0) {};
\node[] (R1) at (-0.4,0) {};
%
\node[] (5) at (2.4,-0.8) {};
\node[] (6) at (1.9,-0.8) {};
\node[] (8) at (1.2,-0.8) {};
\draw (u) edge[latex-,leftblue] (1);
\draw (1) edge[latex-,rightblue] (l1);
\draw (1) edge[latex-,rightblue] (l2);
\draw (1) edge[latex-,rightblue] (r1);
\draw (1) edge[latex-,rightblue] (r2);
\draw (1) edge[latex-,leftblue] (R1);
\draw (1) edge[latex-,leftblue] (R2);
\draw (1) edge[latex-,rightblue] (0);
\draw (0) edge[latex-,rightblue] (5);
\draw (0) edge[latex-,rightblue] (6);
\draw (0) edge[latex-,rightblue] (8);
\end{tikzpicture}}\Ea
\\
&&
+\sum_{I=I_1\sqcup I_2 \atop J=J_1\sqcup J_2}(-1)^{\# I_1(\# I_2+\# J_1)+ \# J_2+1}
\left(
\Ba{c}\resizebox{22mm}{!}
{\begin{tikzpicture}[baseline=-1ex]
\node[] at (-0.9,-0.15) {$_{\underbrace{\ \ \ \ }_{I_1}}$};
\node[] at (0.9,-0.15) {$_{\underbrace{\ \ \ \ }_{J_2}}$};
\node[] at (-0.5,-0.9) {$_{\underbrace{\ }_{I_2}}$};
\node[] at (0.5,-0.9) {$_{\underbrace{\   }_{J_1}}$};
\node[] at (-0.5,-0.6) {$...$};
\node[] at (-0.75,0.2) {$...$};
\node[] at (0.75,0.2) {$...$};
\node[int] (0) at (0,0) {};
\node[int] (1) at (0,0.8) {};
\node[] (u) at (0,1.6) {};
\node[] (l1) at (-1.4,0) {};
\node[] (l2) at (-0.7,0) {};
\node[] (r2) at (1.4,0) {};
\node[] (r1) at (0.7,0) {};
\node[int] (0) at (0,0) {};
\node[] (5) at (-0.9,-0.8) {};
\node[] (6) at (-0.3,-0.8) {};
\node[] (7) at (0.3,-0.8) {};
\node[] (8) at (0.9,-0.8) {};
\draw (u) edge[latex-,leftblue] (1);
\draw (1) edge[latex-,leftblue] (l1);
\draw (1) edge[latex-,leftblue] (l2);
\draw (1) edge[latex-,rightblue] (r1);
\draw (1) edge[latex-,rightblue] (r2);
\draw (1) edge[latex-,leftblue] (0);
\draw (0) edge[latex-,leftblue] (5);
\draw (0) edge[latex-,leftblue] (6);
\draw (0) edge[latex-,rightblue] (7);
\draw (0) edge[latex-,rightblue] (8);
\end{tikzpicture}}\Ea
+
\Ba{c}\resizebox{22mm}{!}
{\begin{tikzpicture}[baseline=-1ex]
\node[] at (-0.9,-0.15) {$_{\underbrace{\ \ \ \ }_{I_1}}$};
\node[] at (0.9,-0.15) {$_{\underbrace{\ \ \ \ }_{J_2}}$};
\node[] at (-0.5,-0.9) {$_{\underbrace{\ }_{I_2}}$};
\node[] at (0.5,-0.9) {$_{\underbrace{\   }_{J_1}}$};
\node[] at (-0.5,-0.6) {$...$};
\node[] at (-0.75,0.2) {$...$};
\node[] at (0.75,0.2) {$...$};
\node[int] (0) at (0,0) {};
\node[int] (1) at (0,0.8) {};
\node[] (u) at (0,1.6) {};
\node[] (l1) at (-1.4,0) {};
\node[] (l2) at (-0.7,0) {};
\node[] (r2) at (1.4,0) {};
\node[] (r1) at (0.7,0) {};
\node[int] (0) at (0,0) {};
\node[] (5) at (-0.9,-0.8) {};
\node[] (6) at (-0.3,-0.8) {};
\node[] (7) at (0.3,-0.8) {};
\node[] (8) at (0.9,-0.8) {};
\draw (u) edge[latex-,rightblue] (1);
\draw (1) edge[latex-,leftblue] (l1);
\draw (1) edge[latex-,leftblue] (l2);
\draw (1) edge[latex-,rightblue] (r1);
\draw (1) edge[latex-,rightblue] (r2);
\draw (1) edge[latex-,rightblue] (0);
\draw (0) edge[latex-,leftblue] (5);
\draw (0) edge[latex-,leftblue] (6);
\draw (0) edge[latex-,rightblue] (7);
\draw (0) edge[latex-,rightblue] (8);
\end{tikzpicture}}\Ea
\right)
\nonumber
\Eeqrn
and similarly for the second corolla. Here the summations run over decompositions of the totally ordered sets into disjoint unions of {\em connected} (with respect to the order) subsets.

\mip

The operad $\cA ss^{(2)}$ is generated by the following planar corollas (in homological degree zero)
$$
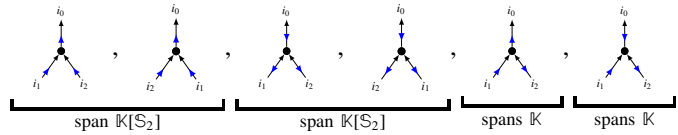
\underbracket{
\Ba{c}\resizebox{12mm}{!}  {
\begin{tikzpicture}[baseline=-1ex]
\node[] at (0,0.8) {$^{i_0}$};
\node[] at (-0.5,-0.71) {$_{i_1}$};
\node[] at (0.5,-0.71) {$_{i_2}$};
\node[int] (a) at (0,0) {};
\node[] (u1) at (0,0.8) {};
\node[int] (a) at (0,0) {};
\node[] (d1) at (-0.5,-0.7) {};
\node[] (d2) at (0.5,-0.7) {};
\draw (u1) edge[latex-, leftblue] (a);
\draw (a) edge[latex-, leftblue] (d1);
\draw (a) edge[latex-, leftblue] (d2);
\end{tikzpicture}
}\Ea
,
\Ba{c}\resizebox{12mm}{!}  {
\begin{tikzpicture}[baseline=-1ex]
\node[] at (0,0.8) {$^{i_0}$};
\node[] at (-0.5,-0.71) {$_{i_2}$};
\node[] at (0.5,-0.71) {$_{i_1}$};
\node[int] (a) at (0,0) {};
\node[] (u1) at (0,0.8) {};
\node[int] (a) at (0,0) {};
\node[] (d1) at (-0.5,-0.7) {};
\node[] (d2) at (0.5,-0.7) {};
\draw (u1) edge[latex-, leftblue] (a);
\draw (a) edge[latex-, leftblue] (d1);
\draw (a) edge[latex-, leftblue] (d2);
\end{tikzpicture}}\Ea}_{\text{span}\ \K[\bS_2]}
,
\underbracket{
\Ba{c}\resizebox{12mm}{!}  {
\begin{tikzpicture}[baseline=-1ex]
\node[] at (0,0.8) {$^{i_0}$};
\node[] at (-0.5,-0.71) {$_{i_1}$};
\node[] at (0.5,-0.71) {$_{i_2}$};
\node[int] (a) at (0,0) {};
\node[] (u1) at (0,0.8) {};
\node[int] (a) at (0,0) {};
\node[] (d1) at (-0.5,-0.7) {};
\node[] (d2) at (0.5,-0.7) {};
\draw (u1) edge[latex-, rightblue] (a);
\draw (a) edge[latex-, rightblue] (d1);
\draw (a) edge[latex-, rightblue] (d2);
\end{tikzpicture}
}\Ea
,
\Ba{c}\resizebox{12mm}{!}  {
\begin{tikzpicture}[baseline=-1ex]
\node[] at (0,0.8) {$^{i_0}$};
\node[] at (-0.5,-0.71) {$_{i_2}$};
\node[] at (0.5,-0.71) {$_{i_1}$};
\node[int] (a) at (0,0) {};
\node[] (u1) at (0,0.8) {};
\node[int] (a) at (0,0) {};
\node[] (d1) at (-0.5,-0.7) {};
\node[] (d2) at (0.5,-0.7) {};
\draw (u1) edge[latex-, rightblue] (a);
\draw (a) edge[latex-, rightblue] (d1);
\draw (a) edge[latex-, rightblue] (d2);
\end{tikzpicture}}\Ea}_{\text{span}\ \K[\bS_2]}
,
\underbracket{
\Ba{c}\resizebox{12mm}{!}  {
\begin{tikzpicture}[baseline=-1ex]
\node[] at (0,0.8) {$^{i_0}$};
\node[] at (-0.5,-0.71) {$_{i_1}$};
\node[] at (0.5,-0.71) {$_{i_2}$};
\node[int] (a) at (0,0) {};
\node[] (u1) at (0,0.8) {};
\node[int] (a) at (0,0) {};
\node[] (d1) at (-0.5,-0.7) {};
\node[] (d2) at (0.5,-0.7) {};
\draw (u1) edge[latex-, leftblue] (a);
\draw (a) edge[latex-, leftblue] (d1);
\draw (a) edge[latex-, rightblue] (d2);
\end{tikzpicture}
}\Ea}_{\text{spans}\ \K}
,
\underbracket{
\Ba{c}\resizebox{12mm}{!}  {
\begin{tikzpicture}[baseline=-1ex]
\node[] at (0,0.8) {$^{i_0}$};
\node[] at (-0.5,-0.71) {$_{i_1}$};
\node[] at (0.5,-0.71) {$_{i_2}$};
\node[int] (a) at (0,0) {};
\node[] (u1) at (0,0.8) {};
\node[int] (a) at (0,0) {};
\node[] (d1) at (-0.5,-0.7) {};
\node[] (d2) at (0.5,-0.7) {};
\draw (u1) edge[latex-, rightblue] (a);
\draw (a) edge[latex-, leftblue] (d1);
\draw (a) edge[latex-, rightblue] (d2);
\end{tikzpicture}
}\Ea}_{\text{spans}\ \K}
$$
while the relations are given by
\Beq\label{3new: Ass^2 relations 1st set}
\Ba{c}\resizebox{14mm}{!}  {
\begin{tikzpicture}[baseline=-1ex]
\node[] at (0.5,1.6) {$^{{i_0}}$};
\node[] at (-0.5,-0.72) {$_{{i_1}}$};
\node[] at (0.57,-0.72) {$_{{i_2}}$};
\node[] at (1.12,-0.02) {$_{{i_3}}$};
\node[int] (a) at (0,0) {};
\node[] (u2) at (0.5,1.6) {};
\node[] (ud2) at (1,0.0) {};
\node[int] (u1) at (0.5,0.8) {};
\node[int] (a) at (0,0) {};
\node[] (d1) at (-0.5,-0.7) {};
\node[] (d2) at (0.5,-0.7) {};
\draw (u1) edge[latex-,leftblue] (a);
\draw (u2) edge[latex-,leftblue] (u1);
\draw (u1) edge[latex-, leftblue] (ud2);
\draw (a) edge[latex-,leftblue] (d1);
\draw (a) edge[latex-,leftblue] (d2);
\end{tikzpicture}
}\Ea
=
\Ba{c}\resizebox{14mm}{!}  {
\begin{tikzpicture}[baseline=-1ex]
\node[] at (0.5,1.6) {$^{{i_0}}$};
\node[] at (1.55,-0.72) {$_{{i_3}}$};
\node[] at (0.57,-0.72) {$_{{i_2}}$};
\node[] at (0.0,-0.02) {$_{{i_1}}$};
\node[] (u2) at (0.5,1.6) {};
\node[] (ud2) at (0,0.0) {};
\node[int] (u1) at (0.5,0.8) {};
\node[int] (a) at (1.0,0) {};
\node[] (d1) at (0.5,-0.7) {};
\node[] (d2) at (1.5,-0.7) {};
\draw (u1) edge[latex-,leftblue] (a);
\draw (u2) edge[latex-,leftblue] (u1);
\draw (u1) edge[latex-,leftblue] (ud2);
\draw (a) edge[latex-,leftblue] (d1);
\draw (a) edge[latex-,leftblue] (d2);
\end{tikzpicture}
}\Ea, \ \ \ \ \ \ \
\Ba{c}\resizebox{14mm}{!}  {
\begin{tikzpicture}[baseline=-1ex]
\node[] at (0.5,1.6) {$^{{i_0}}$};
\node[] at (-0.5,-0.72) {$_{{i_1}}$};
\node[] at (0.57,-0.72) {$_{{i_2}}$};
\node[] at (1.12,-0.02) {$_{{i_3}}$};
\node[int] (a) at (0,0) {};
\node[] (u2) at (0.5,1.6) {};
\node[] (ud2) at (1,0.0) {};
\node[int] (u1) at (0.5,0.8) {};
\node[int] (a) at (0,0) {};
\node[] (d1) at (-0.5,-0.7) {};
\node[] (d2) at (0.5,-0.7) {};
\draw (u1) edge[latex-,leftblue] (a);
\draw (u2) edge[latex-,leftblue] (u1);
\draw (u1) edge[latex-, rightblue] (ud2);
\draw (a) edge[latex-,leftblue] (d1);
\draw (a) edge[latex-,leftblue] (d2);
\end{tikzpicture}
}\Ea
=
\Ba{c}\resizebox{14mm}{!}  {
\begin{tikzpicture}[baseline=-1ex]
\node[] at (0.5,1.6) {$^{{i_0}}$};
\node[] at (1.55,-0.72) {$_{{i_3}}$};
\node[] at (0.57,-0.72) {$_{{i_2}}$};
\node[] at (0.0,-0.02) {$_{{i_1}}$};
\node[] (u2) at (0.5,1.6) {};
\node[] (ud2) at (0,0.0) {};
\node[int] (u1) at (0.5,0.8) {};
\node[int] (a) at (1.0,0) {};
\node[] (d1) at (0.5,-0.7) {};
\node[] (d2) at (1.5,-0.7) {};
\draw (u1) edge[latex-,leftblue] (a);
\draw (u2) edge[latex-,leftblue] (u1);
\draw (u1) edge[latex-,leftblue] (ud2);
\draw (a) edge[latex-,leftblue] (d1);
\draw (a) edge[latex-,rightblue] (d2);
\end{tikzpicture}
}\Ea
+
\Ba{c}\resizebox{14mm}{!}  {
\begin{tikzpicture}[baseline=-1ex]
\node[] at (0.5,1.6) {$^{{i_0}}$};
\node[] at (1.55,-0.72) {$_{{i_3}}$};
\node[] at (0.57,-0.72) {$_{{i_2}}$};
\node[] at (0.0,-0.02) {$_{{i_1}}$};
\node[] (u2) at (0.5,1.6) {};
\node[] (ud2) at (0,0.0) {};
\node[int] (u1) at (0.5,0.8) {};
\node[int] (a) at (1.0,0) {};
\node[] (d1) at (0.5,-0.7) {};
\node[] (d2) at (1.5,-0.7) {};
\draw (u1) edge[latex-,rightblue] (a);
\draw (u2) edge[latex-,leftblue] (u1);
\draw (u1) edge[latex-,leftblue] (ud2);
\draw (a) edge[latex-,leftblue] (d1);
\draw (a) edge[latex-,rightblue] (d2);
\end{tikzpicture}
}\Ea, \ \ \
\Ba{c}\resizebox{14mm}{!}  {
\begin{tikzpicture}[baseline=-1ex]
\node[] at (0.5,1.6) {$^{{i_0}}$};
\node[] at (-0.5,-0.72) {$_{{i_1}}$};
\node[] at (0.57,-0.72) {$_{{i_2}}$};
\node[] at (1.12,-0.02) {$_{{i_3}}$};
\node[int] (a) at (0,0) {};
\node[] (u2) at (0.5,1.6) {};
\node[] (ud2) at (1,0.0) {};
\node[int] (u1) at (0.5,0.8) {};
\node[int] (a) at (0,0) {};
\node[] (d1) at (-0.5,-0.7) {};
\node[] (d2) at (0.5,-0.7) {};
\draw (u1) edge[latex-,leftblue] (a);
\draw (u2) edge[latex-,leftblue] (u1);
\draw (u1) edge[latex-, rightblue] (ud2);
\draw (a) edge[latex-,leftblue] (d1);
\draw (a) edge[latex-,rightblue] (d2);
\end{tikzpicture}
}\Ea
=
\Ba{c}\resizebox{14mm}{!}  {
\begin{tikzpicture}[baseline=-1ex]
\node[] at (0.5,1.6) {$^{{i_0}}$};
\node[] at (1.55,-0.72) {$_{{i_3}}$};
\node[] at (0.57,-0.72) {$_{{i_2}}$};
\node[] at (0.0,-0.02) {$_{{i_1}}$};
\node[] (u2) at (0.5,1.6) {};
\node[] (ud2) at (0,0.0) {};
\node[int] (u1) at (0.5,0.8) {};
\node[int] (a) at (1.0,0) {};
\node[] (d1) at (0.5,-0.7) {};
\node[] (d2) at (1.5,-0.7) {};
\draw (u1) edge[latex-,rightblue] (a);
\draw (u2) edge[latex-,leftblue] (u1);
\draw (u1) edge[latex-,leftblue] (ud2);
\draw (a) edge[latex-,rightblue] (d1);
\draw (a) edge[latex-,rightblue] (d2);
\end{tikzpicture}
}\Ea
\Eeq


\Beq\label{3new: Ass^2 relations 2nd set}
\Ba{c}\resizebox{14mm}{!}  {
\begin{tikzpicture}[baseline=-1ex]
\node[] at (0.5,1.6) {$^{{i_0}}$};
\node[] at (-0.5,-0.72) {$_{{i_1}}$};
\node[] at (0.57,-0.72) {$_{{i_2}}$};
\node[] at (1.12,-0.02) {$_{{i_3}}$};
\node[int] (a) at (0,0) {};
\node[] (u2) at (0.5,1.6) {};
\node[] (ud2) at (1,0.0) {};
\node[int] (u1) at (0.5,0.8) {};
\node[int] (a) at (0,0) {};
\node[] (d1) at (-0.5,-0.7) {};
\node[] (d2) at (0.5,-0.7) {};
\draw (u1) edge[latex-,rightblue] (a);
\draw (u2) edge[latex-,rightblue] (u1);
\draw (u1) edge[latex-, rightblue] (ud2);
\draw (a) edge[latex-,rightblue] (d1);
\draw (a) edge[latex-,rightblue] (d2);
\end{tikzpicture}
}\Ea
=
\Ba{c}\resizebox{14mm}{!}  {
\begin{tikzpicture}[baseline=-1ex]
\node[] at (0.5,1.6) {$^{{i_0}}$};
\node[] at (1.55,-0.72) {$_{{i_3}}$};
\node[] at (0.57,-0.72) {$_{{i_2}}$};
\node[] at (0.0,-0.02) {$_{{i_1}}$};
\node[] (u2) at (0.5,1.6) {};
\node[] (ud2) at (0,0.0) {};
\node[int] (u1) at (0.5,0.8) {};
\node[int] (a) at (1.0,0) {};
\node[] (d1) at (0.5,-0.7) {};
\node[] (d2) at (1.5,-0.7) {};
\draw (u1) edge[latex-,rightblue] (a);
\draw (u2) edge[latex-,rightblue] (u1);
\draw (u1) edge[latex-,rightblue] (ud2);
\draw (a) edge[latex-,rightblue] (d1);
\draw (a) edge[latex-,rightblue] (d2);
\end{tikzpicture}
}\Ea, \ \ \ \ \ \ \
\Ba{c}\resizebox{14mm}{!}  {
\begin{tikzpicture}[baseline=-1ex]
\node[] at (0.5,1.6) {$^{{i_0}}$};
\node[] at (-0.5,-0.72) {$_{{i_1}}$};
\node[] at (0.57,-0.72) {$_{{i_2}}$};
\node[] at (1.12,-0.02) {$_{{i_3}}$};
\node[int] (a) at (0,0) {};
\node[] (u2) at (0.5,1.6) {};
\node[] (ud2) at (1,0.0) {};
\node[int] (u1) at (0.5,0.8) {};
\node[int] (a) at (0,0) {};
\node[] (d1) at (-0.5,-0.7) {};
\node[] (d2) at (0.5,-0.7) {};
\draw (u1) edge[latex-,rightblue] (a);
\draw (u2) edge[latex-,rightblue] (u1);
\draw (u1) edge[latex-, rightblue] (ud2);
\draw (a) edge[latex-,leftblue] (d1);
\draw (a) edge[latex-,rightblue] (d2);
\end{tikzpicture}
}\Ea
+
\Ba{c}\resizebox{14mm}{!}  {
\begin{tikzpicture}[baseline=-1ex]
\node[] at (0.5,1.6) {$^{{i_0}}$};
\node[] at (-0.5,-0.72) {$_{{i_1}}$};
\node[] at (0.57,-0.72) {$_{{i_2}}$};
\node[] at (1.12,-0.02) {$_{{i_3}}$};
\node[int] (a) at (0,0) {};
\node[] (u2) at (0.5,1.6) {};
\node[] (ud2) at (1,0.0) {};
\node[int] (u1) at (0.5,0.8) {};
\node[int] (a) at (0,0) {};
\node[] (d1) at (-0.5,-0.7) {};
\node[] (d2) at (0.5,-0.7) {};
\draw (u1) edge[latex-,leftblue] (a);
\draw (u2) edge[latex-,rightblue] (u1);
\draw (u1) edge[latex-, rightblue] (ud2);
\draw (a) edge[latex-,leftblue] (d1);
\draw (a) edge[latex-,rightblue] (d2);
\end{tikzpicture}
}\Ea
=
\Ba{c}\resizebox{14mm}{!}  {
\begin{tikzpicture}[baseline=-1ex]
\node[] at (0.5,1.6) {$^{{i_0}}$};
\node[] at (1.55,-0.72) {$_{{i_3}}$};
\node[] at (0.57,-0.72) {$_{{i_2}}$};
\node[] at (0.0,-0.02) {$_{{i_1}}$};
\node[] (u2) at (0.5,1.6) {};
\node[] (ud2) at (0,0.0) {};
\node[int] (u1) at (0.5,0.8) {};
\node[int] (a) at (1.0,0) {};
\node[] (d1) at (0.5,-0.7) {};
\node[] (d2) at (1.5,-0.7) {};
\draw (u1) edge[latex-,rightblue] (a);
\draw (u2) edge[latex-,rightblue] (u1);
\draw (u1) edge[latex-,leftblue] (ud2);
\draw (a) edge[latex-,rightblue] (d1);
\draw (a) edge[latex-,rightblue] (d2);
\end{tikzpicture}
}\Ea,
 \ \ \
\Ba{c}\resizebox{14mm}{!}  {
\begin{tikzpicture}[baseline=-1ex]
\node[] at (0.5,1.6) {$^{{i_0}}$};
\node[] at (-0.5,-0.72) {$_{{i_1}}$};
\node[] at (0.57,-0.72) {$_{{i_2}}$};
\node[] at (1.12,-0.02) {$_{{i_3}}$};
%
\node[int] (a) at (0,0) {};
\node[] (u2) at (0.5,1.6) {};
\node[] (ud2) at (1,0.0) {};
\node[int] (u1) at (0.5,0.8) {};
\node[int] (a) at (0,0) {};
\node[] (d1) at (-0.5,-0.7) {};
\node[] (d2) at (0.5,-0.7) {};
\draw (u1) edge[latex-,leftblue] (a);
\draw (u2) edge[latex-,rightblue] (u1);
\draw (u1) edge[latex-, rightblue] (ud2);
\draw (a) edge[latex-,leftblue] (d1);
\draw (a) edge[latex-,leftblue] (d2);
\end{tikzpicture}
}\Ea
=
\Ba{c}\resizebox{14mm}{!}  {
\begin{tikzpicture}[baseline=-1ex]
\node[] at (0.5,1.6) {$^{{i_0}}$};
\node[] at (1.55,-0.72) {$_{{i_3}}$};
\node[] at (0.57,-0.72) {$_{{i_2}}$};
\node[] at (0.0,-0.02) {$_{{i_1}}$};
%
\node[] (u2) at (0.5,1.6) {};
\node[] (ud2) at (0,0.0) {};
\node[int] (u1) at (0.5,0.8) {};
\node[int] (a) at (1.0,0) {};
\node[] (d1) at (0.5,-0.7) {};
\node[] (d2) at (1.5,-0.7) {};
\draw (u1) edge[latex-,rightblue] (a);
\draw (u2) edge[latex-,rightblue] (u1);
\draw (u1) edge[latex-,leftblue] (ud2);
\draw (a) edge[latex-,leftblue] (d1);
\draw (a) edge[latex-,rightblue] (d2);
\end{tikzpicture}
}\Ea
\Eeq

One can describe similarly
 the operad $\cA ss^{(k+1)}$ in terms of generators and relations.

\subsubsection{\bf Theorem}\label{3new: Theorem on Ass_infty^k+1} {\em The natural projection $\cA ss_\infty^{(k+1)}
\lon \cA ss^{(k+1)}$ is a quasi-isomorphism.}

\begin{proof}
We have to show that $H^\bu(\cA ss_\infty^{(k+1)}(I,\fs))=\cA ss^{(k+1)}(I,\fs)$ for any
$(k+1)$-oriented set $(I,\fs)$. In fact, it is enough to show that the cohomology of the operad $H^\bu(\cA ss_\infty^{(k+1)}(I,\fs))$ is concentrated in degree zero because that would imply the required equality due to the fact that the complex $\cA ss_\infty^{(k+1)}(I,\fs)$ is non-positively graded.

 \sip

 We shall prove the claim by induction over $\# I=n+1$, and abbreviate
 the notation $\cA ss_\infty^{(2)}(n):=\cA ss_\infty^{(2)}(I,\fs)$
 and $\cA ss^{(2)}(n):=\cA ss^{(2)}(I,\fs)$.
 When $n=2$, the equality $H^\bu(\cA ss_\infty^{(2)}(n))=\cA ss^{(2)}(n)$
is obvious. Assume
it is true for all multi-oriented sets with $\# I \leq n+1$,
and consider the complex $\cA ss_\infty^{(k+1)}(n+1)$; we can assume without loss of generality that the input (with respect to the basic colour) legs of any graph from $\cA ss_\infty^{(2)}(n+1)$ are labelled from left to right (in accordance with the planar structure) by $1,2,\ldots, n+1$ (while the root vertex by $n+2$).

\sip

 Consider first a filtration of $\cA ss_\infty^{(2)}(n+1)$ by the total number of vertices lying on the  path from the root edge to the leg labelled by $1$ (and call it a {\em special} path), and let $Gr(n+1)$ denote the associated graded. Consider next
a filtration of $Gr(n+1)$ by the total valency of vertices lying on the special path (and denote the set of such vertices by $V_{sp}$), and let $(\cE_r, \delta_r)$ be the associated spectral sequence (converging to $H^\bu(\cA ss_\infty^{(2)}(n+1)$). The initial page $(\cE_0,\delta_0)$ is isomorphic to the direct sum of tensor products of complexes of the form
 $\cA ss_\infty^{(k+1)}(n')$ with all  $n' \leq n$ so that by the induction hypothesis we can easily describe the next page
 of the spectral sequence:
 $$
 \cE_1=H^\bu(\cE_0)\cong \bigoplus_{\text{special\ paths}}\bigoplus_{ n=\sum_{v\in V_{sp}} n_v\atop n_v\geq 1}\bigotimes_{v\in V_{sp}} C_v(n_v)
 $$
 where $C_v(n_v)$ is a complex spanned by planar corollas of the form
$$
 
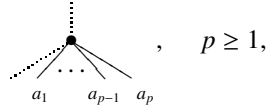
\begin{xy}
 <0mm,0mm>*{\bullet};<0mm,0mm>*{}**@{},
 <0mm,0mm>*{};<-8mm,-5mm>*{}**@{.},
 <0mm,0mm>*{};<-4.5mm,-5mm>*{}**@{-},
 <0mm,0mm>*{};<0mm,-4mm>*{\ldots}**@{},
 <0mm,0mm>*{};<4.5mm,-5mm>*{}**@{-},
 <0mm,0mm>*{};<8mm,-5mm>*{}**@{-},
   <0mm,0mm>*{};<-4mm,-7.9mm>*{^{a_1}}**@{},
    <0mm,0mm>*{};<4.5mm,-7.9mm>*{^{a_{p-1}}}**@{},
   <0mm,0mm>*{};<10.0mm,-7.9mm>*{^{a_p}}**@{},
 <0mm,0mm>*{};<0mm,5mm>*{}**@{.},
 \end{xy}, \ \ \ \ \ p\geq 1,
$$
whose dashed legs belong to the given special path (and are equipped with the induced multi-orientations from that special path) while solid legs  are decorated by arbitrary elements  of the unital extension of the operad $\cA ss^{(k+1)}$,
$$
a_i \in \cA ss^{(k+1)}_u(n_i):=\left\{\Ba{rl}
\cA ss^{(k+1)}(n_i) & \text{if}\ \# n_i\geq 2\\
\K & \text{if}\ \# n_i=1\\
0 & \text{if}\ \# n_i=0\\
\Ea \right. ,\ \ \ \  i \in [p],
 $$
 subject to the condition that
 $$
 \sum_{i=1}^p n_i=n_v,
 $$
The differential on $C_v(n_v)$ is non-trivial only on the root corolla
on which it acts as follows (we suppress some extra orientations in the picture),
\Beq\label{d1+}
\delta_1
 \begin{xy}
 <0mm,0mm>*{\bullet};<0mm,0mm>*{}**@{},
 <0mm,0mm>*{};<-8mm,-5mm>*{}**@{.},
 <0mm,0mm>*{};<-4.5mm,-5mm>*{}**@{-},
 <0mm,0mm>*{};<0mm,-4mm>*{\ldots}**@{},
 <0mm,0mm>*{};<4.5mm,-5mm>*{}**@{-},
 <0mm,0mm>*{};<8mm,-5mm>*{}**@{-},
   <0mm,0mm>*{};<-4mm,-7.9mm>*{^{a_1}}**@{},
    <0mm,0mm>*{};<4.5mm,-7.9mm>*{^{a_{p-1}}}**@{},
   <0mm,0mm>*{};<10.0mm,-7.9mm>*{^{a_p}}**@{},
 <0mm,0mm>*{};<0mm,5mm>*{}**@{.},
 \end{xy}
=\left\{\Ba{lr}
 \sum_{i=0}^{p-2}\sum_{\bar{\fs}\in \f_{k}}(-1)^{p-i+1}
\begin{xy}
 <0mm,0mm>*{\bullet};<0mm,0mm>*{}**@{},
 <0mm,0mm>*{};<-8mm,-5mm>*{}**@{.},
 <0mm,0mm>*{};<-4mm,-5mm>*{}**@{-},
 <0mm,0mm>*{};<-5.6mm,-5mm>*{..}**@{},
 <0mm,0mm>*{};<2.6mm,-5mm>*{..}**@{},
 <0mm,0mm>*{};<4.5mm,-5mm>*{}**@{-},
 <0mm,0mm>*{};<-0.5mm,-5mm>*{}**@{-},
 <0mm,0mm>*{};<8mm,-5mm>*{}**@{-},
   <0mm,0mm>*{};<0.5mm,-3.5mm>*{^{\bar{\fs}}}**@{},
   <0mm,0mm>*{};<-4.5mm,-7.9mm>*{^{a_i}}**@{},
   <0mm,0mm>*{};<10.0mm,-7.9mm>*{^{a_p}}**@{},
 <0mm,0mm>*{};<0mm,5mm>*{}**@{.},
  <-0.5mm,-5.4mm>*{\bullet};
<-0.7mm,-5.4mm>*{};<-3.6mm,-9mm>*{}**@{-},
<-0.7mm,-5.4mm>*{};<3mm,-9mm>*{}**@{-},
<-0.7mm,-5.4mm>*{};<-4mm,-11mm>*{_{a_{i+1}}}**@{},
<-0.7mm,-5.4mm>*{};<4mm,-11mm>*{_{a_{i+2}}}**@{},
 \end{xy} & \mbox{for}\ p\geq 3\\
 -\sum_{\bar{\fs}\in \f_{k}}
 \Ba{c}
\begin{xy}
 <0mm,0mm>*{\bullet};<0mm,0mm>*{}**@{},
 <0mm,0mm>*{};<-4mm,-5mm>*{}**@{.},
 <0mm,0mm>*{};<0mm,-5mm>*{}**@{-},
 <0mm,0mm>*{};<0mm,5mm>*{}**@{.},
 <0mm,0mm>*{};<0.9mm,-3.5mm>*{^{\bar{\fs}}}**@{},
  <0mm,-5.4mm>*{\bullet};
<-0mm,-5.4mm>*{};<-3mm,-9mm>*{}**@{-},
<-0.mm,-5.4mm>*{};<3mm,-9mm>*{}**@{-},
<-0.7mm,-5.4mm>*{};<-4mm,-11mm>*{_{a_1}}**@{},
<-0.7mm,-5.4mm>*{};<4mm,-11mm>*{_{a_2}}**@{},
 \end{xy}\Ea
 & \mbox{for}\ p=2 \vspace{3mm}  \\
0 & \mbox{for}\ p=1.\\
 \Ea
 \right.
\Eeq

{\sc Claim}. {\em The cohomology of the complex $C_v(n_v)$ is concentrated in cohomological degree zero}.

Indeed, consider a one-step filtration of $C_v(n_v)$ by the number
of three-valent vertices of the form $\Ba{c}\resizebox{3mm}{!}{\begin{xy}
 <0mm,0mm>*{\bullet};<0mm,0mm>*{}**@{},
 <0mm,0mm>*{};<-4mm,-5mm>*{}**@{.},
 <0mm,0mm>*{};<0mm,-5mm>*{}**@{-},
 <0mm,0mm>*{};<0mm,5mm>*{}**@{.},
 \end{xy}}\Ea$ and the associated two pages spectral sequence.
 It is easy to see that the complex on the initial page is a direct sum of a trivial complex
 spanned by graphs of the form $\Ba{c}\resizebox{4mm}{!}{\begin{xy}
 <0mm,0mm>*{\bullet};<0mm,0mm>*{}**@{},
 <0mm,0mm>*{};<-4mm,-5mm>*{}**@{.},
 <0mm,0mm>*{};<0mm,-6mm>*{_{a_1}}**@{},
 <0mm,0mm>*{};<0mm,-5mm>*{}**@{-},
 <0mm,0mm>*{};<0mm,5mm>*{}**@{.},
 \end{xy}}\Ea$ with $a_1\in \cA ss_u^{(k+1)}(n_v)$ and a non-trivial complex
 which
 is quasi-isomorphic to the degree shifted (direct summand) subcomplex  of
 $(\cA ss_\infty^{(k+1)}(n_v)[1],\delta)$ spanned by graphs with the orientation of the unique root leg fixed by the multi-orientation of the corresponding dashed edge of the given special path
 (indeed, take a filtration of the latter sub-complex
 by the valency of the root vertex and use the induction assumption). As $n_v\leq n$ we conclude (again by the induction assumption) that  its cohomology is equal to
 $$
A:=  \text{span} \left\langle
 \Ba{c}\resizebox{10mm}{!}{\begin{xy}
 <0mm,0mm>*{\bullet};<0mm,0mm>*{}**@{},
 <0mm,0mm>*{};<-8mm,-5mm>*{}**@{.},
 <0mm,0mm>*{};<-4.5mm,-5mm>*{}**@{-},
 <0mm,0mm>*{};<4.5mm,-5mm>*{}**@{-},
   <0mm,0mm>*{};<-4mm,-7.9mm>*{^{a_1}}**@{},
    <0mm,0mm>*{};<4.5mm,-7.9mm>*{^{a_{2}}}**@{},
 <0mm,0mm>*{};<0mm,5mm>*{}**@{.},
 \end{xy}}\Ea \bmod \cA ss^{(k+1)}\text{-relations},\ a_1\in \cA ss_u^{(k+1)}(n_1),\  a_2\in \cA ss_u^{(k+1)}(n_2),\ n_1+n_2=n_v
 \right\rangle
 $$
 The induced differential on the next (and final) page of the spectral sequence
is an injection
$$
\Ba{rccc}
d: & A  &\lon &  \cA ss^{(k+1)}(n_v)=\text{span}\left\langle \Ba{c}\resizebox{4mm}{!}{\begin{xy}
 <0mm,0mm>*{\bullet};<0mm,0mm>*{}**@{},
 <0mm,0mm>*{};<-4mm,-5mm>*{}**@{.},
 <0mm,0mm>*{};<0mm,-6mm>*{_{a}}**@{},
 <0mm,0mm>*{};<0mm,-5mm>*{}**@{-},
 <0mm,0mm>*{};<0mm,5mm>*{}**@{.},
 \end{xy}}\Ea, \ a\in \cA ss_u^{(k+1)}(n_v) \right\rangle\\
 &  \Ba{c}\resizebox{10mm}{!}{\begin{xy}
 <0mm,0mm>*{\bullet};<0mm,0mm>*{}**@{},
 <0mm,0mm>*{};<-8mm,-5mm>*{}**@{.},
 <0mm,0mm>*{};<-4.5mm,-5mm>*{}**@{-},
 <0mm,0mm>*{};<4.5mm,-5mm>*{}**@{-},
   <0mm,0mm>*{};<-4mm,-7.9mm>*{^{a_1}}**@{},
    <0mm,0mm>*{};<4.5mm,-7.9mm>*{^{a_{2}}}**@{},
 <0mm,0mm>*{};<0mm,5mm>*{}**@{.},
 \end{xy}}\Ea
 &\lon & -\sum_{\bar{\fs}\in \f_{k}}
\Ba{c}\resizebox{8mm}{!}{
\begin{xy}
 <0mm,0mm>*{\bullet};<0mm,0mm>*{}**@{},
 <0mm,0mm>*{};<-4mm,-5mm>*{}**@{.},
 <0mm,0mm>*{};<0mm,-5mm>*{}**@{-},
 <0mm,0mm>*{};<0mm,5mm>*{}**@{.},
 <0mm,0mm>*{};<1.0mm,-3.5mm>*{^{\bar{\fs}}}**@{},
  <0mm,-5.4mm>*{\bullet};
<-0mm,-5.4mm>*{};<-3mm,-9mm>*{}**@{-},
<-0.mm,-5.4mm>*{};<3mm,-9mm>*{}**@{-},
<-0.7mm,-5.4mm>*{};<-4mm,-11mm>*{_{a_1}}**@{},
<-0.7mm,-5.4mm>*{};<4mm,-11mm>*{_{a_2}}**@{},
 \end{xy}}\Ea
 \Ea
$$
which proves the {\sc Claim}.

\sip

We conclude that the cohomology $H^\bu(\cA ss_\infty^{(k+1)}(n+1))$ is generated
by multi-oriented graphs
 of the form (modulo some relations corresponding to the image of the injection $d$)
$$
\Ba{c}
\begin{xy}
 <0mm,-16mm>*{};<0mm,5mm>*{}**@{.},
<0mm,1mm>*{\bullet};
<0mm,1mm>*{};<3mm,-3mm>*{}**@{-},
<5mm,-4.5mm>*{_{a_{v_1}}};
<0mm,-4mm>*{\bullet};
<0mm,-4mm>*{};<3mm,-8mm>*{}**@{-},
<5mm,-9.5mm>*{_{a_{v_2}}};
<3mm,-10mm>*{.};
<3mm,-11mm>*{.};
<3mm,-12mm>*{.};
<0mm,-12mm>*{\bullet};
<0mm,-12mm>*{};<3mm,-16mm>*{}**@{-},
<5mm,-17.5mm>*{_{a_{v_l}}};
 \end{xy}
\Ea\ \ \ \ \ \ \text{where}\ l:=\# V_{sp}, \ \ \ a_{v_i}\in Ass_u^{(k+1)}(n_{v_i}), \ \ \sum_{i=1}^l n_{v_i}=n
$$
which all have cohomological degree zero. Hence $H^\bu(\cA ss_\infty^{(k+1)}(n+1))$ is concentrated in degree zero implying its identification with $\cA ss^{(k+1)}(n+1)$. The induction argument and hence the proof the Theorem are completed.
\end{proof}

In the next subsection we discuss representations of $\cA ss^{(k+1)}$, that is,
{\em associative algebras with $k$ branes}. Rather surprisingly, we recover, in particular, a well-known notion of {\em infinitesimal bialgebra}  as an {\em associative algebra with one (symplectic Lagrangian) brane}. This interesting fact can be seen already now (i.e.\ in purely combinatorial way) as follows.

\subsection{\bf Infinitesimal bialgebras as 2-oriented associative algebras}
Recall that an ordinary (i.e.\ 1-oriented) dioperad of {\em infinitesimal associative bialgebras}\, is, by definition \cite{JR}, the quotient of the  1-oriented free dioperad
$$
\cI \cB:= \cF ree_0^{1\text{-or}}\left\langle B\right\rangle/R
$$
generated
by an $\bS$-bimodule $B=\{B(m,n)\}$
\[
B(m,n):=\left\{
\Ba{rr}
\K[\bS_2]\ot \id_1\equiv\mbox{span}\left\langle
\begin{xy}
 <0mm,-0.55mm>*{};<0mm,-2.5mm>*{}**@{.},
 <0.5mm,0.5mm>*{};<2.2mm,2.2mm>*{}**@{.},
 <-0.48mm,0.48mm>*{};<-2.2mm,2.2mm>*{}**@{.},
 <0mm,0mm>*{\circ};<0mm,0mm>*{}**@{},
 <0mm,-0.55mm>*{};<0mm,-3.8mm>*{_0}**@{},
 <0.5mm,0.5mm>*{};<2.7mm,2.8mm>*{^2}**@{},
 <-0.48mm,0.48mm>*{};<-2.7mm,2.8mm>*{^1}**@{},
 \end{xy}
\,
,\,
\begin{xy}
 <0mm,-0.55mm>*{};<0mm,-2.5mm>*{}**@{.},
 <0.5mm,0.5mm>*{};<2.2mm,2.2mm>*{}**@{.},
 <-0.48mm,0.48mm>*{};<-2.2mm,2.2mm>*{}**@{.},
 <0mm,0mm>*{\circ};<0mm,0mm>*{}**@{},
 <0mm,-0.55mm>*{};<0mm,-3.8mm>*{_0}**@{},
 <0.5mm,0.5mm>*{};<2.7mm,2.8mm>*{^1}**@{},
 <-0.48mm,0.48mm>*{};<-2.7mm,2.8mm>*{^2}**@{},
 \end{xy}
   \right\rangle  & \mbox{if}\ m=2, n=1,\vspace{3mm}\\
\id_1\ot \K[\bS_2]\equiv
\mbox{span}\left\langle
\begin{xy}
 <0mm,0.66mm>*{};<0mm,3mm>*{}**@{.},
 <0.39mm,-0.39mm>*{};<2.2mm,-2.2mm>*{}**@{.},
 <-0.35mm,-0.35mm>*{};<-2.2mm,-2.2mm>*{}**@{.},
 <0mm,0mm>*{\circ};<0mm,0mm>*{}**@{},
   <0mm,0.66mm>*{};<0mm,3.4mm>*{^0}**@{},
   <0.39mm,-0.39mm>*{};<2.9mm,-4mm>*{^2}**@{},
   <-0.35mm,-0.35mm>*{};<-2.8mm,-4mm>*{^1}**@{},
\end{xy}
\,
,\,
\begin{xy}
 <0mm,0.66mm>*{};<0mm,3mm>*{}**@{.},
 <0.39mm,-0.39mm>*{};<2.2mm,-2.2mm>*{}**@{.},
 <-0.35mm,-0.35mm>*{};<-2.2mm,-2.2mm>*{}**@{.},
 <0mm,0mm>*{\circ};<0mm,0mm>*{}**@{},
   <0mm,0.66mm>*{};<0mm,3.4mm>*{^0}**@{},
   <0.39mm,-0.39mm>*{};<2.9mm,-4mm>*{^1}**@{},
   <-0.35mm,-0.35mm>*{};<-2.8mm,-4mm>*{^2}**@{},
\end{xy}
\right\rangle
\ & \mbox{if}\ m=1, n=2, \vspace{3mm}\\
0 & \mbox{otherwise}
\Ea
\right.
\]
by the ideal $R$ generated by the following relations
$$
\Ba{c}
\begin{xy}
 <0mm,0mm>*{\circ};<0mm,0mm>*{}**@{},
 <0mm,-0.49mm>*{};<0mm,-3.0mm>*{}**@{.},
 <0.49mm,0.49mm>*{};<1.9mm,1.9mm>*{}**@{.},
 <-0.5mm,0.5mm>*{};<-1.9mm,1.9mm>*{}**@{.},
 <-2.3mm,2.3mm>*{\circ};<-2.3mm,2.3mm>*{}**@{},
 <-1.8mm,2.8mm>*{};<0mm,4.9mm>*{}**@{.},
 <-2.8mm,2.9mm>*{};<-4.6mm,4.9mm>*{}**@{.},
   <0.49mm,0.49mm>*{};<2.7mm,2.3mm>*{^3}**@{},
   <-1.8mm,2.8mm>*{};<0.4mm,5.3mm>*{^2}**@{},
   <-2.8mm,2.9mm>*{};<-5.1mm,5.3mm>*{^1}**@{},
 \end{xy}\Ea
\ - \
\Ba{c}
\begin{xy}
 <0mm,0mm>*{\circ};<0mm,0mm>*{}**@{},
 <0mm,-0.49mm>*{};<0mm,-3.0mm>*{}**@{.},
 <0.49mm,0.49mm>*{};<1.9mm,1.9mm>*{}**@{.},
 <-0.5mm,0.5mm>*{};<-1.9mm,1.9mm>*{}**@{.},
 <2.3mm,2.3mm>*{\circ};<-2.3mm,2.3mm>*{}**@{},
 <1.8mm,2.8mm>*{};<0mm,4.9mm>*{}**@{.},
 <2.8mm,2.9mm>*{};<4.6mm,4.9mm>*{}**@{.},
   <0.49mm,0.49mm>*{};<-2.7mm,2.3mm>*{^1}**@{},
   <-1.8mm,2.8mm>*{};<0mm,5.3mm>*{^2}**@{},
   <-2.8mm,2.9mm>*{};<5.1mm,5.3mm>*{^3}**@{},
 \end{xy}\Ea=0, \ \ \ \ \
 \Ba{c}\begin{xy}
 <0mm,0mm>*{\circ};<0mm,0mm>*{}**@{},
 <0mm,0.69mm>*{};<0mm,3.0mm>*{}**@{.},
 <0.39mm,-0.39mm>*{};<2.4mm,-2.4mm>*{}**@{.},
 <-0.35mm,-0.35mm>*{};<-1.9mm,-1.9mm>*{}**@{.},
 <-2.4mm,-2.4mm>*{\circ};<-2.4mm,-2.4mm>*{}**@{},
 <-2.0mm,-2.8mm>*{};<0mm,-4.9mm>*{}**@{.},
 <-2.8mm,-2.9mm>*{};<-4.7mm,-4.9mm>*{}**@{.},
    <0.39mm,-0.39mm>*{};<3.3mm,-4.0mm>*{^3}**@{},
    <-2.0mm,-2.8mm>*{};<0.5mm,-6.7mm>*{^2}**@{},
    <-2.8mm,-2.9mm>*{};<-5.2mm,-6.7mm>*{^1}**@{},
 \end{xy}\Ea
\ - \
 \Ba{c}\begin{xy}
 <0mm,0mm>*{\circ};<0mm,0mm>*{}**@{},
 <0mm,0.69mm>*{};<0mm,3.0mm>*{}**@{.},
 <0.39mm,-0.39mm>*{};<2.4mm,-2.4mm>*{}**@{.},
 <-0.35mm,-0.35mm>*{};<-1.9mm,-1.9mm>*{}**@{.},
 <2.4mm,-2.4mm>*{\circ};<-2.4mm,-2.4mm>*{}**@{},
 <2.0mm,-2.8mm>*{};<0mm,-4.9mm>*{}**@{.},
 <2.8mm,-2.9mm>*{};<4.7mm,-4.9mm>*{}**@{.},
    <0.39mm,-0.39mm>*{};<-3mm,-4.0mm>*{^1}**@{},
    <-2.0mm,-2.8mm>*{};<0mm,-6.7mm>*{^2}**@{},
    <-2.8mm,-2.9mm>*{};<5.2mm,-6.7mm>*{^3}**@{},
 \end{xy}\Ea=0,\ \ \ \ \ \
 \ \ \
 \Ba{c} \begin{xy}
 <0mm,2.47mm>*{};<0mm,-0.5mm>*{}**@{.},
 <0.5mm,3.5mm>*{};<2.2mm,5.2mm>*{}**@{.},
 <-0.48mm,3.48mm>*{};<-2.2mm,5.2mm>*{}**@{.},
 <0mm,3mm>*{\circ};<0mm,3mm>*{}**@{},
  <0mm,-0.8mm>*{\circ};<0mm,-0.8mm>*{}**@{},
<0mm,-0.8mm>*{};<-2.2mm,-3.5mm>*{}**@{.},
 <0mm,-0.8mm>*{};<2.2mm,-3.5mm>*{}**@{.},
     <0.5mm,3.5mm>*{};<2.8mm,5.7mm>*{^2}**@{},
     <-0.48mm,3.48mm>*{};<-2.8mm,5.7mm>*{^1}**@{},
   <0mm,-0.8mm>*{};<-2.7mm,-5.2mm>*{^3}**@{},
   <0mm,-0.8mm>*{};<2.7mm,-5.2mm>*{^4}**@{},
\end{xy}\Ea
\  - \
\begin{xy}
<0.39mm,-0.39mm>*{};<-3mm,3.0mm>*{^1}**@{},
<0.39mm,-0.39mm>*{};<2.4mm,6.4mm>*{^2}**@{},
    <-2.0mm,-2.8mm>*{};<0mm,-6.1mm>*{^3}**@{},
    <-2.8mm,-2.9mm>*{};<5.2mm,-2.6mm>*{^4}**@{},
 <0mm,-1.3mm>*{};<0mm,-3.5mm>*{}**@{.},
 <0.38mm,-0.2mm>*{};<2.2mm,2.2mm>*{}**@{.},
 <-0.38mm,-0.2mm>*{};<-2.2mm,2.2mm>*{}**@{.},
<0mm,-0.8mm>*{\circ};<0mm,0.8mm>*{}**@{},
 <2.4mm,2.4mm>*{\circ};<2.4mm,2.4mm>*{}**@{},
 <2.5mm,2.3mm>*{};<4.4mm,-0.8mm>*{}**@{.},
 <2.4mm,2.5mm>*{};<2.4mm,5.2mm>*{}**@{.},
    \end{xy}
\  - \
\begin{xy}
<0.39mm,-0.39mm>*{};<3mm,3.0mm>*{^2}**@{},
<0.39mm,-0.39mm>*{};<-2.4mm,6.4mm>*{^1}**@{},
    <-2.0mm,-2.8mm>*{};<0mm,-6.1mm>*{^4}**@{},
    <-2.8mm,-2.9mm>*{};<-5.2mm,-2.6mm>*{^3}**@{},
 <0mm,-1.3mm>*{};<0mm,-3.5mm>*{}**@{.},
 <0.38mm,-0.2mm>*{};<2.2mm,2.2mm>*{}**@{.},
 <-0.38mm,-0.2mm>*{};<-2.2mm,2.2mm>*{}**@{.},
<0mm,-0.8mm>*{\circ};<0mm,0.8mm>*{}**@{},
 <-2.4mm,2.4mm>*{\circ};<2.4mm,2.4mm>*{}**@{},
 <-2.5mm,2.3mm>*{};<-4.4mm,-0.8mm>*{}**@{.},
 <-2.4mm,2.5mm>*{};<-2.4mm,5.2mm>*{}**@{.},
    \end{xy} \ =\ 0.
$$
Here  all internal edges and legs are assumed to be oriented along the flow running from the bottom of a graph to its top.

\subsubsection{{\bf Proposition (cf.\ \cite{Wa})}}\label{3new: Prop on map from ASS^2 to IB}
{\em There is a (forgetting the basic orientation) morphism of dioperads
$$
\al: \cA ss^{(2)} \lon \cI\cB
$$
given on the generators as follows:
$$
\al\left(\Ba{c}\resizebox{11mm}{!}  {
\begin{tikzpicture}[baseline=-1ex]
\node[] at (0,0.8) {$^{0}$};
\node[] at (-0.5,-0.71) {$_{1}$};
\node[] at (0.5,-0.71) {$_{2}$};
\node[int] (a) at (0,0) {};
\node[] (u1) at (0,0.8) {};
\node[int] (a) at (0,0) {};
\node[] (d1) at (-0.5,-0.7) {};
\node[] (d2) at (0.5,-0.7) {};
\draw (u1) edge[latex-, leftblue] (a);
\draw (a) edge[latex-, leftblue] (d1);
\draw (a) edge[latex-, leftblue] (d2);
\end{tikzpicture}}
\Ea
\right):=
\begin{xy}
 <0mm,0.66mm>*{};<0mm,3mm>*{}**@{.},
 <0.39mm,-0.39mm>*{};<2.2mm,-2.2mm>*{}**@{.},
 <-0.35mm,-0.35mm>*{};<-2.2mm,-2.2mm>*{}**@{.},
 <0mm,0mm>*{\circ};<0mm,0mm>*{}**@{},
   <0mm,0.66mm>*{};<0mm,3.4mm>*{^{^0}}**@{},
   <0.39mm,-0.39mm>*{};<2.9mm,-4mm>*{_{_2}}**@{},
   <-0.35mm,-0.35mm>*{};<-2.8mm,-4mm>*{_{_1}}**@{},
\end{xy}\
,\ \ \ \
\al\left(\Ba{c}\resizebox{11mm}{!}  {
\begin{tikzpicture}[baseline=-1ex]
\node[] at (0,0.8) {$^{0}$};
\node[] at (-0.5,-0.71) {$_{1}$};
\node[] at (0.5,-0.71) {$_{2}$};
\node[int] (a) at (0,0) {};
\node[] (u1) at (0,0.8) {};
\node[int] (a) at (0,0) {};
\node[] (d1) at (-0.5,-0.7) {};
\node[] (d2) at (0.5,-0.7) {};
\draw (u1) edge[latex-, rightblue] (a);
\draw (a) edge[latex-, rightblue] (d1);
\draw (a) edge[latex-, rightblue] (d2);
\end{tikzpicture}
}\Ea\right)
=
\begin{xy}
 <0mm,-0.55mm>*{};<0mm,-2.5mm>*{}**@{.},
 <0.5mm,0.5mm>*{};<2.2mm,2.2mm>*{}**@{.},
 <-0.48mm,0.48mm>*{};<-2.2mm,2.2mm>*{}**@{.},
 <0mm,0mm>*{\circ};<0mm,0mm>*{}**@{},
 <0mm,-0.55mm>*{};<0mm,-3.8mm>*{_{_0}}**@{},
 <0.5mm,0.5mm>*{};<2.7mm,2.8mm>*{^{^1}}**@{},
 <-0.48mm,0.48mm>*{};<-2.7mm,2.8mm>*{^{^2}}**@{},
 \end{xy}
,\ \ \
\al\left(\Ba{c}\resizebox{11mm}{!}  {
\begin{tikzpicture}[baseline=-1ex]
\node[] at (0,0.8) {$^{0}$};
\node[] at (-0.5,-0.71) {$_{1}$};
\node[] at (0.5,-0.71) {$_{2}$};
\node[int] (a) at (0,0) {};
\node[] (u1) at (0,0.8) {};
\node[int] (a) at (0,0) {};
\node[] (d1) at (-0.5,-0.7) {};
\node[] (d2) at (0.5,-0.7) {};
\draw (u1) edge[latex-, leftblue] (a);
\draw (a) edge[latex-, leftblue] (d1);
\draw (a) edge[latex-, rightblue] (d2);
\end{tikzpicture}
}\Ea\right)=
\begin{xy}
 <0mm,-0.55mm>*{};<0mm,-2.5mm>*{}**@{.},
 <0.5mm,0.5mm>*{};<2.2mm,2.2mm>*{}**@{.},
 <-0.48mm,0.48mm>*{};<-2.2mm,2.2mm>*{}**@{.},
 <0mm,0mm>*{\circ};<0mm,0mm>*{}**@{},
 <0mm,-0.55mm>*{};<0mm,-3.8mm>*{_{_1}}**@{},
 <0.5mm,0.5mm>*{};<2.7mm,2.8mm>*{^{^2}}**@{},
 <-0.48mm,0.48mm>*{};<-2.7mm,2.8mm>*{^{^0}}**@{},
 \end{xy} \
,\ \ \
\al\left(\Ba{c}\resizebox{11mm}{!}  {
\begin{tikzpicture}[baseline=-1ex]
\node[] at (0,0.8) {$^{0}$};
\node[] at (-0.5,-0.71) {$_{1}$};
\node[] at (0.5,-0.71) {$_{2}$};
\node[int] (a) at (0,0) {};
\node[] (u1) at (0,0.8) {};
\node[int] (a) at (0,0) {};
\node[] (d1) at (-0.5,-0.7) {};
\node[] (d2) at (0.5,-0.7) {};
\draw (u1) edge[latex-, rightblue] (a);
\draw (a) edge[latex-, leftblue] (d1);
\draw (a) edge[latex-, rightblue] (d2);
\end{tikzpicture}
}\Ea\right)
=
\begin{xy}
 <0mm,0.66mm>*{};<0mm,3mm>*{}**@{.},
 <0.39mm,-0.39mm>*{};<2.2mm,-2.2mm>*{}**@{.},
 <-0.35mm,-0.35mm>*{};<-2.2mm,-2.2mm>*{}**@{.},
 <0mm,0mm>*{\circ};<0mm,0mm>*{}**@{},
   <0mm,0.66mm>*{};<0mm,3.4mm>*{^{^2}}**@{},
   <0.39mm,-0.39mm>*{};<2.9mm,-4mm>*{_{_1}}**@{},
   <-0.35mm,-0.35mm>*{};<-2.8mm,-4mm>*{_{_0}}**@{},
\end{xy}\
$$
}
\begin{proof}
It is straightforward to check that each of the six relations in (\ref{3new: Ass^2 relations 1st set})-(\ref{3new: Ass^2 relations 2nd set}) is mapped under $\al$ into one of the above three relations in $\cI\cB$. Hence the map is well-defined indeed.
\end{proof}

This proposition indicates that the notion of {\em representation}\, of a multi-oriented prop(erad) can not be an immediate  generalization of that notion for ordinary coloured prop(erad)s --- the extra orientations are not really ``colours", and the distinction between operads and dioperads should become
void.

\subsection{Example: Multi-oriented operad of Lie and $\caL ie_\infty$ algebras} Recall that the ordinary operad of strongly homotopy Lie algebras is the free operad
$\caL ie_\infty:=(\cF\! ree^{1\text{-or}}\langle L\rangle, \delta)$ generated by an $\bS$-module
$L=\{L(n)\}_{\geq 2}$ with
$$
L(n):=\sgn_n[n-2]= \mbox{span}
\left(
\Ba{c}\resizebox{19mm}{!}
{\begin{tikzpicture}[baseline=-1ex]
\node[] at (-1.2,-0.9) {$_{\scriptstyle 1}$};
\node[] at (-0.6,-0.9) {$_{\scriptstyle 2}$};
\node[] at (-0.0,-0.5) {$...$};
\node[] at (1.3,-0.9) {$_{\scriptstyle n}$};
\node[] at (0.6,-0.9) {$_{\scriptstyle n-1}$};
\node[int] (0) at (0,0) {};
\node[] (1) at (0,0.8) {};
\node[int] (0) at (0,0) {};
\node[] (5) at (-1.2,-0.8) {};
\node[] (6) at (-0.6,-0.8) {};
\node[] (7) at (0.6,-0.8) {};
\node[] (8) at (1.2,-0.8) {};
\draw (1) edge[latex-] (0);
\draw (0) edge[latex-] (5);
\draw (0) edge[latex-] (6);
\draw (0) edge[latex-] (7);
\draw (0) edge[latex-] (8);
\end{tikzpicture}}\Ea
=(-1)^\sigma
\Ba{c}\resizebox{20mm}{!}
{\begin{tikzpicture}[baseline=-1ex]
\node[] at (-1.2,-0.9) {$_{\scriptstyle \sigma(1)}$};
\node[] at (-0.6,-0.9) {$_{\scriptstyle \sigma(2)}$};
\node[] at (-0.0,-0.5) {$...$};
\node[] at (1.3,-0.9) {$_{\scriptstyle \sigma(n)}$};
\node[] at (0.6,-0.9) {$_{\scriptstyle \sigma(n-1)}$};
\node[int] (0) at (0,0) {};
\node[] (1) at (0,0.8) {};
\node[int] (0) at (0,0) {};
\node[] (5) at (-1.2,-0.8) {};
\node[] (6) at (-0.6,-0.8) {};
\node[] (7) at (0.6,-0.8) {};
\node[] (8) at (1.2,-0.8) {};
\draw (1) edge[latex-] (0);
\draw (0) edge[latex-] (5);
\draw (0) edge[latex-] (6);
\draw (0) edge[latex-] (7);
\draw (0) edge[latex-] (8);
\end{tikzpicture}}\Ea, {\forall\ \sigma\in \bS_n}
\right)
$$
where $\sgn_n$ is the 1-dimensional sign representation of $\bS_n$.
The differential is given on the generators by
$$
\delta \hspace{-5mm}
\Ba{c}\resizebox{19mm}{!}
{\begin{tikzpicture}[baseline=-1ex]
\node[] at (-0.0,-1.0) {$_{\scriptstyle I}$};
\node[] at (-0.0,-0.8) {$\underbrace{\ \ \ \ \ \ \ \ \ \ \ \ \ \ \ \ \ \  \ \ \ \ \ \ \ \ \ \ }$};
\node[] at (-0.0,-0.5) {$...$};
\node[int] (0) at (0,0) {};
\node[] (1) at (0,0.8) {};
\node[int] (0) at (0,0) {};
\node[] (5) at (-1.2,-0.8) {};
\node[] (6) at (-0.5,-0.8) {};
\node[] (7) at (0.5,-0.8) {};
\node[] (8) at (1.2,-0.8) {};
\draw (1) edge[latex-] (0);
\draw (0) edge[latex-] (5);
\draw (0) edge[latex-] (6);
\draw (0) edge[latex-] (7);
\draw (0) edge[latex-] (8);
\end{tikzpicture}}\Ea
=
\sum_{I=I_1\sqcup I_2}
(-1)^{\#I_2 + 1}
\Ba{c}\resizebox{28mm}{!}
{\begin{tikzpicture}[baseline=-1ex]
\node[] at (-1.7,-0.9) {$_{\underbrace{\ \ \ \ \ \ \ \ \ \ \ \ \ \ \ \ }_{I_1}}$};
\node[] at (0.0,-0.15) {$_{\underbrace{\ \ \ \ \ \ \ \ \ \ \ \ \ \ \ \ \ \ \ \ \ }_{I_2}}$};
\node[] at (-2.2,-0.65) {$...$};
\node[] at (-1.16,-0.65) {$...$};
\node[] at (0.6,0.2) {$...$};
\node[] at (-0.55,0.2) {$...$};
%
\node[int] (0) at (-1.7,0) {};
\node[int] (1) at (0,0.8) {};
\node[] (u) at (0,1.6) {};
\node[] (r2) at (-0.4,0) {};
\node[] (r1) at (-1.1,0) {};
\node[] (R2) at (1.2,0) {};
\node[] (R1) at (0.4,0) {};
%
\node[] (5) at (-2.6,-0.8) {};
\node[] (6) at (-2.1,-0.8) {};
\node[] (7) at (-1.3,-0.8) {};
\node[] (8) at (-0.8,-0.8) {};
\draw (u) edge[latex-] (1);
\draw (1) edge[latex-] (r1);
\draw (1) edge[latex-] (r2);
\draw (1) edge[latex-] (R1);
\draw (1) edge[latex-] (R2);
\draw (1) edge[latex-] (0);
\draw (0) edge[latex-] (5);
\draw (0) edge[latex-] (6);
\draw (0) edge[latex-] (7);
\draw (0) edge[latex-] (8);
\end{tikzpicture}}\Ea
$$
It is essentially a skewsymmetrized version of $\cA ss_\infty$ (there is a canonical morphism of operads $\caL ie_\infty \rar \cA ss_\infty$ sending a generator of $\caL ie_\infty$ into a skewsymmetrization of the corresponding generator of $\cA ss_\infty$).
If $I$ is the differential closure of the ideal in $\caL ie_\infty$ generated
by all corollas with negative cohomological degree, then the quotient
$$
\caL ie:= \caL ie_\infty/I
$$
is an operad controlling Lie algebras. It is generated by degree zero planar skewsymmetric corollas
$$
\Ba{c}\resizebox{8mm}{!}  {\begin{tikzpicture}[baseline=-1ex]
\node[] at (0,0.8) {$^{0}$};
\node[] at (-0.5,-0.71) {$_{1}$};
\node[] at (0.5,-0.71) {$_{2}$};
\node[int] (a) at (0,0) {};
\node[] (u1) at (0,0.8) {};
\node[int] (a) at (0,0) {};
\node[] (d1) at (-0.5,-0.7) {};
\node[] (d2) at (0.5,-0.7) {};
\draw (u1) edge[latex-] (a);
\draw (a) edge[latex-] (d1);
\draw (a) edge[latex-] (d2);
\end{tikzpicture}}\Ea
=-
\Ba{c}\resizebox{8mm}{!}  {\begin{tikzpicture}[baseline=-1ex]
\node[] at (0,0.8) {$^{0}$};
\node[] at (-0.5,-0.71) {$_{1}$};
\node[] at (0.5,-0.71) {$_{2}$};
\node[int] (a) at (0,0) {};
\node[] (u1) at (0,0.8) {};
\node[int] (a) at (0,0) {};
\node[] (d1) at (-0.5,-0.7) {};
\node[] (d2) at (0.5,-0.7) {};
\draw (u1) edge[latex-] (a);
\draw (a) edge[latex-] (d1);
\draw (a) edge[latex-] (d2);
\end{tikzpicture}}\Ea
$$
modulo the Jacobi relation,
$$
\Ba{c}\resizebox{10mm}{!}  {
\begin{tikzpicture}[baseline=-1ex]
\node[] at (0.5,1.6) {$^{0}$};
\node[] at (1.55,-0.72) {$_{3}$};
\node[] at (0.57,-0.72) {$_{2}$};
\node[] at (0.0,-0.02) {$_{1}$};
\node[] (u2) at (0.5,1.6) {};
\node[] (ud2) at (0,0.0) {};
\node[int] (u1) at (0.5,0.8) {};
\node[int] (a) at (1.0,0) {};
\node[] (d1) at (0.5,-0.7) {};
\node[] (d2) at (1.5,-0.7) {};
\draw (u1) edge[latex-] (a);
\draw (u2) edge[latex-] (u1);
\draw (u1) edge[latex-] (ud2);
\draw (a) edge[latex-] (d1);
\draw (a) edge[latex-] (d2);
\end{tikzpicture}
}\Ea
+
\Ba{c}\resizebox{10mm}{!}  {
\begin{tikzpicture}[baseline=-1ex]
\node[] at (0.5,1.6) {$^{0}$};
\node[] at (1.55,-0.72) {$_{2}$};
\node[] at (0.57,-0.72) {$_{1}$};
\node[] at (0.0,-0.02) {$_{3}$};
\node[] (u2) at (0.5,1.6) {};
\node[] (ud2) at (0,0.0) {};
\node[int] (u1) at (0.5,0.8) {};
\node[int] (a) at (1.0,0) {};
\node[] (d1) at (0.5,-0.7) {};
\node[] (d2) at (1.5,-0.7) {};
\draw (u1) edge[latex-] (a);
\draw (u2) edge[latex-] (u1);
\draw (u1) edge[latex-] (ud2);
\draw (a) edge[latex-] (d1);
\draw (a) edge[latex-] (d2);
\end{tikzpicture}
}\Ea
+
\Ba{c}\resizebox{10mm}{!}  {
\begin{tikzpicture}[baseline=-1ex]
\node[] at (0.5,1.6) {$^{0}$};
\node[] at (1.55,-0.72) {$_{1}$};
\node[] at (0.57,-0.72) {$_{3}$};
\node[] at (0.0,-0.02) {$_{2}$};
%
\node[] (u2) at (0.5,1.6) {};
\node[] (ud2) at (0,0.0) {};
\node[int] (u1) at (0.5,0.8) {};
\node[int] (a) at (1.0,0) {};
\node[] (d1) at (0.5,-0.7) {};
\node[] (d2) at (1.5,-0.7) {};
\draw (u1) edge[latex-] (a);
\draw (u2) edge[latex-] (u1);
\draw (u1) edge[latex-] (ud2);
\draw (a) edge[latex-] (d1);
\draw (a) edge[latex-] (d2);
\end{tikzpicture}
}\Ea
=0
$$
The natural surjection
$
\caL ie_\infty \lon \caL ie
$
is a quasi-isomorphism.

\sip

An operad of $(k+1)$-oriented strongly homotopy Lie algebras is defined  as an obvious skew-symmetrization of the operad $\cA ss_\infty^{(k+1)}$ introduced in the previous subsection (so that there is again a canonical morphism of dg operads
$\caL ie_\infty^{(k+1)}\rar \cA ss_\infty^{(k+1)}$).
More precisely, the operad $\caL ie_\infty^{(k+1)}$ is a free $(k+1)$-oriented
 operad generated by corollas with the same symmetries and degrees as in the case of  $\caL ie_\infty$, but now with each leg decorated
with extra $k$-orientations,
$$
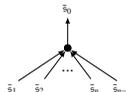
\Ba{c}\resizebox{20mm}{!}
{\begin{tikzpicture}[baseline=-1ex]
\node[] at (0,0.8) {$^{\bar{\fs}_0}$};
\node[] at (-1.2,-0.9) {$_{\bar{\fs}_{1}}$};
\node[] at (-0.6,-0.9) {$_{\bar{\fs}_{2}}$};
\node[] at (1.2,-0.9) {$_{\bar{\fs}_{n-1}}$};
\node[] at (0.6,-0.9) {$_{\bar{\fs}_{n}}$};
\node[] at (0.0,-0.5) {$...$};
\node[int] (0) at (0,0) {};
\node[] (1) at (0,0.8) {};
\node[int] (0) at (0,0) {};
\node[] (5) at (-1.2,-0.8) {};
\node[] (6) at (-0.6,-0.8) {};
\node[] (7) at (0.6,-0.8) {};
\node[] (8) at (1.2,-0.8) {};
\draw (1) edge[latex-] (0);
\draw (0) edge[latex-] (5);
\draw (0) edge[latex-] (6);
\draw (0) edge[latex-] (7);
\draw (0) edge[latex-] (8);
\end{tikzpicture}}\Ea
$$
with the condition that there is at least one ingoing edge and outgoing edge with respect to each of the new directions. The differential is given by the  same formula as in the case of $\caL ie_\infty$ except that now we sum over all possible (and admissible) new orientations attached to the new edge
$$
\delta
\Ba{c}\resizebox{20mm}{!}
{\begin{tikzpicture}[baseline=-1ex]
\node[] at (0,0.8) {$^{\bar{\fs}_0}$};
\node[] at (-1.2,-0.9) {$_{\bar{\fs}_{1}}$};
\node[] at (-0.6,-0.9) {$_{\bar{\fs}_{2}}$};
\node[] at (1.2,-0.9) {$_{\bar{\fs}_{n-1}}$};
\node[] at (0.6,-0.9) {$_{\bar{\fs}_{n}}$};
\node[] at (0.0,-0.5) {$...$};
\node[int] (0) at (0,0) {};
\node[] (1) at (0,0.8) {};
\node[int] (0) at (0,0) {};
\node[] (5) at (-1.2,-0.8) {};
\node[] (6) at (-0.6,-0.8) {};
\node[] (7) at (0.6,-0.8) {};
\node[] (8) at (1.2,-0.8) {};
\draw (1) edge[latex-] (0);
\draw (0) edge[latex-] (5);
\draw (0) edge[latex-] (6);
\draw (0) edge[latex-] (7);
\draw (0) edge[latex-] (8);
\end{tikzpicture}}\Ea
=
 \sum_{[1,\ldots,n]=J_1\sqcup J_2\atop
 {|J_1|\geq 1, |J_2|\geq 1}
} \sum_{\bar{\fs}\in \f_k}
(-1)^{1 + \# J_2  + \cdot \text{sgn}(J_1,J_2)  }
\Ba{c}\resizebox{29mm}{!}
{\begin{tikzpicture}[baseline=-1ex]
\node[] at (1.0,0.74) {$^{ \scriptstyle \bar{\fs}}$};
\node[] at (-0.0,-1.2) {$_{\scriptstyle \bar{\fs}_i,\ i\in J_1}$};
\node[] at (-0.0,-0.9) {$\underbrace{\ \ \ \ \ \ \ \ \ \ \ \ \ \ \ \ \ \ \ \ \  \ \ \ \ \ }$};
\node[] at (1.8,1.9) {$^{\scriptstyle \bar{\fs}_0}$};
\node[] at (1.8,0.0) {$_{\scriptstyle \bar{\fs}_i,\ i\in J_2}$};
\node[] at (1.8,0.27) {$\underbrace{\ \ \ \ \ \ \ \ \ \ \ \ \ \ \ \ \ \ \  }$};
\node[] at (0.0,-0.5) {$...$};
\node[] at (1.8,0.5) {$...$};
\node[int] (0) at (0,0) {};
\node[int] (4) at (1.8,1.1) {};
\node[int] (0) at (0,0) {};
\node[] (5) at (-1.2,-0.8) {};
\node[] (6) at (-0.6,-0.8) {};
\node[] (7) at (0.6,-0.8) {};
\node[] (8) at (1.2,-0.8) {};
\node[] (u1) at (1.8,1.9) {};
\node[] (u2) at (1.4,1.9) {};
\node[] (u3) at (2.1,1.9) {};
\node[] (u4) at (2.6,1.9) {};
\node[] (d1) at (1.0,0.3) {};
\node[] (d2) at (1.4,0.3) {};
\node[] (d3) at (2.1,0.3) {};
\node[] (d4) at (2.6,0.3) {};
\draw (4) edge[latex-] (0);
\draw (0) edge[latex-] (5);
\draw (0) edge[latex-] (6);
\draw (0) edge[latex-] (7);
\draw (0) edge[latex-] (8);
\draw (u1) edge[latex-] (4);
\draw (4) edge[latex-] (d1);
\draw (4) edge[latex-] (d2);
\draw (4) edge[latex-] (d3);
\draw (4) edge[latex-] (d4);
\end{tikzpicture}}\Ea
$$
where the first sum runs over decompositions of the ordered set
$[n]$ into the disjoint union of (not necessarily connected)  ordered subsets, and $\text{sign}(I_1,I_2)$ stands for the parity of the permutation $[n]\rar I_1\sqcup I_2$.

\sip

In more detail,
the operad of $2$-oriented strongly homotopy Lie algebras
is, by definition, a free 2-oriented operad generated by the following skewsymmetric planar corollas
of degree $2-n$, $n\geq 2$,
$$
\underbracket{
\Ba{c}\resizebox{22mm}{!}
{\begin{tikzpicture}[baseline=-1ex]
\node[] at (0,0.8) {$^{\scriptstyle 0}$};
\node[] at (-1.2,-0.9) {$_{\scriptstyle \sigma(1)}$};
\node[] at (-0.3,-0.9) {$_{\scriptstyle \sigma(r)}$};
\node[] at (-0.5,-0.6) {$...$};
\node[] at (0.5,-0.6) {$...$};
\node[] at (1.2,-0.9) {$_{\scriptstyle \tau(r+1)}$};
\node[] at (0.4,-0.9) {$_{\scriptstyle \tau(r+l)}$};
\node[int] (0) at (0,0) {};
\node[] (1) at (0,0.8) {};
\node[int] (0) at (0,0) {};
\node[] (5) at (-1.2,-0.8) {};
\node[] (6) at (-0.3,-0.8) {};
\node[] (7) at (0.3,-0.8) {};
\node[] (8) at (1.2,-0.8) {};
\draw (1) edge[latex-,leftblue] (0);
\draw (0) edge[latex-,leftblue] (5);
\draw (0) edge[latex-,leftblue] (6);
\draw (0) edge[latex-,rightblue] (7);
\draw (0) edge[latex-,rightblue] (8);
\end{tikzpicture}}\Ea
\hspace{-4mm} =
(-1)^{\sigma+\tau}\hspace{-4mm}
\Ba{c}\resizebox{22mm}{!}
{\begin{tikzpicture}[baseline=-1ex]
\node[] at (0,0.8) {$^{\scriptstyle 0}$};
\node[] at (-1.2,-0.9) {$_{\scriptstyle 1}$};
\node[] at (-0.3,-0.9) {$_{\scriptstyle r}$};
\node[] at (-0.5,-0.6) {$...$};
\node[] at (0.5,-0.6) {$...$};
\node[] at (1.2,-0.9) {$_{\scriptstyle r+1}$};
\node[] at (0.4,-0.9) {$_{\scriptstyle r+l}$};
\node[int] (0) at (0,0) {};
\node[] (1) at (0,0.8) {};
\node[int] (0) at (0,0) {};
\node[] (5) at (-1.2,-0.8) {};
\node[] (6) at (-0.3,-0.8) {};
\node[] (7) at (0.3,-0.8) {};
\node[] (8) at (1.2,-0.8) {};
\draw (1) edge[latex-,leftblue] (0);
\draw (0) edge[latex-,leftblue] (5);
\draw (0) edge[latex-,leftblue] (6);
\draw (0) edge[latex-,rightblue] (7);
\draw (0) edge[latex-,rightblue] (8);
\end{tikzpicture}}\Ea}_{
r+l\geq 2,\ r\geq 1,\ l\geq 0}\ , \
\underbracket{
\Ba{c}\resizebox{22mm}{!}
{\begin{tikzpicture}[baseline=-1ex]
\node[] at (0,0.8) {$^{\scriptstyle 0}$};
\node[] at (-1.2,-0.9) {$_{\scriptstyle \sigma(1)}$};
\node[] at (-0.3,-0.9) {$_{\scriptstyle \sigma(r)}$};
\node[] at (-0.5,-0.6) {$...$};
\node[] at (0.5,-0.6) {$...$};
\node[] at (1.2,-0.9) {$_{\scriptstyle \tau(r+1)}$};
\node[] at (0.4,-0.9) {$_{\scriptstyle \tau(r+l)}$};
\node[int] (0) at (0,0) {};
\node[] (1) at (0,0.8) {};
\node[int] (0) at (0,0) {};
\node[] (5) at (-1.2,-0.8) {};
\node[] (6) at (-0.3,-0.8) {};
\node[] (7) at (0.3,-0.8) {};
\node[] (8) at (1.2,-0.8) {};
\draw (1) edge[latex-,rightblue] (0);
\draw (0) edge[latex-,leftblue] (5);
\draw (0) edge[latex-,leftblue] (6);
\draw (0) edge[latex-,rightblue] (7);
\draw (0) edge[latex-,rightblue] (8);
\end{tikzpicture}}\Ea\hspace{-4mm}
=
(-1)^{\sigma+\tau}\hspace{-4mm}
\Ba{c}\resizebox{22mm}{!}
{\begin{tikzpicture}[baseline=-1ex]
\node[] at (0,0.8) {$^{\scriptstyle 0}$};
\node[] at (-1.2,-0.9) {$_{\scriptstyle 1}$};
\node[] at (-0.3,-0.9) {$_{\scriptstyle r}$};
\node[] at (-0.5,-0.6) {$...$};
\node[] at (0.5,-0.6) {$...$};
\node[] at (1.2,-0.9) {$_{\scriptstyle r+1}$};
\node[] at (0.4,-0.9) {$_{\scriptstyle r+l}$};
\node[int] (0) at (0,0) {};
\node[] (1) at (0,0.8) {};
\node[int] (0) at (0,0) {};
\node[] (5) at (-1.2,-0.8) {};
\node[] (6) at (-0.3,-0.8) {};
\node[] (7) at (0.3,-0.8) {};
\node[] (8) at (1.2,-0.8) {};
\draw (1) edge[latex-,rightblue] (0);
\draw (0) edge[latex-,leftblue] (5);
\draw (0) edge[latex-,leftblue] (6);
\draw (0) edge[latex-,rightblue] (7);
\draw (0) edge[latex-,rightblue] (8);
\end{tikzpicture}}\Ea}_{
r+l\geq 2,\ r\geq 0,\ l\geq 1},
$$
for any $\sigma\in \bS_r$ and $\tau\in\bS_l$. As in the case of $\cA ss_\infty^{(k+1)}$  we require that each corolla has at least one ingoing leg and at least one outgoing leg in each extra orientation (in order to avoid curvature terms in representations).

\sip

The differential is given by (and it is easy to check that $\delta$ is a differential indeed)
$$
\delta \hspace{-2mm}
\Ba{c}\resizebox{22mm}{!}
{\begin{tikzpicture}[baseline=-1ex]
\node[] at (-0.6,-1.0) {$_{\scriptstyle I}$};
\node[] at (-0.7,-0.8) {$\underbrace{\  \ \ \ \ \ \ \ \ \ \ }$};
\node[] at (-0.45,-0.5) {$...$};
\node[] at (0.6,-1.0) {$_{\scriptstyle J}$};
\node[] at (0.7,-0.8) {$\underbrace{\  \ \ \ \ \ \ \ \ \ \ }$};
\node[] at (0.45,-0.5) {$...$};
\node[int] (0) at (0,0) {};
\node[] (1) at (0,0.8) {};
\node[int] (0) at (0,0) {};
\node[] (5) at (-1.2,-0.8) {};
\node[] (6) at (-0.3,-0.8) {};
\node[] (7) at (0.3,-0.8) {};
\node[] (8) at (1.2,-0.8) {};
\draw (1) edge[latex-,leftblue] (0);
\draw (0) edge[latex-,leftblue] (5);
\draw (0) edge[latex-,leftblue] (6);
\draw (0) edge[latex-,rightblue] (7);
\draw (0) edge[latex-,rightblue] (8);
\end{tikzpicture}}\Ea
=
\sum_{I=I_1\sqcup I_2 \atop J=J_1\sqcup J_2}(-1)^{\# I_1(\# I_2+\# J_1)+ \# J_2+1+ \text{sign}(I_1,I_2) + \text{sign}(J_1,J_2)}
\left(
\Ba{c}\resizebox{25mm}{!}
{\begin{tikzpicture}[baseline=-1ex]
\node[] at (-0.9,-0.15) {$_{\underbrace{\ \ \ \ }_{I_1}}$};
\node[] at (0.9,-0.15) {$_{\underbrace{\ \ \ \ }_{J_2}}$};
\node[] at (-0.5,-0.9) {$_{\underbrace{\ }_{I_2}}$};
\node[] at (0.5,-0.9) {$_{\underbrace{\   }_{J_1}}$};
\node[] at (-0.5,-0.6) {$...$};
\node[] at (-0.75,0.2) {$...$};
\node[] at (0.75,0.2) {$...$};
\node[int] (0) at (0,0) {};
\node[int] (1) at (0,0.8) {};
\node[] (u) at (0,1.6) {};
\node[] (l1) at (-1.4,0) {};
\node[] (l2) at (-0.7,0) {};
\node[] (r2) at (1.4,0) {};
\node[] (r1) at (0.7,0) {};
\node[int] (0) at (0,0) {};
\node[] (5) at (-0.9,-0.8) {};
\node[] (6) at (-0.3,-0.8) {};
\node[] (7) at (0.3,-0.8) {};
\node[] (8) at (0.9,-0.8) {};
\draw (u) edge[latex-,leftblue] (1);
\draw (1) edge[latex-,leftblue] (l1);
\draw (1) edge[latex-,leftblue] (l2);
\draw (1) edge[latex-,rightblue] (r1);
\draw (1) edge[latex-,rightblue] (r2);
\draw (1) edge[latex-,leftblue] (0);
\draw (0) edge[latex-,leftblue] (5);
\draw (0) edge[latex-,leftblue] (6);
\draw (0) edge[latex-,rightblue] (7);
\draw (0) edge[latex-,rightblue] (8);
\end{tikzpicture}}\Ea
+
\Ba{c}\resizebox{25mm}{!}
{\begin{tikzpicture}[baseline=-1ex]
\node[] at (-0.9,-0.15) {$_{\underbrace{\ \ \ \ }_{I_1}}$};
\node[] at (0.9,-0.15) {$_{\underbrace{\ \ \ \ }_{J_2}}$};
\node[] at (-0.5,-0.9) {$_{\underbrace{\ }_{I_2}}$};
\node[] at (0.5,-0.9) {$_{\underbrace{\   }_{J_1}}$};
\node[] at (-0.5,-0.6) {$...$};
\node[] at (-0.75,0.2) {$...$};
\node[] at (0.75,0.2) {$...$};
\node[int] (0) at (0,0) {};
\node[int] (1) at (0,0.8) {};
\node[] (u) at (0,1.6) {};
\node[] (l1) at (-1.4,0) {};
\node[] (l2) at (-0.7,0) {};
\node[] (r2) at (1.4,0) {};
\node[] (r1) at (0.7,0) {};
\node[int] (0) at (0,0) {};
\node[] (5) at (-0.9,-0.8) {};
\node[] (6) at (-0.3,-0.8) {};
\node[] (7) at (0.3,-0.8) {};
\node[] (8) at (0.9,-0.8) {};
\draw (u) edge[latex-,leftblue] (1);
\draw (1) edge[latex-,leftblue] (l1);
\draw (1) edge[latex-,leftblue] (l2);
\draw (1) edge[latex-,rightblue] (r1);
\draw (1) edge[latex-,rightblue] (r2);
\draw (1) edge[latex-,rightblue] (0);
\draw (0) edge[latex-,leftblue] (5);
\draw (0) edge[latex-,leftblue] (6);
\draw (0) edge[latex-,rightblue] (7);
\draw (0) edge[latex-,rightblue] (8);
\end{tikzpicture}}\Ea
\right)
$$
and similarly for the second class of corollas.
Here the sums run over all admissible decompositions of the ordered sets
$I$ and $J$ into the disjoint unions of (not necessarily connected)  ordered subsets, and $\text{sign}(I_1,I_2)$ (resp., $\text{sign}(J_1,J_2)$) stands for the parity of the permutation $I\rar I_1\sqcup I_2$ (resp., $J\rar J_1\sqcup J_2$).

\sip

If $I$ is the differential closure of the ideal in $\caL ie_\infty^{(k+1)}$ generated
by all corollas with negative cohomological degree, then the quotient
$$
\caL ie^{(k+1)}:= \caL ie_\infty^{(k+1)}/I
$$
is called an operad of multi-oriented Lie algebras.

\subsubsection{\bf Theorem} {\em The natural projection $\caL ie_\infty^{(k+1)}
\lon \caL ie^{(k+1)}$ is a quasi-isomorphism.}
\begin{proof}
It is enough to show that the cohomology of $\caL ie_\infty^{(k+1)}$ is concentrated
in degree zero, and this can be done by the arity induction
in a close analogy to the proof of Theorem {\ref{3new: Theorem on Ass_infty^k+1}}. We omit the details.
\end{proof}

\subsection{\bf The operad of 2-oriented Lie algebras versus the ordinary dioperad of Lie bialgebras}
The operad
$\caL ie^{(2)}$ can be explicitly described as follows: it is generated by the following list of degree 0 corollas,
$$
\Ba{c}\resizebox{12mm}{!}  {
\begin{tikzpicture}[baseline=-1ex]
\node[] at (0,0.8) {$^{0}$};
\node[] at (-0.5,-0.71) {$_{1}$};
\node[] at (0.5,-0.71) {$_{2}$};
\node[int] (a) at (0,0) {};
\node[] (u1) at (0,0.8) {};
\node[int] (a) at (0,0) {};
\node[] (d1) at (-0.5,-0.7) {};
\node[] (d2) at (0.5,-0.7) {};
\draw (u1) edge[latex-, leftblue] (a);
\draw (a) edge[latex-, leftblue] (d1);
\draw (a) edge[latex-, leftblue] (d2);
\end{tikzpicture}
}\Ea =- \Ba{c}\resizebox{12mm}{!}  {
\begin{tikzpicture}[baseline=-1ex]
\node[] at (0,0.8) {$^{0}$};
\node[] at (-0.5,-0.71) {$_{2}$};
\node[] at (0.5,-0.71) {$_{1}$};
\node[int] (a) at (0,0) {};
\node[] (u1) at (0,0.8) {};
\node[int] (a) at (0,0) {};
\node[] (d1) at (-0.5,-0.7) {};
\node[] (d2) at (0.5,-0.7) {};
\draw (u1) edge[latex-, leftblue] (a);
\draw (a) edge[latex-, leftblue] (d1);
\draw (a) edge[latex-, leftblue] (d2);
\end{tikzpicture}}\Ea
,
\Ba{c}\resizebox{12mm}{!}  {
\begin{tikzpicture}[baseline=-1ex]
\node[] at (0,0.8) {$^{0}$};
\node[] at (-0.5,-0.71) {$_{1}$};
\node[] at (0.5,-0.71) {$_{2}$};
\node[int] (a) at (0,0) {};
\node[] (u1) at (0,0.8) {};
\node[int] (a) at (0,0) {};
\node[] (d1) at (-0.5,-0.7) {};
\node[] (d2) at (0.5,-0.7) {};
\draw (u1) edge[latex-, rightblue] (a);
\draw (a) edge[latex-, rightblue] (d1);
\draw (a) edge[latex-, rightblue] (d2);
\end{tikzpicture}
}\Ea
=-
\Ba{c}\resizebox{12mm}{!}  {
\begin{tikzpicture}[baseline=-1ex]
\node[] at (0,0.8) {$^{0}$};
\node[] at (-0.5,-0.71) {$_{2}$};
\node[] at (0.5,-0.71) {$_{1}$};
\node[int] (a) at (0,0) {};
\node[] (u1) at (0,0.8) {};
\node[int] (a) at (0,0) {};
\node[] (d1) at (-0.5,-0.7) {};
\node[] (d2) at (0.5,-0.7) {};
\draw (u1) edge[latex-, rightblue] (a);
\draw (a) edge[latex-, rightblue] (d1);
\draw (a) edge[latex-, rightblue] (d2);
\end{tikzpicture}}\Ea
,
\Ba{c}\resizebox{12mm}{!}  {
\begin{tikzpicture}[baseline=-1ex]
\node[] at (0,0.8) {$^{0}$};
\node[] at (-0.5,-0.71) {$_{1}$};
\node[] at (0.5,-0.71) {$_{2}$};
\node[int] (a) at (0,0) {};
\node[] (u1) at (0,0.8) {};
\node[int] (a) at (0,0) {};
\node[] (d1) at (-0.5,-0.7) {};
\node[] (d2) at (0.5,-0.7) {};
\draw (u1) edge[latex-, leftblue] (a);
\draw (a) edge[latex-, leftblue] (d1);
\draw (a) edge[latex-, rightblue] (d2);
\end{tikzpicture}
}\Ea
,
\Ba{c}\resizebox{12mm}{!}  {
\begin{tikzpicture}[baseline=-1ex]
\node[] at (0,0.8) {$^{0}$};
\node[] at (-0.5,-0.71) {$_{1}$};
\node[] at (0.5,-0.71) {$_{2}$};
\node[int] (a) at (0,0) {};
\node[] (u1) at (0,0.8) {};
\node[int] (a) at (0,0) {};
\node[] (d1) at (-0.5,-0.7) {};
\node[] (d2) at (0.5,-0.7) {};
\draw (u1) edge[latex-, rightblue] (a);
\draw (a) edge[latex-, leftblue] (d1);
\draw (a) edge[latex-, rightblue] (d2);
\end{tikzpicture}
}\Ea
$$
modulo the following relations

\Beq\label{3new: Lie^2 relations 1st set}
\Ba{c}\resizebox{15mm}{!}  {
\begin{tikzpicture}[baseline=-1ex]
\node[] at (0.5,1.6) {$^{{0}}$};
\node[] at (-0.5,-0.72) {$_{{1}}$};
\node[] at (0.57,-0.72) {$_{{2}}$};
\node[] at (1.12,-0.02) {$_{{3}}$};
\node[int] (a) at (0,0) {};
\node[] (u2) at (0.5,1.6) {};
\node[] (ud2) at (1,0.0) {};
\node[int] (u1) at (0.5,0.8) {};
\node[int] (a) at (0,0) {};
\node[] (d1) at (-0.5,-0.7) {};
\node[] (d2) at (0.5,-0.7) {};
\draw (u1) edge[latex-,leftblue] (a);
\draw (u2) edge[latex-,leftblue] (u1);
\draw (u1) edge[latex-, leftblue] (ud2);
\draw (a) edge[latex-,leftblue] (d1);
\draw (a) edge[latex-,leftblue] (d2);
\end{tikzpicture}
}\Ea
+\hspace{-3mm}
\Ba{c}\resizebox{15mm}{!}  {
\begin{tikzpicture}[baseline=-1ex]
\node[] at (0.5,1.6) {$^{{0}}$};
\node[] at (-0.5,-0.72) {$_{{3}}$};
\node[] at (0.57,-0.72) {$_{{1}}$};
\node[] at (1.12,-0.02) {$_{{2}}$};
\node[int] (a) at (0,0) {};
\node[] (u2) at (0.5,1.6) {};
\node[] (ud2) at (1,0.0) {};
\node[int] (u1) at (0.5,0.8) {};
\node[int] (a) at (0,0) {};
\node[] (d1) at (-0.5,-0.7) {};
\node[] (d2) at (0.5,-0.7) {};
\draw (u1) edge[latex-,leftblue] (a);
\draw (u2) edge[latex-,leftblue] (u1);
\draw (u1) edge[latex-, leftblue] (ud2);
\draw (a) edge[latex-,leftblue] (d1);
\draw (a) edge[latex-,leftblue] (d2);
\end{tikzpicture}
}\Ea
+\hspace{-3mm}
\Ba{c}\resizebox{15mm}{!}  {
\begin{tikzpicture}[baseline=-1ex]
\node[] at (0.5,1.6) {$^{{0}}$};
\node[] at (-0.5,-0.72) {$_{{1}}$};
\node[] at (0.57,-0.72) {$_{{2}}$};
\node[] at (1.12,-0.02) {$_{{3}}$};
\node[int] (a) at (0,0) {};
\node[] (u2) at (0.5,1.6) {};
\node[] (ud2) at (1,0.0) {};
\node[int] (u1) at (0.5,0.8) {};
\node[int] (a) at (0,0) {};
\node[] (d1) at (-0.5,-0.7) {};
\node[] (d2) at (0.5,-0.7) {};
\draw (u1) edge[latex-,leftblue] (a);
\draw (u2) edge[latex-,leftblue] (u1);
\draw (u1) edge[latex-, leftblue] (ud2);
\draw (a) edge[latex-,leftblue] (d1);
\draw (a) edge[latex-,leftblue] (d2);
\end{tikzpicture}
}\Ea\hspace{-3mm}
=0
, \ \ \
\Ba{c}\resizebox{15mm}{!}  {
\begin{tikzpicture}[baseline=-1ex]
\node[] at (0.5,1.6) {$^{{0}}$};
\node[] at (-0.5,-0.72) {$_{{1}}$};
\node[] at (0.57,-0.72) {$_{{2}}$};
\node[] at (1.12,-0.02) {$_{{3}}$};
\node[int] (a) at (0,0) {};
\node[] (u2) at (0.5,1.6) {};
\node[] (ud2) at (1,0.0) {};
\node[int] (u1) at (0.5,0.8) {};
\node[int] (a) at (0,0) {};
\node[] (d1) at (-0.5,-0.7) {};
\node[] (d2) at (0.5,-0.7) {};
\draw (u1) edge[latex-,leftblue] (a);
\draw (u2) edge[latex-,leftblue] (u1);
\draw (u1) edge[latex-, rightblue] (ud2);
\draw (a) edge[latex-,leftblue] (d1);
\draw (a) edge[latex-,leftblue] (d2);
\end{tikzpicture}
}\Ea
\hspace{-3mm}
-
\hspace{-2mm}
\Ba{c}\resizebox{15mm}{!}  {
\begin{tikzpicture}[baseline=-1ex]
\node[] at (0.5,1.6) {$^{{0}}$};
\node[] at (1.55,-0.72) {$_{{3}}$};
\node[] at (0.57,-0.72) {$_{{2}}$};
\node[] at (0.0,-0.02) {$_{{1}}$};
\node[] (u2) at (0.5,1.6) {};
\node[] (ud2) at (0,0.0) {};
\node[int] (u1) at (0.5,0.8) {};
\node[int] (a) at (1.0,0) {};
\node[] (d1) at (0.5,-0.7) {};
\node[] (d2) at (1.5,-0.7) {};
\draw (u1) edge[latex-,leftblue] (a);
\draw (u2) edge[latex-,leftblue] (u1);
\draw (u1) edge[latex-,leftblue] (ud2);
\draw (a) edge[latex-,leftblue] (d1);
\draw (a) edge[latex-,rightblue] (d2);
\end{tikzpicture}
}\Ea
\hspace{-3mm}
+
\hspace{-2mm}
\Ba{c}\resizebox{15mm}{!}  {
\begin{tikzpicture}[baseline=-1ex]
\node[] at (0.5,1.6) {$^{{0}}$};
\node[] at (1.55,-0.72) {$_{{3}}$};
\node[] at (0.57,-0.72) {$_{{1}}$};
\node[] at (0.0,-0.02) {$_{{2}}$};
\node[] (u2) at (0.5,1.6) {};
\node[] (ud2) at (0,0.0) {};
\node[int] (u1) at (0.5,0.8) {};
\node[int] (a) at (1.0,0) {};
\node[] (d1) at (0.5,-0.7) {};
\node[] (d2) at (1.5,-0.7) {};
\draw (u1) edge[latex-,leftblue] (a);
\draw (u2) edge[latex-,leftblue] (u1);
\draw (u1) edge[latex-,leftblue] (ud2);
\draw (a) edge[latex-,leftblue] (d1);
\draw (a) edge[latex-,rightblue] (d2);
\end{tikzpicture}
}\Ea
\hspace{-3mm}
-
\hspace{-2mm}
\Ba{c}\resizebox{15mm}{!}  {
\begin{tikzpicture}[baseline=-1ex]
\node[] at (0.5,1.6) {$^{{0}}$};
\node[] at (1.55,-0.72) {$_{{3}}$};
\node[] at (0.57,-0.72) {$_{{2}}$};
\node[] at (0.0,-0.02) {$_{{1}}$};
\node[] (u2) at (0.5,1.6) {};
\node[] (ud2) at (0,0.0) {};
\node[int] (u1) at (0.5,0.8) {};
\node[int] (a) at (1.0,0) {};
\node[] (d1) at (0.5,-0.7) {};
\node[] (d2) at (1.5,-0.7) {};
\draw (u1) edge[latex-,rightblue] (a);
\draw (u2) edge[latex-,leftblue] (u1);
\draw (u1) edge[latex-,leftblue] (ud2);
\draw (a) edge[latex-,leftblue] (d1);
\draw (a) edge[latex-,rightblue] (d2);
\end{tikzpicture}
}\Ea
\hspace{-3mm}
+
\hspace{-2mm}
\Ba{c}\resizebox{15mm}{!}  {
\begin{tikzpicture}[baseline=-1ex]
\node[] at (0.5,1.6) {$^{{0}}$};
\node[] at (1.55,-0.72) {$_{{3}}$};
\node[] at (0.57,-0.72) {$_{{1}}$};
\node[] at (0.0,-0.02) {$_{{2}}$};
\node[] (u2) at (0.5,1.6) {};
\node[] (ud2) at (0,0.0) {};
\node[int] (u1) at (0.5,0.8) {};
\node[int] (a) at (1.0,0) {};
\node[] (d1) at (0.5,-0.7) {};
\node[] (d2) at (1.5,-0.7) {};
\draw (u1) edge[latex-,rightblue] (a);
\draw (u2) edge[latex-,leftblue] (u1);
\draw (u1) edge[latex-,leftblue] (ud2);
\draw (a) edge[latex-,leftblue] (d1);
\draw (a) edge[latex-,rightblue] (d2);
\end{tikzpicture}
}\Ea
\hspace{-4mm}
=0
\Eeq


\Beq\label{3new: Lie^2 relations 2nd set}
\Ba{c}\resizebox{15mm}{!}  {
\begin{tikzpicture}[baseline=-1ex]
\node[] at (0.5,1.6) {$^{{0}}$};
\node[] at (-0.5,-0.72) {$_{{1}}$};
\node[] at (0.57,-0.72) {$_{{2}}$};
\node[] at (1.12,-0.02) {$_{{3}}$};
\node[int] (a) at (0,0) {};
\node[] (u2) at (0.5,1.6) {};
\node[] (ud2) at (1,0.0) {};
\node[int] (u1) at (0.5,0.8) {};
\node[int] (a) at (0,0) {};
\node[] (d1) at (-0.5,-0.7) {};
\node[] (d2) at (0.5,-0.7) {};
\draw (u1) edge[latex-,leftblue] (a);
\draw (u2) edge[latex-,leftblue] (u1);
\draw (u1) edge[latex-, rightblue] (ud2);
\draw (a) edge[latex-,leftblue] (d1);
\draw (a) edge[latex-,rightblue] (d2);
\end{tikzpicture}
}\Ea
\hspace{-3mm}
-
\hspace{-2mm}
\Ba{c}\resizebox{15mm}{!}  {
\begin{tikzpicture}[baseline=-1ex]
\node[] at (0.5,1.6) {$^{{0}}$};
\node[] at (-0.5,-0.72) {$_{{1}}$};
\node[] at (0.57,-0.72) {$_{{3}}$};
\node[] at (1.12,-0.02) {$_{{2}}$};
\node[int] (a) at (0,0) {};
\node[] (u2) at (0.5,1.6) {};
\node[] (ud2) at (1,0.0) {};
\node[int] (u1) at (0.5,0.8) {};
\node[int] (a) at (0,0) {};
\node[] (d1) at (-0.5,-0.7) {};
\node[] (d2) at (0.5,-0.7) {};
\draw (u1) edge[latex-,leftblue] (a);
\draw (u2) edge[latex-,leftblue] (u1);
\draw (u1) edge[latex-, rightblue] (ud2);
\draw (a) edge[latex-,leftblue] (d1);
\draw (a) edge[latex-,rightblue] (d2);
\end{tikzpicture}
}\Ea
\hspace{-3mm}
-
\hspace{-2mm}
\Ba{c}\resizebox{15mm}{!}  {
\begin{tikzpicture}[baseline=-1ex]
\node[] at (0.5,1.6) {$^{{0}}$};
\node[] at (1.55,-0.72) {$_{{3}}$};
\node[] at (0.57,-0.72) {$_{{2}}$};
\node[] at (0.0,-0.02) {$_{{1}}$};
\node[] (u2) at (0.5,1.6) {};
\node[] (ud2) at (0,0.0) {};
\node[int] (u1) at (0.5,0.8) {};
\node[int] (a) at (1.0,0) {};
\node[] (d1) at (0.5,-0.7) {};
\node[] (d2) at (1.5,-0.7) {};
\draw (u1) edge[latex-,rightblue] (a);
\draw (u2) edge[latex-,leftblue] (u1);
\draw (u1) edge[latex-,leftblue] (ud2);
\draw (a) edge[latex-,rightblue] (d1);
\draw (a) edge[latex-,rightblue] (d2);
\end{tikzpicture}
}\Ea
\hspace{-3mm}
=
0, \ \ \ \
\Ba{c}\resizebox{15mm}{!}  {
\begin{tikzpicture}[baseline=-1ex]
\node[] at (0.5,1.6) {$^{{0}}$};
\node[] at (-0.5,-0.72) {$_{{1}}$};
\node[] at (0.57,-0.72) {$_{{2}}$};
\node[] at (1.12,-0.02) {$_{{3}}$};
\node[int] (a) at (0,0) {};
\node[] (u2) at (0.5,1.6) {};
\node[] (ud2) at (1,0.0) {};
\node[int] (u1) at (0.5,0.8) {};
\node[int] (a) at (0,0) {};
\node[] (d1) at (-0.5,-0.7) {};
\node[] (d2) at (0.5,-0.7) {};
\draw (u1) edge[latex-,rightblue] (a);
\draw (u2) edge[latex-,rightblue] (u1);
\draw (u1) edge[latex-, rightblue] (ud2);
\draw (a) edge[latex-,rightblue] (d1);
\draw (a) edge[latex-,rightblue] (d2);
\end{tikzpicture}
}\Ea
\hspace{-3mm}
+
\hspace{-2mm}
\Ba{c}\resizebox{15mm}{!}  {
\begin{tikzpicture}[baseline=-1ex]
\node[] at (0.5,1.6) {$^{{0}}$};
\node[] at (-0.5,-0.72) {$_{{3}}$};
\node[] at (0.57,-0.72) {$_{{1}}$};
\node[] at (1.12,-0.02) {$_{{2}}$};
\node[int] (a) at (0,0) {};
\node[] (u2) at (0.5,1.6) {};
\node[] (ud2) at (1,0.0) {};
\node[int] (u1) at (0.5,0.8) {};
\node[int] (a) at (0,0) {};
\node[] (d1) at (-0.5,-0.7) {};
\node[] (d2) at (0.5,-0.7) {};
\draw (u1) edge[latex-,rightblue] (a);
\draw (u2) edge[latex-,rightblue] (u1);
\draw (u1) edge[latex-, rightblue] (ud2);
\draw (a) edge[latex-,rightblue] (d1);
\draw (a) edge[latex-,rightblue] (d2);
\end{tikzpicture}
}\Ea
\hspace{-3mm}
+
\hspace{-2mm}
\Ba{c}\resizebox{15mm}{!}  {
\begin{tikzpicture}[baseline=-1ex]
\node[] at (0.5,1.6) {$^{{0}}$};
\node[] at (-0.5,-0.72) {$_{{1}}$};
\node[] at (0.57,-0.72) {$_{{2}}$};
\node[] at (1.12,-0.02) {$_{{3}}$};
\node[int] (a) at (0,0) {};
\node[] (u2) at (0.5,1.6) {};
\node[] (ud2) at (1,0.0) {};
\node[int] (u1) at (0.5,0.8) {};
\node[int] (a) at (0,0) {};
\node[] (d1) at (-0.5,-0.7) {};
\node[] (d2) at (0.5,-0.7) {};
\draw (u1) edge[latex-,rightblue] (a);
\draw (u2) edge[latex-,rightblue] (u1);
\draw (u1) edge[latex-, rightblue] (ud2);
\draw (a) edge[latex-,rightblue] (d1);
\draw (a) edge[latex-,rightblue] (d2);
\end{tikzpicture}
}\Ea
\hspace{-3mm}
=0
\Eeq


\Beq\label{3new: Lie^2 relations 3rd set}
\Ba{c}\resizebox{15mm}{!}  {
\begin{tikzpicture}[baseline=-1ex]
\node[] at (0.5,1.6) {$^{{0}}$};
\node[] at (-0.5,-0.72) {$_{{1}}$};
\node[] at (0.57,-0.72) {$_{{2}}$};
\node[] at (1.12,-0.02) {$_{{3}}$};
\node[int] (a) at (0,0) {};
\node[] (u2) at (0.5,1.6) {};
\node[] (ud2) at (1,0.0) {};
\node[int] (u1) at (0.5,0.8) {};
\node[int] (a) at (0,0) {};
\node[] (d1) at (-0.5,-0.7) {};
\node[] (d2) at (0.5,-0.7) {};
\draw (u1) edge[latex-,rightblue] (a);
\draw (u2) edge[latex-,rightblue] (u1);
\draw (u1) edge[latex-, rightblue] (ud2);
\draw (a) edge[latex-,leftblue] (d1);
\draw (a) edge[latex-,rightblue] (d2);
\end{tikzpicture}
}\Ea
\hspace{-2mm}
+
\hspace{-2mm}
\Ba{c}\resizebox{15mm}{!}  {
\begin{tikzpicture}[baseline=-1ex]
\node[] at (0.5,1.6) {$^{{0}}$};
\node[] at (-0.5,-0.72) {$_{{1}}$};
\node[] at (0.57,-0.72) {$_{{3}}$};
\node[] at (1.12,-0.02) {$_{{2}}$};
\node[int] (a) at (0,0) {};
\node[] (u2) at (0.5,1.6) {};
\node[] (ud2) at (1,0.0) {};
\node[int] (u1) at (0.5,0.8) {};
\node[int] (a) at (0,0) {};
\node[] (d1) at (-0.5,-0.7) {};
\node[] (d2) at (0.5,-0.7) {};
\draw (u1) edge[latex-,rightblue] (a);
\draw (u2) edge[latex-,rightblue] (u1);
\draw (u1) edge[latex-, rightblue] (ud2);
\draw (a) edge[latex-,leftblue] (d1);
\draw (a) edge[latex-,rightblue] (d2);
\end{tikzpicture}
}\Ea
\hspace{-2mm}
+
\hspace{-2mm}
\Ba{c}\resizebox{15mm}{!}  {
\begin{tikzpicture}[baseline=-1ex]
\node[] at (0.5,1.6) {$^{{0}}$};
\node[] at (-0.5,-0.72) {$_{{1}}$};
\node[] at (0.57,-0.72) {$_{{2}}$};
\node[] at (1.12,-0.02) {$_{{3}}$};
\node[int] (a) at (0,0) {};
\node[] (u2) at (0.5,1.6) {};
\node[] (ud2) at (1,0.0) {};
\node[int] (u1) at (0.5,0.8) {};
\node[int] (a) at (0,0) {};
\node[] (d1) at (-0.5,-0.7) {};
\node[] (d2) at (0.5,-0.7) {};
\draw (u1) edge[latex-,leftblue] (a);
\draw (u2) edge[latex-,rightblue] (u1);
\draw (u1) edge[latex-, rightblue] (ud2);
\draw (a) edge[latex-,leftblue] (d1);
\draw (a) edge[latex-,rightblue] (d2);
\end{tikzpicture}
}\Ea
\hspace{-2mm}
-
\hspace{-2mm}
\Ba{c}\resizebox{15mm}{!}  {
\begin{tikzpicture}[baseline=-1ex]
\node[] at (0.5,1.6) {$^{{0}}$};
\node[] at (-0.5,-0.72) {$_{{1}}$};
\node[] at (0.57,-0.72) {$_{{3}}$};
\node[] at (1.12,-0.02) {$_{{2}}$};
\node[int] (a) at (0,0) {};
\node[] (u2) at (0.5,1.6) {};
\node[] (ud2) at (1,0.0) {};
\node[int] (u1) at (0.5,0.8) {};
\node[int] (a) at (0,0) {};
\node[] (d1) at (-0.5,-0.7) {};
\node[] (d2) at (0.5,-0.7) {};
\draw (u1) edge[latex-,leftblue] (a);
\draw (u2) edge[latex-,rightblue] (u1);
\draw (u1) edge[latex-, rightblue] (ud2);
\draw (a) edge[latex-,leftblue] (d1);
\draw (a) edge[latex-,rightblue] (d2);
\end{tikzpicture}
}\Ea
\hspace{-2mm}
-
\hspace{-2mm}
\Ba{c}\resizebox{15mm}{!}  {
\begin{tikzpicture}[baseline=-1ex]
\node[] at (0.5,1.6) {$^{{0}}$};
\node[] at (1.55,-0.72) {$_{{3}}$};
\node[] at (0.57,-0.72) {$_{{2}}$};
\node[] at (0.0,-0.02) {$_{{1}}$};
\node[] (u2) at (0.5,1.6) {};
\node[] (ud2) at (0,0.0) {};
\node[int] (u1) at (0.5,0.8) {};
\node[int] (a) at (1.0,0) {};
\node[] (d1) at (0.5,-0.7) {};
\node[] (d2) at (1.5,-0.7) {};
\draw (u1) edge[latex-,rightblue] (a);
\draw (u2) edge[latex-,rightblue] (u1);
\draw (u1) edge[latex-,leftblue] (ud2);
\draw (a) edge[latex-,rightblue] (d1);
\draw (a) edge[latex-,rightblue] (d2);
\end{tikzpicture}
}\Ea
\hspace{-3mm}
=0,
 \
\Ba{c}\resizebox{15mm}{!}  {
\begin{tikzpicture}[baseline=-1ex]
\node[] at (0.5,1.6) {$^{{0}}$};
\node[] at (-0.5,-0.72) {$_{{1}}$};
\node[] at (0.57,-0.72) {$_{{2}}$};
\node[] at (1.12,-0.02) {$_{{3}}$};
\node[int] (a) at (0,0) {};
\node[] (u2) at (0.5,1.6) {};
\node[] (ud2) at (1,0.0) {};
\node[int] (u1) at (0.5,0.8) {};
\node[int] (a) at (0,0) {};
\node[] (d1) at (-0.5,-0.7) {};
\node[] (d2) at (0.5,-0.7) {};
\draw (u1) edge[latex-,leftblue] (a);
\draw (u2) edge[latex-,rightblue] (u1);
\draw (u1) edge[latex-, rightblue] (ud2);
\draw (a) edge[latex-,leftblue] (d1);
\draw (a) edge[latex-,leftblue] (d2);
\end{tikzpicture}
}\Ea
\hspace{-2mm}
-
\hspace{-2mm}
\Ba{c}\resizebox{15mm}{!}  {
\begin{tikzpicture}[baseline=-1ex]
\node[] at (0.5,1.6) {$^{{0}}$};
\node[] at (1.55,-0.72) {$_{{3}}$};
\node[] at (0.57,-0.72) {$_{{2}}$};
\node[] at (0.0,-0.02) {$_{{1}}$};
\node[] (u2) at (0.5,1.6) {};
\node[] (ud2) at (0,0.0) {};
\node[int] (u1) at (0.5,0.8) {};
\node[int] (a) at (1.0,0) {};
\node[] (d1) at (0.5,-0.7) {};
\node[] (d2) at (1.5,-0.7) {};
\draw (u1) edge[latex-,rightblue] (a);
\draw (u2) edge[latex-,rightblue] (u1);
\draw (u1) edge[latex-,leftblue] (ud2);
\draw (a) edge[latex-,leftblue] (d1);
\draw (a) edge[latex-,rightblue] (d2);
\end{tikzpicture}
}\Ea
\hspace{-2mm}
+
\hspace{-2mm}
\Ba{c}\resizebox{15mm}{!}  {
\begin{tikzpicture}[baseline=-1ex]
\node[] at (0.5,1.6) {$^{{0}}$};
\node[] at (1.55,-0.72) {$_{{3}}$};
\node[] at (0.57,-0.72) {$_{{1}}$};
\node[] at (0.0,-0.02) {$_{{2}}$};
\node[] (u2) at (0.5,1.6) {};
\node[] (ud2) at (0,0.0) {};
\node[int] (u1) at (0.5,0.8) {};
\node[int] (a) at (1.0,0) {};
\node[] (d1) at (0.5,-0.7) {};
\node[] (d2) at (1.5,-0.7) {};
\draw (u1) edge[latex-,rightblue] (a);
\draw (u2) edge[latex-,rightblue] (u1);
\draw (u1) edge[latex-,leftblue] (ud2);
\draw (a) edge[latex-,leftblue] (d1);
\draw (a) edge[latex-,rightblue] (d2);
\end{tikzpicture}
}\Ea
\hspace{-3mm}
=0
\Eeq
Representations of the operad of 2-oriented Lie algebras in symplectic vector spaces with one Lagrangian brane are studied in the next section where it is shown that
they can be identified with famous Manin's triples which give us an alternative (and often very useful) characterization of {\em Lie bialgebras}.
Hence the combinatorics of the latter structure must be hidden in the combinatorics of the former one, and our next our purpose to make this inter-relation
explicit.

\sip

Recall that the ordinary (i.e.\ 1-oriented) dioperad  of Lie bialgebras is the quotient
$$
\LB_{diop}:= \cF\! ree_0^{1\text{-or}}\left\langle M\right\rangle/J
$$
of the
 free 1-oriented free dioperad generated by an $\bS$-bimodule $M=\{M(m,n)\}$ with
$$
M(m,n):=\left\{
\Ba{rr}
sgn_2\ot \id_1\equiv\mbox{span}\left\langle
\begin{xy}
 <0mm,-0.55mm>*{};<0mm,-2.5mm>*{}**@{-},
 <0.5mm,0.5mm>*{};<2.2mm,2.2mm>*{}**@{-},
 <-0.48mm,0.48mm>*{};<-2.2mm,2.2mm>*{}**@{-},
 <0mm,0mm>*{\bu};<0mm,0mm>*{}**@{},
 <0mm,-0.55mm>*{};<0mm,-3.8mm>*{_{_0}}**@{},
 <0.5mm,0.5mm>*{};<2.7mm,2.8mm>*{^{^2}}**@{},
 <-0.48mm,0.48mm>*{};<-2.7mm,2.8mm>*{^{^1}}**@{},
 \end{xy}
=-
\begin{xy}
 <0mm,-0.55mm>*{};<0mm,-2.5mm>*{}**@{-},
 <0.5mm,0.5mm>*{};<2.2mm,2.2mm>*{}**@{-},
 <-0.48mm,0.48mm>*{};<-2.2mm,2.2mm>*{}**@{-},
 <0mm,0mm>*{\bu};<0mm,0mm>*{}**@{},
 <0mm,-0.55mm>*{};<0mm,-3.8mm>*{_{_0}}**@{},
 <0.5mm,0.5mm>*{};<2.7mm,2.8mm>*{^{^1}}**@{},
 <-0.48mm,0.48mm>*{};<-2.7mm,2.8mm>*{^{^2}}**@{},
 \end{xy}
   \right\rangle  & \mbox{if}\ m=2, n=1,\vspace{3mm}\\
\id_1\ot sgn_2\equiv
\mbox{span}\left\langle
\begin{xy}
 <0mm,0.66mm>*{};<0mm,3mm>*{}**@{-},
 <0.39mm,-0.39mm>*{};<2.2mm,-2.2mm>*{}**@{-},
 <-0.35mm,-0.35mm>*{};<-2.2mm,-2.2mm>*{}**@{-},
 <0mm,0mm>*{\bu};<0mm,0mm>*{}**@{},
   <0mm,0.66mm>*{};<0mm,3.4mm>*{^{^0}}**@{},
   <0.39mm,-0.39mm>*{};<2.9mm,-4mm>*{^{_2}}**@{},
   <-0.35mm,-0.35mm>*{};<-2.8mm,-4mm>*{^{_1}}**@{},
\end{xy}=-
\begin{xy}
 <0mm,0.66mm>*{};<0mm,3mm>*{}**@{-},
 <0.39mm,-0.39mm>*{};<2.2mm,-2.2mm>*{}**@{-},
 <-0.35mm,-0.35mm>*{};<-2.2mm,-2.2mm>*{}**@{-},
 <0mm,0mm>*{\bu};<0mm,0mm>*{}**@{},
   <0mm,0.66mm>*{};<0mm,3.4mm>*{^{^0}}**@{},
   <0.39mm,-0.39mm>*{};<2.9mm,-4mm>*{^{_1}}**@{},
   <-0.35mm,-0.35mm>*{};<-2.8mm,-4mm>*{^{_2}}**@{},
\end{xy}
\right\rangle
\ & \mbox{if}\ m=1, n=2, \vspace{3mm}\\
0 & \mbox{otherwise}
\Ea
\right.
$$
modulo the ideal $J$ generated by the following relations
$$
\Ba{c}\begin{xy}
 <0mm,0mm>*{\bu};<0mm,0mm>*{}**@{},
 <0mm,-0.49mm>*{};<0mm,-3.0mm>*{}**@{-},
 <0.49mm,0.49mm>*{};<1.9mm,1.9mm>*{}**@{-},
 <-0.5mm,0.5mm>*{};<-1.9mm,1.9mm>*{}**@{-},
 <-2.3mm,2.3mm>*{\bu};<-2.3mm,2.3mm>*{}**@{},
 <-1.8mm,2.8mm>*{};<0mm,4.9mm>*{}**@{-},
 <-2.8mm,2.9mm>*{};<-4.6mm,4.9mm>*{}**@{-},
   <0.49mm,0.49mm>*{};<2.7mm,2.3mm>*{^{^3}}**@{},
   <-1.8mm,2.8mm>*{};<0.4mm,5.3mm>*{^{^2}}**@{},
   <-2.8mm,2.9mm>*{};<-5.1mm,5.3mm>*{^{^1}}**@{},
 \end{xy}
\ + \
\begin{xy}
 <0mm,0mm>*{\bu};<0mm,0mm>*{}**@{},
 <0mm,-0.49mm>*{};<0mm,-3.0mm>*{}**@{-},
 <0.49mm,0.49mm>*{};<1.9mm,1.9mm>*{}**@{-},
 <-0.5mm,0.5mm>*{};<-1.9mm,1.9mm>*{}**@{-},
 <-2.3mm,2.3mm>*{\bu};<-2.3mm,2.3mm>*{}**@{},
 <-1.8mm,2.8mm>*{};<0mm,4.9mm>*{}**@{-},
 <-2.8mm,2.9mm>*{};<-4.6mm,4.9mm>*{}**@{-},
   <0.49mm,0.49mm>*{};<2.7mm,2.3mm>*{^{^2}}**@{},
   <-1.8mm,2.8mm>*{};<0.4mm,5.3mm>*{^{^1}}**@{},
   <-2.8mm,2.9mm>*{};<-5.1mm,5.3mm>*{^{^3}}**@{},
 \end{xy}
\ + \
\begin{xy}
 <0mm,0mm>*{\bu};<0mm,0mm>*{}**@{},
 <0mm,-0.49mm>*{};<0mm,-3.0mm>*{}**@{-},
 <0.49mm,0.49mm>*{};<1.9mm,1.9mm>*{}**@{-},
 <-0.5mm,0.5mm>*{};<-1.9mm,1.9mm>*{}**@{-},
 <-2.3mm,2.3mm>*{\bu};<-2.3mm,2.3mm>*{}**@{},
 <-1.8mm,2.8mm>*{};<0mm,4.9mm>*{}**@{-},
 <-2.8mm,2.9mm>*{};<-4.6mm,4.9mm>*{}**@{-},
   <0.49mm,0.49mm>*{};<2.7mm,2.3mm>*{^{^1}}**@{},
   <-1.8mm,2.8mm>*{};<0.4mm,5.3mm>*{^{^3}}**@{},
   <-2.8mm,2.9mm>*{};<-5.1mm,5.3mm>*{^{^2}}**@{},
 \end{xy}\Ea =0
 \ \ \  , \ \ \
 \Ba{c}\begin{xy}
 <0mm,0mm>*{\bu};<0mm,0mm>*{}**@{},
 <0mm,0.69mm>*{};<0mm,3.0mm>*{}**@{-},
 <0.39mm,-0.39mm>*{};<2.4mm,-2.4mm>*{}**@{-},
 <-0.35mm,-0.35mm>*{};<-1.9mm,-1.9mm>*{}**@{-},
 <-2.4mm,-2.4mm>*{\bu};<-2.4mm,-2.4mm>*{}**@{},
 <-2.0mm,-2.8mm>*{};<0mm,-4.9mm>*{}**@{-},
 <-2.8mm,-2.9mm>*{};<-4.7mm,-4.9mm>*{}**@{-},
    <0.39mm,-0.39mm>*{};<3.3mm,-4.0mm>*{^{_3}}**@{},
    <-2.0mm,-2.8mm>*{};<0.5mm,-6.7mm>*{^{_2}}**@{},
    <-2.8mm,-2.9mm>*{};<-5.2mm,-6.7mm>*{^{_1}}**@{},
 \end{xy}
\ + \
 \begin{xy}
 <0mm,0mm>*{\bu};<0mm,0mm>*{}**@{},
 <0mm,0.69mm>*{};<0mm,3.0mm>*{}**@{-},
 <0.39mm,-0.39mm>*{};<2.4mm,-2.4mm>*{}**@{-},
 <-0.35mm,-0.35mm>*{};<-1.9mm,-1.9mm>*{}**@{-},
 <-2.4mm,-2.4mm>*{\bu};<-2.4mm,-2.4mm>*{}**@{},
 <-2.0mm,-2.8mm>*{};<0mm,-4.9mm>*{}**@{-},
 <-2.8mm,-2.9mm>*{};<-4.7mm,-4.9mm>*{}**@{-},
    <0.39mm,-0.39mm>*{};<3.3mm,-4.0mm>*{^{_2}}**@{},
    <-2.0mm,-2.8mm>*{};<0.5mm,-6.7mm>*{^{_1}}**@{},
    <-2.8mm,-2.9mm>*{};<-5.2mm,-6.7mm>*{^{_3}}**@{},
 \end{xy}
\ + \
 \begin{xy}
 <0mm,0mm>*{\bu};<0mm,0mm>*{}**@{},
 <0mm,0.69mm>*{};<0mm,3.0mm>*{}**@{-},
 <0.39mm,-0.39mm>*{};<2.4mm,-2.4mm>*{}**@{-},
 <-0.35mm,-0.35mm>*{};<-1.9mm,-1.9mm>*{}**@{-},
 <-2.4mm,-2.4mm>*{\bu};<-2.4mm,-2.4mm>*{}**@{},
 <-2.0mm,-2.8mm>*{};<0mm,-4.9mm>*{}**@{-},
 <-2.8mm,-2.9mm>*{};<-4.7mm,-4.9mm>*{}**@{-},
    <0.39mm,-0.39mm>*{};<3.3mm,-4.0mm>*{^{_1}}**@{},
    <-2.0mm,-2.8mm>*{};<0.5mm,-6.7mm>*{^{_3}}**@{},
    <-2.8mm,-2.9mm>*{};<-5.2mm,-6.7mm>*{^{_2}}**@{},
 \end{xy}\Ea =0,
 $$
 $$
 \begin{xy}
 <0mm,2.47mm>*{};<0mm,0.12mm>*{}**@{-},
 <0.5mm,3.5mm>*{};<2.2mm,5.2mm>*{}**@{-},
 <-0.48mm,3.48mm>*{};<-2.2mm,5.2mm>*{}**@{-},
 <0mm,3mm>*{\bu};<0mm,3mm>*{}**@{},
  <0mm,-0.8mm>*{\bu};<0mm,-0.8mm>*{}**@{},
<-0.39mm,-1.2mm>*{};<-2.2mm,-3.5mm>*{}**@{-},
 <0.39mm,-1.2mm>*{};<2.2mm,-3.5mm>*{}**@{-},
     <0.5mm,3.5mm>*{};<2.8mm,5.7mm>*{^{^2}}**@{},
     <-0.48mm,3.48mm>*{};<-2.8mm,5.7mm>*{^{^1}}**@{},
   <0mm,-0.8mm>*{};<-2.7mm,-5.2mm>*{^{_3}}**@{},
   <0mm,-0.8mm>*{};<2.7mm,-5.2mm>*{^{_4}}**@{},
\end{xy}
\  - \
\begin{xy}
 <0mm,-1.3mm>*{};<0mm,-3.5mm>*{}**@{-},
 <0.38mm,-0.2mm>*{};<2.0mm,2.0mm>*{}**@{-},
 <-0.38mm,-0.2mm>*{};<-2.2mm,2.2mm>*{}**@{-},
<0mm,-0.8mm>*{\bu};<0mm,0.8mm>*{}**@{},
 <2.4mm,2.4mm>*{\bu};<2.4mm,2.4mm>*{}**@{},
 <2.77mm,2.0mm>*{};<4.4mm,-0.8mm>*{}**@{-},
 <2.4mm,3mm>*{};<2.4mm,5.2mm>*{}**@{-},
     <0mm,-1.3mm>*{};<0mm,-5.3mm>*{^{_3}}**@{},
     <2.5mm,2.3mm>*{};<5.1mm,-2.6mm>*{^{_4}}**@{},
    <2.4mm,2.5mm>*{};<2.4mm,5.7mm>*{^{^2}}**@{},
    <-0.38mm,-0.2mm>*{};<-2.8mm,2.5mm>*{^{^1}}**@{},
    \end{xy}
\  + \
\begin{xy}
 <0mm,-1.3mm>*{};<0mm,-3.5mm>*{}**@{-},
 <0.38mm,-0.2mm>*{};<2.0mm,2.0mm>*{}**@{-},
 <-0.38mm,-0.2mm>*{};<-2.2mm,2.2mm>*{}**@{-},
<0mm,-0.8mm>*{\bu};<0mm,0.8mm>*{}**@{},
 <2.4mm,2.4mm>*{\bu};<2.4mm,2.4mm>*{}**@{},
 <2.77mm,2.0mm>*{};<4.4mm,-0.8mm>*{}**@{-},
 <2.4mm,3mm>*{};<2.4mm,5.2mm>*{}**@{-},
     <0mm,-1.3mm>*{};<0mm,-5.3mm>*{^{_4}}**@{},
     <2.5mm,2.3mm>*{};<5.1mm,-2.6mm>*{^{_3}}**@{},
    <2.4mm,2.5mm>*{};<2.4mm,5.7mm>*{^{^2}}**@{},
    <-0.38mm,-0.2mm>*{};<-2.8mm,2.5mm>*{^{^1}}**@{},
    \end{xy}
\  - \
\begin{xy}
 <0mm,-1.3mm>*{};<0mm,-3.5mm>*{}**@{-},
 <0.38mm,-0.2mm>*{};<2.0mm,2.0mm>*{}**@{-},
 <-0.38mm,-0.2mm>*{};<-2.2mm,2.2mm>*{}**@{-},
<0mm,-0.8mm>*{\bu};<0mm,0.8mm>*{}**@{},
 <2.4mm,2.4mm>*{\bu};<2.4mm,2.4mm>*{}**@{},
 <2.77mm,2.0mm>*{};<4.4mm,-0.8mm>*{}**@{-},
 <2.4mm,3mm>*{};<2.4mm,5.2mm>*{}**@{-},
     <0mm,-1.3mm>*{};<0mm,-5.3mm>*{^{_4}}**@{},
     <2.5mm,2.3mm>*{};<5.1mm,-2.6mm>*{^{_3}}**@{},
    <2.4mm,2.5mm>*{};<2.4mm,5.7mm>*{^{^1}}**@{},
    <-0.38mm,-0.2mm>*{};<-2.8mm,2.5mm>*{^{^2}}**@{},
    \end{xy}
\ + \
\begin{xy}
 <0mm,-1.3mm>*{};<0mm,-3.5mm>*{}**@{-},
 <0.38mm,-0.2mm>*{};<2.0mm,2.0mm>*{}**@{-},
 <-0.38mm,-0.2mm>*{};<-2.2mm,2.2mm>*{}**@{-},
<0mm,-0.8mm>*{\bu};<0mm,0.8mm>*{}**@{},
 <2.4mm,2.4mm>*{\bu};<2.4mm,2.4mm>*{}**@{},
 <2.77mm,2.0mm>*{};<4.4mm,-0.8mm>*{}**@{-},
 <2.4mm,3mm>*{};<2.4mm,5.2mm>*{}**@{-},
     <0mm,-1.3mm>*{};<0mm,-5.3mm>*{^{_3}}**@{},
     <2.5mm,2.3mm>*{};<5.1mm,-2.6mm>*{^{_4}}**@{},
    <2.4mm,2.5mm>*{};<2.4mm,5.7mm>*{^{^1}}**@{},
    <-0.38mm,-0.2mm>*{};<-2.8mm,2.5mm>*{^{^2}}**@{},
    \end{xy}=0
$$

\subsubsection{ {\bf Proposition (cf.\ \cite{D})}}\label{3new: Prop on map from Lie^2 to Lieb}
{\em There is a (forgetting the basic orientation) morphism of dioperads
$$
\be: \caL ie^{(2)} \lon \LB_{diop}
$$
given on the generators as follows:
$$
\be\left(\Ba{c}\resizebox{12mm}{!}  {
\begin{tikzpicture}[baseline=-1ex]
\node[] at (0,0.8) {$^{0}$};
\node[] at (-0.5,-0.71) {$_{1}$};
\node[] at (0.5,-0.71) {$_{2}$};
\node[int] (a) at (0,0) {};
\node[] (u1) at (0,0.8) {};
\node[int] (a) at (0,0) {};
\node[] (d1) at (-0.5,-0.7) {};
\node[] (d2) at (0.5,-0.7) {};
\draw (u1) edge[latex-, leftblue] (a);
\draw (a) edge[latex-, leftblue] (d1);
\draw (a) edge[latex-, leftblue] (d2);
\end{tikzpicture}}
\Ea
\right):=
\begin{xy}
 <0mm,0.66mm>*{};<0mm,3mm>*{}**@{-},
 <0.39mm,-0.39mm>*{};<2.2mm,-2.2mm>*{}**@{-},
 <-0.35mm,-0.35mm>*{};<-2.2mm,-2.2mm>*{}**@{-},
 <0mm,0mm>*{\bu};<0mm,0mm>*{}**@{},
   <0mm,0.66mm>*{};<0mm,3.4mm>*{^{^0}}**@{},
   <0.39mm,-0.39mm>*{};<2.9mm,-4mm>*{_{_2}}**@{},
   <-0.35mm,-0.35mm>*{};<-2.8mm,-4mm>*{_{_1}}**@{},
\end{xy}\
,\ \ \ \
\be\left(\Ba{c}\resizebox{12mm}{!}  {
\begin{tikzpicture}[baseline=-1ex]
\node[] at (0,0.8) {$^{0}$};
\node[] at (-0.5,-0.71) {$_{1}$};
\node[] at (0.5,-0.71) {$_{2}$};
\node[int] (a) at (0,0) {};
\node[] (u1) at (0,0.8) {};
\node[int] (a) at (0,0) {};
\node[] (d1) at (-0.5,-0.7) {};
\node[] (d2) at (0.5,-0.7) {};
\draw (u1) edge[latex-, rightblue] (a);
\draw (a) edge[latex-, rightblue] (d1);
\draw (a) edge[latex-, rightblue] (d2);
\end{tikzpicture}
}\Ea\right)
=
\begin{xy}
 <0mm,-0.55mm>*{};<0mm,-2.5mm>*{}**@{-},
 <0.5mm,0.5mm>*{};<2.2mm,2.2mm>*{}**@{-},
 <-0.48mm,0.48mm>*{};<-2.2mm,2.2mm>*{}**@{-},
 <0mm,0mm>*{\bu};<0mm,0mm>*{}**@{},
 <0mm,-0.55mm>*{};<0mm,-3.8mm>*{_{_0}}**@{},
 <0.5mm,0.5mm>*{};<2.7mm,2.8mm>*{^{^1}}**@{},
 <-0.48mm,0.48mm>*{};<-2.7mm,2.8mm>*{^{^2}}**@{},
 \end{xy}
,\ \ \
\be\left(\Ba{c}\resizebox{13mm}{!}  {
\begin{tikzpicture}[baseline=-1ex]
\node[] at (0,0.8) {$^{0}$};
\node[] at (-0.5,-0.71) {$_{1}$};
\node[] at (0.5,-0.71) {$_{2}$};
\node[int] (a) at (0,0) {};
\node[] (u1) at (0,0.8) {};
\node[int] (a) at (0,0) {};
\node[] (d1) at (-0.5,-0.7) {};
\node[] (d2) at (0.5,-0.7) {};
\draw (u1) edge[latex-, leftblue] (a);
\draw (a) edge[latex-, leftblue] (d1);
\draw (a) edge[latex-, rightblue] (d2);
\end{tikzpicture}
}\Ea\right)=
\begin{xy}
 <0mm,-0.55mm>*{};<0mm,-2.5mm>*{}**@{-},
 <0.5mm,0.5mm>*{};<2.2mm,2.2mm>*{}**@{-},
 <-0.48mm,0.48mm>*{};<-2.2mm,2.2mm>*{}**@{-},
 <0mm,0mm>*{\bu};<0mm,0mm>*{}**@{},
 <0mm,-0.55mm>*{};<0mm,-3.8mm>*{_{_1}}**@{},
 <0.5mm,0.5mm>*{};<2.7mm,2.8mm>*{^{^2}}**@{},
 <-0.48mm,0.48mm>*{};<-2.7mm,2.8mm>*{^{^0}}**@{},
 \end{xy} \
,\ \ \
\be\left(\Ba{c}\resizebox{12mm}{!}  {
\begin{tikzpicture}[baseline=-1ex]
\node[] at (0,0.8) {$^{0}$};
\node[] at (-0.5,-0.71) {$_{1}$};
\node[] at (0.5,-0.71) {$_{2}$};
\node[int] (a) at (0,0) {};
\node[] (u1) at (0,0.8) {};
\node[int] (a) at (0,0) {};
\node[] (d1) at (-0.5,-0.7) {};
\node[] (d2) at (0.5,-0.7) {};
\draw (u1) edge[latex-, rightblue] (a);
\draw (a) edge[latex-, leftblue] (d1);
\draw (a) edge[latex-, rightblue] (d2);
\end{tikzpicture}
}\Ea\right)
=
\begin{xy}
 <0mm,0.66mm>*{};<0mm,3mm>*{}**@{-},
 <0.39mm,-0.39mm>*{};<2.2mm,-2.2mm>*{}**@{-},
 <-0.35mm,-0.35mm>*{};<-2.2mm,-2.2mm>*{}**@{-},
 <0mm,0mm>*{\bu};<0mm,0mm>*{}**@{},
   <0mm,0.66mm>*{};<0mm,3.4mm>*{^{^2}}**@{},
   <0.39mm,-0.39mm>*{};<2.9mm,-4mm>*{_{_1}}**@{},
   <-0.35mm,-0.35mm>*{};<-2.8mm,-4mm>*{_{_0}}**@{},
\end{xy}\
$$
}
\begin{proof}
It is straightforward to check that each of the eight relations in (\ref{3new: Lie^2 relations 1st set})-(\ref{3new: Lie^2 relations 3rd set}) is mapped under $\be$ into one of the above three relations for $\LB_{diop}$. Hence the map is well-defined indeed.
\end{proof}

This result gives us a purely combinatorial interpretation of the famous Manin triple construction \cite{D}.

\subsection{Multi-oriented prop of homotopy Lie bialgebras} Let us recall a graded generalization of the classical prop of Lie bialgebras depending on two integer parameters $c,d\in \Z$.
By definition \cite{Dr} (one can also see \cite{MW1} for more details)  $\LBcd$ is a quadratic properad given as the quotient,
$$
\LB_{c,d}:=\cF\! ree^{1\text{-or}}\langle Q\rangle/\langle\cR\rangle,
$$
of the free properad generated by an  $\bS$-bimodule $Q=\{Q(m,n)\}_{m,n\geq 1}$ with
 all $Q(m,n)=0$ except
$$
Q(2,1):=\id_1\ot \sgn_2^{\ot c}[c-1]=\mbox{span}\left\langle
\Ba{c}\begin{xy}
 <0mm,-0.55mm>*{};<0mm,-2.5mm>*{}**@{-},
 <0.5mm,0.5mm>*{};<2.2mm,2.2mm>*{}**@{-},
 <-0.48mm,0.48mm>*{};<-2.2mm,2.2mm>*{}**@{-},
 <0mm,0mm>*{\bu};<0mm,0mm>*{}**@{},
 <0.5mm,0.5mm>*{};<2.7mm,2.8mm>*{^{_2}}**@{},
 <-0.48mm,0.48mm>*{};<-2.7mm,2.8mm>*{^{_1}}**@{},
 \end{xy}\Ea
=(-1)^{c}
\Ba{c}\begin{xy}
 <0mm,-0.55mm>*{};<0mm,-2.5mm>*{}**@{-},
 <0.5mm,0.5mm>*{};<2.2mm,2.2mm>*{}**@{-},
 <-0.48mm,0.48mm>*{};<-2.2mm,2.2mm>*{}**@{-},
 <0mm,0mm>*{\bu};<0mm,0mm>*{}**@{},
 <0.5mm,0.5mm>*{};<2.7mm,2.8mm>*{^{_1}}**@{},
 <-0.48mm,0.48mm>*{};<-2.7mm,2.8mm>*{^{_2}}**@{},
 \end{xy}\Ea
   \right\rangle
$$
$$
Q(1,2):= \sgn_2^{\ot d}\ot \id_1[d-1]=\mbox{span}\left\langle
\Ba{c}\begin{xy}
 <0mm,0.66mm>*{};<0mm,3mm>*{}**@{-},
 <0.39mm,-0.39mm>*{};<2.2mm,-2.2mm>*{}**@{-},
 <-0.35mm,-0.35mm>*{};<-2.2mm,-2.2mm>*{}**@{-},
 <0mm,0mm>*{\bu};<0mm,0mm>*{}**@{},
   <0.39mm,-0.39mm>*{};<2.9mm,-4mm>*{^{_2}}**@{},
   <-0.35mm,-0.35mm>*{};<-2.8mm,-4mm>*{^{_1}}**@{},
\end{xy}\Ea
=(-1)^{d}
\Ba{c}\begin{xy}
 <0mm,0.66mm>*{};<0mm,3mm>*{}**@{-},
 <0.39mm,-0.39mm>*{};<2.2mm,-2.2mm>*{}**@{-},
 <-0.35mm,-0.35mm>*{};<-2.2mm,-2.2mm>*{}**@{-},
 <0mm,0mm>*{\bu};<0mm,0mm>*{}**@{},
   <0.39mm,-0.39mm>*{};<2.9mm,-4mm>*{^{_1}}**@{},
   <-0.35mm,-0.35mm>*{};<-2.8mm,-4mm>*{^{_2}}**@{},
\end{xy}\Ea
\right\rangle
$$
by the ideal generated by the following relations
$$
\Ba{c}\resizebox{7mm}{!}{
\begin{xy}
 <0mm,0mm>*{\bu};<0mm,0mm>*{}**@{},
 <0mm,-0.49mm>*{};<0mm,-3.0mm>*{}**@{-},
 <0.49mm,0.49mm>*{};<1.9mm,1.9mm>*{}**@{-},
 <-0.5mm,0.5mm>*{};<-1.9mm,1.9mm>*{}**@{-},
 <-2.3mm,2.3mm>*{\bu};<-2.3mm,2.3mm>*{}**@{},
 <-1.8mm,2.8mm>*{};<0mm,4.9mm>*{}**@{-},
 <-2.8mm,2.9mm>*{};<-4.6mm,4.9mm>*{}**@{-},
   <0.49mm,0.49mm>*{};<2.7mm,2.3mm>*{^3}**@{},
   <-1.8mm,2.8mm>*{};<0.4mm,5.3mm>*{^2}**@{},
   <-2.8mm,2.9mm>*{};<-5.1mm,5.3mm>*{^1}**@{},
 \end{xy}}\Ea
 +
\Ba{c}\resizebox{7mm}{!}{\begin{xy}
 <0mm,0mm>*{\bu};<0mm,0mm>*{}**@{},
 <0mm,-0.49mm>*{};<0mm,-3.0mm>*{}**@{-},
 <0.49mm,0.49mm>*{};<1.9mm,1.9mm>*{}**@{-},
 <-0.5mm,0.5mm>*{};<-1.9mm,1.9mm>*{}**@{-},
 <-2.3mm,2.3mm>*{\bu};<-2.3mm,2.3mm>*{}**@{},
 <-1.8mm,2.8mm>*{};<0mm,4.9mm>*{}**@{-},
 <-2.8mm,2.9mm>*{};<-4.6mm,4.9mm>*{}**@{-},
   <0.49mm,0.49mm>*{};<2.7mm,2.3mm>*{^2}**@{},
   <-1.8mm,2.8mm>*{};<0.4mm,5.3mm>*{^1}**@{},
   <-2.8mm,2.9mm>*{};<-5.1mm,5.3mm>*{^3}**@{},
 \end{xy}}\Ea
 +
\Ba{c}\resizebox{7mm}{!}{\begin{xy}
 <0mm,0mm>*{\bu};<0mm,0mm>*{}**@{},
 <0mm,-0.49mm>*{};<0mm,-3.0mm>*{}**@{-},
 <0.49mm,0.49mm>*{};<1.9mm,1.9mm>*{}**@{-},
 <-0.5mm,0.5mm>*{};<-1.9mm,1.9mm>*{}**@{-},
 <-2.3mm,2.3mm>*{\bu};<-2.3mm,2.3mm>*{}**@{},
 <-1.8mm,2.8mm>*{};<0mm,4.9mm>*{}**@{-},
 <-2.8mm,2.9mm>*{};<-4.6mm,4.9mm>*{}**@{-},
   <0.49mm,0.49mm>*{};<2.7mm,2.3mm>*{^1}**@{},
   <-1.8mm,2.8mm>*{};<0.4mm,5.3mm>*{^3}**@{},
   <-2.8mm,2.9mm>*{};<-5.1mm,5.3mm>*{^2}**@{},
 \end{xy}}\Ea
 \ \ , \ \
\Ba{c}\resizebox{8.4mm}{!}{ \begin{xy}
 <0mm,0mm>*{\bu};<0mm,0mm>*{}**@{},
 <0mm,0.69mm>*{};<0mm,3.0mm>*{}**@{-},
 <0.39mm,-0.39mm>*{};<2.4mm,-2.4mm>*{}**@{-},
 <-0.35mm,-0.35mm>*{};<-1.9mm,-1.9mm>*{}**@{-},
 <-2.4mm,-2.4mm>*{\bu};<-2.4mm,-2.4mm>*{}**@{},
 <-2.0mm,-2.8mm>*{};<0mm,-4.9mm>*{}**@{-},
 <-2.8mm,-2.9mm>*{};<-4.7mm,-4.9mm>*{}**@{-},
    <0.39mm,-0.39mm>*{};<3.3mm,-4.0mm>*{^3}**@{},
    <-2.0mm,-2.8mm>*{};<0.5mm,-6.7mm>*{^2}**@{},
    <-2.8mm,-2.9mm>*{};<-5.2mm,-6.7mm>*{^1}**@{},
 \end{xy}}\Ea
 +
\Ba{c}\resizebox{8.4mm}{!}{ \begin{xy}
 <0mm,0mm>*{\bu};<0mm,0mm>*{}**@{},
 <0mm,0.69mm>*{};<0mm,3.0mm>*{}**@{-},
 <0.39mm,-0.39mm>*{};<2.4mm,-2.4mm>*{}**@{-},
 <-0.35mm,-0.35mm>*{};<-1.9mm,-1.9mm>*{}**@{-},
 <-2.4mm,-2.4mm>*{\bu};<-2.4mm,-2.4mm>*{}**@{},
 <-2.0mm,-2.8mm>*{};<0mm,-4.9mm>*{}**@{-},
 <-2.8mm,-2.9mm>*{};<-4.7mm,-4.9mm>*{}**@{-},
    <0.39mm,-0.39mm>*{};<3.3mm,-4.0mm>*{^2}**@{},
    <-2.0mm,-2.8mm>*{};<0.5mm,-6.7mm>*{^1}**@{},
    <-2.8mm,-2.9mm>*{};<-5.2mm,-6.7mm>*{^3}**@{},
 \end{xy}}\Ea
 +
\Ba{c}\resizebox{8.4mm}{!}{ \begin{xy}
 <0mm,0mm>*{\bu};<0mm,0mm>*{}**@{},
 <0mm,0.69mm>*{};<0mm,3.0mm>*{}**@{-},
 <0.39mm,-0.39mm>*{};<2.4mm,-2.4mm>*{}**@{-},
 <-0.35mm,-0.35mm>*{};<-1.9mm,-1.9mm>*{}**@{-},
 <-2.4mm,-2.4mm>*{\bu};<-2.4mm,-2.4mm>*{}**@{},
 <-2.0mm,-2.8mm>*{};<0mm,-4.9mm>*{}**@{-},
 <-2.8mm,-2.9mm>*{};<-4.7mm,-4.9mm>*{}**@{-},
    <0.39mm,-0.39mm>*{};<3.3mm,-4.0mm>*{^1}**@{},
    <-2.0mm,-2.8mm>*{};<0.5mm,-6.7mm>*{^3}**@{},
    <-2.8mm,-2.9mm>*{};<-5.2mm,-6.7mm>*{^2}**@{},
 \end{xy}}\Ea
 $$
 $$
 \Ba{c}\resizebox{5mm}{!}{\begin{xy}
 <0mm,2.47mm>*{};<0mm,0.12mm>*{}**@{-},
 <0.5mm,3.5mm>*{};<2.2mm,5.2mm>*{}**@{-},
 <-0.48mm,3.48mm>*{};<-2.2mm,5.2mm>*{}**@{-},
 <0mm,3mm>*{\bu};<0mm,3mm>*{}**@{},
  <0mm,-0.8mm>*{\bu};<0mm,-0.8mm>*{}**@{},
<-0.39mm,-1.2mm>*{};<-2.2mm,-3.5mm>*{}**@{-},
 <0.39mm,-1.2mm>*{};<2.2mm,-3.5mm>*{}**@{-},
     <0.5mm,3.5mm>*{};<2.8mm,5.7mm>*{^2}**@{},
     <-0.48mm,3.48mm>*{};<-2.8mm,5.7mm>*{^1}**@{},
   <0mm,-0.8mm>*{};<-2.7mm,-5.2mm>*{^1}**@{},
   <0mm,-0.8mm>*{};<2.7mm,-5.2mm>*{^2}**@{},
\end{xy}}\Ea
  -
\Ba{c}\resizebox{7mm}{!}{\begin{xy}
 <0mm,-1.3mm>*{};<0mm,-3.5mm>*{}**@{-},
 <0.38mm,-0.2mm>*{};<2.0mm,2.0mm>*{}**@{-},
 <-0.38mm,-0.2mm>*{};<-2.2mm,2.2mm>*{}**@{-},
<0mm,-0.8mm>*{\bu};<0mm,0.8mm>*{}**@{},
 <2.4mm,2.4mm>*{\bu};<2.4mm,2.4mm>*{}**@{},
 <2.77mm,2.0mm>*{};<4.4mm,-0.8mm>*{}**@{-},
 <2.4mm,3mm>*{};<2.4mm,5.2mm>*{}**@{-},
     <0mm,-1.3mm>*{};<0mm,-5.3mm>*{^1}**@{},
     <2.5mm,2.3mm>*{};<5.1mm,-2.6mm>*{^2}**@{},
    <2.4mm,2.5mm>*{};<2.4mm,5.7mm>*{^2}**@{},
    <-0.38mm,-0.2mm>*{};<-2.8mm,2.5mm>*{^1}**@{},
    \end{xy}}\Ea
  - (-1)^{d}
\Ba{c}\resizebox{7mm}{!}{\begin{xy}
 <0mm,-1.3mm>*{};<0mm,-3.5mm>*{}**@{-},
 <0.38mm,-0.2mm>*{};<2.0mm,2.0mm>*{}**@{-},
 <-0.38mm,-0.2mm>*{};<-2.2mm,2.2mm>*{}**@{-},
<0mm,-0.8mm>*{\bu};<0mm,0.8mm>*{}**@{},
 <2.4mm,2.4mm>*{\bu};<2.4mm,2.4mm>*{}**@{},
 <2.77mm,2.0mm>*{};<4.4mm,-0.8mm>*{}**@{-},
 <2.4mm,3mm>*{};<2.4mm,5.2mm>*{}**@{-},
     <0mm,-1.3mm>*{};<0mm,-5.3mm>*{^2}**@{},
     <2.5mm,2.3mm>*{};<5.1mm,-2.6mm>*{^1}**@{},
    <2.4mm,2.5mm>*{};<2.4mm,5.7mm>*{^2}**@{},
    <-0.38mm,-0.2mm>*{};<-2.8mm,2.5mm>*{^1}**@{},
    \end{xy}}\Ea
  - (-1)^{d+c}
\Ba{c}\resizebox{7mm}{!}{\begin{xy}
 <0mm,-1.3mm>*{};<0mm,-3.5mm>*{}**@{-},
 <0.38mm,-0.2mm>*{};<2.0mm,2.0mm>*{}**@{-},
 <-0.38mm,-0.2mm>*{};<-2.2mm,2.2mm>*{}**@{-},
<0mm,-0.8mm>*{\bu};<0mm,0.8mm>*{}**@{},
 <2.4mm,2.4mm>*{\bu};<2.4mm,2.4mm>*{}**@{},
 <2.77mm,2.0mm>*{};<4.4mm,-0.8mm>*{}**@{-},
 <2.4mm,3mm>*{};<2.4mm,5.2mm>*{}**@{-},
     <0mm,-1.3mm>*{};<0mm,-5.3mm>*{^2}**@{},
     <2.5mm,2.3mm>*{};<5.1mm,-2.6mm>*{^1}**@{},
    <2.4mm,2.5mm>*{};<2.4mm,5.7mm>*{^1}**@{},
    <-0.38mm,-0.2mm>*{};<-2.8mm,2.5mm>*{^2}**@{},
    \end{xy}}\Ea
 - (-1)^{c}
\Ba{c}\resizebox{7mm}{!}{\begin{xy}
 <0mm,-1.3mm>*{};<0mm,-3.5mm>*{}**@{-},
 <0.38mm,-0.2mm>*{};<2.0mm,2.0mm>*{}**@{-},
 <-0.38mm,-0.2mm>*{};<-2.2mm,2.2mm>*{}**@{-},
<0mm,-0.8mm>*{\bu};<0mm,0.8mm>*{}**@{},
 <2.4mm,2.4mm>*{\bu};<2.4mm,2.4mm>*{}**@{},
 <2.77mm,2.0mm>*{};<4.4mm,-0.8mm>*{}**@{-},
 <2.4mm,3mm>*{};<2.4mm,5.2mm>*{}**@{-},
     <0mm,-1.3mm>*{};<0mm,-5.3mm>*{^1}**@{},
     <2.5mm,2.3mm>*{};<5.1mm,-2.6mm>*{^2}**@{},
    <2.4mm,2.5mm>*{};<2.4mm,5.7mm>*{^1}**@{},
    <-0.38mm,-0.2mm>*{};<-2.8mm,2.5mm>*{^2}**@{},
    \end{xy}}\Ea
$$
Its minimal resolution $\HoLBcd$ is a 1-oriented dg free properad generated by the following (skew)symmetric corollas of degree $1 +c(1-m)+d(1-n)$
\Beq\label{2: symmetries of HoLiebcd corollas}
\Ba{c}\resizebox{17mm}{!}{\begin{xy}
 <0mm,0mm>*{\bu};<0mm,0mm>*{}**@{},
 <-0.6mm,0.44mm>*{};<-8mm,5mm>*{}**@{-},
 <-0.4mm,0.7mm>*{};<-4.5mm,5mm>*{}**@{-},
 <0mm,0mm>*{};<1mm,5mm>*{\ldots}**@{},
 <0.4mm,0.7mm>*{};<4.5mm,5mm>*{}**@{-},
 <0.6mm,0.44mm>*{};<8mm,5mm>*{}**@{-},
   <0mm,0mm>*{};<-10.5mm,5.9mm>*{^{\sigma(1)}}**@{},
   <0mm,0mm>*{};<-4mm,5.9mm>*{^{\sigma(2)}}**@{},
   <0mm,0mm>*{};<10.0mm,5.9mm>*{^{\sigma(m)}}**@{},
 <-0.6mm,-0.44mm>*{};<-8mm,-5mm>*{}**@{-},
 <-0.4mm,-0.7mm>*{};<-4.5mm,-5mm>*{}**@{-},
 <0mm,0mm>*{};<1mm,-5mm>*{\ldots}**@{},
 <0.4mm,-0.7mm>*{};<4.5mm,-5mm>*{}**@{-},
 <0.6mm,-0.44mm>*{};<8mm,-5mm>*{}**@{-},
   <0mm,0mm>*{};<-10.5mm,-6.9mm>*{^{\tau(1)}}**@{},
   <0mm,0mm>*{};<-4mm,-6.9mm>*{^{\tau(2)}}**@{},
   <0mm,0mm>*{};<10.0mm,-6.9mm>*{^{\tau(n)}}**@{},
 \end{xy}}\Ea
=(-1)^{c|\sigma|+d|\tau|}
\Ba{c}\resizebox{14mm}{!}{\begin{xy}
 <0mm,0mm>*{\bu};<0mm,0mm>*{}**@{},
 <-0.6mm,0.44mm>*{};<-8mm,5mm>*{}**@{-},
 <-0.4mm,0.7mm>*{};<-4.5mm,5mm>*{}**@{-},
 <0mm,0mm>*{};<-1mm,5mm>*{\ldots}**@{},
 <0.4mm,0.7mm>*{};<4.5mm,5mm>*{}**@{-},
 <0.6mm,0.44mm>*{};<8mm,5mm>*{}**@{-},
   <0mm,0mm>*{};<-8.5mm,5.5mm>*{^1}**@{},
   <0mm,0mm>*{};<-5mm,5.5mm>*{^2}**@{},
   <0mm,0mm>*{};<4.5mm,5.5mm>*{^{m\hspace{-0.5mm}-\hspace{-0.5mm}1}}**@{},
   <0mm,0mm>*{};<9.0mm,5.5mm>*{^m}**@{},
 <-0.6mm,-0.44mm>*{};<-8mm,-5mm>*{}**@{-},
 <-0.4mm,-0.7mm>*{};<-4.5mm,-5mm>*{}**@{-},
 <0mm,0mm>*{};<-1mm,-5mm>*{\ldots}**@{},
 <0.4mm,-0.7mm>*{};<4.5mm,-5mm>*{}**@{-},
 <0.6mm,-0.44mm>*{};<8mm,-5mm>*{}**@{-},
   <0mm,0mm>*{};<-8.5mm,-6.9mm>*{^1}**@{},
   <0mm,0mm>*{};<-5mm,-6.9mm>*{^2}**@{},
   <0mm,0mm>*{};<4.5mm,-6.9mm>*{^{n\hspace{-0.5mm}-\hspace{-0.5mm}1}}**@{},
   <0mm,0mm>*{};<9.0mm,-6.9mm>*{^n}**@{},
 \end{xy}}\Ea \ \ \forall \sigma\in \bS_m, \forall\tau\in \bS_n
\Eeq
and has the differential
given on the generators by
\Beq\label{LBk_infty}
\delta
\Ba{c}\resizebox{14mm}{!}{\begin{xy}
 <0mm,0mm>*{\bu};<0mm,0mm>*{}**@{},
 <-0.6mm,0.44mm>*{};<-8mm,5mm>*{}**@{-},
 <-0.4mm,0.7mm>*{};<-4.5mm,5mm>*{}**@{-},
 <0mm,0mm>*{};<-1mm,5mm>*{\ldots}**@{},
 <0.4mm,0.7mm>*{};<4.5mm,5mm>*{}**@{-},
 <0.6mm,0.44mm>*{};<8mm,5mm>*{}**@{-},
   <0mm,0mm>*{};<-8.5mm,5.5mm>*{^1}**@{},
   <0mm,0mm>*{};<-5mm,5.5mm>*{^2}**@{},
   <0mm,0mm>*{};<4.5mm,5.5mm>*{^{m\hspace{-0.5mm}-\hspace{-0.5mm}1}}**@{},
   <0mm,0mm>*{};<9.0mm,5.5mm>*{^m}**@{},
 <-0.6mm,-0.44mm>*{};<-8mm,-5mm>*{}**@{-},
 <-0.4mm,-0.7mm>*{};<-4.5mm,-5mm>*{}**@{-},
 <0mm,0mm>*{};<-1mm,-5mm>*{\ldots}**@{},
 <0.4mm,-0.7mm>*{};<4.5mm,-5mm>*{}**@{-},
 <0.6mm,-0.44mm>*{};<8mm,-5mm>*{}**@{-},
   <0mm,0mm>*{};<-8.5mm,-6.9mm>*{^1}**@{},
   <0mm,0mm>*{};<-5mm,-6.9mm>*{^2}**@{},
   <0mm,0mm>*{};<4.5mm,-6.9mm>*{^{n\hspace{-0.5mm}-\hspace{-0.5mm}1}}**@{},
   <0mm,0mm>*{};<9.0mm,-6.9mm>*{^n}**@{},
 \end{xy}}\Ea
\ \ = \ \
 \sum_{[1,\ldots,m]=I_1\sqcup I_2\atop
 {|I_1|\geq 0, |I_2|\geq 1}}
 \sum_{[1,\ldots,n]=J_1\sqcup J_2\atop
 {|J_1|\geq 1, |J_2|\geq 1}
}\hspace{0mm}
(-1)^{\# I_1\# J_2 + \# I_1 + \# J_2 + c\text{sgn}(I_1,I_2) + d\text{sgn}(J_1,J_2)  }
\Ba{c}\resizebox{22mm}{!}{ \begin{xy}
 <0mm,0mm>*{\bu};<0mm,0mm>*{}**@{},
 <-0.6mm,0.44mm>*{};<-8mm,5mm>*{}**@{-},
 <-0.4mm,0.7mm>*{};<-4.5mm,5mm>*{}**@{-},
 <0mm,0mm>*{};<0mm,5mm>*{\ldots}**@{},
 <0.4mm,0.7mm>*{};<4.5mm,5mm>*{}**@{-},
 <0.6mm,0.44mm>*{};<12.4mm,4.8mm>*{}**@{-},
     <0mm,0mm>*{};<-2mm,7mm>*{\overbrace{\ \ \ \ \ \ \ \ \ \ \ \ }}**@{},
     <0mm,0mm>*{};<-2mm,9mm>*{^{I_1}}**@{},
 <-0.6mm,-0.44mm>*{};<-8mm,-5mm>*{}**@{-},
 <-0.4mm,-0.7mm>*{};<-4.5mm,-5mm>*{}**@{-},
 <0mm,0mm>*{};<-1mm,-5mm>*{\ldots}**@{},
 <0.4mm,-0.7mm>*{};<4.5mm,-5mm>*{}**@{-},
 <0.6mm,-0.44mm>*{};<8mm,-5mm>*{}**@{-},
      <0mm,0mm>*{};<0mm,-7mm>*{\underbrace{\ \ \ \ \ \ \ \ \ \ \ \ \ \ \
      }}**@{},
      <0mm,0mm>*{};<0mm,-10.6mm>*{_{J_1}}**@{},
 <13mm,5mm>*{};<13mm,5mm>*{\bu}**@{},
 <12.6mm,5.44mm>*{};<5mm,10mm>*{}**@{-},
 <12.6mm,5.7mm>*{};<8.5mm,10mm>*{}**@{-},
 <13mm,5mm>*{};<13mm,10mm>*{\ldots}**@{},
 <13.4mm,5.7mm>*{};<16.5mm,10mm>*{}**@{-},
 <13.6mm,5.44mm>*{};<20mm,10mm>*{}**@{-},
      <13mm,5mm>*{};<13mm,12mm>*{\overbrace{\ \ \ \ \ \ \ \ \ \ \ \ \ \ }}**@{},
      <13mm,5mm>*{};<13mm,14mm>*{^{I_2}}**@{},
 <12.4mm,4.3mm>*{};<8mm,0mm>*{}**@{-},
 <12.6mm,4.3mm>*{};<12mm,0mm>*{\ldots}**@{},
 <13.4mm,4.5mm>*{};<16.5mm,0mm>*{}**@{-},
 <13.6mm,4.8mm>*{};<20mm,0mm>*{}**@{-},
     <13mm,5mm>*{};<14.3mm,-2mm>*{\underbrace{\ \ \ \ \ \ \ \ \ \ \ }}**@{},
     <13mm,5mm>*{};<14.3mm,-4.5mm>*{_{J_2}}**@{},
 \end{xy}}\Ea
\Eeq
where $\text{sign}(I_1,I_2)$ (resp., $\text{sign}(J_1,J_2)$) stands for the parity of the permutation $I\rar I_1\sqcup I_2$ (resp., $J\rar J_1\sqcup J_2$). The case $c=d=0$ corresponds to ordinary strong homotopy Lie bialgebras \cite{MaVo,Va}, while the case $c=1$, $d=0$ to formal Poisson structures on graded vector spaces $V$ viewed as linear manifolds \cite{Me1}.

\sip

A $(k+1)$-oriented generalization of $\HoLBcd$ is quite straightforward: the prop $\HoLBcd^{(k+1)\text{-or}}$ is a free $(k+1)$-oriented
prop generated by corollas with the same symmetries and degrees as in the case of  $\HoLBcd$, but now with each leg decorated
with $k$ extra orientations,
$$
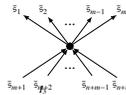
\Ba{c}\resizebox{18mm}{!}
{\begin{tikzpicture}[baseline=-1ex]
\node[] at (-1.2,0.8) {$^{\bar{\fs}_1}$};
\node[] at (-0.6,0.8) {$^{\bar{\fs}_2}$};
\node[] at (1.2,0.8) {$^{\bar{\fs}_m}$};
\node[] at (0.6,0.8) {$^{\bar{\fs}_{m-1}}$};
\node[] at (-1.2,-0.9) {$_{\bar{\fs}_{m+1}}$};
\node[] at (-0.6,-0.9) {$_{\bar{\fs}_{m+2}}$};
\node[] at (1.2,-0.9) {$_{\bar{\fs}_{n+m}}$};
\node[] at (0.6,-0.9) {$_{\bar{\fs}_{n+m-1}}$};
\node[] at (0.0,0.5) {$...$};
\node[] at (0.0,-0.5) {$...$};
\node[] at (-0.6,-1.0) {$_{\scriptstyle I_3}$};
\node[int] (0) at (0,0) {};
\node[] (1) at (-1.2,0.8) {};
\node[] (2) at (-0.6,0.8) {};
\node[] (3) at (0.6,0.8) {};
\node[] (4) at (1.2,0.8) {};
\node[int] (0) at (0,0) {};
\node[] (5) at (-1.2,-0.8) {};
\node[] (6) at (-0.6,-0.8) {};
\node[] (7) at (0.6,-0.8) {};
\node[] (8) at (1.2,-0.8) {};
\draw (1) edge[latex-] (0);
\draw (2) edge[latex-] (0);
\draw (3) edge[latex-] (0);
\draw (4) edge[latex-] (0);
\draw (0) edge[latex-] (5);
\draw (0) edge[latex-] (6);
\draw (0) edge[latex-] (7);
\draw (0) edge[latex-] (8);
\end{tikzpicture}}\Ea
$$
 subject to the condition that there is at least one ingoing edge and one outgoing edge in each of the new directions. The differential is given by the  same formula as in the case of $\HoLBcd$ except that now we sum over all possible (and admissible) new orientations attached to the new edge
$$
\delta
\Ba{c}\resizebox{18mm}{!}
{\begin{tikzpicture}[baseline=-1ex]
\node[] at (-1.2,0.8) {$^{\bar{\fs}_1}$};
\node[] at (-0.6,0.8) {$^{\bar{\fs}_2}$};
\node[] at (1.2,0.8) {$^{\bar{\fs}_m}$};
\node[] at (0.6,0.8) {$^{\bar{\fs}_{m-1}}$};
\node[] at (-1.2,-0.9) {$_{\bar{\fs}_{m+1}}$};
\node[] at (-0.6,-0.9) {$_{\bar{\fs}_{m+2}}$};
\node[] at (1.2,-0.9) {$_{\bar{\fs}_{n+m}}$};
\node[] at (0.6,-0.9) {$_{\bar{\fs}_{n+m-1}}$};
\node[] at (0.0,0.5) {$...$};
\node[] at (0.0,-0.5) {$...$};
\node[int] (0) at (0,0) {};
\node[] (1) at (-1.2,0.8) {};
\node[] (2) at (-0.6,0.8) {};
\node[] (3) at (0.6,0.8) {};
\node[] (4) at (1.2,0.8) {};
\node[int] (0) at (0,0) {};
\node[] (5) at (-1.2,-0.8) {};
\node[] (6) at (-0.6,-0.8) {};
\node[] (7) at (0.6,-0.8) {};
\node[] (8) at (1.2,-0.8) {};
\draw (1) edge[latex-] (0);
\draw (2) edge[latex-] (0);
\draw (3) edge[latex-] (0);
\draw (4) edge[latex-] (0);
\draw (0) edge[latex-] (5);
\draw (0) edge[latex-] (6);
\draw (0) edge[latex-] (7);
\draw (0) edge[latex-] (8);
\end{tikzpicture}}\Ea
=
 \sum_{[1,\ldots,m]=I_1\sqcup I_2\atop
 {|I_1|\geq 0, |I_2|\geq 1}}
 \sum_{[m+1,\ldots,m+n]=J_1\sqcup J_2\atop
 {|J_1|\geq 1, |J_2|\geq 1}
} \sum_{\bar{\fs}\in \f_k}
(-1)^{\# I_1\# J_2 + \# I_2 + \# J_2 + c\cdot \text{sgn}(I_1,I_2) + d\cdot \text{sgn}(J_1,J_2)  }
\Ba{c}\resizebox{24mm}{!}
{\begin{tikzpicture}[baseline=-1ex]
\node[] at (1.0,0.74) {$^{ \scriptstyle \bar{\fs}}$};
\node[] at (-0.2,1.1) {$^{ \scriptstyle \bar{\fs}_i,\ i\in I_1}$};
\node[] at (-0.3,0.8) {$\overbrace{\ \ \  \ \ \ \ \ \ \ \ \ \ \  \ \ \ \ \ }$};
\node[] at (-0.0,-1.2) {$_{\scriptstyle \bar{\fs}_i,\ i\in J_1}$};
\node[] at (-0.0,-0.9) {$\underbrace{\ \ \ \ \ \ \ \ \ \ \ \ \ \ \ \ \ \ \ \ \  \ \ \ \ \ }$};
\node[] at (1.8,2.1) {$^{\scriptstyle \bar{\fs}_i,\ i\in I_2}$};
\node[] at (1.8,1.9) {$\overbrace{\ \ \  \ \ \ \ \ \ \ \ \ \ \  \ \ \ \ \ }$};
\node[] at (1.8,0.0) {$_{\scriptstyle \bar{\fs}_i,\ i\in J_2}$};
\node[] at (1.8,0.27) {$\underbrace{\ \ \ \ \ \ \ \ \ \ \ \ \ \ \ \ \ \ \  }$};
\node[] at (0.0,0.5) {$...$};
\node[] at (0.0,-0.5) {$...$};
\node[] at (1.8,1.7) {$...$};
\node[] at (1.8,0.5) {$...$};
\node[int] (0) at (0,0) {};
\node[] (1) at (-1.2,0.8) {};
\node[] (2) at (-0.6,0.8) {};
\node[] (3) at (0.6,0.8) {};
\node[int] (4) at (1.8,1.1) {};
\node[int] (0) at (0,0) {};
\node[] (5) at (-1.2,-0.8) {};
\node[] (6) at (-0.6,-0.8) {};
\node[] (7) at (0.6,-0.8) {};
\node[] (8) at (1.2,-0.8) {};
\node[] (u1) at (1.0,1.9) {};
\node[] (u2) at (1.4,1.9) {};
\node[] (u3) at (2.1,1.9) {};
\node[] (u4) at (2.6,1.9) {};
\node[] (d1) at (1.0,0.3) {};
\node[] (d2) at (1.4,0.3) {};
\node[] (d3) at (2.1,0.3) {};
\node[] (d4) at (2.6,0.3) {};
\draw (1) edge[latex-] (0);
\draw (2) edge[latex-] (0);
\draw (3) edge[latex-] (0);
\draw (4) edge[latex-] (0);
\draw (0) edge[latex-] (5);
\draw (0) edge[latex-] (6);
\draw (0) edge[latex-] (7);
\draw (0) edge[latex-] (8);
\draw (u1) edge[latex-] (4);
\draw (u2) edge[latex-] (4);
\draw (u3) edge[latex-] (4);
\draw (u4) edge[latex-] (4);
\draw (4) edge[latex-] (d1);
\draw (4) edge[latex-] (d2);
\draw (4) edge[latex-] (d3);
\draw (4) edge[latex-] (d4);
\end{tikzpicture}}\Ea
$$
The homotopy theory of such props can be highly non-trivial.
As we discuss in more detail below in \S 5,
the automorphism group of $\HoLB_{c,d}^{(c+d-1)\text{-or}}$
is equal to  the Grothendieck-Teichm\"uller group $GRT=GRT_1\ltimes \K^*$ (for any $c,d\in \Z$
with $c+d\geq 2$) and hence this prop
can be a foundation for a rich deformation quantization theory in {\em every}\, geometric dimension $c+d\geq 2$.
This fact was one of our main motivations to introduce and study the multi-oriented props.

\sip

Let $K$ be the differential closure of the ideal in $\HoLBcd^{(k+1)\text{-or}}$
generated by all corollas with total arity $\geq 4$ and denote the quotient by
\Beq\label{3new: p from Holib^k+1 to Lieb^k+1}
\LBcd^{(k+1)\text{-or}}=\HoLBcd^{(k+1)\text{-or}}/K.
\Eeq
This prop gives us a $(k+1)$ oriented version of $\LBcd$.  It is quite easy to see that the automorphism group of the prop
$\LBcd^{(c+d-1)\text{-or}}$ with $c,d\geq 1$ and $c+d\geq 3$ is almost trivial --- it is equal to $\K^*$ acting by rescalings of the generators. In particular, the induced action of $GRT_1$ on $\LBcd^{(k+1)\text{-or}}$ with $c,d\geq 1$ and $c+d\geq 3$ is trivial; this is in sharp contrast to the 1-oriented case $c=d=1$ where that action remains highly non-trivial.

\sip
Our main purpose in this paper is to introduce the notion of a representation of a multi-oriented prop in the category of dg vector spaces (with branes) which is done in the next section.

\bip

\bip

{\Large
\section{\bf Multidirected endomorphism  prop and homotopy algebras with branes}
}

\bip



\subsection{Tensor algebra of infinite-dimensional vector spaces}\label{App: subsec on inf dim vector spaces}
By a {\em countably  infinite-dimensional}\, graded vector space $V$ we understand in this paper
any {\em direct}\, limit $\displaystyle V:=\lim_{\longrightarrow} V_p$ of a direct system of finite dimensional vector spaces $V_p$, $p\geq 1$,
\Beq\label{App: injections i_n}
V_0\lon  V_1 \stackrel{i_1}{\lon} V_2 \stackrel{i_2}{\lon} \ldots \stackrel{i_{p-1}}{\lon} V_p \stackrel{i_{p}}{\lon} V_{p+1} \stackrel{i_{p+1}}{\lon} \ldots.
\Eeq
where  all arrows $i_p$ are proper injections. For example, $\displaystyle \K^\infty=\lim_{\longrightarrow\atop n} \K^n$ with
$$
0\lon  \K \stackrel{i_1}{\lon} \K^2 \stackrel{i_2}{\lon} \ldots \stackrel{i_p}{\lon} \K^p \stackrel{i_{p+1}}{\lon} \K^{p+1} \stackrel{i_{p+2}}{\lon} \ldots, \ \ \ i_p(a_1,\ldots,a_p):=(a_1,\ldots,a_p,0),
$$
is an example of a countably infinite-dimensional vector space.

\sip

Next we define  (non-countably) infinite-dimensional vector spaces
$$
\Hom(\ot^r V, \ot^l V)
= \lim_{\longleftarrow\atop (p_1,\ldots,p_r)}\left( V_{p_1}^*\ot \ldots \ot V_{p_r}^* \ot V^{\ot^l}\right).
$$
which are equipped with the standard projective limit topology.
An element $f\in \Hom(\ot^r V, \ot^l V)$ is called a {\em linear map}
$$
f: \ot^r V \lon \ot^l V
$$
from $\ot^r V$ to  $\ot^l V$. Such maps can be composed (no divergences) so that one has a well-defined
endomorphism prop
$$
\cE nd_V=\left\{ \cE nd_V(l,r):=\Hom(\ot^r V, \ot^l V)\right\}
$$
associated to $V$ and hence talk about representations of ordinary props in $V$.

\sip

Note that $\Hom(V,V)$ can contain infinite sums of the form
$$
\sum_{n,m=1}^\infty a_m^n e^m\ot e_n, \ \ e^m\in V_m^*, e_n\in V_n, a_m^n\in \K, \ \ \text{all}\ a_m^n\neq 0,
$$
so that the trace operation on $\Hom(V,V)$ (and hence on $\Hom(\ot^k V, \ot^l V)$) is not well-defined in general so that $V$ can not be used for representations of {\em wheeled}\, props.

\sip

If one has a collection of $k$ infinite-dimensional vector spaces,
$\displaystyle V_\tau=\lim_{\lon} V_{\tau,p}$, $\tau\in [k]$, one can define similarly a $k$-coloured
endomorphism prop $\cE nd_{V_1,\ldots,V_k}$ based on topological $\bS$-modules
\Beq\label{App: formula for End V_1..V_k}
\cE nd_{V_1,\ldots,V_k}=\left\{ \Hom(V_{1}^{\ot r_1}\ot
 \ldots \ot V_{k}^{\ot r_k},V_{1}^{\ot l_1}\ot
 \ldots \ot V_{k}^{\ot l_k}):=
\lim_{\longleftarrow\atop (n_1,\ldots,n_r)}\left( (V_{1,p_1}^*)^{\ot r_1}\ot
 \ldots \ot (V_{k, p_k}^*)^{\ot r_k} \ot  V_{1}^{\ot l_1}\ot
 \ldots \ot V_{k}^{\ot l_k}    \right)\right\}.
\Eeq
Its elements give us linear maps
$
V_{1}^{\ot r_1}\ot
 \ldots \ot V_{k}^{\ot r_k}  \lon   V_{1}^{\ot l_1}\ot
 \ldots \ot V_{k}^{\ot l_k}
$
and hence can be used to define a representation of a  $k$-coloured prop.

\subsection{\bf An infinite-dimensional graded vector space with $k$ branes}\label{App: subsec on vector space with branes} Let $\displaystyle V:=\lim_{\longrightarrow\atop p} V_p$ be a  countably infinite-dimensional graded vector space. It is called a
 vector space {\em  with $k$
branes}\, (and denoted by $(V,W_1,\ldots,W_k)$ or simply by $V^{k\text{-br}}$) if the following conditions hold:
\Bi
\item[(i)] $V$
comes equipped with a descending filtration
$$
V=F^0V \supset F^1V \supset F^2V \supset\ldots \supset F^pV\supset F^{p+1}V \supset \ldots
$$
such that each quotient vector space $V/F^{p}V$ is finite-dimensional and is isomorphic to $V_p$ for any $p\in \N$,
\item[(ii)] For any $p\geq 0$ we have $k$ different non trivial direct sum decompositions $V_p=W_{\tau,p}^+ \oplus W_{\tau,p}^-$,  $\tau\in [k]$, which are compatible with the given injections $i_p: V_p\rar V_{p+1}$,
    $$
    i_p(W_{\tau,p}^\pm)\subset W_{\tau,p+1}^\pm
    $$
\Ei
Note that the inclusion $F^{p+1}V\subset F^{p}V$ induces a projection
$$
\pi_{p+1}: V_{p+1}\equiv F^{p+1}V/F^{p+2}V \lon V_p:=F^{p}V/F^{p+1}V
$$
providing  us with an inverse system of finite-dimensional vector spaces,
$$
\ldots \stackrel{\pi_{p+1}}{\lon} W_{\tau,p}^\pm \stackrel{\pi_{p}}{\lon} W_{\tau,p-1}^\pm \stackrel{\pi_{p-2}}{\lon}\ldots \lon 0
$$
Hence can consider two limits for branes (and their intersections, see below), the direct  and projective ones,
$$
W_{\tau}^\pm := \lim_{\longrightarrow\atop p} W_{\tau,p}^\pm \subset  \hat{W}_{\tau}^\pm := \lim_{\longleftarrow\atop p} W_{\tau,p}^\pm,\ \ \ \
(W_{\tau}^\pm)^* = \lim_{\longleftarrow\atop p} (W_{\tau,p}^\pm)^* \supset (\hat{W}_{\tau}^\pm)^* =  \lim_{\longrightarrow\atop p} (W_{\tau,p}^\pm)^*, \ \ \ \forall \tau\in [k].
$$
Note that the spaces $W_{\tau}^\pm$ and $(\hat{W}_{\tau}^\pm)^*$ are always countably dimensional, while $(W_{\tau}^\pm)^*$ and $\hat{W}_{\tau}^\pm$ are, in general, not
(but as a compensation they come equipped with nice topologies). Note also that
$$
((\hat{W}_{\tau}^\pm)^*)^*=\hat{W}_{\tau}^\pm, \ \ \ \ (({W}_{\tau}^\pm)^*)^*={W}_{\tau}^\pm.
$$

\sip

To define a suitable multi-oriented endomorphism prop out of an infinite-dimensional vector space with $k$ branes,
one has to work with both types of completions simultaneously. This fact motivates the extra filtration condition  (i) in the definition of $V^{k\text{-br}}$ above.

\subsubsection{\bf Basic example}\label{App: Example of an inf dim space with branes} Let $\{x_1,x_2,\ldots,\}$ be a countably infinite set of formal variables of some homological degrees $|x_i|\in \Z$, $i\in \N_{\geq 1}$.
The graded vector space
$$
V= \text{span}\langle x_1, x_2, \ldots\ \rangle
$$
 is a typical example
of an infinite-dimensional vector space satisfying conditions (i) and (ii) above
with
$$
F^pV=\text{span}\langle x_i \rangle_{i\geq p+1}, \ \ \  V_p=\text{span}\langle x_1, x_2, \ldots, x_p\rangle.
$$
Let us choose $k$ injections of countably infinite sets (i.e.\ $k$ pairs of disjoint countably infinite subsets
of $\N_{\geq 1}$)
$$
f_\tau: \N_{\geq 1} \oplus \N_{\geq 1} \lon \N_{\geq 1}, \ \ \ \tau\in [k],
$$
and define a family of finite-dimensional vector spaces (equipped with the $\Z$-grading induced in the obvious way from the homological grading of the formal variables $(x_1,x_2,\ldots, x_p)$)
$$
W_{\tau,p}^+:=\text{span}\left\langle \pi_+\circ f_\tau^{-1}\{1,2,\ldots,p\}  \right \rangle,
\ \ \ \ \ \
W_{\tau,p}^-:=\text{span}\left\langle \pi_-\circ f_\tau^{-1}\{1,2,\ldots,p\}  \right \rangle
$$
where $\pi_\pm: \N_{\geq 1} \oplus \N_{\geq 1} \rar \N_{\geq 1}$ is the projection to the first/second summand.
The resulting data $(V, W_1^{\pm},\ldots, W_k^\pm)$ is an example of an infinite-dimensional graded vector space with $k$ branes.

\subsection{Finite-dimensional case} {\em A finite-dimensional vector space
with $k$ branes}\, is simply a finite-dimensional vector space $V$ equipped with $k$ direct sum decompositions $V=W^+_\tau\oplus W^-_\tau$, $\tau\in [k]$. We shall be most interested in the infinite-dimensional case and hence use all the time use the direct/projective limit notation introduced in the previous section. The finite dimensional version
fits into that notation as a special case when $V_p=V$, $W_{\tau,p}^\pm=W^\pm_\tau$ for all $p$.

\subsection{The (simplest non-trivial) case of a 2-directed endomorphism prop}
Let $V^{1\text{-br}}$ be an infinite-dimensional vector space with one brane, and consider
$$
W^+=\lim_{\lon \atop p} W_p^+
, \ \ \
\hat{W}^-=\lim_{\longleftarrow \atop p} W_p^-, \ (W^+)^*=\lim_{\longleftarrow \atop p} (W_p^+)^*
, \ \ \
(\hat{W}^-)^*=\lim_{\longrightarrow \atop p} (W_p^-)^*
$$
Note that out of these four spaces only $W^+$ and $(\hat{W}^{-})^*$ are always countably dimensional.
Introduce next an $\cS^{(2)}$-module, that is a functor
$$
\Ba{rccl}
\cE nd_{V^{1\text{-br}}}: & \cS^{(2)} & \lon & \text{Category of graded vector spaces} \\
              & (I, \fs) & \lon &  \cE nd_{V^{2\text{-br}}}(I,\fs)
              \Ea,
$$
as follows. First notice that one can identify any 2-oriented set\footnote{Here we denote the elements of $[1^+]$ by $\bar{0}$ and $\bar{1}$ so that the value
$\fs_i$ of the map $\fs$ on an element $i\in I$ is itself a map of sets $\fs_i:\{\bar{0},\bar{1}\} \rar \{out,in\}$.}
$$
\left(I,\fs: I\rar \f_{1^+}\right)\equiv \left(I, \fs_{\bar{0}}: I\rar \{out,in\}, \fs_{\bar{1}}: I\rar \{out,in\}\ \text{such\ that}\ \fs_{\bar{0}}(i):= \fs_i(\bar{0}),  \fs_{\bar{1}}(i):= \fs_i(\bar{1})\ \forall i\in I \right)
$$
 with a 2-directed corolla (the subscript $0$ indicates the basic orientation)
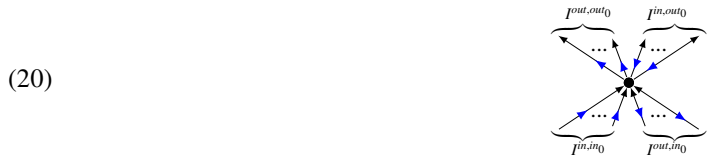
\Beq\label{6: 2-dir corolla}
\Ba{c}\resizebox{23mm}{!}
{\begin{tikzpicture}[baseline=-1ex]
\node[] at (-0.6,1.0) {$^{\scriptstyle I^{out,out_{{0}}}}$};
\node[] at (-0.7,0.8) {$\overbrace{\  \ \ \ \ \ \ \ \ \ \ }$};
\node[] at (-0.45,0.5) {$...$};
\node[] at (0.6,1.0) {$^{\scriptstyle I^{in,out_0}}$};
\node[] at (0.7,0.8) {$\overbrace{\  \ \ \ \ \ \ \ \ \ \ }$};
\node[] at (0.45,0.5) {$...$};
\node[] at (-0.6,-1.0) {$_{\scriptstyle I^{in,in_0}}$};
\node[] at (-0.7,-0.8) {$\underbrace{\  \ \ \ \ \ \ \ \ \ \ }$};
\node[] at (-0.45,-0.5) {$...$};
\node[] at (0.6,-1.0) {$_{\scriptstyle I^{out,in_0}}$};
\node[] at (0.7,-0.8) {$\underbrace{\  \ \ \ \ \ \ \ \ \ \ }$};
\node[] at (0.45,-0.5) {$...$};
\node[int] (0) at (0,0) {};
\node[] (1) at (-1.2,0.8) {};
\node[] (2) at (-0.3,0.8) {};
\node[] (3) at (0.3,0.8) {};
\node[] (4) at (1.2,0.8) {};
\node[int] (0) at (0,0) {};
\node[] (5) at (-1.2,-0.8) {};
\node[] (6) at (-0.3,-0.8) {};
\node[] (7) at (0.3,-0.8) {};
\node[] (8) at (1.2,-0.8) {};
\draw (1) edge[latex-,leftblue] (0);
\draw (2) edge[latex-,leftblue] (0);
\draw (3) edge[latex-,rightblue] (0);
\draw (4) edge[latex-,rightblue] (0);
\draw (0) edge[latex-,leftblue] (5);
\draw (0) edge[latex-,leftblue] (6);
\draw (0) edge[latex-,rightblue] (7);
\draw (0) edge[latex-,rightblue] (8);
\end{tikzpicture}}\Ea
\Eeq
where
$$
 I^{out,out_0}:= \fs_{\bar{1}}^{-1}(out) \sqcup \fs_{\bar{0}}^{-1}(out),
 \ \ \
 I^{in,out_0}:= \fs_{\bar{1}}^{-1}(in) \sqcup \fs_{\bar{0}}^{-1}(out), \ \ \
 I^{out,in_0}:= \fs_{\bar{1}}^{-1}(out) \sqcup \fs_{\bar{0}}^{-1}(in), \ \ \
 I^{in,in_0}:= \fs_{\bar{1}}^{-1}(in) \sqcup \fs_{\bar{0}}^{-1}(in), \ \ \
$$
Set
$$
\# I^{out,out_0}=:m_1, \ \  \#I^{in,out_0}=:m_2,\ \ \  \#I^{in,in_0}=:n_1, \ \
\#I^{out,in_0}=:n_2.
$$

Next we define\footnote{Here we use the facts that for any vector space $M$ and any inverse system of finite-dimensional vector spaces $\{N_i\}$ one has $\displaystyle \lim_{\longleftarrow} \Hom(N_i,M)\cong \Hom(\lim_{\lon}N_i, M)$ and
 $\displaystyle \lim_{\longleftarrow }\Hom(M, N_i) \cong \Hom(M, \lim_{\longleftarrow }N_i)$, while  $\displaystyle \lim_{\longleftarrow} (N_i\ot M)\cong (\lim_{\longleftarrow }N_i)\ot M$ only if $M$ is finite-dimensional. On the other hand, for {\em any}\, direct system $\{N_i\}$ the equality
 $\displaystyle  \lim_{\longrightarrow }(M \ot N_i) \cong M\ot \lim_{\longrightarrow }N_i$ holds true for any $M$, while the equality
 $\displaystyle  \lim_{\longrightarrow }\Hom(M, N_i) \cong \Hom(M,\lim_{\longrightarrow }N_i) $ is true if and only if $M$ is finite-dimensional. } (cf.\ (\ref{App: formula for End V_1..V_k}))
\Beqrn
\cE nd_{V^{1\text{-br}}}(I,\fs)&:=&\lim_{\longleftarrow  \atop {p_a \ \text{for}\atop a\in I_\tau^{in,in_0}}}
\lim_{\longleftarrow  \atop {p_b \ \text{for}\atop b\in I_\tau^{in,out_0}}}
\left(\lim_{\longrightarrow  \atop {p_c \ \text{for}\atop c\in I_\tau^{out,out_0}}}
\lim_{\longrightarrow  \atop {p_e \ \text{for}\atop e\in I_\tau^{out,in_0}}}
 \bigotimes_{a\in I_\tau^{in,in_0}} (W^+_{p_a})^*
 \bigotimes_{b\in I_\tau^{in,out_0}} W^-_{p_b}
  \bigotimes_{c\in I_\tau^{out,out_0}} W^+_{p_c}
  \bigotimes_{e\in I_\tau^{out,in_0}} (W^-_{p_e})^*
 \right)\\
 &=&
 \lim_{\longleftarrow  \atop {p_a \ \text{for}\atop a\in I_\tau^{in,in_0}}}
\lim_{\longleftarrow  \atop {p_b \ \text{for}\atop b\in I_\tau^{in,out_0}}}
\Hom\left(
 \bigotimes_{a\in I_\tau^{in,in_0}} W^+_{p_a}
 \bigotimes_{b\in I_\tau^{in,out_0}} (W^-_{p_b})^*\ , \   \ot^{m_1} W^+ \bigotimes \ot^{n_2} (\hat{W}^-)^*\right)\\
&=& \Hom\left(\ot^{n_1} W^+ \bigotimes \ot^{m_2} (\hat{W}^-)^*\ ,\  \ot^{m_1} W^+ \bigotimes \ot^{n_2} (\hat{W}^-)^*\right)
\Eeqrn
 Thus an element $f\in \cE nd_{V^{1\text{-br}}}(I,\fs)$ gives us a well-defined map
$$
f: \ot^{n_1} W^+\bigotimes  \ot^{m_2} (\hat{W}^-)^* \lon \ot^{m_1} W^+ \bigotimes \ot^{n_2} (\hat{W}^-)^*
$$
between countably dimensional vector spaces, and hence such elements can be composed along
the ``blue direction". What about the basic direction? We can try rearranging tensor factors
in   $\cE nd_{V^{1\text{-br}}}(I,\fs)$  as follows
\Beqrn
\cE nd_{V^{1\text{-br}}}(I,\fs)&:=&\lim_{\longleftarrow  \atop {p_a}}
\lim_{\longleftarrow  \atop {p_b}}
\left(\lim_{\longrightarrow  \atop {p_c}}
\lim_{\longrightarrow  \atop {p_e}}
 \bigotimes_{a\in I_\tau^{in,in_0}} (W^+_{p_a})^* \bigotimes_{e\in I_\tau^{out,in_0}} (W^-_{p_e})^*
\bigotimes_{c\in I_\tau^{out,out_0}} W^+_{p_c} \bigotimes_{b\in I_\tau^{in,out_0}} W^-_{p_b}
 \right)\\
 &=&\lim_{\longleftarrow  \atop {p_a}}
\lim_{\longleftarrow  \atop {p_b}}
\left(\lim_{\longrightarrow  \atop {p_c}}
\lim_{\longrightarrow  \atop {p_e}}
\Hom\left(
 \bigotimes_{a\in I_\tau^{in,in_0}} W^+_{p_a} \bigotimes_{e\in I_\tau^{out,in_0}} W^-_{p_e}\ , \
\bigotimes_{c\in I_\tau^{out,out_0}} W^+_{p_c} \bigotimes_{b\in I_\tau^{in,out_0}} W^-_{p_b}
 \right)\right)\\
 &= &\
 \Hom\left(\lim_{\longrightarrow  \atop {p_a}}
 \lim_{\longleftarrow  \atop {p_e}}
 \bigotimes_{a\in I_\tau^{in,in_0}} W^+_{p_a} \bigotimes_{e\in I_\tau^{out,in_0}} W^-_{p_e}\
 ,\
 \lim_{\longleftarrow  \atop {p_b}}
 \lim_{\longrightarrow  \atop {p_c}}
 \bigotimes_{c\in I_\tau^{out,out_0}} W^+_{p_c}\bigotimes_{b\in I_\tau^{in,out_0}} W^-_{p_b}
 \right)\\
\Eeqrn
However, in general,
$$
 \ot^{n_1} W^+ \bigotimes \hat{\ot}^{n_2} \hat{W^-}:=
\lim_{\longrightarrow  \atop {p_a}}
 \lim_{\longleftarrow  \atop {p_e}}
 \bigotimes_{a\in I_\tau^{in,in_0}} W^+_{p_a} \bigotimes_{e\in I_\tau^{out,in_0}} W^-_{p_e}
 \ \ \neq \ \
 \lim_{\longleftarrow  \atop {p_e}}
 \lim_{\longrightarrow  \atop {p_a}}
 \bigotimes_{a\in I_\tau^{in,in_0}} W^+_{p_a} \bigotimes_{e\in I_\tau^{out,in_0}} W^-_{p_e}=: \ot^{m_1} W^+ \widehat{\bigotimes} \hat{\ot}^{m_2} \hat{W^-}
 $$
 with
the l.h.s. being a (proper, in general!) subspace of the r.h.s. Hence elements
of the $\Hom$-spaces
$$
\cE nd_{V^{1\text{-br}}}(I,\fs)\cong\Hom\left(
\ot^{n_1} W^+ \bigotimes \hat{\ot}^{n_2} \hat{W^-}\ , \
\ot^{m_1} W^+ \widehat{\bigotimes} \hat{\ot}^{m_2} \hat{W^-}
\right)
$$
can not be composed, in general, along graphs
of the type shown in (\ref{2new: wrong compos in 2-or prop}). (Nevertheless
the latest formula shows that any element of $\cE nd_{V^{1\text{-br}}}(I,\fs)$ can be understood  as some linear map along the basic direction.)

\sip

 We conclude that the $\cS^{(2)}$-module $\cE nd_{V^{1\text{-br}}}$ admits nice compositions $\mu_\Ga$ along any graphs $\Ga$ {\em not containing closed paths of directed edges in blue color}\, as in (\ref{2new: good compos in 2-or prop}) (with the ``associativity" axioms are obviously satisfied) and hence gives us an example of 2-oriented prop. We call it
the {\em endomorphism prop}\, of $V^{1\text{-br}}$.

\sip

Note that if $V^{1\text{-br}}$ is finite-dimensional (or at least if $W^-$ is finite-dimensional), then
$$
\cE nd_{V^{1\text{-br}}}(I,\fs)\cong\Hom\left(
\ot^{n_1} W^+ \bigotimes {\ot}^{n_2} {W^-}\, ,
\ot^{m_1} W^+ {\bigotimes} {\ot}^{m_2} {W^-}
\right)\cong \Hom\left(\ot^{n_1} W^+ \bigotimes \ot^{m_2} (\hat{W}^-)^*\ ,\  \ot^{m_1} W^+ \bigotimes \ot^{n_2} (\hat{W}^-)^*\right)
$$
and the compositions $\mu_\Ga$ (in the definition of a multi-directed prop) make sense for any graphs $\Ga\in G^{0\uparrow k+1}$.

\subsubsection{\bf Definition} Let $\cP^{2\text{-or}}$ be a 2-oriented prop(erad). A morphism of 2-oriented prop(erad)s
$$
\rho:\cP^{2\text{-or}} \lon \cE nd_{V^{1\text{-br}}}
$$
is called a {\em representation}\, of $\cP^{2\text{-or}}$ in the vector space $V$ with one brane.

\subsubsection{\bf Example} A representation of a 2-oriented operad $\cA ss^{(2)}$ (resp.,  $\caL ie^{(2)}$) in $V^{1\text{-br}}$ is given by a collection of linear maps
$$
\Ba{c}\resizebox{10mm}{!}  {
\begin{tikzpicture}[baseline=-1ex]
\node[int] (a) at (0,0) {};
\node[] (u1) at (0,0.8) {};
\node[int] (a) at (0,0) {};
\node[] (d1) at (-0.5,-0.7) {};
\node[] (d2) at (0.5,-0.7) {};
\draw (u1) edge[latex-, leftblue] (a);
\draw (a) edge[latex-, leftblue] (d1);
\draw (a) edge[latex-, leftblue] (d2);
\end{tikzpicture}
}\Ea\hspace{-5mm}: W^+ \ot W^+ \rar W^+, \Ba{c}\resizebox{10mm}{!}  {
\begin{tikzpicture}[baseline=-1ex]
\node[int] (a) at (0,0) {};
\node[] (u1) at (0,0.8) {};
\node[int] (a) at (0,0) {};
\node[] (d1) at (-0.5,-0.7) {};
\node[] (d2) at (0.5,-0.7) {};
\draw (u1) edge[latex-, rightblue] (a);
\draw (a) edge[latex-, rightblue] (d1);
\draw (a) edge[latex-, rightblue] (d2);
\end{tikzpicture}
}\Ea\hspace{-5mm}: (\hat{W}^-)^* \rar (\hat{W}^-)^*\ot (\hat{W}^-)^*,
\Ba{c}\resizebox{10mm}{!}  {
\begin{tikzpicture}[baseline=-1ex]
\node[int] (a) at (0,0) {};
\node[] (u1) at (0,0.8) {};
\node[int] (a) at (0,0) {};
\node[] (d1) at (-0.5,-0.7) {};
\node[] (d2) at (0.5,-0.7) {};
\draw (u1) edge[latex-, leftblue] (a);
\draw (a) edge[latex-, leftblue] (d1);
\draw (a) edge[latex-, rightblue] (d2);
\end{tikzpicture}
}\Ea\hspace{-5mm}: W^+ \rar W^+\ot (W^{-})^*, \Ba{c}\resizebox{10mm}{!}  {
\begin{tikzpicture}[baseline=-1ex]
\node[int] (a) at (0,0) {};
\node[] (u1) at (0,0.8) {};
\node[int] (a) at (0,0) {};
\node[] (d1) at (-0.5,-0.7) {};
\node[] (d2) at (0.5,-0.7) {};
\draw (u1) edge[latex-, rightblue] (a);
\draw (a) edge[latex-, leftblue] (d1);
\draw (a) edge[latex-, rightblue] (d2);
\end{tikzpicture}}\Ea\hspace{-5mm}: W^+ {\ot} (\hat{W}^-)^* \rar (\hat{W}^-)^*
$$
such that their compositions satisfy relations (\ref{3new: Ass^2 relations 1st set})-(\ref{3new: Ass^2 relations 2nd set}) (respectively, (\ref{3new: Lie^2 relations 1st set})-(\ref{3new: Lie^2 relations 3rd set})).

\sip

To illustrate the divergency phenomenon let us consider a generic representation
of a prop in which compositions along graphs with {\em non-trivial genus}\, make sense. For example, a representation of, say,  the prop of 2-oriented Lie bialgebras
$\LB^{2\text{-or}}$ in $V^{1\text{-br}}$ is given by maps as above plus the following
ones
$$
\Ba{c}\resizebox{10mm}{!}  {
\begin{tikzpicture}[baseline=-1ex]
\node[int] (a) at (0,0) {};
\node[] (u1) at (0,-0.8) {};
\node[int] (a) at (0,0) {};
\node[] (d1) at (-0.5,0.7) {};
\node[] (d2) at (0.5,0.7) {};
\draw (a) edge[latex-, leftblue] (u1);
\draw (d1) edge[latex-,leftblue] (a);
\draw (d2) edge[latex-, leftblue] (a);
\end{tikzpicture}
}\Ea\hspace{-5mm}: W^+ \rar W^+ {\ot} W^+,
\Ba{c}\resizebox{10mm}{!}  {
\begin{tikzpicture}[baseline=-1ex]
\node[int] (a) at (0,0) {};
\node[] (u1) at (0,-0.8) {};
\node[int] (a) at (0,0) {};
\node[] (d1) at (-0.5,0.7) {};
\node[] (d2) at (0.5,0.7) {};
\draw (a) edge[latex-, rightblue] (u1);
\draw (d1) edge[latex-, leftblue] (a);
\draw (d2) edge[latex-, rightblue] (a);
\end{tikzpicture}
}\Ea\hspace{-5mm}: (\hat{W}^-)^* \rar W^+ {\ot} (\hat{W}^-)^*
 \Ba{c}\resizebox{10mm}{!}  {
\begin{tikzpicture}[baseline=-1ex]
\node[int] (a) at (0,0) {};
\node[] (u1) at (0,-0.8) {};
\node[int] (a) at (0,0) {};
\node[] (d1) at (-0.5,0.7) {};
\node[] (d2) at (0.5,0.7) {};
\draw (a) edge[latex-, leftblue] (u1);
\draw (d1) edge[latex-,leftblue] (a);
\draw (d2) edge[latex-, rightblue] (a);
\end{tikzpicture}
}\Ea\hspace{-5mm}: W^+\ot (\hat{W}^-)^* \rar W^+,
\Ba{c}\resizebox{12mm}{!}  {
\begin{tikzpicture}[baseline=-1ex]
\node[int] (a) at (0,0) {};
\node[] (u1) at (0,-0.8) {};
\node[int] (a) at (0,0) {};
\node[] (d1) at (-0.5,0.7) {};
\node[] (d2) at (0.5,0.7) {};
\draw (a) edge[latex-, rightblue] (u1);
\draw (d1) edge[latex-,rightblue] (a);
\draw (d2) edge[latex-, rightblue] (a);
\end{tikzpicture}
}\Ea\hspace{-5mm}: (\hat{W}^-)^*\ot (\hat{W}^-)^* \rar  (\hat{W}^-)^*
$$
satisfying certain quadratic relations. If $\{x_{a^+}\}_{a^+\in N_{\geq 1}}$ is a countably infinite basis of $W^+$, $\{x_{a^-}\}_{a^-\in N_{\geq 1}}$  a countably infinite basis of ${W}^-$, and $\{y^{a_+}\}_{a_+\in N_{\geq 1}}$ and  $\{y^{a_-}\}_{a_-\in N_{\geq 1}}$ the associated
countably infinite dual bases of  $(\hat{W}^+)^*$ and, respectively, $(\hat{W}^-)^*$
then the corresponding maps, say the following ones
$$
\Ba{c}\resizebox{10mm}{!}  {
\begin{tikzpicture}[baseline=-1ex]
\node[int] (a) at (0,0) {};
\node[] (u1) at (0,0.8) {};
\node[int] (a) at (0,0) {};
\node[] (d1) at (-0.5,-0.7) {};
\node[] (d2) at (0.5,-0.7) {};
\draw (u1) edge[latex-, leftblue] (a);
\draw (a) edge[latex-, leftblue] (d1);
\draw (a) edge[latex-, leftblue] (d2);
\end{tikzpicture}
}\Ea\hspace{-5mm}: W^+ \ot W^+ \stackrel{\mu_1}{\rar} W^+,
\Ba{c}\resizebox{10mm}{!}  {
\begin{tikzpicture}[baseline=-1ex]
\node[int] (a) at (0,0) {};
\node[] (u1) at (0,0.8) {};
\node[int] (a) at (0,0) {};
\node[] (d1) at (-0.5,-0.7) {};
\node[] (d2) at (0.5,-0.7) {};
\draw (u1) edge[latex-, leftblue] (a);
\draw (a) edge[latex-, leftblue] (d1);
\draw (a) edge[latex-, rightblue] (d2);
\end{tikzpicture}
}\Ea\hspace{-5mm}: W^+ \stackrel{\mu_2}{\rar} W^+\ot (\hat{W}^{-})^*,
\Ba{c}\resizebox{10mm}{!}  {
\begin{tikzpicture}[baseline=-1ex]
\node[int] (a) at (0,0) {};
\node[] (u1) at (0,-0.8) {};
\node[int] (a) at (0,0) {};
\node[] (d1) at (-0.5,0.7) {};
\node[] (d2) at (0.5,0.7) {};
\draw (a) edge[latex-, leftblue] (u1);
\draw (d1) edge[latex-,leftblue] (a);
\draw (d2) edge[latex-, leftblue] (a);
\end{tikzpicture}
}\Ea\hspace{-5mm}: W^+ \stackrel{\mu_3}{\rar} W^+ {\ot} W^+,
 \Ba{c}\resizebox{10mm}{!}  {
\begin{tikzpicture}[baseline=-1ex]
\node[int] (a) at (0,0) {};
\node[] (u1) at (0,-0.8) {};
\node[int] (a) at (0,0) {};
\node[] (d1) at (-0.5,0.7) {};
\node[] (d2) at (0.5,0.7) {};
\draw (a) edge[latex-, leftblue] (u1);
\draw (d1) edge[latex-,leftblue] (a);
\draw (d2) edge[latex-, rightblue] (a);
\end{tikzpicture}
}\Ea\hspace{-5mm}: W^+\ot (\hat{W}^-)^* \stackrel{\mu_4}{\rar} W^+,
$$
can be represented by the following infinite, in general, sums
$$
\mu_1=\sum_{a^+,b^+,c^+\in \N_{\geq 1}}\Phi_{a^+b^+}^{c^+} y^{a^+}\ot y^{b^+}\ot x_{c^+}, \ \ \ \
\mu_2=\sum_{a^+,b_-,c^+\in \N_{\geq 1}}\Phi^{a^+}_{c^+b_-} x^{c^+}\ot x_{a^+}\ot y^{b_-},\ \ \ \ \ \ \ \Phi_{a^+b^+}^{c^+}, \Phi^{a^+}_{c^+b_-}\in \K
$$
$$
\mu_3=\sum_{a_+,b_+,c_+\in \N_{\geq 1}}\Psi^{a_+b_+}_{c_+} y^{c^+}\ot x_{a^+}\ot x_{b^+}  \ \ \ \
\mu_4=\sum_{a^+,b_-,c^+\in \N_{\geq 1}}\Psi^{a^+b_-}_{c^+}  y^{c^+}\ot x_{b^-}\ot x_{a^+}, \ \ \ \ \ \ \ \ \Psi^{a_+b_+}_{c_+}, \Psi^{a^+b_-}_{c^+}\in \K.
$$
where the coefficients satisfy the conditions:
\Bi
\item for fixed $a_+,b_-$ only finitely many $\Phi_{a^+b^+}^{c^+}\neq 0$; for fixed $c_+$ only finitely many $\Phi^{a^+}_{c^+b_-}\neq 0$;
\item for fixed $c_+$ only finitely many $\Psi^{a_+b_+}_{c_+}\neq 0$;
for fixed $a_+,b_-$ only finitely many $\Psi_{a^+b_-}^{c^+} \neq 0$.
\Ei
Then the element
$$
\Ba{c}\resizebox{6mm}{!}  {
\begin{tikzpicture}[baseline=-.65ex]
\node[] (u) at (0,1.2) {};
 \node[int] (a) at (0,-0.4) {};
 \node[int] (b) at (0,0.4) {};
\node[] (d) at (0,-1.2) {};
\draw (u) edge[latex-, leftblue] (b);
\draw (a) edge[latex-, leftblue] (d);
 \draw (b) edge[latex-, leftblue, bend left=60] (a);
 \draw (b) edge[latex-, leftblue, bend right=60] (a);
\end{tikzpicture}}\Ea\in \LB^{2\text{-or}}
$$
gets represented in $V^{1\text{-br}}$ as a linear map
$$
\sum_{c^+,d^+\in \N_{\geq 1}}\left(\underbracket{\sum_{a^+,b^+\in \N_{\geq 1}} \Phi_{a^+b^+}^{c^+} \Psi^{a_+b_+}_{d_+}}_{\mathrm{only\ finitely\ many\atop terms\ non-zero}}   \right) y^{d^+} \ot x_{c^+}: W^+\rar W^+
$$
which is always well defined, while  the element
$$
\Ba{c}\resizebox{6mm}{!}  {
\begin{tikzpicture}[baseline=-.65ex]
\node[] (u) at (0,1.2) {};
 \node[int] (a) at (0,-0.4) {};
 \node[int] (b) at (0,0.4) {};
\node[] (d) at (0,-1.2) {};
\draw (u) edge[latex-, leftblue] (b);
\draw (a) edge[latex-, leftblue] (d);
 \draw (b) edge[latex-, leftblue, bend left=60] (a);
 \draw (b) edge[latex-, rightblue, bend right=60] (a);
\end{tikzpicture}}\Ea\in \LB^{1\uparrow 2}
$$
gets represented in $V^{1\text{-br}}$ as a formal sum of linear maps
$$
\sum_{c^+,d^+\in \N_{\geq 1}}\left(\underbracket{\sum_{a^+,b^+\in \N_{\geq 1}} \Phi^{c^+}_{a^+b_-}\Psi^{a^+b_-}_{d_+}}_{\mathrm{infinitely\ many\ terms\atop can\ be\ non-zero\ in \ general}}   \right) y^{d^+} \ot x_{c^+}: W^+\rar W^+
$$
which in general diverges. Such a map makes sense in general in general only
in the case $\dim W^+< \infty$ or $\dim W^-<\infty$.

\subsubsection{\bf Symplectic vector space with Lagrangian branes}
Let $(V,\om: \wedge^2 V \rar \K)$ be a finite-dimensional vector space
equipped with a symplectic form; in general $\dim V=2n$ for some $n\in \N_{\geq 1}$. A subspace $W\subset V$ is called isotropic
if
$$
W\subset W^\bot:=\left\{v\in V\ \mid \om(v,w)=0\ \forall w\in W     \right\}.
$$
Such a subspace is called {\em Lagrangian}\, if $\dim W=n$. It is well-known
that a Lagrangian subspace $W^+\subset V$ always admits a complement $W^-\subset V$ which is also Lagrangian. Moreover, in this case the symplectic form
induces a canonical isomorphism
$$
\om: (W^-)^* \lon W^+.
$$
The data $(V,W^+,W^-)$ is called a {\em finite-dimensional symplectic  vector space
with one Lagrangian brane}. Similarly one defines  {\em finite-dimensional symplectic  vector space
with $k$ Lagrangian branes}, $(V,W^+_\tau,W_\tau^-)_{\tau\in [k]}$.
We generalize this notion to infinite dimensions as follows.

\sip

Let $(V, W_1,\ldots, W_k)$ be a countably infinite dimensional vector space
with $k$ branes such that for each $p$ and each $\tau \in [k]$ the vector space $V_p=W^+_{\tau, p} \oplus W^-_{\tau, p}$ is a finite dimensional symplectic vector space with $k$ Lagrangian branes. Then  the symplectic form induces a linear map
$$
\om: (\hat{W}_\tau^-)^* \lon W_\tau^+
$$
which is an isomorphism for each $\tau\in [k]$. The resulting datum
$(V,W_1,\ldots,W_k,\om)$ is called an {\em infinite-dimensional symplectic
vector space with $k$ Lagrangian branes}\, and denoted by $V^{k\text{-br}}_{symp}$.

\sip

If we consider now a generic representation $\rho$ of, say, $\cA ss^{(2)}$ or $\caL ie^{(2)}$ in $V^{1\text{-br}}_{symp}$, then, due to the canonical isomorphism
$(\hat{W}_\tau^-)^*= W_\tau^+$, we see that multioriented generators
which differ only in the basic orientation stand for linear maps of the same type,
for example
$$
\Ba{c}\resizebox{10mm}{!}  {
\begin{tikzpicture}[baseline=-1ex]
\node[int] (a) at (0,0) {};
\node[] (u1) at (0,0.8) {};
\node[int] (a) at (0,0) {};
\node[] (d1) at (-0.5,-0.7) {};
\node[] (d2) at (0.5,-0.7) {};
\draw (u1) edge[latex-, leftblue] (a);
\draw (a) edge[latex-, leftblue] (d1);
\draw (a) edge[latex-, leftblue] (d2);
\end{tikzpicture}
}\Ea\hspace{-5mm}: W^+ \ot W^+ \rar W^+, \ \  \Ba{c}\resizebox{10mm}{!}  {
\begin{tikzpicture}[baseline=-1ex]
\node[int] (a) at (0,0) {};
\node[] (u1) at (0,0.8) {};
\node[int] (a) at (0,0) {};
\node[] (d1) at (-0.5,-0.7) {};
\node[] (d2) at (0.5,-0.7) {};
\draw (u1) edge[latex-, rightblue] (a);
\draw (a) edge[latex-, leftblue] (d1);
\draw (a) edge[latex-, rightblue] (d2);
\end{tikzpicture}}\Ea\hspace{-5mm}: W^+ {\ot} W^+ \rar W^+,
$$
and hence it makes sense to identify them. We call a representation $\rho$
in $V^{1\text{-br}}_{symp}$ {\em reduced symplectic Lagrangian}\, if $\rho$
takes identical values on all generating corollas of $\cP^{2\text{-or}}$ which
become identical (as $k$-oriented graphs) once the { basic} directions on edges are forgotten.
Then we can reformulate observations {\ref{3new: Prop on map from ASS^2 to IB}} and {\ref{3new: Prop on map from Lie^2 to Lieb}}
as follows.

\subsubsection{\bf Proposition} {\em
(i) There is a one-to-one correspondence between reduced symplectic Lagrangian
representations of $\cA ss^{(2)}$ in $V^{1\text{-br}}_{symp}$ and infinitesimal bialgebra structures in the Lagrangian subspace $W^+$.

\sip

(ii) There is a one-to-one correspondence between reduced symplectic Lagrangian
representations of $\caL ie^{(2)}$ in $V^{1\text{-br}}_{symp}$ and Lie bialgebra structures in the Lagrangian subspace $W^+$.
}

 \subsection{\bf Remark} In principle one can use symplectic
structures on $V$ of homological degree $q\neq 0$ so that the
induced isomorphism takes the form  $\om: (\hat{W}_\tau^-)^* \lon W_\tau^+[q]$
 but then the basic direction  can not be forgotten completely in representations as it stands now for a degree shift of linear maps.

\subsection{\bf Multidirected endomorphism  prop  of a graded vector space with $k$ branes}\label{App: subsubsec on tensor algebra of Ws}
Let $(I,\fs: I\rar \f_{k^+})$ be a multi-oriented set. Recall that for any fixed $i\in I$ there is an associated  map
$$
\fs_i: [k^+]\stackrel{\fs(i)}{\lon} \{out,in\}.
$$
while for any fixed $\tau\in [k^+]$ there is a map
$$
\Ba{rccc}
\fs_\tau: & I & \lon & \{out,in\}\\
            & i & \lon & \fs_\tau(i):= \fs_i(\tau).
\Ea
$$
The latter map can be used to decompose $I$ into two disjoint subsets
$$
I=\fs_\tau^{-1}(out) \sqcup \fs_\tau^{-1}(in)
$$
The basic direction $\tau=0$ plays a special role. For any $\tau\neq 0$, i.e.\ for any $\tau\in [k]$ we can further decompose the set $I$ into four disjoint subsets
$$
I=(\fs_\tau^{-1}(out) \sqcup \fs_\tau^{-1}(in))\cap (\fs_0^{-1}(out) \sqcup \fs_0^{-1}(in)):= I_{\tau}^{out,out_0}\sqcup  I_{\tau}^{out,in_0}\sqcup  I_{\tau}^{in,out_0}  \sqcup I_{\tau}^{in,in_0}
$$
where
$$
 I_{\tau}^{out,out_0}:=\fs_\tau^{-1}(out)\cap \fs_0^{-1}(out), \ \ \
 I_{\tau}^{out,in_0}:=\fs_\tau^{-1}(out)\cap \fs_0^{-1}(in), \ \ \
  I_{\tau}^{in,out_0}:=\fs_\tau^{-1}(in)\cap \fs_0^{-1}(out), \ \ \
 I_{\tau}^{in,in_0}:=\fs_\tau^{-1}(in)\cap \fs_0^{-1}(in), \ \ \
$$

\sip

\sip

Given a graded vector space with $k$  branes, $V^{k\text{-br}}=(V=\displaystyle \lim_{\lon} V_p,W_1,\ldots, W_k)$, consider a collection of linear subspaces
for each $p\in \N$,
$$
 W^{\fm}_{p}:=  W_{1,p}^{\fm(1)}\cap W_{2,p}^{\fm(2)} \cap \ldots \cap W_{k,p}^{\fm(k)},
$$
one for each multidirection $\fm:[k^+]\rar  \{out,in\}$ from $\f_{k^+}$, where we set for each $\tau\in [k]$,
$$
 W_{\tau,p}^{\fm(\tau)}:=
 \left\{\Ba{ll}
 W_{\tau,p}^+ & \text{if}\ \fm(0)=\fm(\tau)=out\\
  (W_{\tau,p}^+)^* & \text{if}\ \fm(0)=\fm(\tau)=in \\
 W_{\tau,p}^- & \text{if}\ \fm(0)=out,\ \fm(\tau)=in\\
  (W_{\tau,p}^-)^* & \text{if}\ \fm(0)=in,\  \fm(\tau)=out\\
\Ea
 \right.
$$
Note that $(W_{\tau,p}^{\fm})^*= W_{\tau,p}^{\fm^{opp}}$.
For example,
$$
\text{for}\ \ \fm= \Ba{c}
\begin{tikzpicture}[baseline=-2ex]
\node[int] at (0,0) {};
\node[] at (2.75,0) {};
\draw (0,0) edge[leftgreen] node[above] { \ \ \ $\scriptstyle \fm(1)$}(0.9,0);
\draw (0.9,0) edge[rightblue] node[above] {\ $\scriptstyle \fm(2)$ }(1.5,0);
\draw (1.5,0) edge[] node[above] {$...$\ \ \ }(1.9,0);
\draw (1.9,0) edge[leftred] node[above] { $\scriptstyle \fm(k)$ }(2.3,0);
\draw (2.2,0) edge[-latex] node[above] {\ \ \ \ \ $\scriptstyle \fm(0)$ }(2.8,0);
\end{tikzpicture}
\Ea
\ \ \ \ \text{one has}\ \ \ \  W^{\fm}_{p}=W^{-}_{1,p}\cap W^{+}_{2,p} \cap \ldots \cap W^{-}_{k,p}
$$
while
$$
\hspace{13mm}
\text{for}\ \ \ \fm^{opp}= \Ba{c}
\begin{tikzpicture}[baseline=-2ex]
\node[] at (0,0) {};
\node[int] at (2.85,0) {};
\draw (0,0) edge[leftgreen] node[above] { \ \ \ $\scriptstyle \fm(1)$}(0.9,0);
\draw (0.9,0) edge[rightblue] node[above] {\ $\scriptstyle \fm(2)$ }(1.5,0);
\draw (1.5,0) edge[] node[above] {$...$\ \ \ }(1.9,0);
\draw (1.9,0) edge[leftred] node[above] { $\scriptstyle \fm(k)$ }(2.3,0);
\draw (2.2,0) edge[-latex] node[above] {\ \ \ \ \ $\scriptstyle \fm(0)$ }(2.8,0);
\end{tikzpicture}
\Ea
\ \ \ \ \text{one has}\ \ \ \  W^{\fm^{opp}}_{p}=(W^{-}_{1,p})^*\cap (W^{+}_{2,p})^* \cap \ldots \cap (W^{-}_{k,p})^*.
$$
We define a countably dimensional vector space
$$
W^{\fm}:=\lim_{\lon\atop p}   W^{\fm}_{p},
$$

\sip

 Define an $\cS^{(k+1)}$-module $\cE nd_{V^{k\text{-br}}}$, that is, a functor
$$
\Ba{rccl}
\cE nd_{V^{k\text{-br}}}: & \cS^{(k+1)} & \lon & \text{Category of dg vector spaces} \\
              & (I, \fs) & \lon &  \cE nd_{V^{k\text{-br}}}(I,\fs)
              \Ea,
$$
by setting
$$
\cE nd_{V^{k\text{-br}}}(I,\fs):=\bigcap_{\tau\in [k]} \Hom_\tau(\fs, I),
$$
where (cf.\ (\ref{App: formula for End V_1..V_k}))
\Beqrn
\Hom_\tau(\fs, I)&:=&\lim_{\longleftarrow  \atop {p_i \ \text{for}\atop i\in I_\tau^{in,out_0} \cup I_\tau^{in,in_0}}}  \left(\lim_{\longrightarrow  \atop {p_j \ \text{for}\atop i\in I_\tau^{out,out_0} \cup I_\tau^{out,in_0}}} \bigotimes_{i\in I} W^{\fs_i}_{p_i}\right)\\
&=& \lim_{\longleftarrow  \atop p_i}\left(
 \bigotimes_{i\in i\in I_\tau^{in,out_0} \cup I_\tau^{in,in_0}}W^{\fs_i}_{p_i}\ \ \ot \bigotimes_{j\in I_\tau^{out,out_0} \cup I_\tau^{out,in_0}}  W^{\fs_j}\right)\\
  &=& \Hom\left(
 \bigotimes_{i\in I_\tau^{in,out_0} \cup I_\tau^{in,in_0}}\lim_{\lon\atop p_i} W^{\fs^{opp}_i}_{p_i}, \bigotimes_{j\in I_\tau^{out,out_0} \cup I_\tau^{out,in_0}}  W^{\fs_j}\right)\\
 &=& \Hom\left(
 \bigotimes_{i\in i\in I_\tau^{in,out_0} \cup I_\tau^{in,in_0}}W^{\fs^{opp}_i}, \bigotimes_{j\in I_\tau^{out,out_0} \cup I_\tau^{out,in_0}}  W^{\fs_j}\right).
\Eeqrn
Thus a single element $f\in \cE nd_{V^{k\text{-br}}}(I,\fs)$ has
$k$ incarnations as a linear map, one for each ``coloured direction" $\tau\in [k]$.
Note that all the $k$ spaces $\Hom_\tau(\fs, I)$, $\tau\in [k]$,  belong to one and the same vector space
\Beq\label{4new: Hom(I, s)}
\Hom(\fs, I):=\lim_{\longleftarrow  \atop p_i } \bigotimes_{i\in I} W^{\fs_i}_{p_i}
\Eeq
so that it makes sense to talk about their intersection. If $V$ is finite-dimensional, then, of course, $\Hom_\tau(\fs, I)=\Hom(\fs, I)$ for any $\tau\in [k]$.
\sip

 Therefore elements of  $\cE nd_{V^{k\text{-br}}}$  can be composed (when it makes sense) along each of the ``coloured" direction, but in general they can {\em not}\, be composed along the basic direction (i.e.\ compositions of the type
(\ref{2new: wrong compos in 2-or prop}) have no sense in general).

\sip

Let $\cP^{(k+1)\text{-or}}$ be a $(k+1)$-oriented prop(erad). A morphism of $(k+1)$-oriented prop(erad)s
$$
\rho:\cP^{(k+1)\text{-or}} \lon \cE nd_{V^{k\text{-br}}}
$$
is called a {\em representation}\, of $\cP^{(k+1)\text{-or}}$ in a vector space with one $k$ branes $V^{k\text{-br}}$. If $V^{k\text{-br}}$ happens to be a symplectic vector
space with $k$ Lagrangian branes, then a representation $\rho$ is called
{\em reduced symplectic Lagrangian}\, if $\rho$
takes identical values on all those generating corollas of $\cP^{(k+1)\text{-or}}$ which
become identical (as $k$-oriented decorated graphs) once the {\em basic}\, directions on edges are forgotten.

\subsection{Example: 3-oriented endomorphism prop} Let $V^{2\text{-or}}=(V,W_1^+,W_2^+)$ be a countably dimensional graded vector space with $2$ branes. We would like to describe in more details the structure of the associated endomorphism prop which is a functor
$$
\Ba{rccc}
\cE nd_{V^{2\text{-or}}}: & \cS^{(3)} & \lon & \text{Category\ of\ dg\   vector\ spaces}\\
& (I,\fs)=
\Ba{c}\resizebox{48mm}{!}
{\begin{tikzpicture}[baseline=-1ex]
\node[] at (-0.6,1.3) {$^{\scriptstyle I_2}$};
\node[] at (-0.7,1.1) {$\overbrace{\  \ \ \ \ \ \ \ \ \ \ }$};
\node[] at (-0.6,0.9) {$...$};
\node[] at (-2.1,1.3) {$^{\scriptstyle I_1}$};
\node[] at (-2.1,1.1) {$\overbrace{\  \ \ \ \ \ \ \ \ \ \ }$};
\node[] at (-2.0,0.9) {$...$};
\node[] at (0.6,1.3) {$^{\scriptstyle I_3}$};
\node[] at (0.7,1.1) {$\overbrace{\  \ \ \ \ \ \ \ \ \ \ }$};
\node[] at (0.6,0.9) {$...$};
\node[] at (2.1,1.3) {$^{\scriptstyle I_4}$};
\node[] at (2.1,1.1) {$\overbrace{\  \ \ \ \ \ \ \ \ \ \ }$};
\node[] at (2.0,0.9) {$...$};
\node[] at (-0.6,-1.9) {$_{\scriptstyle I_6}$};
\node[] at (-0.7,-1.1) {$\underbrace{\  \ \ \ \ \ \ \ \ \ \ }$};
\node[] at (-0.6,-0.9) {$...$};
\node[] at (-2.1,-1.3) {$_{\scriptstyle I_5}$};
\node[] at (-2.1,-1.1) {$\underbrace{\  \ \ \ \ \ \ \ \ \ \ }$};
\node[] at (-2.0,-0.9) {$...$};
\node[] at (0.6,-1.3) {$_{\scriptstyle I_7}$};
\node[] at (0.7,-1.1) {$\underbrace{\  \ \ \ \ \ \ \ \ \ \ }$};
\node[] at (0.6,-0.9) {$...$};
\node[] at (2.1,-1.3) {$_{\scriptstyle I_8}$};
\node[] at (2.1,-1.1) {$\underbrace{\  \ \ \ \ \ \ \ \ \ \ }$};
\node[] at (2.0,-0.9) {$...$};
\node[int] (0) at (0,0) {};
\node[] (1) at (-2.7,1.0) {};
\node[] (2) at (-1.8,1.0) {};
\node[] (3) at (-1.2,1.0) {};
\node[] (4) at (-0.3,1.0) {};
\node[] (1') at (2.7,1.0) {};
\node[] (2') at (1.8,1.0) {};
\node[] (3') at (1.2,1.0) {};
\node[] (4') at (0.3,1.0) {};
\node[] (5) at (-2.7,-1.0) {};
\node[] (6) at (-1.8,-1.0) {};
\node[] (7) at (-1.2,-1.0) {};
\node[] (8) at (-0.3,-1.0) {};
\node[int] (0) at (0,0) {};
\node[] (5') at (2.7,-1.0) {};
\node[] (6') at (1.8,-1.0) {};
\node[] (7') at (1.2,-1.0) {};
\node[] (8') at (0.3,-1.0) {};
\draw (1) edge[latex-,leftblue,leftred] (0);
\draw (2) edge[latex-,leftblue, leftred] (0);
\draw (3) edge[latex-,rightblue,leftred] (0);
\draw (4) edge[latex-,rightblue,leftred] (0);
\draw (1') edge[latex-,rightblue,rightred] (0);
\draw (2') edge[latex-,rightblue, rightred] (0);
\draw (3') edge[latex-,leftblue,rightred] (0);
\draw (4') edge[latex-,leftblue,rightred] (0);
\draw (0) edge[latex-,leftblue,leftred] (5);
\draw (0) edge[latex-,leftblue, leftred] (6);
\draw (0) edge[latex-,rightblue,leftred] (7);
\draw (0) edge[latex-,rightblue,leftred] (8);
\draw (0) edge[latex-,rightblue,rightred] (5');
\draw (0) edge[latex-,rightblue, rightred] (6');
\draw (0) edge[latex-,leftblue,rightred] (7');
\draw (0) edge[latex-,leftblue,rightred] (8');
\end{tikzpicture}}\Ea
                         &\lon & \cE nd_{V^{2\text{-or}}}(I,\fs)
\Ea
$$
The multi-orientation $\fs$ defines (and can be reconstructed from) the decomposition $I=I_1 \sqcup I_2\sqcup  \ldots \sqcup I_8$ as explained in the picture. If
$$
\# I_i=m_i\ \ \text{for}\ \ i\in \{1,2,3,4\}\ , \  \# I_i=n_{i-4}\ \ \text{for}\ \ i\in \{5,6,7,8\}.
$$
then, by definition, $\cE nd_{V^{2\text{-or}}}(I,\fs)$ is the intersection
in (\ref{4new: Hom(I, s)}) of two graded vector spaces
$$
{\color{blue}{\Hom}(I,\fs)}:= \Hom\left(\ot^{n_1} W^{++} \ot^{n_3} W^{+-}  \ot^{m_2} (\hat{W}^{-+})^*
  \ot^{m_4} (\hat{W}^{--})^*\ , \
\ot^{m_1} W^{++}  \ot^{m_3} W^{+-} \ot^{n_2} (\hat{W}^{-+})^*
  \ot^{n_4} (\hat{W}^{--})^*\right)
$$
and
$$
{\color{red}{\Hom}(I,\fs)}:= \Hom\left(\ot^{n_1} W^{++} \ot^{n_2} W^{-+}  \ot^{m_3} (\hat{W}^{+-})^*
  \ot^{m_4} (\hat{W}^{--})^*\ , \
\ot^{m_1} W^{++}  \ot^{m_2} W^{-+} \ot^{n_3} (\hat{W}^{+-})^*
  \ot^{n_4} (\hat{W}^{--})^*\right)
$$
where we set
$$
W^{++}:=\lim_{\lon\atop p} W_{1,p}^+\cap W_{2,p}^+,\ \
W^{-+}:=\lim_{\lon\atop p} W_{1,p}^-\cap W_{2,p}^+, \ \
W^{+-}:=\lim_{\lon\atop p} W_{1,p}^+\cap W_{2,p}^-, \ \
W^{--}:=\lim_{\lon\atop p} W_{1,p}^-\cap W_{2,p}^-,
$$
$$
\hat{W}^{++}:=\lim_{\longrightarrow \atop p} W_{1,p}^-\cap W_{2,p}^-,\ \
\hat{W}^{-+}:=\lim_{\longleftarrow \atop p} W_{1,p}^-\cap W_{2,p}^+, \ \
\hat{W}^{--}:=\lim_{\longrightarrow \atop p} W_{1,p}^-\cap W_{2,p}^-,\ \
\hat{W}^{+-}:=\lim_{\longleftarrow \atop p} W_{1,p}^+\cap W_{2,p}^-,
$$
All the tensor factors shown in the above formulae for ${\color{blue}{\Hom}(I,\fs)}$
and ${\color{red}{\Hom}(I,\fs)}$ are countably dimensional vector spaces. Let
$\{x_{A_{++}}\}$, $\{x_{A_{+-}}\}$, $\{x_{A_{-+}}\}$,  $\{x_{A_{--}}\}$ be bases for the (direct limit) vector spaces $W^{++}$, $W^{+-}$,  $W^{-+}$ and
 $W^{--}$, while $y^{A^{++}}$, $y^{A^{+-}}$, $y^{A^{-+}}$ and $y^{A^{--}}$  be the associated dual bases for (also direct limit) vector spaces $(\hat{W}^{++})^*$, $(\hat{W}^{+-})^*$, $(\hat{W}^{+-})^*$   and $(\hat{W}^{--})^*$. Then the ``big" vector space (\ref{4new: Hom(I, s)}) consists of all formal power series of the form
$$
\sum_{A_{\bu,\bu}, B^{\bu\bu}} F^{A_{++}A_{+-}A_{-+}A_{--}}_{B^{++}B^{+-}B^{-+}B^{--}}\ x_{A_{-+}} \ot x_{A_{-+}} \ot x_{A_{-+}} \ot x_{A_{-+}} \ot y^{B^{++}} \ot y^{B^{+-}} \ot
y^{B^{-+}} \ot y^{B^{--}}, \ \ \ \ \ F^{A_{++}A_{+-}A_{-+}A_{--}}_{B^{++}B^{+-}B^{-+}B^{--}}\in \K,
$$
its subspace $\color{blue}{\Hom}(I,\fs)$ is spanned by those formal series whose coefficients
satisfy the condition
\Bi
\item for any fixed values of indices  $A_{++}$, $A_{+-}$, $B^{-+}$ and $B^{--}$ only {\em finitely many}\, $F^{A_{++}A_{+-}A_{-+}A_{--}}_{B^{++}B^{+-}B^{-+}B^{--}}\neq 0$,
\Ei
while the subspace $\color{red}{\Hom}(I,\fs)$ is characterized by
\Bi
\item for any fixed values of indices  $A_{++}$, $A_{-+}$, $B^{+-}$ and $B^{--}$ only {\em finitely many}\, $F^{A_{++}A_{+-}A_{-+}A_{--}}_{B^{++}B^{+-}B^{-+}B^{--}}\neq 0$,
\Ei
This gives us a  ``down to earth" characterization of the endomorphism
prop $\cE nd_{V^{2\text{-br}}}\cong \{ {\color{blue}{\Hom}(I,\fs)}\cap {\color{red}{\Hom}(I,\fs)}\}$.

\bip

\bip

{\Large
\section{\bf Action of the Grothendieck-Tiechm\"uller group on
some multi-oriented props}}

\bip

\subsection{An operad of multi-oriented graphs}\label{2: subsec on DRGC}
For any $l\geq -1$ and $k\geq 0$
let $G^{l+1\uparrow k+1}_{n,p}$ be a set of $(l+1)$-oriented $(k+1)$-directed (see \S {\ref{2new: subsec on
multid graphs}}) graphs $\Ga$ with $n$ vertices and $p$  edges such that
some bijections $V(\Ga)\rar [n]$ and $E(\Ga)\rar [p]$ are fixed, i.e.\ every vertex and every edge of $\Ga$
has a numerical label.
There is
a natural right action of the group $\bS_n \times  \bS_p$ on the set $G^{l+1\uparrow k+1}_{n,p}$ with $\bS_n$
acting by relabeling the vertices and  $\bS_p$ by relabeling the
edges.
For each fixed integer $d$, consider a collection of $\bS_n$-modules
$\cG ra^{l+1\uparrow k+1}_{d}=\{\cG ra^{l+1\uparrow k+1}_d(n)\}_{n\geq 1}$, where
$$
\cG ra_d^{l+1\uparrow k+1}(n):= \prod_{p\geq 0} \K \langle G_{n,p}^{l+1\uparrow k+1}\rangle \ot_{ \bS_p}
\sgn_p^{\ot |d-1|}  [p(d-1)].
$$
where  $\sgn_p$ is the 1-dimensional sign representation of $\bS_p$.
It has an (ordinary, i.e $1$-oriented!) operad structure with the composition rules
$$
\Ba{rccc}
\circ_i: &  \cG ra_d^{l+1\uparrow k+1}(n) \times \cG ra^{l+1\uparrow k+1}_d(m) &\lon & \cG ra_d^{l+1\uparrow
k+1}(n+m-1),  \ \ \forall\ i\in [n]\\
         &       (\Ga_1, \Ga_2) &\lon &      \Ga_1\circ_i \Ga_2,
\Ea
$$
given by substituting the graph $\Ga_2$ into the $i$-labeled vertex $v_i$ of $\Ga_1$ and taking the sum over
re-attachments of dangling edges (attached before to $v_i$) to vertices of $\Ga_2$
in all possible ways. If $l=k$ we abbreviate  $\cG ra_{d}^{(k+1)\text{-or}}=\cG ra_d^{k+1\uparrow k+1}$ and
call it the {\em operad of $(k+1)$-oriented graphs}.

 \sip

 Note also that for $l>l'$  the operad $\cG ra_{d}^{l+1\uparrow k+1}$ is a suboperad of $\cG
 ra_{d}^{l'+1\uparrow k+1}$.

\sip

There is a canonical injection
$$
\cG ra_{d}^{l+1\uparrow k+1} \lon \cG ra_{d}^{l+1\uparrow k+2}
$$
sending a $(k+1)$-directed graph $\Ga$ into a sum of $(k+2)$-directed graphs obtained from $\Ga$ by attaching
the new $(k+2)$nd direction to each edge in  two possible ways.

\sip

Let $\caL ie_d$ be a (degree shifted) ordinary operad of Lie algebras whose representations are  graded Lie
algebras equipped with the Lie bracket in degree
$1-d$, and consider the standard (cf.\ \cite{MV}) deformation complex
of the trivial morphism of operads,
\Beq\label{2: def of fGC_k,d^l}
\fGC_{d}^{l+1\uparrow k+1}:= \Def\left(\caL ie_d \stackrel{0}{\lon} \cG ra_{d}^{l+1\uparrow k+1}\right)\simeq
\prod_{n\geq 1}\cG ra_{d}^{l+1\uparrow k+1}(n)^{\bS_n}[d(1-n)]\ \ \ \ \forall\ k\geq 0,\ l\in
\{-1,0,1,\ldots, k\}.
\Eeq
This is a Lie algebra. Moreover, it admits  a non-trivial Maurer-Cartan element $\ga_0$
which corresponds to a morphism
$$
\ga_0: \caL ie_d\ \lon \  \cG ra_{d}^{l+1\uparrow k+1}
$$
 given explicitly on the generator (of homological degree $1-d$)
 $$
 \Ba{c}\begin{xy}
 <0mm,0.66mm>*{};<0mm,3mm>*{}**@{-},
 <0.39mm,-0.39mm>*{};<2.2mm,-2.2mm>*{}**@{-},
 <-0.35mm,-0.35mm>*{};<-2.2mm,-2.2mm>*{}**@{-},
 <0mm,0mm>*{\circ};<0mm,0mm>*{}**@{},
   <0.39mm,-0.39mm>*{};<2.9mm,-4mm>*{^{_2}}**@{},
   <-0.35mm,-0.35mm>*{};<-2.8mm,-4mm>*{^{_1}}**@{},
\end{xy}\Ea = (-1)^{d} \Ba{c}\begin{xy}
 <0mm,0.66mm>*{};<0mm,3mm>*{}**@{-},
 <0.39mm,-0.39mm>*{};<2.2mm,-2.2mm>*{}**@{-},
 <-0.35mm,-0.35mm>*{};<-2.2mm,-2.2mm>*{}**@{-},
 <0mm,0mm>*{\circ};<0mm,0mm>*{}**@{},
   <0.39mm,-0.39mm>*{};<2.9mm,-4mm>*{^{_1}}**@{},
   <-0.35mm,-0.35mm>*{};<-2.8mm,-4mm>*{^{_2}}**@{},
\end{xy} \Ea \in \caL ie_d(2)
 $$
  by the following explicit formula (cf.\ \cite{Wi1,Wi2})
\Beq\label{2:  map from Lie to dgra}
\ga_0 \left(\Ba{c}\begin{xy}
 <0mm,0.66mm>*{};<0mm,3mm>*{}**@{-},
 <0.39mm,-0.39mm>*{};<2.2mm,-2.2mm>*{}**@{-},
 <-0.35mm,-0.35mm>*{};<-2.2mm,-2.2mm>*{}**@{-},
 <0mm,0mm>*{\circ};<0mm,0mm>*{}**@{},
   <0.39mm,-0.39mm>*{};<2.9mm,-4mm>*{^{_2}}**@{},
   <-0.35mm,-0.35mm>*{};<-2.8mm,-4mm>*{^{_1}}**@{},
\end{xy}\Ea\right)=\sum_{\fa\in \f_k}\left(
\Ba{c}\resizebox{8.3mm}{!}{\xy
(0,1)*+{_1}*\cir{}="b",
(14,1)*+{_2}*\cir{}="c",
(7,3)*+{\fa}
\ar @{->} "b";"c" <0pt>
\endxy}
\Ea  + (-1)^d
\Ba{c}\resizebox{8.3mm}{!}{\xy
(0,1)*+{_2}*\cir{}="b",
(14,1)*+{_1}*\cir{}="c",
(7,3)*+{\fa}
\ar @{->} "b";"c" <0pt>
\endxy}
\Ea\right)=:\xy
 (0,0)*{\bullet}="a",
(5,0)*{\bu}="b",
\ar @{->} "a";"b" <0pt>
\endxy
\Eeq
where the summation runs over all possible ways to attach extra $k$ directions to the 1-oriented edge.
Note that elements of $\fGC_{d}^{l+1\uparrow k+1}$ can be identified with graphs from $\cG
ra_{d}^{l+1\uparrow k+1}$ whose vertices' labels are symmetrized (for $d$ even) or skew-symmetrized (for $d$
odd) so that in pictures we can forget about labels of vertices  and denote them by unlabelled black bullets
as in the formula above. Note also that graphs from  $\cG ra_{d}^{l+1\uparrow k+1}$ come equipped with a {\em
orientation}\,  which is a choice of ordering of edges (for $d$ even) or a choice of ordering of vertices
(for $d$ odd) up to an even permutation in both cases. Thus every graph $\Ga\in \fGC_{d}^{l+1\uparrow k+1}$
has at most two different orientations, $or$ and $or^{opp}$, and one has
the standard relation, $(\Ga, or)=-(\Ga, or^{opp})$; as usual, the data $(\Ga, or)$ is abbreviated to $\Ga$
(with some choice of orientation implicitly assumed).  Note that the homological degree of graph $\Ga$ from
$\fGC_{d}^{l+1\uparrow k+1}$ is given by
$$
|\Ga|=d(\# V(\Ga) -1) + (1-d) \# E(\Ga).
$$
We show in \cite{Me5} some other explicit Maurer-Cartan elements in the Lie algebra $\fGC_{d}^{l+1\uparrow
k+1}$ given by transcendental formulae; in this paper we need only $\ga_0$.
\sip

The above Maurer-Cartan element (\ref{2:  map from Lie to dgra}) makes
 $(\fGC_{d}^{l+1\uparrow k+1}, [\ ,\ ])$ into a {\em differential}\, Lie algebra with the differential
 \Beq\label{2: delta_0}
 \delta_0:= [\xy
 (0,0)*{\bullet}="a",
(5,0)*{\bu}="b",
\ar @{->} "a";"b" <0pt>
\endxy ,\ ].
\Eeq
 This dg Lie algebra contains a  dg subalgebra $\fcGC_{d}^{l+1\uparrow k+1}$ spanned by {\em connected}\,
 graphs which in turn contains   a dg Lie subalgebra
 $\GC_{d}^{l+1\uparrow k+1}$ spanned by connected graphs
with at least bivalent vertices. It was proven in \cite{Wi1,Wi2} (for the case $k=0$ and $l=-1,0$ but the
arguments works in greater generality) that the latter two subalgebras are quasi-isomorphic,
$$
H^\bu(\fcGC_{d}^{l+1\uparrow k+1}, \delta_0)=H^\bu(\GC_{d}^{l+1\uparrow k+1}, \delta_0)
$$

It was also proven in \cite{Wi1,Wi2} (in the cases $k=0$ and $l\in\{-1,0\}$ but the arguments work in greater
generality) that
$$
H^\bu(\fGC_{d}^{l+1\uparrow k+1}, \delta_0)= \odot^{\bu\geq 1}\left(H^\bu(\GC_{d}^{l+1\uparrow k+1},\delta_0)[-d]\right)[d]
$$
so that there is no loss of information to working solely with $\GC_{d}^{l+1\uparrow k+1}$  instead of the full graph
complex
$\fGC_{d}^{l+1\uparrow k+1}$.
There is a remarkable isomorphism of Lie algebras \cite{Wi1},
$$
H^0({\GC}_{2}^{0\uparrow 1},\delta_0)=\fg\fr\ft_1,
$$
where $\fg\fr\ft_1$ is the Lie algebra of the Grothendieck-Teichm\"uller group $GRT_1$ introduced by
Drinfeld in the context of the deformation quantization of Lie bialgebras. Nowadays, this group plays an
important role in many  areas of mathematics.

\sip

The multidirected graph complexes  have been introduced and studied in \cite{Z}; more precisely, Marko \v
Zivkovi\' c studied fully oriented graph complexes which are dual to the complexes
$
\left(\GC_{d}^{k+1\uparrow k+1},\delta_0\right), k\geq 0.
$
We often abbreviate  $\GC_{d}^{(k+1)\text{-or}}:= \GC_{d}^{k+1\uparrow k+1}$ for $k\geq 0$ and
$\GC_{d}^{0\text{-or}}:= \GC_{d}^{0\uparrow 1}$

\sip

Note that for $l'<l$ the Lie algebra $\GC_{d}^{l+1\uparrow k+1}$ is a Lie subalgebra of
$\GC_{d}^{l'+1\uparrow k+1}$.

\subsection{Cohomology of (partially) oriented multi-directed graph complexes}\label{2: subsec on cohom of
multi graph complexes}
For any $k\geq 0$ and any $-1\leq l\leq k$ there is an obvious map of graph complexes
 $$
i: \GC_{d}^{l+1\uparrow k+1} \lon \GC_{d}^{l+1\uparrow k+2}
 $$
which sends an $(l+1)$-oriented graph  with $k+1$ directions to  an $(l+1)$-oriented graph with $k+2$
directions by taking a sum over all possible ways to attach a new $(k+2)$-nd direction to each
$(k+1)$-directed edge.

\subsubsection{\bf Theorem \cite{Wi1}}\label{2: W theorem on l,k multi} {\em The injection $i:
\GC_{d}^{l+1\uparrow k+1} \lon \GC_{d}^{l+1\uparrow k+2}$ is a quasi-isomorphism of dg Lie algebras.}

\mip

This theorem was proved  by Thomas Willwacher in \cite{Wi1} in the case $k=0$, $l\in \{-1,0\}$, but the
argument works in greater generality. This result implies
$$
H^\bu(\GC_{d}^{l+1\uparrow k+1},\delta_0)=
H^\bu(\GC_{d}^{(l+1)\text{-or}}, \delta_0)
  \ \ \
\forall k\geq 0, \ \ -1\leq l\leq k.
$$
Put another way, multidirections which are not {\em oriented}\, can be forgotten, they do not give us
something really new.

\sip

Thomas Willwacher also proved the following
\subsubsection{\bf Theorem \cite{Wi2}}$\label{2: Thomas theorem on GC and GC-or}
H^\bu(\GC_{d}^{0\text{-or}},\delta_0)=H^\bu(\GC^{1\text{-or}}_{d+1}; \delta_0)$\ \
{\em  for any $d\in \Z$. }

\sip

In particular, one has an isomorphism
$$
H^0(\GC_{3}^{1\text{-or}}, \delta_0)=H^0(\GC_2^{0\text{-or}},\delta_0)=\fg\fr\ft_1
$$
which plays an important role in the homotopy theory of (involutive) Lie bialgebras \cite{MW1}.

\mip

This Theorem has been recently generalized to $(k+1)$-oriented graphs by  Marko \v Zivkovi\' c.

\subsubsection{\bf Theorem \cite{Z}}\label{2: Zivkovic theorem} {
$H^\bu(\GC_{d}^{(k+1)\text{-or}},
\delta_0)=H^\bu(\GC_{d+1}^{(k+2)\text{-or}},\delta_0) $
\ \em  for any $d\in \Z$ and any $k\geq 0$.}

\mip

Theorems {\ref{2: W theorem on l,k multi}} and {\ref{2: Zivkovic theorem}} imply the equalities
\Beq\label{GC l+1,k+1 = GC l+2 k+2}
H^\bu(\GC_{d}^{l+1\uparrow k+1},\delta_0)=
H^\bu(\GC_{d+1}^{l+2\uparrow k+2},\delta_0)=
  \ \ \
\forall d\in \Z,\ \ k\geq 0, \ \ -1\leq l\leq k.
\Eeq
In particular we have  isomorphisms of Lie algebras,
\Beq\label{GCd+2 d-or and grt}
H^0({\GC}_{d+2}^{d\text{-or}},\delta_0)= H^0({\GC}_{2}^{0\uparrow 1},\delta_0)=\fg\fr\ft_1,
\Eeq
for any $d\geq 0$. For $d=2$ and $d=3$ the algebro-geometric meanings of the associated graph complex
incarnations of the Grothendieck-Teichm\"uller group $GRT_1$ are clear: the $d=2$ case corresponds to the
action of $GRT_1$  (through cocycle representatives in ${\GC}_{2}^{0\uparrow 1}$) on universal Kontsevich
formality maps
associated with the deformation quantization of Poisson structures (given explicitly with the help of suitable
configuration spaces in the {\em two}\, dimensional upper half-plane \cite{Ko2}), while the case $d=3$
corresponds to the action of $GRT_1$
 (through cocycle representatives in ${\GC}_{3}^{1\text{-or}}$) on universal formality maps associated with
 the deformation quantization of Lie bialgebras (see \cite{MW2} where compactified configuration spaces in
 {\em three}\, dimensions have been used).

\sip

The above results tell us that the Grothendieck-Teichm\"uller group
survives in any geometric dimension $\geq 4$ but now in the multi-oriented graph complex incarnation. What can
the associated to $\fg\fr\ft_1$ degree zero cocycles in ${\GC}_{d+2}^{d\text{-or}}$  act on? It is an attempt to answer this question
which motivated much of the present work. In the first approximation the answer is that it acts
on
the multi-oriented props $\HoLBcd^{(c+d-1)\text{-or}}$ (more precisely, on their genus
completed versions $\wHoLBcd^{(c+d-1)\text{-or}}$), and it is not hard to see how. Recall the main result of
\cite{MW1} which says that there is a morphism
of dg Lie algebras
$$
F \colon \GC_{c+d+1}^{\text{1-or}}\to \Der(\wHoLBcd)
$$
where $\wHoLBcd$ is the genus completion of $\HoLBcd$ and $\Der(\wHoLBcd)$
is the Lie algebra of continuous derivations of $\wHoLBcd$ (see \cite{MW1} for some subtlety in its
definition). This map is a quasi-isomorphism (up to one rescaling class), and it can be given by a simple formula: for any $\Ga\in
\GC_{c+d+1}^{\text{1-or}}$ one has
$$
F(\Ga) =
\sum_{m,n\geq 1} \sum_{s:[n]\rar V(\Ga)\atop \hat{s}:[m]\rar V(\Ga)}  \Ba{c}\resizebox{11mm}{!}  {\xy
 (-6,7)*{^1},
(-3,7)*{^2},
(2.5,7)*{},
(7,7)*{^m},
(-3,-8)*{_2},
(3,-6)*{},
(7,-8)*{_n},
(-6,-8)*{_1},
(0,4.5)*+{...},
(0,-4.5)*+{...},
(0,0)*+{\Ga}="o",
(-6,6)*{}="1",
(-3,6)*{}="2",
(3,6)*{}="3",
(6,6)*{}="4",
(-3,-6)*{}="5",
(3,-6)*{}="6",
(6,-6)*{}="7",
(-6,-6)*{}="8",
\ar @{-} "o";"1" <0pt>
\ar @{-} "o";"2" <0pt>
\ar @{-} "o";"3" <0pt>
\ar @{-} "o";"4" <0pt>
\ar @{-} "o";"5" <0pt>
\ar @{-} "o";"6" <0pt>
\ar @{-} "o";"7" <0pt>
\ar @{-} "o";"8" <0pt>
\endxy}\Ea
$$
 where the second sum in  taken over all ways, $s$ and $\hat{s}$, of attaching the in- and outgoing legs to
 the graph $\Ga$, and then setting to zero every  graph containing a vertex with valency $\leq 2$ or with no
 input legs or no output legs
    (there is an implicit rule of signs in-built into this formula).
    In a complete analogy one can define an action of the dg Lie algebra $
\GC_{c+d+1}^{\text{(k+1)-or}}$  as  derivations on the multi-oriented dg prop $\wHoLBcd^{(k+1)\text{-or}}$, that is, a morphism of dg Lie algebras
$$
F \colon \GC_{c+d+1}^{(k+1)\text{-or}}\to \Der(\wHoLBcd^{(k+1)\text{-or}}).
$$
It was proven by Assar Andersson in \cite{A} that this map is a quasi-isomorphism (up to one rescaling class). This result together with  equality
(\ref{GCd+2 d-or and grt}) imply a highly non-trivial action of $GRT_1$
on the infinite family of the multi-oriented props $\wHoLBcd^{(c+d-1)\text{-or}}$, $c+d\geq 3$.

\bip

\bip

\def\cprime{$'$}

\end{document}